\documentclass[preprint,11pt]{elsarticle}
\usepackage[bookmarks=true,colorlinks=true,linkcolor=blue]{hyperref}

\pdfoutput=1

\biboptions{sort&compress}% compress bibtex entries

\usepackage{color}

\usepackage[usenames,dvipsnames,svgnames,table]{xcolor}

\usepackage{amssymb,amsmath}
\usepackage{amsthm}
\usepackage{multirow}

\usepackage{calc}
\usepackage[margin=1.in]{geometry}

\renewcommand{\url}[1]{}% ******** Remove URL's from bibtex entries ********
\newcommand{\citeCount}[1]{}% for citation counts papers for WDH citations
\newcommand{\f}[2]{\frac{#1}{#2}}% \frac short form
\newtheorem{theorem}{Theorem}

\newtheorem*{AMPInterfaceCondition}{AMP Interface Condition}

\newcommand{\esup}{{(e)}}
\newcommand{\psup}{{(p)}}

\newcommand{\dx}{\Delta x}

\newcommand{\dt}{\Delta t}

\newcommand{\dn}{\Delta n}

\newcommand{\bogus}[1]{{}}

\newenvironment{myIndent}%
 {\list{}{\leftmargin=0.1in\rightmargin=0.1in}\item[]}%
  {\endlist}

{\noindent\textbf{Procedure~}{#1}\begin{myIndent}\em}% start 
{\end{myIndent}}% end

\newcommand{\eqdef}{\overset{{\rm def}}{=}}% definition = 

\newcommand{\av}{\mathbf{ a}}
\newcommand{\bv}{\mathbf{ b}}

\newcommand{\ev}{\mathbf{ e}}
\newcommand{\fv}{\mathbf{ f}}
\newcommand{\gv}{\mathbf{ g}}

\newcommand{\iv}{\mathbf{ i}}

\newcommand{\nv}{\mathbf{ n}}

\newcommand{\rv}{\mathbf{ r}}

\newcommand{\tv}{\mathbf{ t}}
\newcommand{\uv}{\mathbf{ u}}
\newcommand{\vv}{\mathbf{ v}}
\newcommand{\wv}{\mathbf{ w}}
\newcommand{\xv}{\mathbf{ x}}

\newcommand{\Dv}{\mathbf{ D}}

\newcommand{\Fv}{\mathbf{ F}}
\newcommand{\Gv}{\mathbf{ G}}

\newcommand{\Iv}{\mathbf{ I}}

\newcommand{\Mv}{\mathbf{ M}}
\newcommand{\Nv}{\mathbf{ N}}
\newcommand{\Ov}{\mathbf{ O}}
\newcommand{\Pv}{\mathbf{ P}}
\newcommand{\Qv}{\mathbf{ Q}}
\newcommand{\Rv}{\mathbf{ R}}

\newcommand{\Real}{{\mathbb R}}

\newcommand{\zerov}{\mathbf{0}}

\newcommand{\Dc}{{\mathcal D}}

\newcommand{\Fc}{{\mathcal F}}
\newcommand{\Gc}{{\mathcal G}}

\newcommand{\Nc}{{\mathcal N}}

\newcommand{\omegav}{\boldsymbol{\omega}}
\newcommand{\tauv}{\boldsymbol{\tau}}
\newcommand{\phiv}{\boldsymbol{\phi}}

\newcommand{\sigmav}{\boldsymbol{\sigma}}

\newcommand{\grad}{\nabla}

\newcommand{\tableFont}{\scriptsize}

\newcommand{\num}[2]{#1e#2} % Use this macro to define the format of the numbers in the table

\clearpage

\newcommand{\rhos}{\rho_b}% body density 
\newcommand{\rhob}{\rho_b}% body density 

\newcommand{\normalss}{\sffamily}

\newcommand{\ds}{\Delta s}

\newcommand{\uHat}{\hat{u}}
\newcommand{\vHat}{\hat{v}}
\newcommand{\pHat}{\hat{p}}

\newcommand{\OmegaF}{\Omega}% fluid domain
\newcommand{\OmegaB}{{\Omega_b}}% domain of the body
\newcommand{\partialB}{{\Gamma_b}}% body boundary
\newcommand{\GammaB}{\Gamma_b}% body boundary

\newcommand{\strutt}{\rule{0pt}{9pt}}% strutt to make table column height bigger
\newcommand{\Gchv}{\Gc_{\rm{hv}}}% grid for the heart valve 
\newcommand{\Gcts}{\Gc_{\rm{ts}}}% grid for the 2d heart valve (two sticks)

\newcommand\tnv{{\mathbf{t}}}

\renewcommand\tv{\tnv}

\newcommand{\mrb}{m_{b}}% mass of the body
\newcommand{\xvb}{\xv_b}% position of the center of mass
\newcommand{\vvb}{\vv_b}% position of the center of mass
\newcommand{\avb}{\av_b}% acceleration of the center of mass
\newcommand{\bvb}{\bv_b}% angular acceleration of the center of mass
\newcommand{\xvcm}{\xv_b}% position of the center of mass
\newcommand{\vvcm}{\vv_b}% velocity of the center of mass
\newcommand{\dotxvcm}{\dot\xv_b}% position of the center of mass
\newcommand{\dotvvcm}{\dot\vv_b}% velocity of the center of mass

\newcommand{\omegavb}{\omegav_b}% body angular acceleration
\newcommand{\Eb}{E_b}% axes of inertial
\newcommand{\fvbe}{\fv_e}% extranal force on body
\newcommand{\gvbe}{\gv_e}% external torque on body

\newcommand{\mb}{{m_b}}

\newcommand{\Da}{\Dc}% added damping matrix
\newcommand{\Dvv}{\Da^{v v}}
\newcommand{\Dvw}{\Da^{v \omega}}
\newcommand{\Dwv}{\Da^{\omega v}}
\newcommand{\Dww}{\Da^{\omega \omega}}
\newcommand{\Dvva}{\tilde\Da^{v v}}
\newcommand{\Dvwa}{\tilde\Da^{v \omega}}

\newcommand{\Dwwa}{\tilde\Da^{\omega \omega}}
\newcommand{\Ib}{\Iv_b}% moment of inertia matrix
\newcommand{\dArea}{dS}% area element for surface integrals

\newcommand{\ampRB}{AMP-RB}% name of scheme

\newcommand{\ADtensor}{\Dv}% added damping composite tensor 
\newcommand{\dr}{\Delta r}

\newlength{\ycbTop}% For colour bar
\newlength{\ycbMid}%

\newcommand{\bodyStepIComment}{Preliminary body evolution step}

\newcommand{\wHat}{{\hat w}}
\def\aa{r_1}
\def\bb{r_2}

\newcommand{\NAv}{\Nc_v}
\newcommand{\IbBar}{{\bar I_b}}

\newcommand{\fap}{\fv}% advection and pressure terms

\newcommand{\diskRadius}{{R_b}}

\newcommand{\hj}{h^{(j)}}

\newcommand{\channelRight}{L}% right end of the channel 
\newcommand{\channelHeight}{H}% right end of the channel 
\newcommand{\Gcp}{\Gc_p}% piston grid 
\newcommand{\Amplitude}{\alpha_b}

\usepackage{tikz}
\usetikzlibrary{calc,patterns,angles,quotes}

\newlength{\tfwidth}
\newlength{\tfheight}
\newlength{\tfxa}
\newlength{\tfxb}
\newlength{\tfya}
\newlength{\tfyb}
\newcommand{\trimFigWithBox}[6]{%
\setlength\fboxsep{0pt}%
\setlength\fboxrule{1.0pt}% border thickness
\fbox{\includegraphics[width=#2, clip, trim=#3 #4 #5 #6]{#1}}%
}
\newcommand{\trimFigNoBox}[6]{%
\setlength\fboxsep{1pt}% note: make this 1pt and rule thickness zero so box size matches that below
\setlength\fboxrule{0.0pt}% border thickness
\fbox{\includegraphics[width=#2, clip, trim=#3 #4 #5 #6]{#1}}%
}
\newcommand{\trimFigHeightWithBox}[6]{%
\setlength\fboxsep{0pt}%
\setlength\fboxrule{1.0pt}% border thickness
\fbox{\includegraphics[height=#2, clip, trim=#3 #4 #5 #6]{#1}}%
}
\newcommand{\trimFigHeightNoBox}[6]{%
\setlength\fboxsep{1pt}% note: make this 1pt and rule thickness zero so box size matches that below
\setlength\fboxrule{0.0pt}% border thickness
\fbox{\includegraphics[height=#2, clip, trim=#3 #4 #5 #6]{#1}}%
}
\newcommand{\trimFig}[6]{%
\setlength{\tfwidth}{(#2+#2*\real{#3})+#2*\real{#4}}%   % width of un-clipped fig
\setlength{\tfheight}{(#2+#2*\real{#5})+#2*\real{#6}}%
\setlength{\tfxa}{\tfwidth*\real{#3}}%
\setlength{\tfxb}{\tfwidth*\real{#4}}%
\setlength{\tfya}{\tfheight*\real{#5}}%
\setlength{\tfyb}{\tfheight*\real{#6}}%
\trimFigNoBox{#1}{#2}{\tfxa}{\tfya}{\tfxb}{\tfyb}%
}

\newsavebox\figBox

\newcommand{\trimw}[6]{%
\sbox\figBox{\includegraphics{#1}}
\setlength{\tfwidth}{\the\wd\figBox}
\setlength{\tfheight}{\the\ht\figBox}
\setlength{\tfxa}{\tfwidth*\real{#3}}%
\setlength{\tfxb}{\tfwidth*\real{#4}}%
\setlength{\tfya}{\tfheight*\real{#5}}%
\setlength{\tfyb}{\tfheight*\real{#6}}%
\trimFigNoBox{#1}{#2}{\tfxa}{\tfya}{\tfxb}{\tfyb}%
}

\newcommand{\trimwb}[6]{%
\sbox\figBox{\includegraphics{#1}}
\setlength{\tfwidth}{\the\wd\figBox}
\setlength{\tfheight}{\the\ht\figBox}
\setlength{\tfxa}{\tfwidth*\real{#3}}%
\setlength{\tfxb}{\tfwidth*\real{#4}}%
\setlength{\tfya}{\tfheight*\real{#5}}%
\setlength{\tfyb}{\tfheight*\real{#6}}%
\trimFigWithBox{#1}{#2}{\tfxa}{\tfya}{\tfxb}{\tfyb}%
}

\newcommand{\trimh}[6]{%
\sbox\figBox{\includegraphics{#1}}
\setlength{\tfwidth}{\the\wd\figBox}
\setlength{\tfheight}{\the\ht\figBox}
\setlength{\tfxa}{\tfwidth*\real{#3}}%
\setlength{\tfxb}{\tfwidth*\real{#4}}%
\setlength{\tfya}{\tfheight*\real{#5}}%
\setlength{\tfyb}{\tfheight*\real{#6}}%
\trimFigHeightNoBox{#1}{#2}{\tfxa}{\tfya}{\tfxb}{\tfyb}%
}

\newcommand{\trimhb}[6]{%
\sbox\figBox{\includegraphics{#1}}
\setlength{\tfwidth}{\the\wd\figBox}
\setlength{\tfheight}{\the\ht\figBox}
\setlength{\tfxa}{\tfwidth*\real{#3}}%
\setlength{\tfxb}{\tfwidth*\real{#4}}%
\setlength{\tfya}{\tfheight*\real{#5}}%
\setlength{\tfyb}{\tfheight*\real{#6}}%
\trimFigHeightWithBox{#1}{#2}{\tfxa}{\tfya}{\tfxb}{\tfyb}%
}
\usepackage{algorithm} 
\usepackage{algpseudocode} % algorithmcx 

\begin{document}

\small

\begin{frontmatter}
\title{
A stable partitioned FSI algorithm for rigid bodies and incompressible flow in three dimensions
   }

\author[rpi]{J.~W.~Banks\fnref{DOEThanks,PECASEThanks}}
\ead{banksj3@rpi.edu}

\author[rpi]{W.~D.~Henshaw\fnref{DOEThanks,NSFgrantNew}}
\ead{henshw@rpi.edu}

\author[rpi]{D.~W.~Schwendeman\fnref{DOEThanks,NSFgrantNew}}
\ead{schwed@rpi.edu}

\author[rpi]{Qi Tang \corref{cor1}\fnref{QiThanks}}
\ead{tangq3@rpi.edu}

\address[rpi]{Department of Mathematical Sciences, Rensselaer Polytechnic Institute, Troy, NY 12180, USA.}

\cortext[cor1]{Department of Mathematical Sciences, Rensselaer Polytechnic Institute, 110 8th Street, Troy, NY 12180, USA.}

\fntext[DOEThanks]{This work was supported by contracts from the U.S. Department of Energy ASCR Applied Math Program.}

\fntext[PECASEThanks]{Research supported by a U.S. Presidential Early Career Award for Scientists and Engineers.}

% New NSF
\fntext[NSFgrantNew]{Research supported by the National Science Foundation under grant DMS-1519934.}

\fntext[QiThanks]{Research supported by the Eliza Ricketts Postdoctoral Fellowship.}

\begin{abstract}
This paper describes a novel partitioned algorithm for fluid-structure interaction (FSI) problems that couples the motion of 
rigid bodies and incompressible flow. This is the first partitioned algorithm that remains stable and second-order accurate,
without sub-time-step iterations, for very light, and even zero-mass, bodies in three dimensions. This new 
{\em added-mass partitioned} (AMP) algorithm extends the previous developments in~\cite{rbinsmp2017, rbins2017}
by generalizing the added-damping tensors to account for arbitrary three-dimensional rotations, 
and by employing a general quadrature for the surface integral over a rigid body to derive the discrete AMP 
interface condition for the fluid pressure. Stability analyses for two three-dimensional model problems show 
that the algorithm remains stable for bodies of any mass {when applied to the relevant model problems}. The resulting AMP algorithm is implemented in 
parallel using a moving composite grid framework to treat one or more rigid bodies in complex three-dimensional 
configurations.  The new three-dimensional algorithm is verified and validated though several benchmark problems, 
including the motion of a sphere in a viscous incompressible fluid and the interaction of a bi-leaflet mechanical 
heart valve and a pulsating fluid. Numerical simulations confirm the predictions of the stability analysis {even for complex problems}, and show that 
the AMP algorithm remains stable, without sub-iterations, for light and even zero-mass three-dimensional rigid bodies 
of general shape. These benchmark problems are further used to examine the parallel performance of the algorithm and to  
investigate the conditioning of the linear system for the pressure including the newly derived AMP interface conditions.

\end{abstract}


\end{frontmatter}

% ------------- Table of contents
\clearpage
\tableofcontents

\clearpage

\section{Introduction}

Fluid-structure interaction (FSI) problems involving the motion of rigid bodies in a viscous
incompressible fluid arise in many important scientific and engineering applications, and the
accurate and stable simulation of such problems remains a significant challenge.  Many available
numerical approaches to simulate FSI problems of this type are {\em partitioned}
methods.  For such methods, the equations governing the fluid (the incompressible Navier-Stokes equations)
and the rigid solids (the Newton-Euler equations) are advanced separately in each time-step of the
algorithm with the matching conditions at the fluid-solid interfaces providing a coupling of the
solutions.  For example, in a standard traditional-partitioned scheme, the coupling is implemented
sequentially by first using the motion of the rigid bodies to provide the position and velocity as boundary conditions for the fluid
solver, and then evaluating the forces and torques on the bodies from the computed fluid stress to advance
the bodies.  While this approach is appealing due to its ease of adapting available
solvers for the fluid and rigid bodies, it is known that traditional partitioned methods can posses
severe stability issues due to their choice of fluid-solid interface coupling conditions. These so-called
added-mass instabilities are especially challenging for the case of light solids where added-mass effects are large.

In recent work~\cite{rbinsmp2017, rbins2017}, a new partitioned scheme was developed and shown, in two-dimensions, to be
accurate and stable for any mass of the solid, including the extreme case of a solid with zero mass,
thus overcoming instabilities due to added-mass effects.  The key development for this new scheme,
referred to as the~\ampRB~scheme, was the treatment of the interface coupling.  In contrast to a
Dirichlet-Neumann-type coupling as in the aforementioned traditional (TP-RB) scheme,
the~\ampRB~scheme employs a new coupling strategy in which the translational and rotational {\em
accelerations} of the rigid solid are balanced with the fluid accelerations, leading to an AMP~interface condition for the fluid
pressure that addresses the added-mass instabilities.
This interface condition also incorporates an {\em added-damping} tensor designed to suppress a
secondary instability associated with viscous fluid stresses.  The work
in~\cite{rbinsmp2017, rbins2017} provides a detailed derivation of the AMP~interface conditions, as
well as an analysis of the instabilities associated with added-mass and added-damping effects for
two-dimensional FSI problems. 

The motion of rigid bodies in three dimensions is significantly more complicated than in two
dimensions, involving one additional translational and two additional rotational degrees of freedom.
A question arises
as to whether the proposed AMP interface conditions, which were founded on one-dimensional model
problems, can be extended to complex three-dimensional motions of general-shaped bodies.
The added-damping tensor
in three dimensions, for example, consists of four $3\times 3$ component tensors (whose entries have a complicated dependence on
the viscosity, time-step and grid spacing), and the~\ampRB~interface conditions must then incorporate 
surface integrals of these tensors into the linear system for pressure using accurate numerical approximations 
of the surface integrals for complex geometries. The aim of the present work is therefore to extend the 
formulation and implementation to three dimensions,
and to evaluate this new scheme using careful analyses and numerical tests.
To that end, a stability analysis of the three-dimensional time-stepping scheme is performed
for an FSI model problem involving a rigid spherical solid coupled to a viscous fluid in a
spherical shell surrounding the solid. This problem naturally decomposes into two 
independent model problems; the first related to translational motion of
the solid and the second to rotational motion of the solid.  The first sub-problem reveals that 
the~\ampRB~scheme is stable, without sub-time-step
relaxation iterations of the matching conditions, even for a light solid where added-mass effects
are large. The second sub-problem is used to prove that the scheme is stable when added-damping effects are 
important.

The implementation of the~\ampRB~scheme is carried out using the moving overlapping grid
framework~\cite{mog2006}.  In this approach, a collection of overlapping structured grids are used
to cover the domain of the fluid with moving interface-fitted grids attached to the surfaces of the
evolving rigid solids. The evolution equations for the rigid solids, as well as the evaluation of
added-damping tensors for the AMP-interface conditions, requires accurate approximations of
integrals over the surface of the solids, including the general case of surfaces covered by multiple overlapping grids.
An effective approach to deriving these approximate surface quadrature weights using the left-null vector of 
an overlapping-grid elliptic boundary-value problem is described. The fluid solver used here is based on the fractional-step
scheme described in~\cite{ICNS,splitStep2003}, with changes to incorporate the new AMP coupling scheme.
In addition, the full~\ampRB~scheme is implemented in parallel using a domain-decomposition approach based on the work
described in~\cite{max2006b, pog2008a}.

In addition to the development, analysis, and implementation of the~\ampRB~scheme, the accuracy and stability
of the algorithm is verified by computing solutions of several benchmark FSI problems in three
dimensions.  As a first check of the~\ampRB~scheme, we consider the one-dimensional motion of a
piston in a three-dimensional fluid chamber.  An exact solution is available for this problem, and
thus it provides a useful benchmark problem to verify the accuracy and stability of the method for a
range of solid densities.  A second benchmark problem involves the settling of a light solid in a
fluid-filled container.  A self-convergence test of solutions of the~\ampRB~scheme is performed, and
the results are also compared with experimental and numerical data available in the literature.  As
a third problem, we examine a solid particle rising or falling under gravity in a long fluid
chamber, and compare the solutions given by the~\ampRB~scheme with the results given by other
algorithms discussed in the literature.  We use this problem, in particular, to demonstrate the
stability of the~\ampRB~scheme for a problem involving a very light solid (and even zero-mass solid)
where added-mass effects are large.  As a final demonstration of the~\ampRB~scheme for a complex FSI
problem in three dimensions, we consider the dynamical behavior of two flapping leaflets in a model
of a bi-leaflet mechanical heart valve.  This problem is chosen to illustrate the effectiveness of
the new scheme for a difficult FSI problem that has received significant attention in the recent
literature.  For example, the problem has been simulated using various FSI algorithms, such as those
discussed in~\cite{cheng2004three, dumont2007comparison, tai2007numerical, borazjani2008curvilinear,
  nobili2008numerical, de2009direct, borazjani2010high}, and the review
paper~\cite{sotiropoulos2009review} summarizes the practical and numerical challenges of the
problem.  Our particular attention is on the strong added-mass effects in this simulation due to a
very small moment of inertia of the leaflets.  Numerical instabilities have been
reported~\cite{borazjani2008curvilinear, de2009direct, borazjani2010high}, and sub-time-step
iterations of the matching conditions are typically used by others to stabilize their FSI
algorithms.  Apart from difficulties associated with added-mass effects, this problem is challenging
since the two leaflets are hinged and thus follow a constrained motion.  Thus, this final FSI
problem provides a significant test of the accuracy and stability of the~\ampRB~scheme in three
dimensions, as well as a excellent check of the various implementation issues associated with the
scheme as noted above.

It is worth noting that the FSI regime considered in this paper has drawn significant attention due to its importance for many applications.  Besides the moving overlapping grid approach employed here, a variety  of classical numerical techniques for moving complex geometries have been extended to this regime, including arbitrary Lagrangian--Eulerian (ALE) methods~\cite{takashi1992, HuPatankarZhu2001, vierendeels2005analysis}, level-set methods~\cite{coquerelle2008vortex, gibou2012efficient}, fictitious domain methods~\cite{glowinski1999distributed}, embedded boundary methods~\cite{costarelli2016embedded} and immersed boundary methods~\cite{kajishima2002interaction, uhlmann2005immersed, kim2006immersed, lee2008immersed, borazjani2008curvilinear, breugem2012second, kempe2012improved, yang2012simple, bhalla2013unified, yang2015non, wang2015strongly, kim2016penalty, lacis2016stable}.  Recently, new approaches to handle moving geometries have also been developed for this regime, such as methods based on boundary-integral equations~\cite{corona2017integral} and implicit mesh discontinuous Galerkin methods~\cite{Saye2017a, Saye2017b}.

As noted above, a significant issue for any partitioned scheme is the numerical treatment of the interface conditions coupling the evolution of the fluid and the solid.  The~\ampRB~scheme, first introduced in~\cite{rbinsmp2017, rbins2017}, uses a new AMP~interface condition designed and shown to suppress added-mass and added-damping instabilities.  The standard TP-RB~scheme requires costly sub-time-step relaxation iterations of the matching conditions to suppress these instabilities, see~\cite{vierendeels2005analysis, borazjani2008curvilinear, wang2015strongly, Koblitz2016, kadapa2017} for example.  Apart from these two coupling approaches, other strategies have been proposed in the literature with the aim of suppressing added-mass instabilities, see~\cite{kempe2012improved, yang2015non, lacis2016stable} for example.  Through this work, the stability bound for the ratio of the density of the solid to that of the fluid has improved over the years.  For instance, the immersed boundary method described in~\cite{lacis2016stable} remains stable for a density ratio as low as~$10^{-4}$.  In this approach, a set of equations that couples the pressure with the motion for the body is exposed through an approximate block LU~factorization.  This treatment bears some similarity with the handling of added-mass effects in the~\ampRB~scheme except that the AMP~interface condition used here is derived at a continuous level, which can facilitate the development of high-order accurate schemes.  Recently, Kempe and Fr{\"o}hlich, and their collaborators, proposed two other coupling strategies for immersed-boundary methods applied to spherical particles, see~\cite{schwarz2015temporal, tschisgale2017non}.  Their numerical results demonstrate that both of the approaches are stable for spherical particles of any mass, including the case of zero mass, without sub-time-step iterations.  In~\cite{schwarz2015temporal}, a virtual mass with an ad~hoc parameter is incorporated into the equations of motion of the rigid body.  This approach appears to have a similar formulation as the added-damping tensor in the AMP~interface condition.  However, the added-damping tensor in our approach is again derived at a continuous level for bodies of general shapes, and thus it is not limited to spherical bodies and no ad~hoc parameters are needed.  Another difference between the approach in~\cite{schwarz2015temporal} and one discussed here is that the added-damping tensor is only used to handle added-damping effects, while the virtual mass in~\cite{schwarz2015temporal} is used to handle both added-mass and added-damping effects.  In their more recent work~\cite{tschisgale2017non}, the mass and inertia of the fluid in a Lagrangian layer surrounding the particle surface are exposed at a continuous level, and this is incorporated into the equations of motions of the rigid body.  This approach results in a strongly coupled FSI scheme that only needs a semi-implicit update and avoids sub-time-step iterations completely.  Although work in~\cite{tschisgale2017non} only considered spherical particles, the approach appears to be promising and seems extendible to bodies of general shape.

The~\ampRB~scheme has been implemented in serial and parallel within the {\em Overture} object-oriented framework. This implementation, as well as the numerical examples presented in this work, is freely available at~{\tt overtureFramework.org}.  The parallel implementation follows a domain-decomposition approach on overlapping grids discussed in~\cite{pog2008a}.  Note that the~\ampRB~scheme discussed herein follows the line of research on AMP schemes for various FSI regimes, including inviscid compressible fluids coupled to rigid solids~\cite{lrb2013}, linearly elastic solids~\cite{fsi2012}, nonlinear elastic solids~\cite{flunsi2016}, as well as incompressible fluids coupled to elastic bulk solids~\cite{fib2014}, elastic structural shells~\cite{fis2014} and deforming beams~\cite{beamins2016}. 

The remaining sections of the paper are organized as follows.  The governing equations are given in
Section~\ref{sec:governing}.  The details of the~\ampRB~algorithm are outlined in
Section~\ref{sec:algorithm}, including a brief review of the AMP interface condition and the
time-stepping algorithm.  A stability analysis of the time-stepping scheme for two FSI model
problems in a three-dimensional geometry is also presented in Section~\ref{sec:algorithm}.
Section~\ref{sec:numericalApproach} provides a description of the numerical approach. This
  section includes a brief outline of the moving overlapping grid approach, a description of the
  approach for computing surface integrals on bodies covered by multiple overlapping grids, and a
  discussion of the parallel implementation.  Numerical results for the various benchmark problems
  are presented in Section~\ref{sec:numericalResults}, and concluding remarks as well as future
  directions are given in Section~\ref{sec:conclusions}.

% end reva
%A collision model designed specifically for the heart-valve problem is discussed in~\ref{sec:collision}.
%A collision model designed specifically for the heart valve simulation is discussed in Section~\ref{sec:collision}.

% In those numerical schemes, a variety of coupling strategies between rigid body and incompressible fluids have been proposed. 
%The introduction should  focus on low density/mass/inertia, added-mass/damping and 3D work in the literature 
%The focus: 3D, heart valve, light rigid body, new reference.
%To successfully implement the algorithm, we rely on several numerical approaches,
%including the moving three-dimensional overlapping grids, an accurate evaluation of the surface integral on bodies covered by multiple overlapping grids and parallel implementation.

\section{Governing equations} \label{sec:governing}

We consider an FSI problem involving one or more rigid bodies fully immersed in an incompressible fluid as illustrated in Figure~\ref{fig:heartValve3Dintro}.  The fluid occupies the domain, $\xv\in\OmegaF(t)$ at a time~$t$, and is governed by the incompressible Navier-Stokes equations,
\begin{alignat}{3}
  &  \frac{\partial\vv}{\partial t} + (\vv\cdot\grad)\vv 
                 =  \frac{1}{\rho} \grad\cdot\sigmav  , \qquad&& \xv\in\OmegaF(t) ,  \label{eq:fluidMomentum}  \\
  & \grad\cdot\vv =0,  \quad&& \xv\in\OmegaF(t) ,  \label{eq:fluidDiv3d}
\end{alignat}
where $\vv=\vv(\xv,t)$ is the fluid velocity, $\sigmav=\sigmav(\xv,t)$ is the stress tensor, and $\rho$ is the (constant) fluid density.  The fluid stress tensor is given by 
\begin{equation}
  \sigmav = -p \Iv + \tauv,\qquad \tauv = \mu \left[ \grad\vv + (\grad\vv)^T \right],
  \label{eq:fluidStress}
\end{equation}
where $p=p(\xv,t)$ is the pressure, $\Iv$ is the identity tensor, $\tauv$ is the viscous stress tensor and $\mu$ is the (constant) dynamic viscosity. 
Equations~\eqref{eq:fluidMomentum}--\eqref{eq:fluidStress} can be used to derive an elliptic equation for the pressure satisfying
\begin{alignat}{3}
  &  \Delta p &= - \rho \grad\vv:(\grad\vv)^T ,  \qquad&& \xv\in\OmegaF(t) , \label{eq:fluidPressure}
\end{alignat}
where
\begin{align*}
\grad\vv:(\grad\vv)^T 
      ~ \eqdef \sum_{i=1}^{3}\sum_{j=1}^{3} \frac{\partial v_i}{\partial x_j}\frac{\partial v_j}{\partial x_i},
\end{align*}
with an additional boundary condition, $\grad\cdot\vv=0$ for $\xv\in\partial\OmegaF(t)$, see~\cite{ICNS} for more details.
In practice, our time-stepping scheme uses the velocity-pressure form of the incompressible Navier-Stokes equations, in which the continuity equation~\eqref{eq:fluidDiv3d} is replaced by the Poisson equation~\eqref{eq:fluidPressure}. 

{% ------ STREAMLINES ------
% ============================ DRAW ===================================
\newcommand{\figWidth}{4.5cm}
\newcommand{\trimfig}[2]{\trimw{#1}{#2}{.02}{.02}{.1}{.1}}
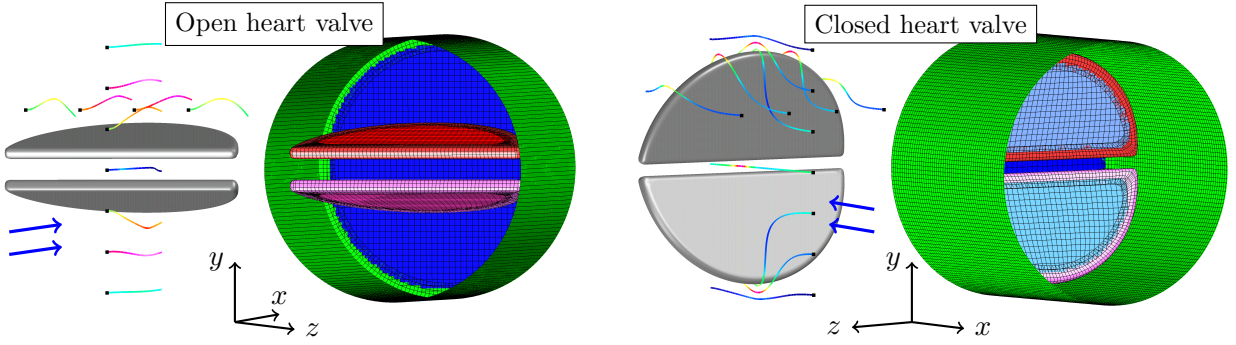
\begin{figure}[htb]
\begin{center}
%\resizebox{16cm}{!}{% START resize box
\begin{tikzpicture}[scale=1]
  \useasboundingbox (0.0,.5) rectangle (16.,4.3);  % set the bounding box (so we have less surrounding white space)
 \draw(-1,0) node[anchor=south west,xshift=-4pt,yshift=+0pt] {\trimfig{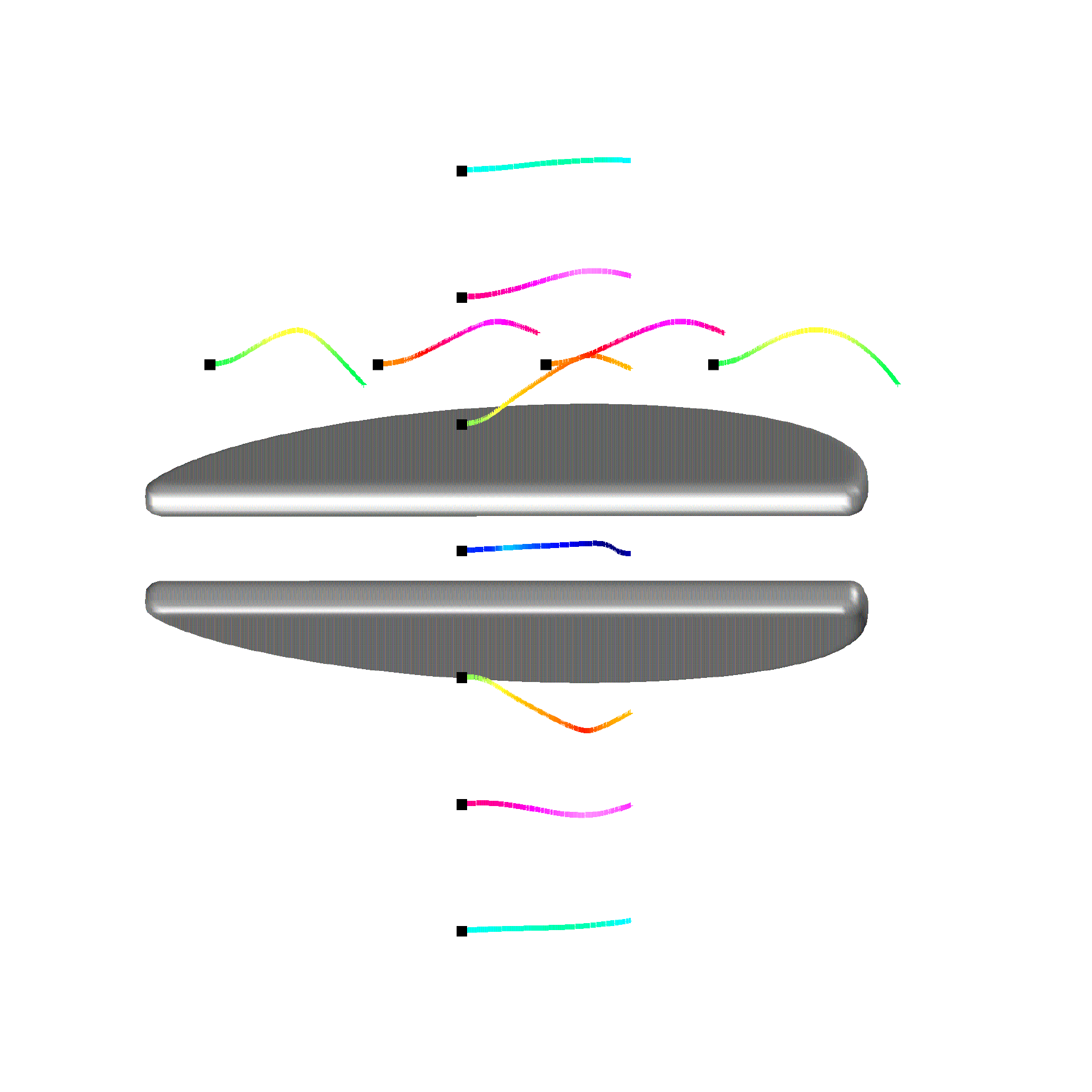}{\figWidth}};
 \draw(2.8,0) node[anchor=south west,xshift=-4pt,yshift=+0pt] {\trimfig{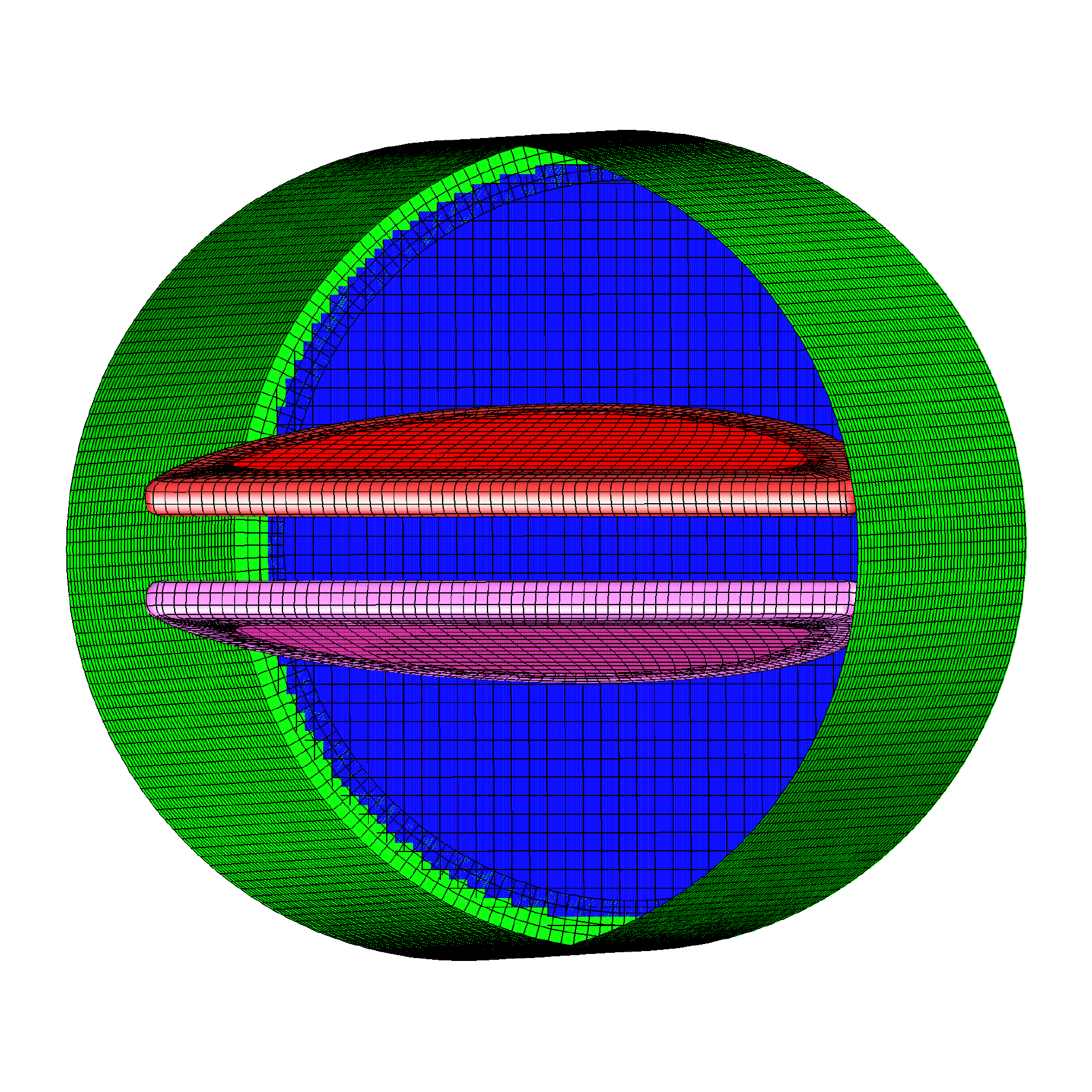}{\figWidth}};
 \draw(7.4,0) node[anchor=south west,xshift=-4pt,yshift=+0pt] {\trimfig{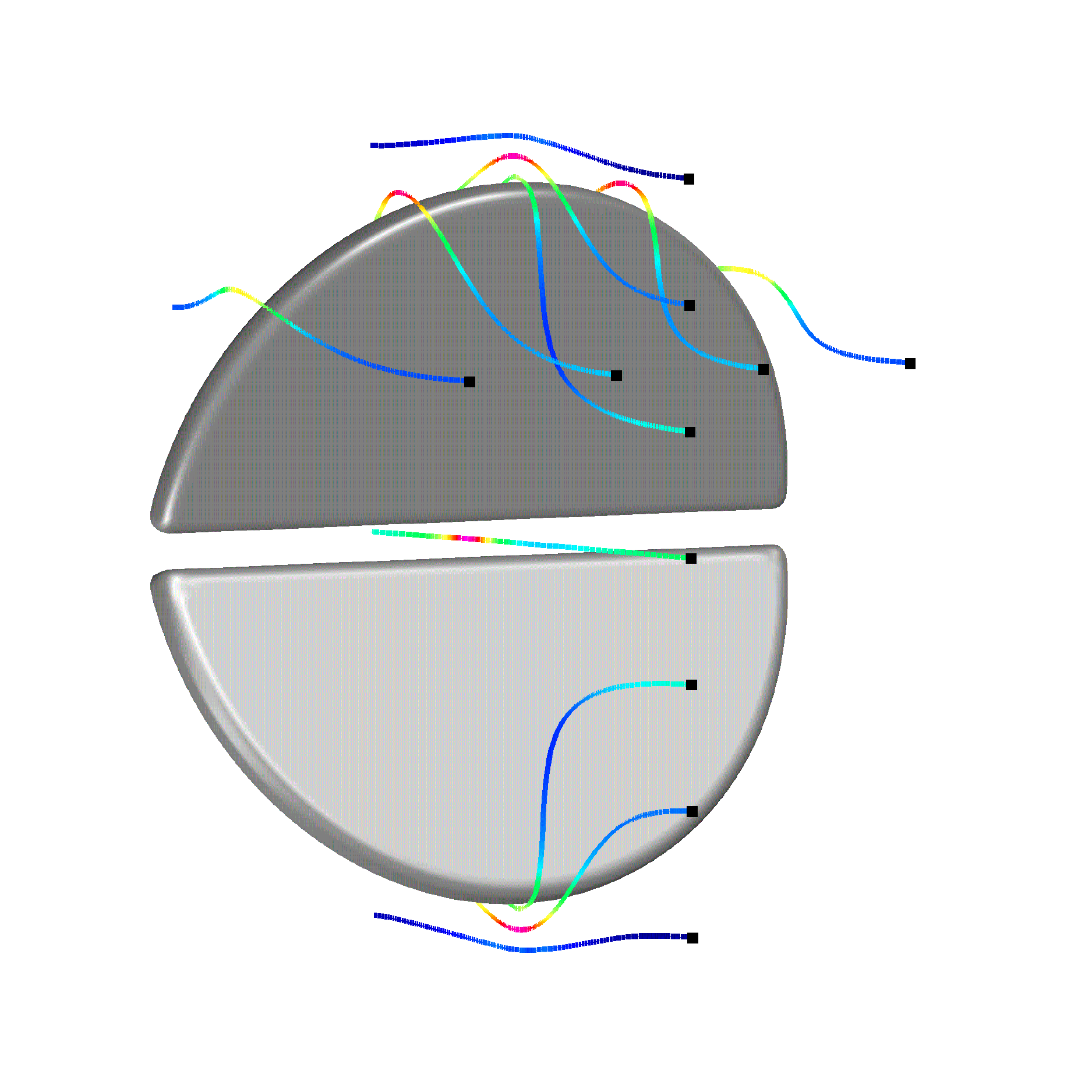}{\figWidth}};
 \draw(11.3,0) node[anchor=south west,xshift=-4pt,yshift=+0pt] {\trimfig{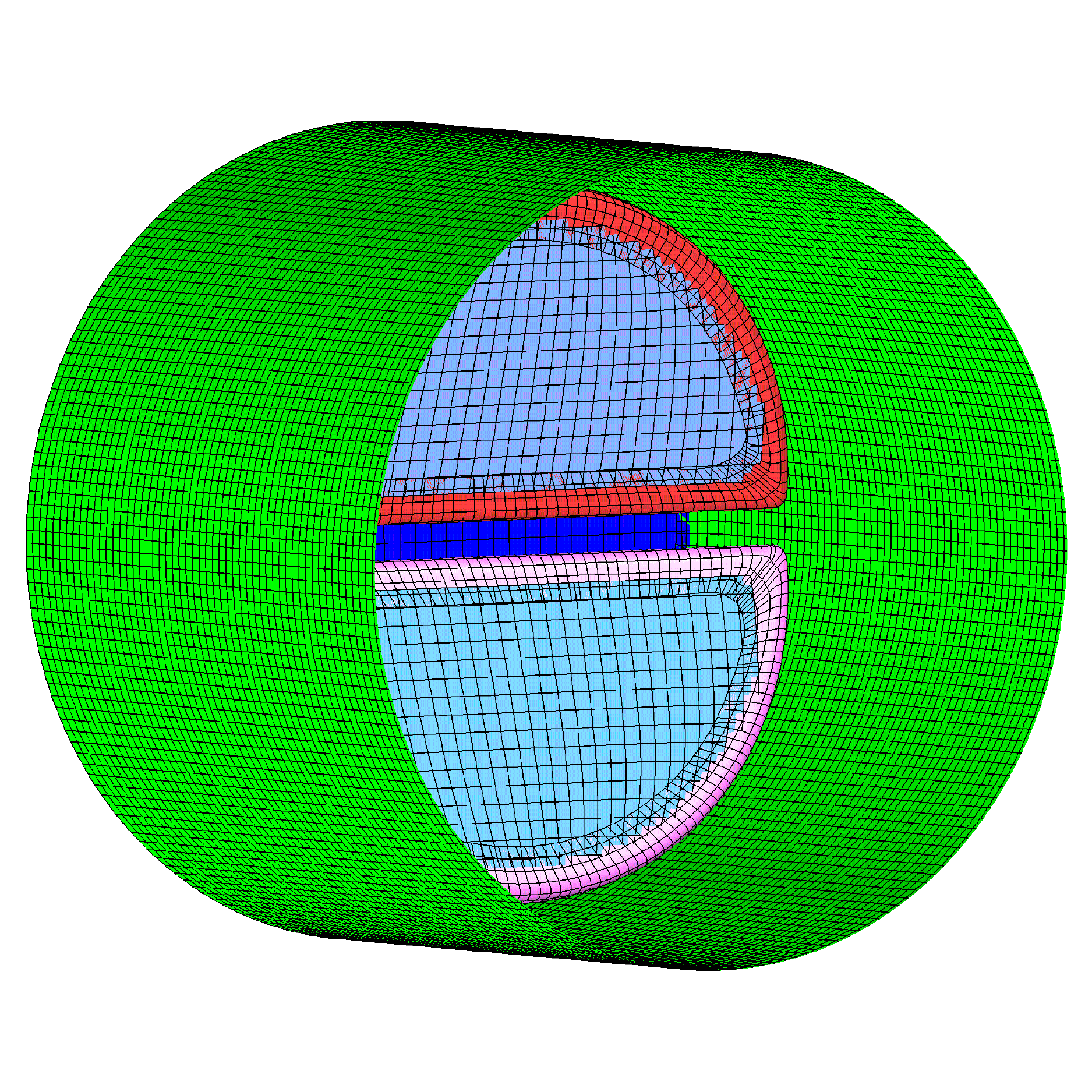}{\figWidth}};
\draw(2.,4) node[draw,fill=white,anchor=west,xshift=2pt,yshift=-1pt] {\small Open heart valve};
\draw(10.5,4) node[draw,fill=white,anchor=west,xshift=2pt,yshift=-1pt] {\small Closed heart valve};
%open inflow
\draw[->,very thick,blue] (0,.9) -- ( .7,1.); 
\draw[->,very thick,blue] (0,1.2) -- ( .7,1.3); 
%close inflow
\draw[->,very thick,blue] (11.5,1.5) --(10.9 ,1.6); 
\draw[->,very thick,blue] (11.5,1.2) --(10.9,1.3); 
{
\def\thetax{10}
\def\thetay{-7}
\def\thetaz{90}
\def\radx{.6}
\def\rady{.8}
\def\radz{.8}
%plot coordiantes
\begin{scope}[xshift=3.cm,yshift=0.cm]
\draw[thick,->] (0,0) -- ({\radx*cos(\thetax)},{\radx*sin(\thetax)}) node[anchor=south]{$x$};
\draw[thick,->] (0,0) -- ({\rady*cos(\thetay)},{\rady*sin(\thetay)}) node[anchor=west]{$z$}; 
\draw[thick,->] (0,0) -- ({\radz*cos(\thetaz)},{\radz*sin(\thetaz)}) node[anchor=east]{$y$};
\end{scope}
}
{
\def\thetax{-7}
\def\thetay{185}
\def\thetaz{90}
\def\radx{.7}
\def\rady{.8}
\def\radz{.8}
%plot coordiantes
\begin{scope}[xshift=12cm,yshift=0.cm]
\draw[thick,->] (0,0) -- ({\radx*cos(\thetax)},{\radx*sin(\thetax)}) node[anchor=west]{$x$};
\draw[thick,->] (0,0) -- ({\rady*cos(\thetay)},{\rady*sin(\thetay)}) node[anchor=east]{$z$}; 
\draw[thick,->] (0,0) -- ({\radz*cos(\thetaz)},{\radz*sin(\thetaz)}) node[anchor=east]{$y$};
\end{scope}
}
%
% ---------------
% grid:
%\draw[step=1cm,gray] (0,0) grid (16,4.3);
\end{tikzpicture}
%}% end resize box
\end{center}
  \caption{Flow past the bi-leaflet mechanical heart valve at the open and closed states. 
      Note the change of the flow directions at the two positions. 
      Computed streamlines (colored by the flow speed) and the composite grids are presented.  
}
  \label{fig:heartValve3Dintro}
\end{figure}
}

Within the FSI problem, each body has mass $m_b$ and moment of inertia $\Ib$ and occupies the domain $\OmegaB(t)$.
The dynamics of the body are characterized by
\begin{alignat*}{3}
  & \avb(t)\in\Real^3         ~&&: \text{linear acceleration of the centre of mass},\\ 
  & \bvb(t)\in\Real^3         ~&&: \text{angular acceleration of the centre of mass}, \\
  & \vvb(t)\in\Real^3         ~&&: \text{velocity of the centre of mass},\\ 
  & \omegavb(t)\in\Real^3     ~&&: \text{angular velocity},\\ 
  & \xvb(t)\in\Real^3         ~&&: \text{position of the centre of mass},\\ 
  & \Eb(t)\in\Real^{3\times 3} ~&&: \text{matrix with columns being the principle axes of inertia}.
\end{alignat*}
The motion of each rigid body is determined by a surface integral of the fluid traction~$\sigmav\nv$ along the interface~$\GammaB(t)=\bar\OmegaF(t)\cap {\bar\Omega}_b(t)$, where $\sigmav$ is the fluid stress and $\nv$ is the unit outward normal at the interface.
In three dimensions, the rigid body has six degrees of freedom and
its motion satisfies the Newton-Euler equations,
\begin{align}
    \mrb \avb & = \int_{\GammaB} \sigmav \nv \, \dArea + \fvbe  ,   \label{eq:linearAccelerationEquation} \\
    \Ib \bvb &= - \omegavb\times \Ib \omegavb + \int_{\GammaB} (\rv-\xvb)\times \sigmav \nv \, \dArea + \gvbe,
                 \label{eq:angularAccelerationEquation}\\
   \dotvvcm &= \avb , \label{eq:centerOfMassVelocityEquation} \\
    \dot\omegav_b &= \bvb , \label{eq:angularVelocityEquation} \\
   \dotxvcm &= \vvcm , \label{eq:centerOfMassPositionEquation} \\
   \dot \Eb &= \omegavb\times \Eb ,  \label{eq:axesOfInertiaEquation}
\end{align}
where $\rv$ denotes a point on the surface of the body, and $\fvbe$ and $\gvbe$ are an applied external force and torque on the body, respectively.

The velocity of the fluid is coupled to that of each rigid body on its surface $\GammaB$ by the matching condition
\begin{align}
 &  \vv(\rv(t),t)  = \dot{\rv}(t), \qquad \rv(t)\in\GammaB. \label{eq:RBsurfaceVelocity}
\end{align}
%For a rigid body that has six degrees of freedom, 
The motion of a point $\rv(t)$ on the surface of the body is given by a translation together with a rotation about the
initial centre of mass,
\begin{align}
    \rv(t) &= \xvcm(t) + R(t) (\rv(0)-\xvcm(0)) , \label{eq:bodySurface}
\end{align}
where $R(t)$ is the rotation matrix given by
\begin{align}
   R(t) &= \Eb(t) \Eb^{T}(0).    \label{eq:rotationMatrix}
\end{align}
The velocity $\dot{\rv}(t)$ in~\eqref{eq:RBsurfaceVelocity} is determined by differentiating~\eqref{eq:bodySurface}, and using~\eqref{eq:centerOfMassPositionEquation}, \eqref{eq:axesOfInertiaEquation} and~\eqref{eq:rotationMatrix} to find  \begin{alignat}{3}
 &  \dot{\rv}(t) = \vvb(t) +\omegavb(t)\times(\rv(t)-\xvb(t)), \qquad&& \rv(t)\in\GammaB. \label{eq:BodyPointVelocity}
\end{alignat}
We also use the acceleration of a point on the surface of the body for the~\ampRB~algorithm, and this is given by
\begin{alignat}{3}
 & \ddot{\rv}(t) = \avb(t) + \bvb(t)\times(\rv(t)-\xvb(t))  + \omegavb(t)\times\big[ \omegavb(t)\times(\rv(t)-\xvcm(t))\big], \qquad&& \rv(t)\in\GammaB.\label{eq:BodyPointAcceleration} 
\end{alignat}

% -------------------------------------------------------------------------------
\section{The \ampRB~algorithm and stability analysis} \label{sec:algorithm}

In this section, we describe the interface conditions used in the~\ampRB~scheme to handle added-mass and added-damping effects arising from the fluid stress on the rigid body.  These AMP interface conditions are a key ingredient in the~\ampRB~time-stepping scheme, which is also discussed in this section.  Lastly, we perform a stability analysis of the time-stepping scheme applied to two FSI model problems in three dimensions.

\newcommand{\Fvt}{\Fc_\mu}% viscous traction force
\newcommand{\Gvt}{\Gc_\mu}% viscous traction force

% ----------------------------------------------------------------------------------------------------
\subsection{The added-mass and added-damping interface conditions}
\label{sec:AMPDPinterfaceConditions}

To handle added-mass effects, the fluid and body accelerations on the interface are matched.
Taking a time derivative of the matching conditions in~\eqref{eq:RBsurfaceVelocity} and using the acceleration of a point on the body given in~\eqref{eq:BodyPointAcceleration}, we obtain
\[
\partial_t\vv+(\dot{\rv}\cdot\grad)\vv=  \avb + \bvb\times(\rv-\xvb)  + \omegavb\times\big[ \omegavb\times(\rv-\xvcm)\big], \qquad \rv\in\GammaB,
\]
and using the momentum equation~\eqref{eq:fluidMomentum} gives
\begin{align}
-\frac{1}{\rho}\Big[ \grad p -\mu\Delta\vv \Big] 
      =  \avb + \bvb\times(\rv-\xvb)  + \omegavb\times\big[ \omegavb\times(\rv-\xvcm)\big], \qquad \rv\in\GammaB. \label{eq:vectorComp}
  \end{align}
The normal component of~\eqref{eq:vectorComp} satisfies 
\begin{align}
    \partial_n p
              & = -\rho\nv^T \Big( \avb + \bvb\times(\rv-\xvb)  + \omegavb\times\big[ \omegavb\times(\rv-\xvcm)\big]  \Big)  
                + \mu \nv^T\Delta \vv ,  \qquad \rv\in\GammaB, \label{eq:pressureBC} 
\end{align}
where $\nv=\nv(\rv,t)$ is the outward normal to the body. 
When the pressure and body accelerations are coupled strongly in the~\ampRB~scheme, the condition in~\eqref{eq:pressureBC} is equivalent to a generalized Robin boundary condition for the pressure.

%added-damping
Added-damping effects are treated in the~\ampRB~scheme by considering the contributions to the force and torque on the body in~\eqref{eq:linearAccelerationEquation} and~\eqref{eq:angularAccelerationEquation}, respectively, due to the integrals,
\begin{equation*}
   \Fvt(\vv,\vvb,\omegavb) = \int_{\GammaB} \tauv\nv\,\dArea, \qquad 
   \Gvt(\vv,\vvb,\omegavb) = \int_{\GammaB}  (\rv-\xvb)\times (\tauv\nv)\,\dArea,
%\label{eq:surfaceIntegrals}
\end{equation*}
involving the viscous shear stress $\tauv\nv$ of the fluid.  The integrals depend on the fluid velocity $\vv$, and on the velocity $\vvb$ and angular velocity $\omegavb$ of the body due to its motion.  The linearized forms
\begin{align}
   \Fvt(\vv,\vvb,\omegavb) \approx \Fvt(\vv^p,\vvb^p,\omegavb^p) - \Dvv\, (\vvb-\vvb^p) - \Dvw\,(\omegavb-\omegavb^p) ,  \label{eq:addedDampingFdws} \\
   \Gvt(\vv,\vvb,\omegavb) \approx \Gvt(\vv^p,\vvb^p,\omegavb^p) - \Dwv\,(\vvb-\vvb^p) - \Dww\,(\omegavb-\omegavb^p), \label{eq:addedDampingGdws}
\end{align}
about the predicted states $\vv^p$, $\vvb^p$ and $\omegavb^p$ are used to expose the dependence of the integrals on the motion of the body.
  In~\cite{rbins2017} it was postulated that the added-damping tensors can be taken from the approximate solution to
  a discrete potential problem. The approximate forms, written here using surface integrals over the rigid-body (in~\cite{rbins2017} these
  were written as quadratures), are
\begin{alignat}{3} 
  & \Dvv \eqdef \mu \int_{\GammaB} \frac{1}{\dn} \, (\Iv - \nv \nv^T) \,d S , \label{eq:ADvv} \\
  & \Dvw \eqdef \mu \int_{\GammaB} \frac{1}{\dn} \, (\Iv - \nv \nv^T)  [\rv-\xv_b]_\times^T \, d S , \label{eq:ADvw} \\
  & \Dwv \eqdef \mu \int_{\GammaB} \frac{1}{\dn} \,  [\rv-\xv_b]_\times  (\Iv - \nv \nv^T)   \, d S , \label{eq:ADwv} \\
  & \Dww \eqdef \mu \int_{\GammaB} \frac{1}{\dn} \, [\rv-\xv_b]_\times (\Iv - \nv \nv^T)  [\rv-\xv]_\times^T \, d S  ,\label{eq:ADww}
\end{alignat}
where $[\rv-\xvb]_\times$ denotes the matrix representation of the cross product involving the vector $\rv-\xvb$, 
and $\dn=\dn(\rv)$ is the {\em added-damping length-scale} parameter given by
\begin{equation}
\dn \eqdef \frac{\ds_n(\rv)}{1-e\sp{-\delta}},\qquad \delta \eqdef \frac{\ds_n(\rv)}{\sqrt{\alpha\nu\dt}}, \qquad \rv\in\GammaB.
\label{eq:dn}
\end{equation}
Here, $\ds_n(\rv)$ is the mesh spacing in the normal direction to the surface, $\alpha=1/2$ is the coefficient in the implicit
time-stepping scheme, and $\nu=\mu/\rho$ is the kinematic viscosity.
Note that $\dn$ varies with $\delta$, which is the ratio of the mesh spacing in the normal direction to a {\em viscous} grid spacing, and it takes the limiting values
\[
  \dn \sim \begin{cases}
                 \sqrt{\alpha \nu \dt} & \text{for $\delta\rightarrow 0$} , \\
                 \ds_n         & \text{for $\delta\rightarrow \infty$} .
           \end{cases}
\] 
In two-dimensions it was found to be very important to use the correct form for $\dn$ to ensure stability. Here
we will show that the same form is valid in three dimensions.

For convenience, a larger composite tensor, $\ADtensor(t) \in\Real^{6\times 6}$, is defined as 
\begin{align}
  \ADtensor \eqdef   \begin{bmatrix}
    \Dvv & \Dvw  \\
    \Dwv & \Dww
  \end{bmatrix} ,  \label{eq:AddedDampingTensor} 
\end{align}
and we note that the component tensors in $\ADtensor(t)$ can be computed based on their initial states using the transformation
\begin{align}
   {\cal D}^{\alpha\beta}(t) = R(t)\, {\cal D}^{\alpha\beta}(0) \, R^T(t) ,  \label{eq:addedDampingInTime}
\end{align}
where $R(t)$ is the rotation matrix from~\eqref{eq:rotationMatrix}.

Note that the linearizations~\eqref{eq:addedDampingFdws} and~\eqref{eq:addedDampingGdws} can be expressed in terms of the body accelerations $\avb$ and $\bvb$ as follows:
\begin{align}
   \Fvt(\vv,\vvb,\omegavb) \approx \Fvt(\vv^p,\vvb^p,\omegavb^p) - \dt \, \Dvv\, (\avb-\avb^p) 
                 - \dt \,\Dvw\,(\bvb-\bvb^p) ,  \label{eq:addedDampingFdot} \\
   \Gvt(\vv,\vvb,\omegavb) \approx \Gvt(\vv^p,\vvb^p,\omegavb^p) - \dt \,\Dwv\,(\avb-\avb^p) 
                 - \dt \,\Dww\,(\bvb-\bvb^p). \label{eq:addedDampingGdot}
\end{align}
where $\dt$ is the discrete time-step in the~\ampRB~time-stepping scheme.

The pressure boundary condition in~\eqref{eq:pressureBC} and the linearizations in~\eqref{eq:addedDampingFdot} and~\eqref{eq:addedDampingGdot} are used in the following primary AMP interface condition:
% ---------------------------------------------------------------------------------------------------------------
\begin{AMPInterfaceCondition} The AMP interface conditions on the surface of a rigid body $\rv\in\GammaB$ for the pressure equation~\eqref{eq:fluidPressure} are 
\begin{align}
    \partial_n p + \rho\nv^T \Big( \avb + \bvb\times(\rv-\xvb)   \Big)  
               &= -\rho\nv^T \Big( \omegavb\times\big[ \omegavb\times(\rv-\xvcm)\big]  \Big)  
                + \mu \nv^T\Delta \vv , \label{eq:AMPDPpressure} \\
 \left\{ 
\begin{bmatrix} \mrb\, I_{3\times 3} & 0 \smallskip\\  0 & \Ib  \end{bmatrix}
+ \dt\, \ADtensor 
%  \begin{bmatrix}
%     \Dvv & \Dvw  \\
%     \Dwv & \Dww
%   \end{bmatrix}
\right\} 
\begin{bmatrix} \avb \smallskip\\ \bvb  \end{bmatrix}
+
\Fv(p) 
& = 
-  \begin{bmatrix} 0 \smallskip\\ \omegavb\times \Ib \omegavb\end{bmatrix} + 
         \Gv(\vv) + \dt\,\ADtensor \begin{bmatrix} \avb \smallskip\\ \bvb  \end{bmatrix},  \label{eq:AMPDPrigidBodyAcceleration}
\end{align}
where
\begin{equation}
\Fv(p) \eqdef \begin{bmatrix}
                 \displaystyle \int_{\partialB} p\, \nv \, \dArea \\ 
                 \displaystyle \int_{\partialB} (\rv-\xvb)\times (p\, \nv ) \, \dArea 
             \end{bmatrix} ,
\qquad
\Gv(\vv) \eqdef 
\begin{bmatrix}
                \displaystyle \Fvt \, + \fvbe \\
                \displaystyle \Gvt \, + \gvbe
              \end{bmatrix} \,=\,  
\begin{bmatrix}
                \displaystyle \int_{\partialB} \tauv\nv \, \dArea + \fvbe \\
                \displaystyle \int_{\partialB} (\rv-\xvb)\times (\tauv\nv) \, \dArea + \gvbe
              \end{bmatrix}.  
       \label{eq:Pdef}
\end{equation}
Here $\dt$ is the time-step used in the~\ampRB~time-stepping scheme.
\end{AMPInterfaceCondition}
% -----------------------------------------------------------------------------------------------------------------------
The AMP interface conditions given here differ slightly from the ones derived in~\cite{rbins2017}.  First, the signs on the integrals of the two vector components of $\Fv(p)$ in~\eqref{eq:Pdef} were incorrect in our previous paper and have been fixed here.  Second, an added-damping parameter appeared in the earlier AMP interface conditions for the purpose of studying the stability of the~\ampRB~scheme.  As mentioned previously, the results of the stability analysis discussed in~\cite{rbinsmp2017} suggest that this parameter can be taken to be one.  This value has been used in the form of the interface conditions given here and so the added-damping parameter does not appear. {Note also that one component of the AMP interface conditions is a nonlocal Robin boundary condition on the pressure, which is the same condition at the core of the method discussed recently in~\cite{yoon2018}. }

\subsection{Time-stepping algorithm}
Algorithm~\ref{alg:ampRB} provides the details of the~\ampRB~time-stepping scheme used for the calculations presented in this paper.
The algorithm given here follows the ones discussed in our previous work~\cite{rbinsmp2017, rbins2017}, but extended to three dimensions.
It is a predictor-corrector-type fractional-step scheme that advances the solution from the time step $t^n$ to $t^{n+1}$.
The fluid is advanced with an IMEX scheme in time that combines an Adams-Bashforth scheme with the trapezoidal rule,
 while the solid equations use a Leapfrog predictor together with a trapezoidal-rule corrector.
The algorithm starts with a preliminary body evolution stage consisting of Steps~1 and~2, 
in which the rigid-body variables and the moving composite grid at the new time step are predicted 
from the analogous quantities at previous time steps.
 The moving grid $\Gc$ in Step~2 is determined by the predicted values for $\xvb$ and $\Eb$ using equations~\eqref{eq:bodySurface} and~\eqref{eq:rotationMatrix}, see Section~\ref{sec:comositeGrids} for more details.  
The preliminary stage is followed by the predictor stage consisting of Steps~3--6.
In Step 3, the fluid velocity, $\vv_{\iv}\sp\star$, is computed implicitly using the rigid-body velocities predicted in Step 1 
for the boundary condition on the surface of the body.
Note that the advection and pressure terms in the momentum equation~\eqref{eq:fluidMomentum} 
are written as a single term,
\begin{align*}
  \fap_{\iv}^n &\eqdef \big((\vv_{\iv}^n-\wv_{\iv}^n\,)\cdot\grad_h\big) \vv_{\iv}^n + \grad_h  p_\iv^{n}, %\label{eq:advPressure}
\end{align*}
in Step~3, where the {\em grid velocity}, $\wv_{\iv}^n$, is included as required to express the equations in a moving coordinate system.  The ghost points for the velocity calculation at this step are determined by the
divergence of the velocity applied on the boundary together with an extrapolation of the velocity in the tangential directions $\tv_1$ and $\tv_2$.
The pressure, $p^{\star}_\iv$, and rigid-body body accelerations, $\avb^{\star}$ and $\bvb^{\star}$, are then determined in Step~4, by solving a discrete Poisson equation with the AMP interface conditions.
Steps~5 and~6 update the variables for the rigid body and the composite grid at $t_{n+1}$, and these steps are similar to the ones in the preliminary stage.
The corrector stage, consisting of Steps~7--10, follows the steps in the predictor stage, and is usually applied once 
in order to increase the stability region of the scheme to allow a larger advection time-step.
Finally, we have included an optional velocity correction in Step~11 as was done in~\cite{rbinsmp2017, rbins2017},
which has been shown to significantly improve the numerical stability to overcome added-damping effects.

In Algorithm~\ref{alg:ampRB}, the symbols $\grad_h$ and $\Delta_h$ denote discrete approximations to the gradient and Laplacian operators, respectively, and these operators are implemented using standard second-order finite differences.
In the implementation of the AMP interface conditions in Steps~4 and~8, the term $\Delta\vv$ in the boundary condition~\eqref{eq:AMPDPpressure} has been replaced 
by the equivalent form $-\grad\times\grad\times\vv$ to improve the time-step restriction. 
(see~\cite{Petersson00} for more details of its effects).
The interface conditions also use discrete approximations of the added-damping tensors as discussed in more detail in Section~\ref{sec:surfaceIntegral}.
Finally, the time step for the~\ampRB~scheme is then given by
\[
    \dt = \lambda_\text{CFL} \min_g \dt_g,
\]
where $\lambda_\text{CFL}$ is a CFL number, which is taken to be $0.9$ in the current work unless otherwise noted, $\dt_g$ is the time step associated with a component grid $g$, and the minimum is taken over all grids in the overlapping grid (see Section~\ref{sec:comositeGrids}).  A stable value for $\dt_g$ is estimated by
\[
    \dt_g = \bigg(\max_\iv \sum_{j=1}^{n_d} \frac{ |v_{j,\iv}|}{\dx_{j,\iv}}\bigg)^{-1},
\]
where $\dx_{j,\iv}$ is the grid spacing at grid-point $\iv$ in the $j\sp{{\rm th}}$
coordinate direction and $v_{j,\iv}$ is the component of the fluid velocity in the coordinate direction $j$.

% ================================ ALGORITHM ============================================
{
\def\alignspace{\hspace{1.2em}}
\begin{algorithm}\caption{\rm Added-mass partitioned (\ampRB) scheme}\scriptsize
\[
\begin{array}{l}
\hbox{// \textsl{\bodyStepIComment}}\smallskip\\
1.\text{ Predict the rigid-body degrees of freedom:}\\ 
\alignspace \avb^{(e)} = 2 \avb^n - \avb^{n-1} , \quad \bvb^{(e)} = 2 \bvb^n - \bvb^{n-1},  \quad
    \vvb^{(e)} = \vvb^{n-1} + 2\dt\, \avb^{n} , \\
\alignspace \omegavb^{(e)} = \omegavb^{n-1} + 2\dt\, \bvb^{n} ,
\quad \xvb^{(e)} = \xvb^{n-1} + 2\dt\, \vvb^{n} , \quad
\Eb^{(e)} = \Eb^{n-1} + 2\dt\, \omegavb^{n}\times\Eb^{n} .\\
2.\text{ Predict the moving grid }\Gc^{(e)} \text{ by $\xvb^{(e)}$ and $\Eb^{(e)}$.}
\medskip\\
\hbox{// \textsl{Prediction steps}}\smallskip\\
3. \text{ Advance the fluid velocity $\vv^{\star}_{\iv}$:}\\
\alignspace
\begin{cases}
    \rho\Big({\vv^{\star}_{\iv} -\vv^n_{\iv}}\Big)/{\dt} + \Big( 3\fap_\iv^{n}-\fap_\iv^{n-1} \Big)/2 
= {\mu}\Big( \Delta_h \vv^{\star}_{\iv} + \Delta_h \vv_\iv^n \Big)/2  , & \quad \iv\in\OmegaF_h, \\
\vv^{\star}_\iv = \vvb^{(e)} + \omegavb^{(e)}\times( \rv_\iv^{(e)}-\xvb^{(e)}),\quad \grad_h\cdot\vv^{\star}_\iv=0, \quad \text{Extrapolate ghost: $\tv_m^T\vv^{\star}_\iv$}, & \quad \iv\in\Gamma_h,\\
\text{Velocity boundary conditions on $\partial \Omega_h\backslash \Gamma_h$}.
\end{cases}
  \smallskip\\

  4.\text{ Update the pressure $p^{\star}_\iv$ and body accelerations $\avb^{\star}$, $\bvb^{\star}$:}\\
  \alignspace
\begin{cases}
    \Delta_h p^{\star}_\iv = -\grad_h\vv^{\star}_\iv : {(\grad_h\vv^{\star}_\iv)^T} , & \iv \in \Omega_h, \\
   \nv_\iv^T\grad_h p^{\star}_\iv + \rho{\nv_\iv^T} \Big( \avb^{\star} + \bvb^{\star}\times(\rv_\iv^{(e)}-\xvb^{(e)} )   \Big)  
   = -\rho\nv_\iv^T \Big( \omegavb^{(e)} \times\big[ \omegavb^{(e)}\times(\rv^{(e)}_\iv-\xvcm^{(e)})\big]  \Big)   - \mu \nv_\iv^T\, (\grad_h\times\grad_h\times \vv_\iv^{\star}) , & \iv\in\Gamma_h,\\
 \left\{ 
\begin{bmatrix} \mrb\, I_{3\times 3} & 0 \\  0 & \Ib  \end{bmatrix}
+ \dt\,\ADtensor^{(e)}
 \right\} 
\begin{bmatrix} \avb^{\star} \\ \bvb^{\star}  \end{bmatrix}
+
\Fv(p^{\star}_\iv) 
 = 
-  \begin{bmatrix} 0 \\ \omegavb^{(e)}\times \Ib \omegavb^{(e)}\end{bmatrix} + 
         \Gv(\vv_\iv^{\star}) + \dt\,\ADtensor^{(e)} \begin{bmatrix} \avb^{(e)} \\ \bvb^{(e)}  \end{bmatrix}, &\\
 
\text{Pressure boundary conditions on $\partial \Omega_h\backslash \Gamma_h$}.
\end{cases}
\smallskip\\

5.\text{ Update the rigid-body degrees of freedom:}\\
 \alignspace  \vvb^{\star} = \vvb^{n} + \dt\, (\avb^{\star}+ \avb^{n})/2 , \quad
\omegavb^{\star} = \omegavb^{n} + \dt\,( \bvb^{\star} + \bvb^{n})/2 , \\
  \alignspace
    \xvb^{\star} = \xvb^{n} + \dt\, (\vvb^{\star}+\vvb^{n})/2 , \quad
    \Eb^{\star} = \Eb^{n} + \dt\,( \omegavb^{\star}\times\Eb^{(e)} +\omegavb^{n}\times\Eb^{n})/2 .\\
6.\text{ Update the moving grid }\Gc^{\star} \text{ by $\xvb^{\star}$ and $\Eb^{\star}$.}
  \medskip\\
\hbox{// \textsl{Correction steps}}\smallskip\\
7. \text{ Advance the fluid velocity $\vv_\iv^{n+1}$:}\\
  \alignspace
\begin{cases}
    \rho\Big({\vv^{n+1}_{\iv} -\vv^n_{\iv}}\Big)/{\dt} + \Big( \fap_\iv^{\star}   +  \fap_{\iv}^n \Big)/2 
= {\mu}\Big( \Delta_h \vv^{n+1}_{\iv} + \Delta_h \vv_\iv^n \Big)/2  , & \quad \iv\in\OmegaF_h, \\
     \vv^{n+1}_\iv = \vvb^{\star} + \omegavb^{\star}\times( \rv_\iv^{\star}-\xvb^{\star}),\quad \grad_h\cdot\vv^{n+1}_\iv=0, \quad \text{Extrapolate ghost: $\tv_m^T\vv^{n+1}_\iv$}, & \quad \iv\in\Gamma_h, \\
\text{Velocity boundary conditions on $\partial \Omega_h\backslash \Gamma_h$}.
\end{cases}
  \smallskip\\

8.\text{ Update the pressure $p^{n+1}_\iv$ and body accelerations $\avb^{n+1}$, $\bvb^{n+1}$:}\\
  \alignspace
\begin{cases}
    \Delta_h p^{n+1}_\iv = -\grad_h\vv^{n+1}_\iv : {(\grad_h\vv^{n+1}_\iv)^T} , & \iv \in \Omega_h, \\
   \nv_\iv^T\grad_h p^{n+1}_\iv + \rho{\nv_\iv^T} \Big( \avb^{n+1} + \bvb^{n+1}\times(\rv_\iv^{\star}-\xvb^{\star} )   \Big)  
   = -\rho\nv_\iv^T \Big( \omegavb^{\star} \times\big[ \omegavb^{\star}\times(\rv^{\star}_\iv-\xvcm^{\star})\big]  \Big)   - \mu \nv_\iv^T\, (\grad_h\times\grad_h\times \vv_\iv^{n+1}) , & \iv\in\Gamma_h,\\
 \left\{ 
\begin{bmatrix} \mrb\, I_{3\times 3} & 0 \\  0 & \Ib  \end{bmatrix}
+ \dt\,\ADtensor^{\star}
 \right\} 
\begin{bmatrix} \avb^{n+1} \\ \bvb^{n+1}  \end{bmatrix}
+
\Fv(p^{n+1}_\iv) 
 = 
-  \begin{bmatrix} 0 \\ \omegavb^{\star}\times \Ib \omegavb^{\star}\end{bmatrix} + 
         \Gv(\vv_\iv^{n+1}) + \dt\,\ADtensor^{\star} \begin{bmatrix} \avb^{\star} \\ \bvb^{\star}  \end{bmatrix}, &\\
 
\text{Pressure boundary conditions on $\partial \Omega_h\backslash \Gamma_h$}.
\end{cases}
\smallskip\\

9.\text{ Update the rigid-body degrees of freedom:}\\
\alignspace
\vvb^{(n+1)} = \vvb^{n} + \dt\, (\avb^{(n+1)}+ \avb^{n})/2 , \quad
\omegavb^{(n+1)} = \omegavb^{n} + \dt\,( \bvb^{(n+1)} + \bvb^{n})/2 , \\
  \alignspace
    \xvb^{(n+1)} = \xvb^{n} + \dt\, (\vvb^{(n+1)}+\vvb^{n})/2 , \quad
    \Eb^{(n+1)} = \Eb^{n} + \dt\,( \omegavb^{(n+1)}\times\Eb^{\star} +\omegavb^{n}\times\Eb^{n})/2 .\\
10.\text{ Update the moving grid }\Gc^{(n+1)} \text{ by $\xvb^{(n+1)}$ and $\Eb^{(n+1)}$.} 
  \medskip\\

\hbox{// \textsl{Fluid-velocity correction step (optional)}}\smallskip\\
11. \text{ Advance the fluid velocity $\vv_\iv^{n+1}$:}\\
  \alignspace
\begin{cases}
    \rho\Big({\vv^{n+1}_{\iv} -\vv^n_{\iv}}\Big)/{\dt} + \Big( \fap_\iv^{n+1}   +  \fap_{\iv}^n \Big)/2 
= {\mu}\Big( \Delta_h \vv^{n+1}_{\iv} + \Delta_h \vv_\iv^n \Big)/2  , & \qquad \iv\in\OmegaF_h, \\
     \vv^{n+1}_\iv = \vvb^{n+1} + \omegavb^{n+1}\times( \rv_\iv^{n+1}-\xvb^{n+1}),\quad \grad_h\cdot\vv^{n+1}_\iv=0, \quad \text{Extrapolate ghost: $\tv_m^T\vv^{n+1}_\iv$}, & \qquad \iv\in\Gamma_h,\\
\text{Velocity boundary conditions on $\partial \Omega_h\backslash \Gamma_h$}.
\end{cases}
\end{array}
\]
\label{alg:ampRB}
\end{algorithm}
}

\newcommand{\ama}{{\cal M}_a}
\def\linop{{\mathcal L}}

\subsection{Stability analysis of two FSI model problems in a spherical geometry}

We consider the stability of the time-stepping algorithm for two FSI model problems in a spherical geometry.
For the model problems, the fluid domain is a spherical shell 
between an inner radius $r_1$ and an outer sphere of radius $r_2$.
The fluid surrounds a spherical rigid body of radius $r_1$ with mass $m_b$ and moment of inertia $I_b$.
The interface between the body and fluid is denoted as $\partialB$ with a unit normal vector of $\nv$.
Following our previous work~\cite{rbinsmp2017}, we consider small displacements of the sphere 
and linearize the motion of the body and fluid.
Therefore, it is sufficient to consider Stokes flow in the spherical shell for the purpose of the analysis.

The Stokes flow is described using spherical coordinates
$(r, \theta, \phi)$, where $r$ is the radial distance of range $(r_1, r_2)$, 
$\theta$ is the polar angle from the positive $z$-axis of range $[0, \pi]$, and 
$\phi$ is the azimuthal angle of range $[0, 2\pi)$.
The velocity of the flow is then described by the components $(u,v,w)$ in the directions of the unit vectors
$(\hat{\bf r},\hat{\pmb{\theta}},\hat{\pmb\phi})$.
As for the analysis of the FSI problems in annular geometry in~\cite{rbinsmp2017},
it is sufficient to consider translational motion of the rigid sphere parallel to the $z$-direction 
and a rotation of the sphere about the $z$-axis.
The governing equations corresponding to these two motions can be decoupled 
and are considered in the following two model problems.

\subsubsection{Added-mass model problem}
We first consider the model problem corresponding to the rigid sphere translating parallel to the $z$-direction.
The motion of the sphere is dominated by added-mass effects, and therefore it is sufficient to drop
the viscous terms in the fluid.
The flow is axisymmetric, and thus there is no $\phi$-dependence, and $w$ can be set to zero without loss. 
The governing equations for the FSI model problem are
\begin{equation*}
\left\{ 
  \begin{alignedat}{3}
      & \rho \frac{\partial u}{\partial t} +\frac{\partial p}{\partial r} = 0 , \quad&&  r\in(r_1,r_2),\quad \theta\in[0,\pi],   \\
      & \rho \frac{\partial v}{\partial t} +\frac 1 r\frac{\partial p}{\partial \theta} = 0 , \quad&&  r\in(r_1,r_2), \quad\theta\in[0,\pi], \\
      &  \grad \cdot \uv = 0   , \quad&&  r\in(r_1,r_2), \quad\theta\in[0,\pi], \\
 &   m_b \frac{d w_b}{dt} = \int_\partialB \ev^z \cdot[(\rv-\xvb)\times (-p \nv)] \, dS + f_e(t), \\
 &  u=w_b(t)\cos \theta, && r=r_1, \\
 &  u = 0, && r=r_2,
  \end{alignedat}  \right. 
\end{equation*}
where $w_b$ is the velocity of the rigid sphere in the $z$-direction with unit vector $\ev^z$.
Note that the boundary condition of the fluid at $r=r_1$ is only related to 
the component of the solid velocity along $\hat\rv$-direction, and there is a homogeneous
boundary condition on the fluid at $r=r_2$ assuming a zero-normal-flow condition there.
We consider solutions of the form
\[
    u = \uHat(r,t)\cos\theta, \qquad v = - \vHat(r,t) \sin\theta, \qquad p=\pHat(r,t) \cos\theta,
\]
which gives the added-mass model problem (MP-AM) of the form
\begin{equation}
\text{{MP-AM}:}\;  
\left\{ 
  \begin{alignedat}{3}
 & \rho \frac{\partial \uHat}{\partial t} +\frac{\partial \pHat}{\partial r} = 0 , \quad&&  r\in(r_1,r_2),   \\
 & \rho \frac{\partial \vHat}{\partial t} +\frac {\pHat} r = 0 , \quad&&  r\in(r_1,r_2), \\
 &  \frac 1 {r^2} \frac{\partial}{\partial r} (r^2 \uHat) - \frac 2 r \vHat  = 0   , \quad&&  r\in(r_1,r_2),\\
 &   m_b \frac{d w_b}{dt} = -\frac{4\pi  r_1^2}3   \pHat(r_1,t) + f_e(t), \\
 &  \uHat= w_b(t), && r=r_1, \\
 &  \uHat = 0, && r=r_2.
  \end{alignedat}  \right. 
\end{equation}
The exact solution is given by
\begin{align}
w_b(t)&={1\over m_b+\ama}\int_0\sp{t}f_e(\tau)\,d\tau+w_b(0), \label{eq:ubsolna} \\
\pHat(r,t)&={\rho\aa a_w(t)\over (\bb/\aa)^2-\aa/\bb}\left({r\over\bb}+{\bb^2\over2 r^2}\right), \label{eq:pHatsolna} \\
\uHat(r,t)&={\bb\sp3-r\sp3\over\bb\sp3-\aa\sp3}\left({\aa\over r}\right)\sp3w_b(t), \label{eq:uHatsolna} \\
\vHat(r,t)&={2r\sp3+\bb\sp3\over2\aa\sp3-2\bb\sp3}\left({\aa\over r}\right)\sp3w_b(t), \label{eq:vHatsolna}
\end{align}
where
\begin{equation}
\ama=\frac{4}{3} \rho\pi\aa\sp3\left[{1+2\bigl(\aa/\bb\bigr)\sp3\over2-2\bigl(\aa/\bb\bigr)\sp3}\right],
\label{eq:addedmassannular}
\end{equation}
is the added-mass for this model problem, and the acceleration of the rigid body is given by
\begin{equation}
    a_w(t)={d{w_b}\over dt}={f_e(t)\over m_b+\ama}.
\label{eq:ausolna}
\end{equation}
Note that in the limit, $r_2/r_1\rightarrow\infty$, the added-mass in~\eqref{eq:addedmassannular} becomes
\begin{equation*}
\ama= \frac{2}{3} \rho\pi\aa\sp3 = \frac1 2 \rho V_b,
\end{equation*}
which is the classical added-mass of a sphere immersed in a fluid of infinite extent.
This added-mass formula can be found in many fluid dynamics textbooks, see~\cite{YihFluidMechanics} for example.

To study the stability of the {\ampRB}~algorithm for the MP-AM model problem,
we consider a semi-discrete (in time) version of Algorithm~\ref{alg:ampRB}. 
Applying Step~4 in Algorithm~\ref{alg:ampRB} to the MP-AM model problem leads to a system of the form
\begin{equation*}
\left\{ 
  \begin{alignedat}{2}
    & \partial_r (r^2\partial_r \pHat^\star) - 2 \pHat^\star = 0,  & \qquad  r\in (r_1, r_2), \\
    & \partial_r \pHat^\star (r_1) + \rho a_w^\star  = 0,\\
    & m_b a_w^\star + \frac{4\pi r_1^2}3 \pHat(r_1)  = f_e(t^{n+1}), \\
    & \partial_r \pHat^\star(r_2) = 0.
  \end{alignedat} \right.
\end{equation*} 
Provided such a system is well-posed (i.e.~$m_b+\ama$ is bounded away from zero), 
its solution is 
\[
\pHat\sp\star(r)=\pHat(r,t\sp{n+1}),\qquad a_w\sp\star=a_w(t\sp{n+1}),
\]
where $\pHat(r,t)$ and $a_w(t)$ are the exact fluid pressure and acceleration of the body, respectively.
The updates of the velocity of the body in Step~5 and the components of the fluid velocity in Step~7 use trapezoidal-rule quadratures
of the exact body acceleration and fluid pressure,
 which implies the prediction step is unconditionally stable. 
Similarly, it can be shown that the correction step is unconditionally stable.
Therefore,  the stability analysis based on the semi-discrete version of the~\ampRB~algorithm
gives the following result:

\begin{theorem} The \ampRB~algorithm given in Algorithm~\ref{alg:ampRB} for the MP-AM model problem
is unconditionally stable provided there exists a constant $K>0$ such that 
\begin{align*}
    m_b+\ama \ge K.
\end{align*}
\end{theorem}

Note that the stability analysis and its result are analogous to that discussed in~\cite{rbinsmp2017} for the MP-AMA model problem involving a two-dimensional annular geometry.

\subsubsection{Added-damping model problem}

We now consider a second model problem corresponding to the rigid sphere rotating about the $z$-axis. 
For this problem, the only nonzero component of the fluid velocity is the azimuthal component, and there are only added-damping effects.  As in the previous model problem, the flow is axisymmetric so that there is no $\phi$-dependence.
The rotational motion of the rigid sphere thus satisfies
\begin{align*}
    I_b \frac{d \omega_b}{dt} & = \int_\partialB  \ev^z \cdot[ (\rv-\xvb)\times\sigmav \nv] \, dS + g_e(t) \\
    & = \int_\partialB \sin \theta \, r_1 \tau_{r \phi} \, dS + g_e(t)  \\
    & = 2\pi r_1^3  \int_0^\pi \mu \left[ r \frac \partial {\partial r} \left( \frac {w} r \right) \right]_{r=r_1} \sin^2  \theta \, d\theta +g_e(t).
\end{align*}
Therefore, the equations governing the model problem are given by
\begin{equation*}
\left\{ 
  \begin{alignedat}{3}
      & \rho \frac{\partial w}{\partial t} = \mu \left[ \frac 1 {r^2} \frac{\partial}{\partial r} \left( r^2 \frac{\partial w}{\partial r} \right) + \frac{1}{r^2 \sin \theta}
\frac{\partial }{\partial \theta}  \left( \sin \theta \frac{\partial w}{\partial \theta}  \right)
- \frac {w}{r^2\sin^2 \theta} \right] , \quad&&  r\in(r_1,r_2), \quad \theta\in[0,\pi],  \\
 &   I_b \frac{d \omega_b}{dt} = 2\pi r_1^3  \int_0^\pi \mu \left[ r \frac \partial {\partial r} \left( \frac {w} r \right) \right]_{r=r_1} \sin^2  \theta \, d\theta+ g_e(t), \\
 &  w =r_1 w_b\sin \theta  , && r=r_1, \\
 &  w  = 0, && r=r_2.
  \end{alignedat}  \right. 
 \label{eq:MP-RF}
\end{equation*}
Assume the azimuthal component of the fluid velocity satisfies
\[
    w(r, \theta, t) = \wHat(r, t)\sin\theta .
\]
The $\theta$-dependence can be separated and the problem is simplified to an added-damping model  problem
(MP-AD) satisfying
\begin{equation*}
\text{{MP-AD}:}\;  
\left\{ 
  \begin{alignedat}{3}
      & \frac{\partial \wHat}{\partial t} = \nu \linop \wHat , \quad&&  r\in(r_1,r_2),\\
&   I_b \frac{d \omega_b}{dt} = \frac{8\pi r_1^3 }3  \mu \left[ r \frac \partial {\partial r} \left( \frac {\wHat} r \right) \right]_{r=r_1} + g_e(t), \\
 &  \wHat = r_1 \omega_b, && r=r_1, \\
 &  \wHat = 0, && r=r_2,
  \end{alignedat}  \right. 
 \label{eq:MP-RF}
\end{equation*}
where $\nu=\mu/\rho$ is the kinematic viscosity and the differential operator in the equation for $\wHat(r,t)$ is defined by
\[
    \linop \wHat = \frac 1 {r^2} \frac{\partial}{\partial r} \left( r^2 \frac{\partial \wHat}{\partial r} \right) 
      - \frac {2 \wHat}{r^2}.
\]

To study the stability of the {\ampRB}~algorithm for the MP-AD model problem,
we still consider a semi-discrete (in time) version of Algorithm~\ref{alg:ampRB},
which is rewritten as Algorithm~\ref{alg:ampRBad}.
The added-damping tensor in Algorithm~\ref{alg:ampRB} reduces to the scalar $\Dww$, which for a sphere has been given in~\cite{rbins2017} as
\begin{equation}
    \Dww = \mu \frac{8}3 \pi r_1^4 \frac{1-e^{-\delta}}{\dr}, \qquad \delta \eqdef \frac{\dr}{\sqrt{\nu\dt/2}}	.
\label{eq:addedDampingTerm}
\end{equation}
Also, the added-damping parameter $\beta_d$ is introduced
in Algorithm~\ref{alg:ampRBad} as an adjustable parameter for the purpose of the stability analysis (see~\cite{rbinsmp2017}).
The difference operator, $D_{rh}$, in Steps~3 and~6 of the algorithm for MP-AD is defined by
\[
D_{rh}f(r)\eqdef {-3f(r)+4f(r+\Delta r)-f(r+2\Delta r)\over 2\Delta r}.
\]
The goal here is to prove that the algorithm is stable for any body inertia, $I_b$, 
if and only if the added-damping term is involved in Steps~3 and~6.

% -------------------------------- ADDED-DAMPING -------------------------------------------------------------
\begin{algorithm}\caption{\ampRB~scheme for the MP-AD model problem.}\scriptsize
\[
\begin{array}{l}
\hbox{// \textsl{{\bodyStepIComment}}}\smallskip\\
1.\quad b_\omega^\esup = 2 b_\omega^n - b_\omega^{n-1} , \qquad \omega_b^\esup = \omega_b^{n-1} + 2\dt\, b_\omega^n,\medskip\\
\hbox{// \textsl{Prediction steps}}\smallskip\\
2.\quad \wHat\sp\psup = \wHat^n + \frac{\nu\dt}{2}\bigl(\linop\wHat\sp\psup+\linop\wHat\sp{n}\bigr), \quad r\in(\aa,\bb), \qquad \wHat^\psup(\aa) = \aa\omega_b^\esup, \quad \wHat\sp\psup(\bb)=0, \smallskip\\
3.\quad (I_b+\beta_d\dt\,\Dww) b_\omega^\psup=2V_b\mu \bigl[rD_{rh}(\wHat\sp\psup/r)\bigr]_{r=\aa}+\beta_d\dt\,\Dww b_\omega\sp\esup+g_e(t\sp{n+1}), \medskip\\
4.\quad \omega_b^\psup = \omega_b^{n} + \frac{\dt}{2}( b_\omega^\psup + b_\omega^n ),\medskip\\
\hbox{// \textsl{Correction steps}}\smallskip\\
5.\quad \wHat\sp{(c)} = \wHat^n + \frac{\nu\dt}{2}\bigl(\linop\wHat\sp{n+1}+\linop\wHat\sp{n}\bigr), \quad r\in(\aa,\bb), \qquad \wHat^{(c)}(\aa) = \aa\omega_b^\psup, \quad \wHat\sp{(c)}(\bb)=0, \smallskip\\
6.\quad (I_b+\beta_d\dt\,\Dww) b_\omega^{n+1}=2V_b\mu \bigl[rD_{rh}(\wHat\sp{(c)}/r)\bigr]_{r=\aa}+\beta_d\dt\,\Dww b_\omega\sp\psup+g_e(t\sp{n+1}), \medskip\\
7.\quad \omega_b^{n+1} = \omega_b^{n} + \frac{\dt}{2}( b_\omega^{n+1} + b_\omega^n ), \medskip\\
\hbox{// \textsl{Fluid-velocity correction step}}\smallskip\\
8.\quad \wHat\sp{n+1} = \wHat^n + \frac{\nu\dt}{2}\bigl(\linop\wHat\sp{n+1}+\linop\wHat\sp{n}\bigr), \quad r\in(\aa,\bb), \qquad \wHat^{n+1}(\aa) = \aa\omega_b^{n+1}, \quad \wHat\sp{n+1}(\bb)=0.
\end{array}
\]
\label{alg:ampRBad}
\end{algorithm}

To examine the stability of Algorithm~\ref{alg:ampRBad}, we consider the homogeneous MP-AD problem, i.e.~$g_e(t)=0$, and set $\wHat\sp{n}=A\sp{n}\bar w\sp{0}$, $\omega_b\sp{n}=A\sp{n}\omega_b\sp{0}$ and $b_\omega\sp{n}=A\sp{n}b_\omega\sp{0}$, where $A$ is an amplification factor. 
The intermediate values in the prediction and correction steps can be defined similarly, for instance, 
$\wHat\sp\psup=A\sp{n+1}\bar w\sp\psup$.
The solutions for $\bar w\sp\psup$, $\bar w\sp{(c)}$ and $\bar w\sp{0}$ of the boundary-value problems implied by Steps~2, 5 and~8, respectively, take the form
\begin{align}
\bar w\sp\psup(r) &=\aa\bigl(\bar\omega_b\sp\esup-\bar\omega_b\sp{0}\bigr)\phi_1(r)+\bar w\sp{0}(r),\label{eq:bdySolution1} \\
\bar w\sp{(c)}(r) &=\aa\bigl(\bar\omega_b\sp\psup-\bar\omega_b\sp{0}\bigr)\phi_1(r)+\bar w\sp{0}(r), \label{eq:bdySolution2} \\
\bar w\sp{0}(r) &=\aa\bar\omega_b\sp{0} \,\phi_2(r),
\label{eq:bdySolution3}
\end{align}
where
\[
\phi_m(r)={i_1(\zeta_mr)k_1(\zeta_m\bb)-i_1(\zeta_m\bb)k_1(\zeta_mr)\over i_1(\zeta_m\aa)k_1(\zeta_m\bb)-i_1(\zeta_m\bb)k_1(\zeta_m\aa)},\qquad\hbox{$m=1$ or 2,}
\]
and
\[
\zeta_1=\sqrt{{2\over\nu\dt}},\qquad \zeta_2=\zeta_1\sqrt{{A-1\over A+1}}.
\]
Here $i_1(z)$ and $k_1(z)$ are modified spherical Bessel functions of the first and second kind of order one,
satisfying
\[
i_1(z)=\frac{z\cosh z -\sinh z}{z^2},\qquad k_1(z)=\frac{e^{-z} (z +1) }{z^2}.
\]
The two equations in Step~1 can be collected in the form
\begin{equation}
\Bigg[\begin{array}{c}
\bar \omega_b\sp{(e)} \medskip\\
\bar b_\omega\sp{(e)}
\end{array}\Bigg]=
\Mv\Bigg[\begin{array}{c}
\bar \omega_b\sp{0} \medskip\\
\bar b_\omega\sp{0}
\end{array}\Bigg],
\label{eq:sysE}
\end{equation}
where the $2\times2$ matrix $\Mv$ involves the time-step $\dt$ and the amplification factor $A$.
Similarly, Steps~3 and~4, along with the solution for $\bar w\sp\psup$, lead to a system of the form
\begin{equation}
\Nv\Bigg[\begin{array}{c}
\bar \omega_b\sp{(p)} \medskip\\
\bar b_\omega\sp{(p)}
\end{array}\Bigg]=
\Ov\Bigg[\begin{array}{c}
\bar \omega_b\sp{(e)} \medskip\\
\bar b_\omega\sp{(e)}
\end{array}\Bigg] + \Pv\Bigg[\begin{array}{c}
\bar \omega_b\sp{0} \medskip\\
\bar b_\omega\sp{0}
\end{array}\Bigg],
\label{eq:sysP}
\end{equation}
where $\Nv$, $\Ov$ and $\Pv$ are $2\times2$ matrices.
Finally, Steps~6 and~7, and the solution for $\bar w\sp{(c)}$, give a system of the form
\begin{equation}
\Qv\Bigg[\begin{array}{c}
\bar \omega_b\sp{0} \medskip\\
\bar b_\omega\sp{0}
\end{array}\Bigg]=
\Rv\Bigg[\begin{array}{c}
\bar \omega_b\sp{(p)} \medskip\\
\bar b_\omega\sp{(p)}
\end{array}\Bigg],
\label{eq:sysC}
\end{equation}
where $\Qv$ and $\Rv$ are $2\times2$ matrices.
The intermediate states, $\bigl[\bar \omega_b\sp{(e)},\bar b_\omega\sp{(e)}\bigr]\sp{T}$ and $\bigl[\bar \omega_b\sp{(p)},\bar b_\omega\sp{(p)}\bigr]\sp{T}$, can be eliminated from the systems in~\eqref{eq:sysE}, \eqref{eq:sysP} and~\eqref{eq:sysC} to give a homogeneous system of the form
\newcommand{\Mt}{{\widetilde M}}
\begin{equation}
\big[\Qv-\Rv \Nv^{-1} (\Ov \Mv+\Pv)\big]\left[\begin{array}{c}
\bar \omega_b\sp{0} \medskip\\
\bar b_\omega\sp{0}
\end{array}\right]=0.
\label{eq:sysDa}
\end{equation}
Nontrivial solutions of~\eqref{eq:sysDa} are obtained if the determinant of the coefficient matrix is zero.
This determinant condition leads to a transcendental equation involving the amplification factor $A$. 
The time-stepping scheme of Algorithm~\ref{alg:ampRBad}  is considered to be stable if none of the solutions of this transcendental equation have $|A|>1$.
This requirement leads to the following theorem:
% ----------------------------------- BEGIN THEOREM ---------------------------------------------------------
\begin{theorem}
   The \ampRB~algorithm given in Algorithm~\ref{alg:ampRBad} for the MP-AD model problem is stable in the 
sense of Godunov and Ryabenkii if and only if there are no roots with $|A|>1$ to the following equation
\begin{equation}
  \NAv  \eqdef \gamma_v(A-1)\sp3+\gamma_0\bigl[\IbBar(A-1)+C_2(A+1)\bigr]A\sp2=0,
\label{eq:advconstrainta}
\end{equation}
where the coefficients in~\eqref{eq:advconstrainta} are given by
\[
\gamma_v\eqdef (2\beta_d\bar{\cal D}\sp{\omega\omega}-C_1)\sp2,\qquad 
\gamma_0 \eqdef \IbBar+4\beta_d\bar{\cal D}\sp{\omega\omega}-C_1.
\]
Here, $\bar{\cal D}\sp{\omega\omega}$ and $\IbBar$ are a dimensionless added-damping coefficient and moment of inertia of the body, respectively, given by
\begin{align*}
 \bar{\cal D}\sp{\omega\omega} \eqdef{\Delta r{\cal D}\sp{\omega\omega}\over\mu(2V_b\aa)}=1-e\sp{-\delta}, \qquad
  \IbBar \eqdef {I_b\over\rho(2V_b\aa)\dr } \, {\delta\sp2}, 
%%  & \tilde\delta \eqdef {\Delta r\over\sqrt{\nu\dt/2}},
\end{align*}
and $C_1$ and $C_2$ are dimensionless Dirichlet-to-Neumann transfer coefficients given by
\[
C_m \eqdef - \Delta r \left[rD_{rh}\left({\phi_m(r)\over r}\right)\right]_{r=\aa}={\aa\over2}\left({3\phi_m(\aa)\over\aa}-{4\phi_m(\aa+\Delta r)\over\aa+\Delta r}+{\phi_m(\aa+2\Delta r)\over\aa+2\Delta r}\right),\quad\text{$m=1$ or $2$}.
\]
\end{theorem}
% ----------------------------------- END THEOREM ---------------------------------------------------------

The number of roots to~\eqref{eq:advconstrainta} is not know a-priori due to the appearance of a transcendental function of $\phi_m(r)$.
Similar to the previous analysis in~\cite{rbinsmp2017}, 
we rely on a simple numerical continuation approach to trace out the stability boundaries in
the parameter space of $(\delta, \beta_d)$, where $\delta$ was defined in~\eqref{eq:addedDampingTerm} and $\beta_d$ is the added-damping parameter.
The stability region is presented in Figure~\ref{fig:stabilityRegion} for the case of $\bar I_b = 0$.
Here we only present the stability region of the massless sphere since it is the most challenging case
for the algorithm.
The important conclusion drawn from Figure~\ref{fig:stabilityRegion} is that there is a large range for $\beta_d$,
near $\beta_d=1$, where the~\ampRB~scheme is stable for any choise of $\delta$.  The scaling of the added-damping term for the full {\ampRB}~algorithm in Algorithm~\ref{alg:ampRB} is equivalent to the choice $\beta_d=1$.

{% ------ MATLAB CURVES  ------
\newcommand{\labelFont}{\scriptsize}
\newcommand{\figWidth}{8.cm}
\newcommand{\trimfig}[2]{\trimw{#1}{#2}{.0}{.0}{.0}{.0}}
\begin{figure}[htb]
\begin{center}
\resizebox{16cm}{!}{% START resize box
\begin{tikzpicture}[scale=1]
  \useasboundingbox (0.0,1.) rectangle (16,6.2);  % set the bounding box (so we have less surrounding white space)
  \draw(3.5,0) node[anchor=south west,xshift=-4pt,yshift=+0pt] {\trimfig{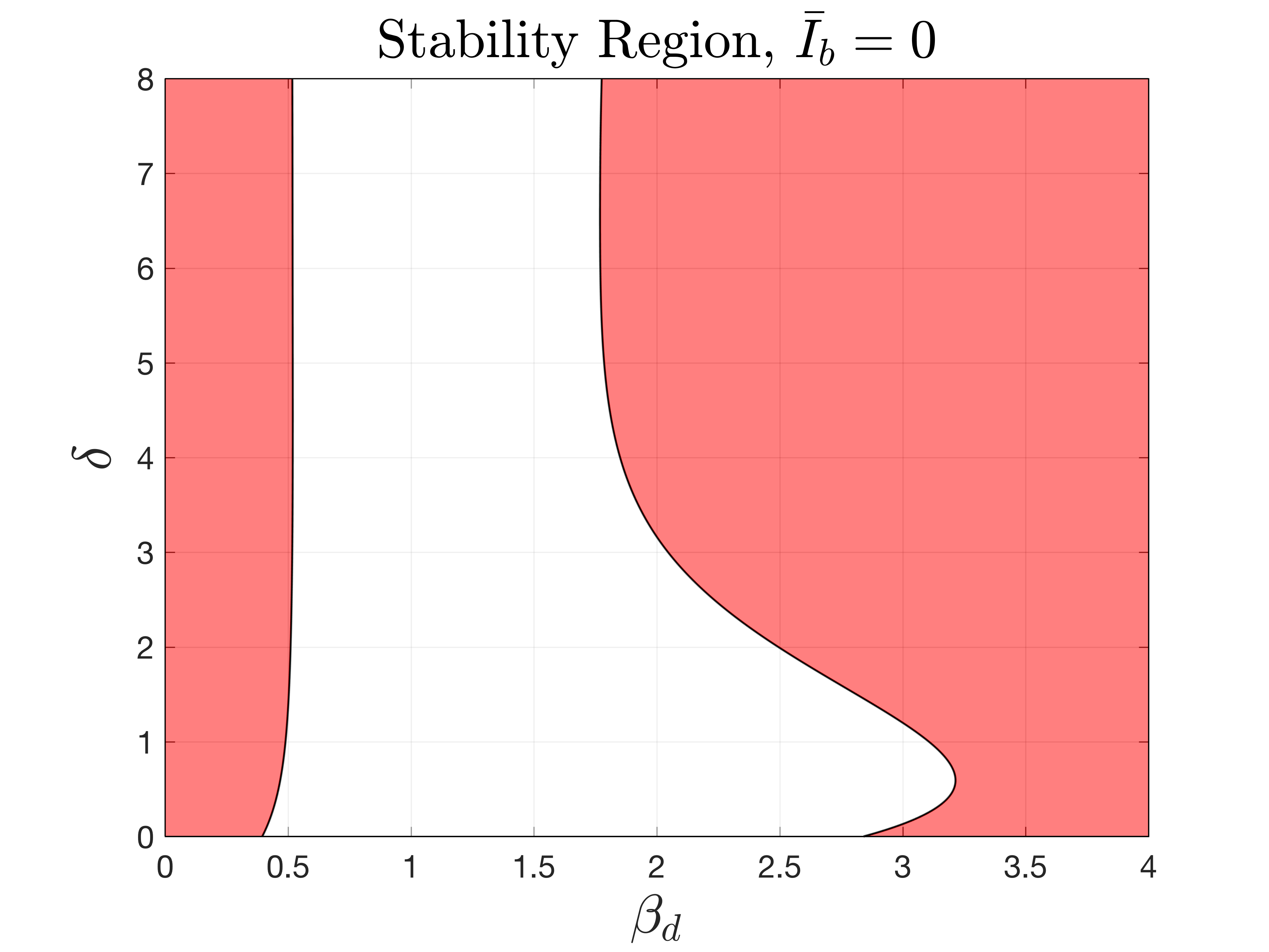}{\figWidth}};
  \draw(7,2.5) node[draw,fill=white] {\labelFont{Stable region}};
% grid:
%\draw[step=1cm,gray] (0,0) grid (16,6.);
\end{tikzpicture}
}% end resize box
\end{center}
\caption{Stability region for the AMP-RB scheme (with velocity correction) for the MP-AD model problem.  The stability region is presented in the $\beta_d$--$\delta$ plane for the zero-mass body. Here $r_1 =20 \sqrt{\nu \dt/2}$ and $r_2 =2 r_1$.}
  \label{fig:stabilityRegion}
\end{figure}
}

Note that the stability analysis and its result are analogous to those of the MP-ADA model problem
(a similar problem but in an annular geometry) in the previous work~\cite{rbinsmp2017}.
In particular, the steps in the analysis and the derived theorem are in parallel to the steps and theorem 
therein. Note that there is another theorem corresponding to the~\ampRB~scheme without the velocity correction 
step (i.e., without Step 8 in Algorithm~\ref{alg:ampRBad}) in~\cite{rbinsmp2017}. 
Its parallel version for the current model problem is similar,
and therefore is not shown for conciseness.

\section{Numerical approach} \label{sec:numericalApproach}

\subsection{Moving composite grids}
\label{sec:comositeGrids}

{
% To change the size:
%    scale the numbers on the first 3 lines : 
% -------------------- Small version ------------------
 %\psset{xunit=.5cm,yunit=.5cm,runit=.5cm}
 \newcommand{\figWidth}{5.5cm}  % 8.25*.5/.75 
 \newcommand{\figWidtha}{3.75cm}   % 5.625*.5/.75 
 \newcommand{\labelSize}{\small}
% -------------------- Medium version -----------------
% \psset{xunit=.6cm,yunit=.6cm,runit=.6cm}
% \newcommand{\figWidth}{6.6cm}  % 8.25*.6/.75 
% \newcommand{\figWidtha}{4.5cm}   % 5.625*.6/.75 
% \newcommand{\labelSize}{\normalss}
% -------------------- Big version --------------------
% \psset{xunit=.75cm,yunit=.75cm,runit=.75cm}
% \newcommand{\figWidth}{8.25cm}
% \newcommand{\figWidtha}{5.625cm}
% \newcommand{\labelSize}{\largess}
% ------------------------------------------------------
% 
\begin{figure}[hbt]
\newcommand{\trimfiga}[2]{\trimFig{#1}{#2}{.05}{.025}{.1}{.02}}
\newcommand{\trimfigb}[2]{\trimFig{#1}{#2}{.055}{.05}{.1}{.04}}
 \newcommand{\trimfig}[2]{\trimFig{#1}{#2}{.05}{.05}{.0}{.0}}
 \begin{center}
 %\begin{pspicture}(-1,-3.2)(14.5,13.6)
%\resizebox{14cm}{!}{% START resize box
\begin{tikzpicture}[scale=.5]
  \useasboundingbox (-1,-2.8) rectangle (14.5,13.8);  % set the bounding box (so we have less surrounding white space)
  \draw(9,8.5) node[]{\trimfig{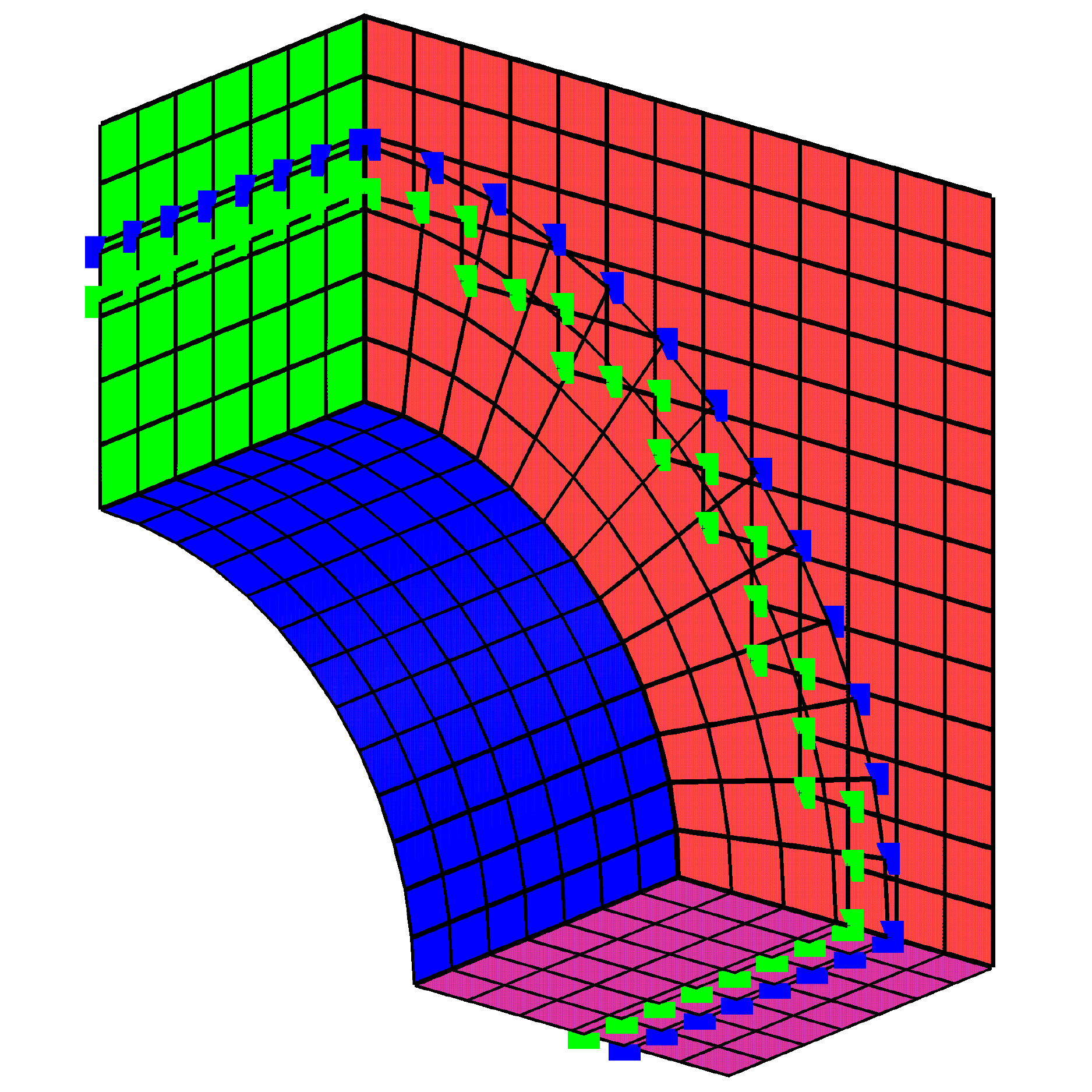}{\figWidth}};
  \draw(2.5,3.8) node[]{\trimfigb{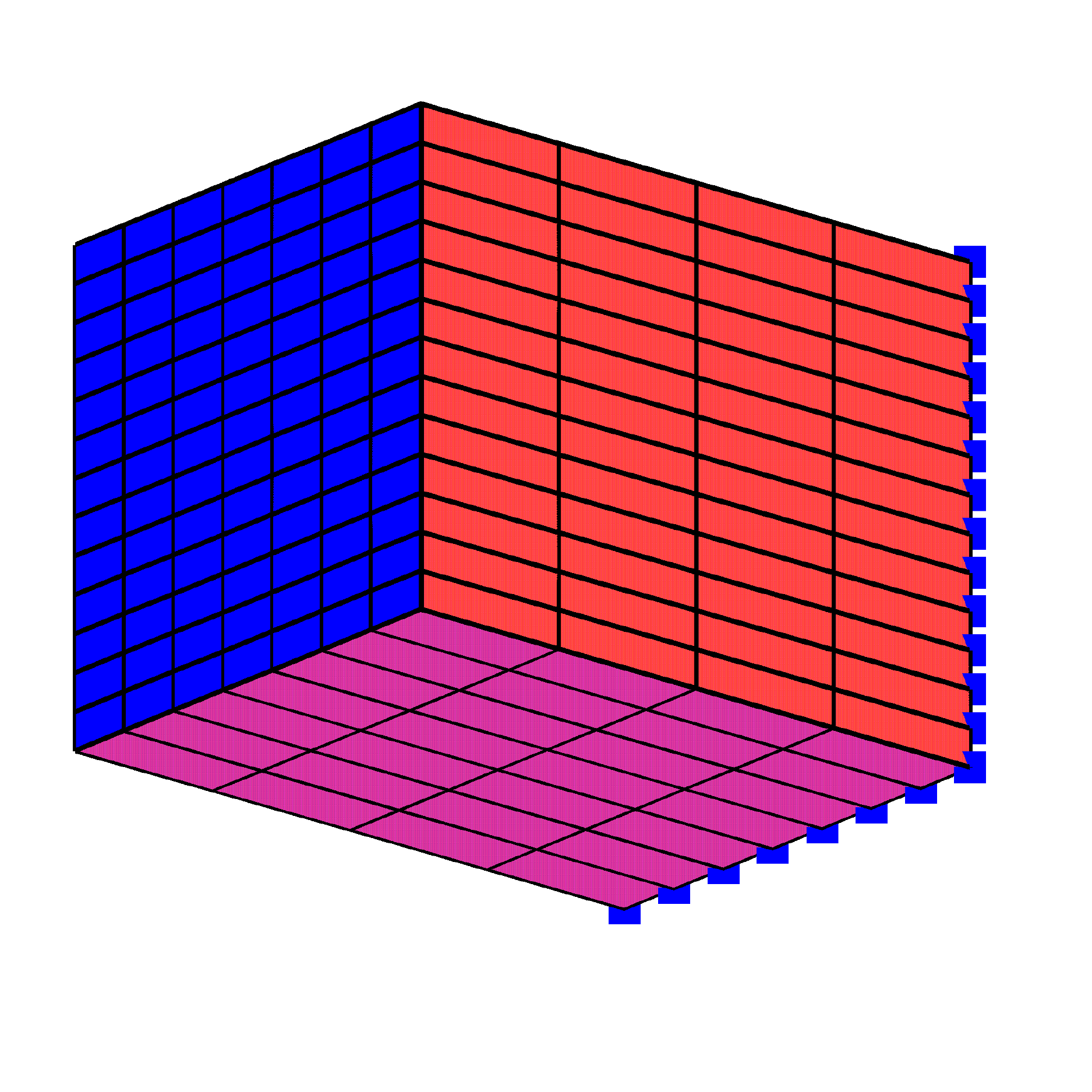}{\figWidtha}};
  \draw(9.6,-.8) node[]{\trimfiga{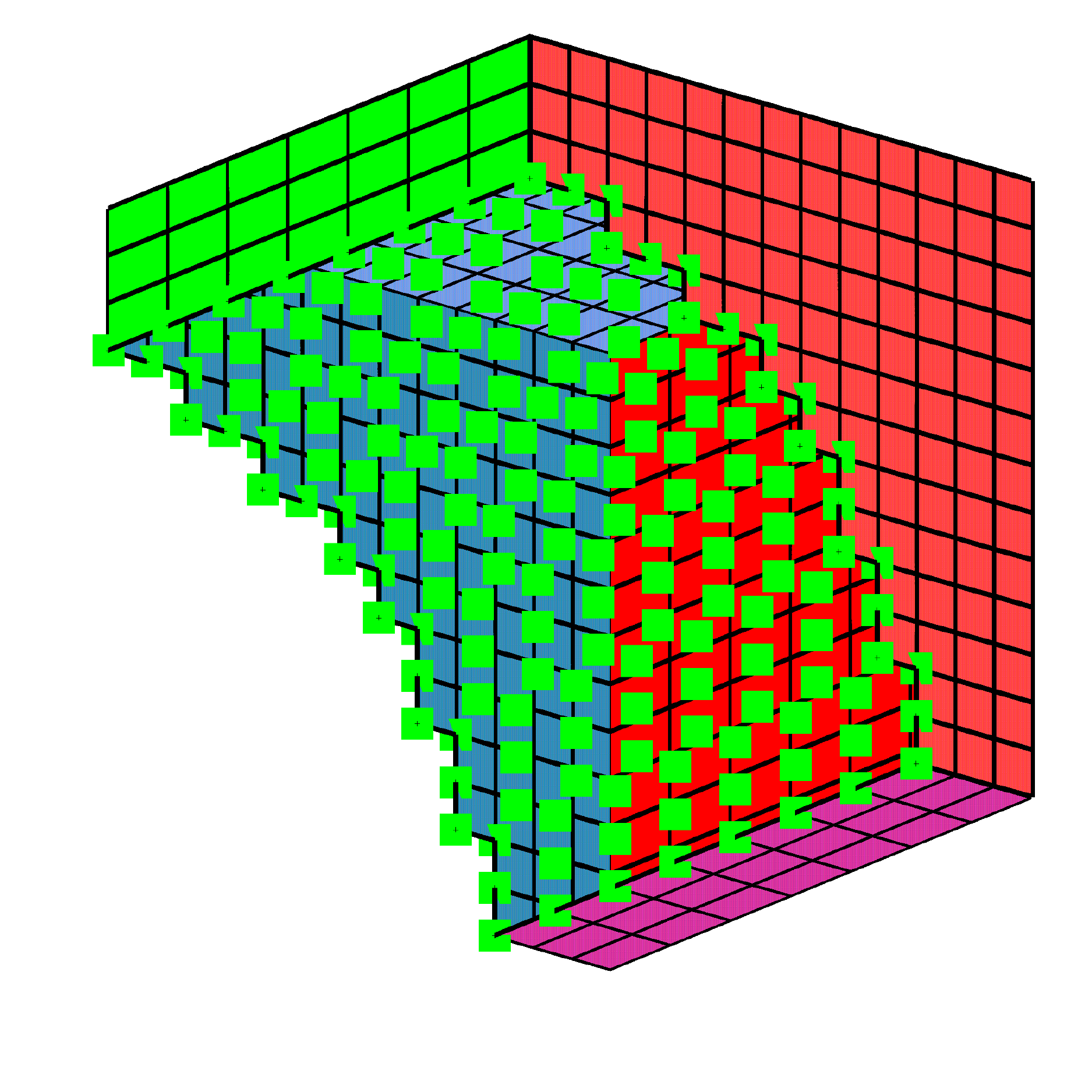}{\figWidtha}};
\draw(12.5,11.5) node[draw,fill=white,anchor=east,xshift=2pt,yshift=-1pt]{\labelSize box};
\draw(8.9,8) node[draw,fill=white,anchor=east,xshift=2pt,yshift=-1pt]{\labelSize cylinder};
\draw(3.,-1.)node[draw,fill=white,anchor=west,xshift=2pt,yshift=-1pt]{\labelSize box unit-cube};
\draw(-1.1,.55)node[draw,fill=white,anchor=west,xshift=2pt,yshift=-1pt]{\labelSize cylinder unit-cube};
\begin{scope}[xshift=-.7cm]
\draw(3.6,10) node[draw,fill=white,anchor=east,xshift=2pt,yshift=-1pt]{\labelSize interpolation points};
\draw[fill=green,draw=green](3.9,10) rectangle ++(.3,.3);
\draw[fill=blue,draw=blue](3.9,9.6) rectangle ++(.3,.3);
\end{scope}
%
% \psgrid[subgriddiv=2]
%\draw[step=1cm,gray] (0,0) grid (14.5,13.8);
\end{tikzpicture}
%}
\end{center}
\caption{Top: a three-dimensional overlapping grid for a quarter-cylinder in a box.
    Bottom left and right: component grids for the cylindrical and box grids in the unit cube parameter space.
Interpolation points at the grid overlap are marked and
color-coded for each component grid. }
\label{fig:overlap3dCartoon}
\end{figure}
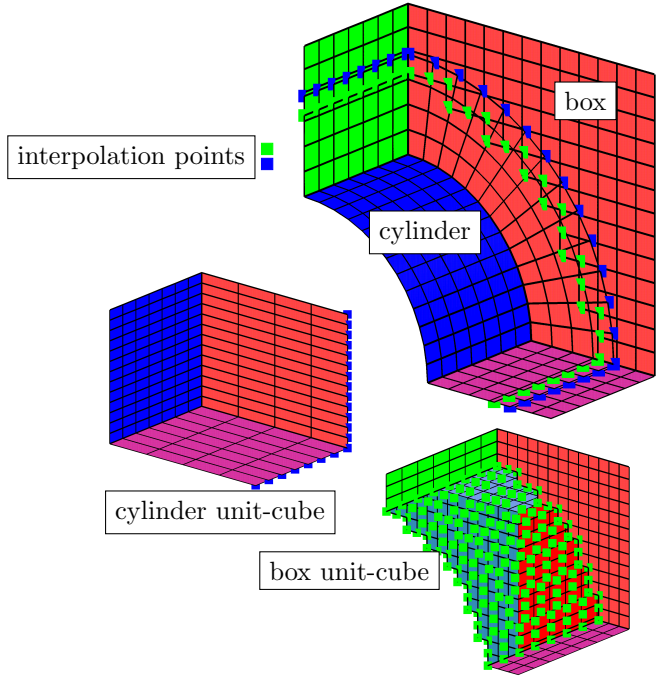
}

\newcommand{\unitvariable}{\boldsymbol{\xi}}

%summary
The moving composite grid approach described originally in~\cite{mog2006} within the context of inviscid, compressible flow is adapted here to handle the moving geometry in three space dimensions associated with the surfaces of the moving rigid bodies in the current FSI regime.
As illustrated in Figure~\ref{fig:overlap3dCartoon},  an overlapping grid, denoted as $\Gc$, 
covers the entire fluid domain and consists of a set of component grids $\{G_g\}$, $g=1,\ldots,{\mathcal N}$.
In three dimensions, each component grid, $G_g$, is a logically rectangular, curvilinear grid
defined by a smooth mapping from a unit cube parameter space~$\unitvariable$ to physical
space~$\xv$,
\begin{equation}
  \xv = \Gv_g(\unitvariable,t),\qquad \unitvariable\in[0,1]^3,\qquad \xv\in\Real^3.
\label{eq:gridMapping}
\end{equation}
All grid points in $\Gc$ are classified as discretization, interpolation or unused points.
The overlapping grid generator {\bf Ogen}~\cite{ogen} from the {\em Overture} framework of codes is used to provide this information.
%In the beginnning of a grid generating process, multiple raw component grids are generated through its definition.
In a typical composite grid for the current FSI application, one or more boundary-fitted curvilinear grids in the fluid are used to represent the surface of each rigid body.  These grids move with the body according to its computed location information $\xvb(t)$ and $\Eb(t)$ following the formulas in~\eqref{eq:bodySurface} and~\eqref{eq:rotationMatrix}.
The remainder of the fluid domain is covered by one or more static Cartesian grids.
{\bf Ogen} cuts holes in the appropriate component grids by locating the physical boundary, which includes the interface between the rigid bodies and the fluid,
and thus determines the unused points based on their location.
For instance, the \lq\lq cylinder'' grid displayed in the upper right image of Figure~\ref{fig:overlap3dCartoon} cuts a hole in the Cartesian \lq\lq box'' grid so that the latter grid has many unused points (those not being plotted in the lower right image).
{\bf Ogen} also provides the interpolation information for all interpolation points in the overlap region between component grids.
Two types of interpolation, called {\em explicit} and {\em implicit}, are supported in {\bf Ogen} and both types are used in the current work.
For explicit interpolation, the grid points used in the interpolation stencil of each donor component grid must be discretization points.  The solution values at interpolation points for explicit interpolation are thus computed explicitly from the values at discretization points on the donor grids.
In contrast, the grid points needed on donor component grids for implicit interpolation may themselves be interpolation points, and thus the solution values at interpolation points are coupled between component grids and must be computed implicitly.
Explicit interpolation is typically used in our simulations because of its faster overall performance.  However, wider overlaps between component grids are needed for explicit interpolation, and thus for some difficult problems, especially when the grids are relatively coarse, implicit interpolation is more robust in producing valid composite grids.
For example, the coarsest grids in the calculations of the settling particle and three-dimensional heart valve use implicit interpolation.
Further discussion of issues related to three-dimensional composite grids can be found in~\cite{pog2008a}.
On discretization points of each component grid, the equations governing the fluid are discretized with second-order accurate approximations on the unit cube parameter space following an exact transformation of the equations using the metrics of the known mapping in~\eqref{eq:gridMapping}.
See~\cite{ICNS, mog2006, pog2008a} for more details on the discretization approach on overlapping grids.

\subsection{Computing surface integrals on composite grids}
\label{sec:surfaceIntegral}

The computation of the forces and torques on the body~\eqref{eq:Pdef} as well as the evaluation
%of the added damping tensors~\eqref{eq:disceteAddedDampingTensors}
of discrete approximations to
the approximate added-damping tensors~\eqref{eq:ADvv}-\eqref{eq:ADww}
% \old{the exact added damping tensors~\eqref{eq:addedDampingTensors}}
requires discrete approximations to integrals over the surface of the rigid body.
This surface quadrature must account for regions on the surface covered by all boundary-fitted
component grids, and must avoid over counting where two or more component grids
overlap. One common way to avoid the over counting issue is to construct an unstructured surface triangulation
in the areas of overlap~\cite{ChanHybridMesh2009}. The solution from the overlapping grid can be interpolated to
the vertices or cell centers of the surface triangles from which a discrete surface integral can be defined.
In this paper we use an alternative approach
to find the surface quadrature that does not require
a surface triangulation\footnote{This approach was discovered by one of the authors
in the course of forming conservative approximations for composite grids~\cite{CGCN,CGCN94}.}. The approach is based on the fact that the Neumann problem 
\begin{subequations}
\label{eq:neumannProblem}
\begin{alignat}{3}
  \Delta \phi &= f,                     \quad&&  \xv\in\Omega, \\
  \f{\partial \phi}{\partial n} &= g    \quad&&  \xv\in\partial\Omega, 
\end{alignat}
\end{subequations}
only has a solution (and that solution is determined up to an additive constant) when
\begin{align}
   \iint_{\Omega} f(\xv) \, dV = \int_{\partial\Omega} g (\xv) \, dS. \label{eq:continuousConstraint} 
\end{align} 
Similarly a discrete approximation to~\eqref{eq:neumannProblem} on an overlapping (or non-overlapping) grid,
\[
      A \phiv = \fv 
\]
only has a solution when 
\begin{align}
  \wv^T \fv=0,  \label{eq:discreteConstraint} 
\end{align}
where $\wv$ is the left null vector of $A$, $\wv^TA=0$.  (As before the solution is determined up to an additive constant.)
The constraint in~\eqref{eq:discreteConstraint} is a discrete counter-part of the integral constraint in~\eqref{eq:continuousConstraint}, and
thus $\wv$ holds the coefficients of a discrete approximation to the volume and surface integrals in~\eqref{eq:continuousConstraint}.
To determine these coefficients, $\wv$ must be scaled by an appropriate constant, and this constant
can be determined by matching the coefficients in $\wv$ to the expected quadrature weights in regions away from
the overlap.  We note that the quadrature weights need only be computed once at the start of the simulation in
view of the transformation in~\eqref{eq:addedDampingInTime}.

\subsection{Parallel implementation and performance}
\label{sec:parallel}

%\todo{add parallel performance here. I think this is a better place to discuss issues and the current performance.
%    I need put choices of linear solver and other information in this section.
%}

%Follow~\cite{pog2008a}. Things to discuss: parallel implementations, the parallel grid figure.

Both serial and parallel versions of the three-dimensional~\ampRB~time-stepping scheme have been implemented based on the incompressible Navier-Stokes solver~\textbf{Cgins}~\cite{CginsUserGuide} within the {\em Overture} framework.  The parallel version of the code follows the work in~\cite{max2006b, pog2008a}, and is based on a domain-decomposition approach in which each component grid $G_g$ is partitioned across different processors  of a distributed-memory parallel computer.
The component grids are partitioned using a load-balancing algorithm that
  depends on the number of the grid points on each grid and the number of processors available.
In the current FSI application, the static background Cartesian grids contain the majority of the grid points in each simulation, and these grids are generally partitioned across all of the available processors.  The boundary-fitted grids are usually narrow in the direction normal to the surface of the rigid bodies and so these grids contain a relatively small number of grid points.  As a result, these narrow grids are usually partitioned across a subset of the available processors, or may even reside on a single processor.
For more details on the parallel implementation of PDE solvers on overlapping grids, see~\cite{max2006b, pog2008a}.

%On each component grid, the field data, such as the fluid velocity and pressure,  is store in a {\em grid function}.
%The partitioning of grid functions follows the decomposition of the component grids and is implemented through the P++ array class~\cite{A++}.
%The P++ array is a C++ class representing distributed multi-dimensional arrays and there is one serial multi-dimensional array on each processor in a P++ array.
%The P++ array can be manipulated through high-level array operations or through low-level Fortran subroutine after each serial array is passed manually.
%Most expensive operations in the algorithm are implemented by the latter approach for efficiency.

%The operations on the P++ array can be m
%Operations of the P++ array is typically implemented by using the local serial array for efficiency. 

In the present work, special attention has been given to the implementation of the~\ampRB~interface conditions in parallel.  These conditions provide additional boundary conditions for the pressure Poisson equation in Steps~4 and~8 of Algorithm~\ref{alg:ampRB}.
The extra equations due to the interface conditions have the same number as the degrees of freedom for the rigid bodies, and these equations are appended to the end of the linear system for the pressure equation.
In our implementation, the extra equations reside in the last available processor in the parallel sparse linear solver.
However, the grid points corresponding to the entries of $\Fv(p)$ in \eqref{eq:Pdef} may reside over several processors in general and this requires some communication across processors.
This communication is very fast as compared to the linear solver since the surface integrals in $\Fv(p)$ are one dimension lower than the full three dimensions of the FSI problem.
All of the parallel communication, including the communication for the interface conditions, are performed  using the Message Passing Interface (MPI).

{
\newcommand{\figWidth}{7.75cm}
\newcommand{\trimfig}[2]{\trimFig{#1}{#2}{.0}{.0}{.0}{.0}}
\begin{figure}[hbt]
\begin{center}
\begin{minipage}[b]{0.44\linewidth}
\begin{tikzpicture}[scale=1]
  \useasboundingbox (0,0.) rectangle (8.,6);  % set the bounding box (so we have less surrounding white space)
  \draw(0,-0.5) node[anchor=south west,xshift=0pt,yshift=+0pt] {\trimfig{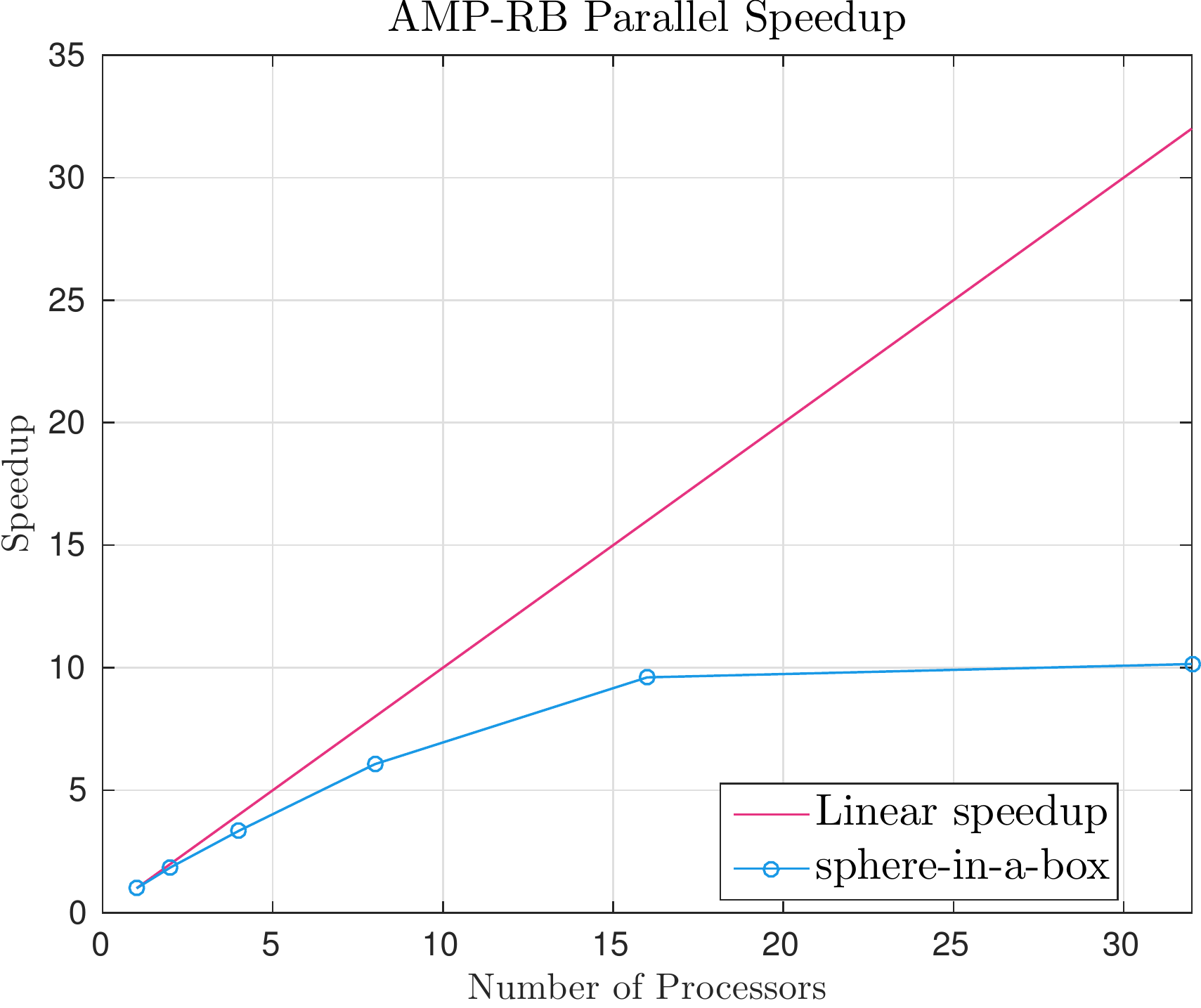}{\figWidth}};
  % grid:
% \draw[step=1cm,gray] (0,0) grid (8,6);
 \end{tikzpicture}
\end{minipage}
\hfill
\begin{minipage}[b]{0.52\linewidth}
    \centering
\tableFont
\begin{tabular}{| c | c | c | c| c | c| c|}
\hline
\multicolumn{7}{|c|}{Parallel results for a sphere-in-a-box grid } \\ \hline
& \multicolumn{2}{c|}{1 precessor} & \multicolumn{2}{c|}{4 processors} &\multicolumn{2}{c|}{16 processors}  \\ \hline
 & s/step & \% & s/step & \% & s/step & \%   \\ \hline
\hline
grid    &   18      &   3.3     & 5     &   3.0     &  2    & 3.2   \\ \hline
psolver &   305     &   56.7    & 95    &   59.4    &  33   & 58.7  \\ \hline
vsolver &   185     &   34.4    &  52   &   32.7    &  19   & 33.2  \\ \hline
others  &   30      &   5.6     &  9    &   4.9     &  2    & 4.9   \\ \hline
total   &   538     & 100       & 161   &  100      &  56   & 100   \\ \hline
%G32 results:
%2    & 4.1 
%34   & 56.3
%20   & 33.5
%4    & 6.1 
%60   & 100 
\end{tabular}
\vspace{1.5cm}
\end{minipage}
\end{center}
  \caption{
Left: strong scaling results for the sphere-in-a-box grid with approximately 13 million grid points. 
Right: CPU time (in seconds) for various parts of the~\ampRB~scheme, and their percentages of the total CPU time per step.
``Grid'' stands for the moving composite grid generator, while
``psolver'' and ``vsolver'' for the Krylov solvers for the pressure and velocity respectively.
The results of 1, 4 and 16 processors are presented.
}
\label{fig:parallelPerformance}
\end{figure}
}

%parallel linear solver
%\paragraph{Sparse Linear Solver} 
All linear systems in the~\ampRB~scheme are solved in parallel using the Krylov methods available in the PETSc library~\cite{petsc-user-ref}.  These linear systems include the implicit equations for the fluid velocity in Steps~3 and~7 (and optionally Step~11) of Algorithm~\ref{alg:ampRB} as well as the aforementioned pressure Poisson equations.
The stabilized bi-conjugate-gradient (Bi-CGSTAB) method with a block-Jacobi preconditioning is used for most of the simulations with an incomplete LU (ILU) preconditioner for each individual block.
The velocity equation uses ILU(1) as the preconditioner, while the pressure equation uses ILU(1) for easy problems and ILU(3) for some hard problem that require 20 to 50 iterations to converge.
We found that these choices give the best performance amongst the Krylov solvers tested.
For one three-dimensional heart valve simulation, the parallel Bi-CGSTAB solver failed to converge in the first pressure solve, presumably from a poor initial guess. 
For this case, a parallel GMRES solver was used for a few time steps before switching back to the Bi-CGSTAB solver for improved efficiency.

%\dws{I suggest that we move this next sentence to the Conclusions section since it mentions what we plan to do next as an extension of the present paper.} 
%was found to give the best results amongst the Krylov methods that were tested. 
%The Krylov solvers used here are from the PETSc library.
%Algebraic multigrid preconditioners do not perform well in our case.
%In Future work, we are interested in using Geometric Multigrid solver to solve the new system.
%Simple problem we used ILU(1). Hard problem ILU(3) with GMRES in the first few time steps.

%\paragraph{Parallel Performance} 
Figure~\ref{fig:parallelPerformance} presents the parallel performance of the~\ampRB~scheme on a distributed memory Linux cluster with 2.5 GHz Xeon processors.
The figure gives strong scaling results for the settling-particle problem described later in Section~\ref{sec:settlingParticle} using a sphere-in-a-box grid  with $1.3 \times 10^{7}$ grid points. 
The graph and the CPU time table show that the parallel scaling is reasonable up to 16 processors.
%This scaling performance increases as the grid is being refined.
The breakdown of the CPU times in the table shows that most of the computational cost is used in the Krylov solves of the pressure and velocity equations, and thus the parallel performance depends heavily on the scaling of the Krylov solvers for this problem.
The drop-off in parallel efficiency in the Krylov solvers is
partially due to the preconditioners being precomputed at every time-step.
We note that the computational domain for this problem is relatively simple, and hence the {\bf Ogen} grid generator is very efficient and scales well in parallel.
However, the grid generator has not yet been optimized for moving grid problems and can be more costly for complex problems.  For example, it has been observed that the grid generator can use up to half of the computational cost for the heart valve simulations described later in Section~\ref{sec:heartValve3D}.

The results displayed in Figure~\ref{fig:parallelPerformance} should be taken as a baseline for the initial parallel implementation as much work remains to improve the parallel scaling. Possible improvements include (i)~reusing the preconditioner (which is currently recomputed at every step),
(ii)~using a multigrid solver, and (iii)~optimizing the parallel grid generation.
%\QT{
%G4 seems acceptable. G1 and G2 are bad. The overhead is too high in the coarse grid.
%
%In another case, the scaling of all periodic boundary condition is too good to be true. Blame the dense equation fix.
%
%%one strong scaling for G4 of settling particles
%%one table for each pieces in parallel. np=1 vs np=8 vs np=32
%
%%cpu info: cgpi, a distributed memory Linux cluster (with 2.5 GHz Xeon processors).
%}

\section{Numerical results}
\label{sec:numericalResults}

Several benchmark problems are considered in this section to demonstrate the stability and accuracy of the~\ampRB~scheme.
We begin in Section~\ref{sec:piston3D} by revisiting the one-dimensional motion of a piston in a fluid chamber discussed in~\cite{rbinsmp2017, rbins2017}.  Here, we solve the problem using the~\ampRB~scheme for a three-dimensional fluid chamber and compare the numerical results with the exact solution.  This benchmark problem provides a good first test of both the accuracy of the scheme as well as its added-mass stability properties.
%An additional grid is intentionally involved in the composite grid to confirm the accuracy of the surface integral on bodies covered by multiple overset grids.
%
In Section~\ref{sec:settlingParticle}, we consider a second benchmark problem involving the settling of a moderately light particle in a small container.  For this problem, the numerical solution given by the~\ampRB~scheme is compared with experimental data and with other numerical results from the literature, and a self-convergence study is conducted to demonstrate the accuracy of the scheme.
We also use this second problem to study the added-damping stability properties of the~\ampRB~scheme.
Section~\ref{sec:compareSchemes} is devoted to a comparison of the~\ampRB~scheme and a traditional partitioned (TP-RB) scheme, which is performed using the previous settling particle problem.  Both schemes require the solution of a Poisson problem for the pressure at each time step.  The pressure solve for the TP-RB scheme uses a boundary condition for the pressure alone whereas the~\ampRB~scheme uses interface conditions involving the pressure and the accelerations of the rigid body.  The corresponding linear systems are different but have similar conditioning.
%, which demonstrates that there is no loss of accuracy in the solution of the linear system for the~\ampRB~scheme.  
We are also interested in a comparison of the performance of the two schemes, and this demonstrates the efficiency of the present scheme.
In Section~\ref{sec:risingParticle}, a third benchmark problem involving a particle rising or falling in a long container is considered.
The numerical solution given by the~\ampRB~scheme is compared with the results of various algorithms discussed in the literature.
Importantly, this problem is used to confirm that the~\ampRB~scheme is stable for rigid-body FSI problems involving different density ratios of the solid bodies and the fluid, including the difficult case of a body with small, or even zero, mass.  This difficult case has proved to be challenging, or impossible, for other available schemes due to the severe added-mass instability.
The last problem we consider is an FSI problem generalized from a bi-leaflet mechanical heart valve.
Simulations of this problem in two dimensions (Section~\ref{sec:heartValve2D}) and three dimensions (Section~\ref{sec:heartValve3D}) are discussed.
For this application, two leaflets are hinged  and allowed to rotate along the hinge axes.
An oscillating pressure gradient is applied to drive a fluid flow through the value, and repulsive torques and damping terms are implemented to restrict the rotation of the leaflets to the designed range.
This example demonstrates the behavior of the~\ampRB~scheme when multiple light bodies are involved,
and further confirms that the numerical instabilities have been suppressed in the present scheme for such cases.

\subsection{One-dimensional motion of a piston}
\label{sec:piston3D}

%general
In this section, we consider the piston problem discussed previously in~\cite{rbinsmp2017, rbins2017}, but extended here to three space dimensions.
In this problem a rigid body is adjacent to one end of a fluid channel with rectangular cross section.  Two problems were considered previously, one in which the body moved normal to the channel and the other where the body moves tangentially, and exact solutions were found for each case.
%The problem has been solved in two dimensions therein and served as a good test problem to verify added-mass effects of the~\ampRB~scheme.
Here, we focus on the problem from~\cite{rbins2017} in which the body moves normal to the channel, and we use this problem to verify the stability and accuracy of the~\ampRB~scheme. 
Unlike the previous test of the two-dimensional version of the scheme where a single body-fitted rectangular grid is used to represent the surface of the body, the current test uses multiple body-fitted component grids to verify the implementation and accuracy of the calculation of surface integrals on the body for a case when the body-fitted grids overlap.

%exact solutions:
For this problem, a rectangular rigid body of size $L_b \times \channelHeight \times W$ and mass $\mb$ is adjacent to the fluid domain, $\OmegaF(t)=[x_I(t),\channelRight]\times[0,\channelHeight]\times[0,W]$, where $x = x_I(t)$ is the body-fitted end of the fluid channel.  The fluid in the channel is governed by the incompressible Navier-Stokes equations with a density and viscosity taken to be $\rho=1$ and $\mu=0.1$, respectively.
The boundary conditions on the fluid at the fixed end of the channel, $x=\channelRight$, are 
\[
    p(\channelRight,y,z,t) = p_\channelRight(t), \qquad v_2(\channelRight,y,z,t) = v_3(\channelRight,y,z,t) = 0,
    \qquad (y,z)\in[0,\channelHeight]\times[0, W],\quad t>0,
\]
where $p_\channelRight(t)$ is an applied pressure to be determined.  
The velocity matching condition in~\eqref{eq:RBsurfaceVelocity} is applied at $x = x_I(t)$, and slip-wall boundary conditions are applied on the remaining four sides of the fluid channel.
Exact solutions for the velocities of the body, $\vvb=(v_{1,b},v_{2,b}, v_{3,b})\sp{T}$, and fluid, $\vv=(v_{1},v_{2},v_3)\sp{T}$, and for the fluid pressure $p$ have the form
\begin{align}
v_{1,b}(t) &= \int_0\sp{t}{-H W p_\channelRight(\tau)\,d\tau\over \mb+M_a(\tau)}+v_{1,b}(0), \label{eq:horizontalVelocity} \\
p(x,t) &= p_\channelRight(t) + \left({\channelRight-x\over\channelRight-x_I(t)}\right){-p_\channelRight(t)\over \mb/M_a(t)+1}, \label{eq:channelpressure} \\
v_1(t) &= v_{1,b}(t), \label{eq:channelVelocity}
\end{align}
and $v_2=v_3=v_{2,b}=v_{3,b}=0$.  The {\em added-mass} for this problem found analytically in~\cite{rbins2017} to be
\begin{equation}
M_a(t)=\rho\channelHeight W \bigl(\channelRight-x_I(t)\bigr).
\label{eq:addedMass}
\end{equation}
%The horizontal position of the body center is given by
Because the horizontal velocity is $v_{1,b}(t)$, the horizontal position is calculated as
\begin{equation}
x_b(t) = \int_0\sp{t} v_{1,b}(\tau)\,d\tau + x_b(0).
\label{eq:horizontalPosition}
\end{equation}
The applied pressure at the fixed end of the channel given by $p_\channelRight(t)$ can be specified by choosing a motion for the rigid body.  For the present tests, we take the interface $x_{I}(t)$ to be

\begin{equation}
x_I(t) =x_b(t)+L_b/2 = \Amplitude\sin(2 \pi t),\qquad \Amplitude=1/4,
\label{eq:prescribedPistonMotion}
\end{equation}
which is sufficient to specify the exact solution.
The geometric configurations based on the given interface equation are presented in Figure~\ref{fig:pistonGrids}
for two different times.

\input texFiles/piston3DGridFigure

%composite grids
Numerical solutions are computed using a composite grid, denoted by $\Gcp^{(j)}$, with resolution factor~$j$, as shown in Figure~\ref{fig:pistonGrids}.
The composite grid consists of three component grids, the first of which is a body-fitted Cartesian grid (green in the figure) of fixed length $1/2$ in $x$-direction, and of width $W = 1$ and height $\channelHeight=1$ in $y$ and $z$-directions, respectively.
Inside this first body-fitted grid, there is an extra body-fitted cylindrical grid (red in the figure) of fixed length equal to $0.4$ in the $x$-direction, and with inner radius $0.2$ and outer radius $0.4$.
These two body-fitted grids overlap and move in time as the body moves. 
Note that the cylindrical grid is not needed to handle the geometry of the problem, but is intentionally added to test the calculation of surface integrals on the body as discussed in Section~\ref{sec:surfaceIntegral}.
In addition to these body-fitted grids, there is a static Cartesian background grid (blue in the figure) covering the domain 
$[-3/4,\channelRight]\times[0,\channelHeight]\times[0,W]$ with $\channelRight=3/2$.
The approximate grid spacing for all of the grids is $\hj=1/(10 j)$ in each of their coordinate directions.  
The length of the rigid body is chosen to be $L_b=1$ and hence $\mb=\rho_b$, which will be varied in different cases.
%The density of the fluid is taken to be $\rho=1$ and its viscosity is $\mu=0.1$ for all cases.  
Figure~\ref{fig:pistonGrids} also shows the computed pressure at $t = 0.8$ (using $\dt=0.01$) for the case of a very light body with $\rhob=0.001$ along with a wire frame of parts the composite grid $\Gcp^{(4)}$.

{% ------ CURVES ------
\newcommand{\figWidtha}{6cm}% height 
\newcommand{\figWidths}{5.95cm}
\newcommand{\trimfiga}[2]{\trimh{#1}{#2}{.0}{.0}{.0}{.0}}
\begin{figure}[htb]
\begin{center}
\resizebox{14cm}{!}{% START resize box 
\begin{tikzpicture}[scale=1]
  \useasboundingbox (0,0) rectangle (16.,11.5);  % set the bounding box (so we have less surrounding white space)
%  --- curves:
  \draw(-.5,5.4) node[anchor=south west,xshift=0pt,yshift=+0pt] {\trimfiga{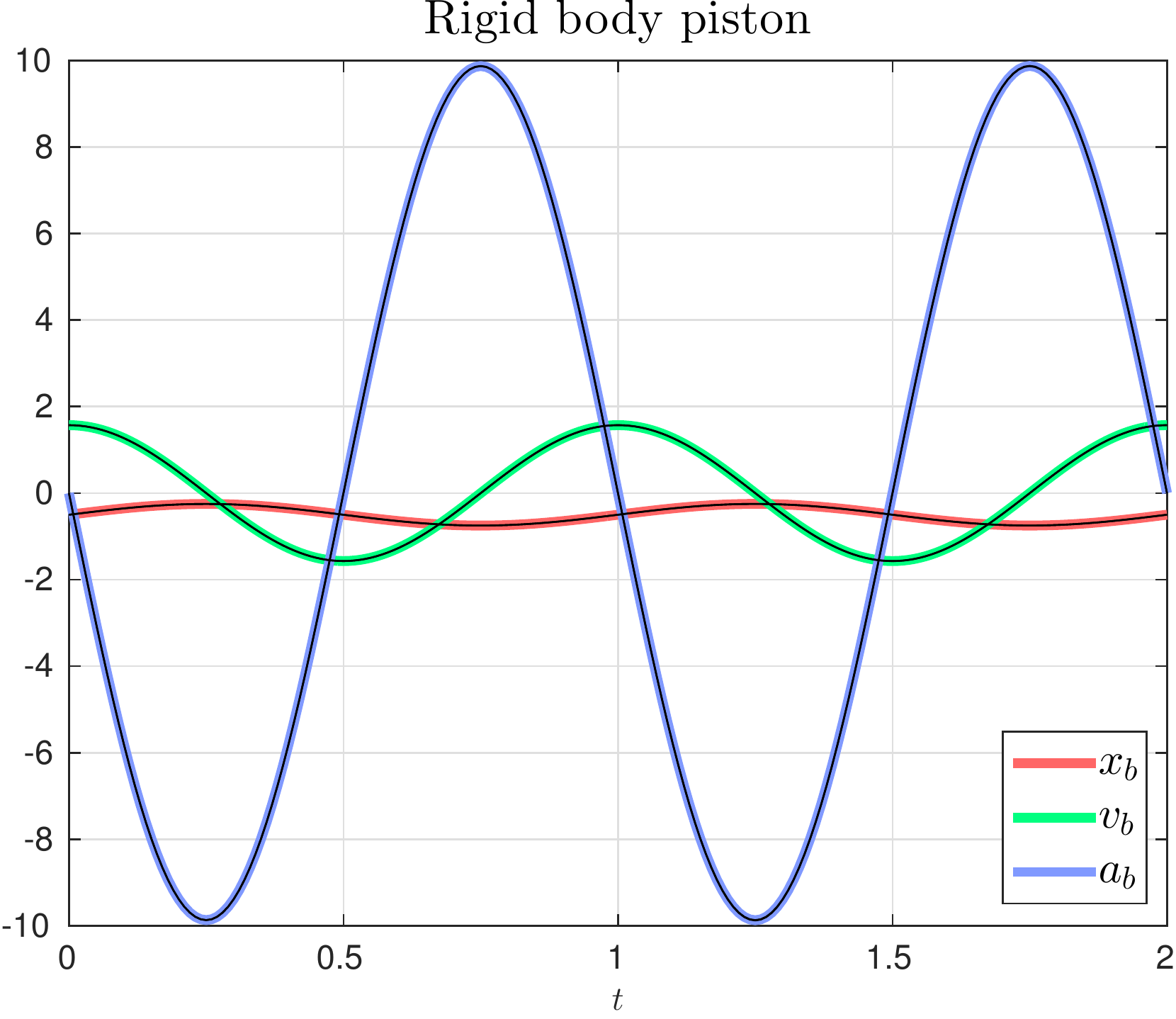}{\figWidtha}};
  \draw(8.0,5.4) node[anchor=south west,xshift=0pt,yshift=+0pt] {\trimfiga{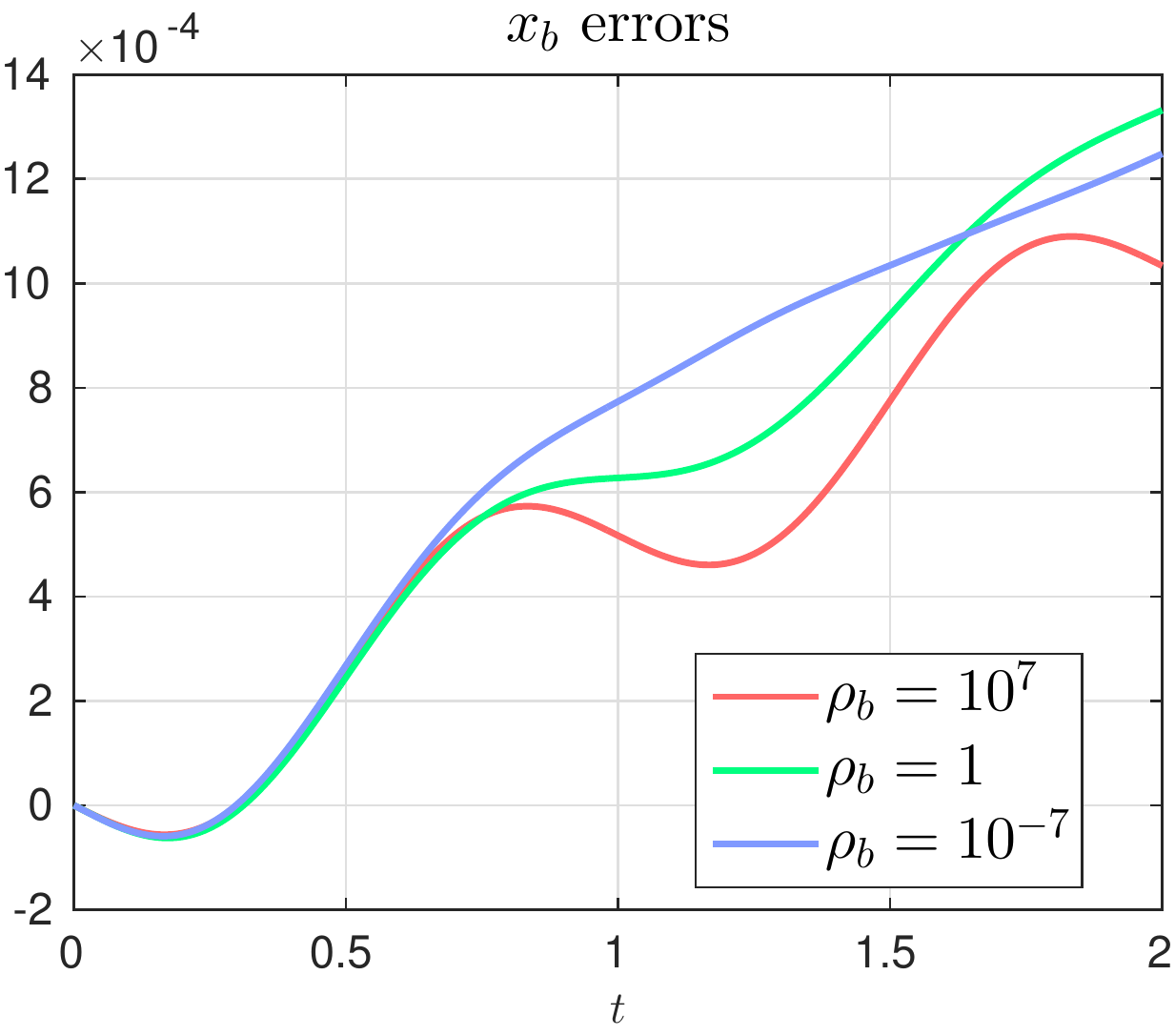}{\figWidths}};
  \draw(-.5,-.95) node[anchor=south west,xshift=0pt,yshift=+0pt] {\trimfiga{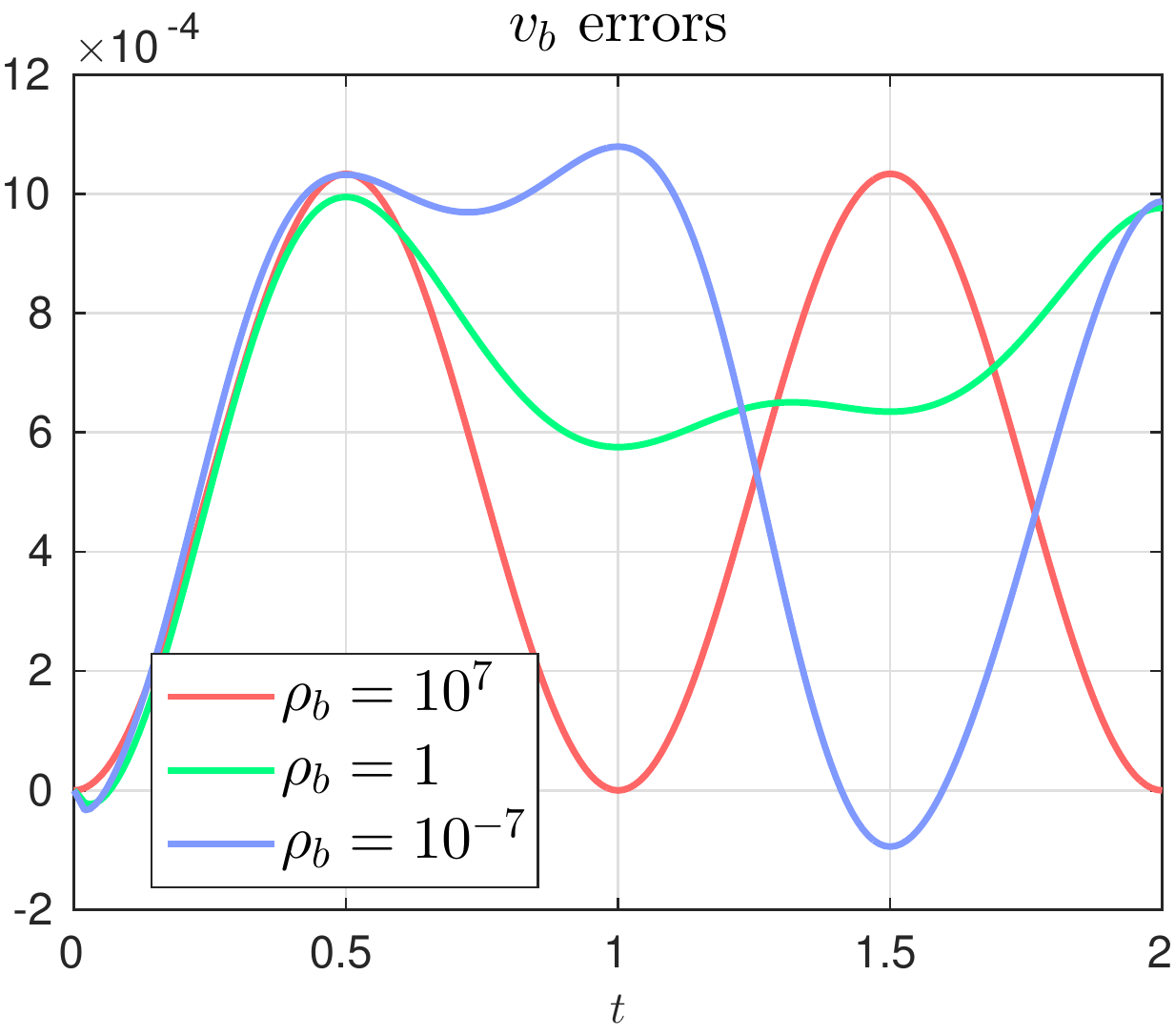}{\figWidths}};
  \draw(8,-.95) node[anchor=south west,xshift=0pt,yshift=+0pt] {\trimfiga{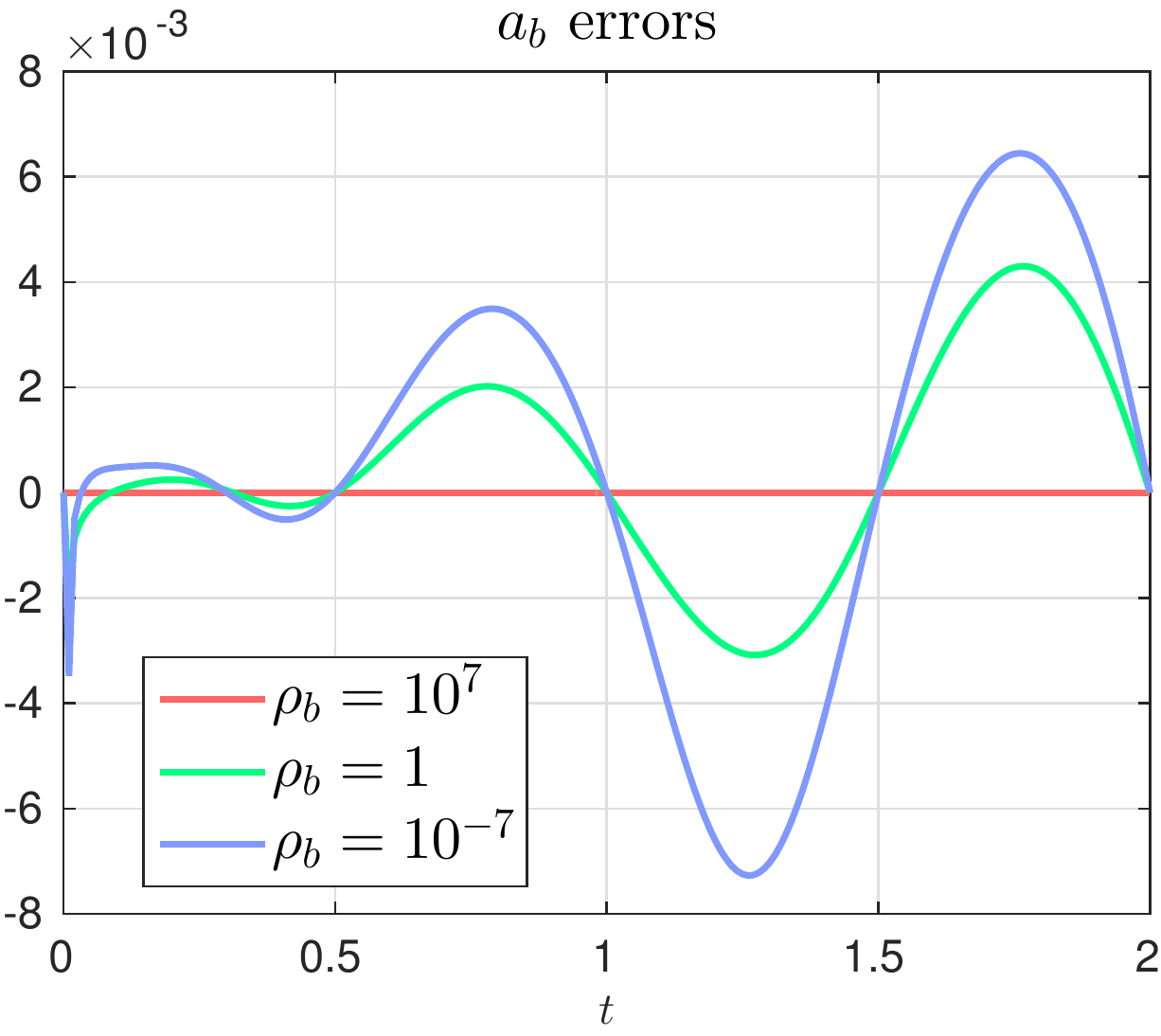}{\figWidths}};
%
% grid:
 %\draw[step=1cm,gray] (0,0) grid (16,11.5);
\end{tikzpicture}
}% end resize box
\end{center}
\caption{Computed piston motion and errors for $\rhos=10^7$, $\rhos=1$  and $\rhos=10^{-7}$. The solutions are computed on grid $\Gcp^{(4)}$.
   }
  \label{fig:pistonCurves}
\end{figure}
}

%added mass issue
As a test, the~\ampRB~scheme has been run for this piston problem using $\Gcp^{(4)}$ and $\dt=0.01$ for a wide range of densities of the body, from very light {($\rhob = 10^{-7}$)} to very heavy {($\rhob = 10^{7}$)}, and for all cases we observe that the scheme remains stable and gives accurate results.  For example, the time history of the position, velocity and acceleration of the rigid body, 
as well as their errors, are presented in Figure~\ref{fig:pistonCurves} for a very heavy body with $\rhob=10^7$, a medium body with $\rhob=1$ and a very light body $\rhob=10^{-7}$.  These results illustrate the stability and accuracy of the simulations.  As a comparison, we have also used the TP-RB scheme without sub-timestep-iterations to compute numerical solutions of this problem and have found that it is unstable for calculations using $\Gcp^{(4)}$ and $\dt=0.01$ if $\rhob =\mb \approx 2.0$ or smaller.  
This observation is in reasonable agreement with the stability results in~\cite{rbinsmp2017, rbins2017}, where additional details are discussed.

% See runs/ins/conv3/memo
{ % =================================== START TABLE ===============================================
\newcommand{\errp}{$E_j^{p}$}
\newcommand{\erruv}{$E_j^{\vv}$}
\newcommand{\errxA}{$E_j^{\xvb}$}
\newcommand{\errvA}{$E_j^{\vvb}$}
\newcommand{\erraA}{$E_j^{\avb}$}
%
% -------------------------------------------------------------------------------------------------------
\newcommand{\pistonTableI}{%
\begin{tabular}{|c|c|c|c|c|c|c|c|c|c|c|} \hline 
  \multicolumn{11}{|c|}{Piston motion, $\rhob=10$} \\ \hline 
\strutt$\hj$&     \errp     &  r   &     \erruv     &  r   &     \errxA     &  r   &     \errvA     &  r   &     \erraA     &  r    \\[3pt] \hline 
  1/10  & \num{1.4}{-2} &      & \num{6.9}{-3} &      & \num{9.2}{-3} &      & \num{6.9}{-3} &      & \num{7.5}{-3} &     \\ \hline
  1/20  & \num{3.2}{-3} & 4.3  & \num{1.7}{-3} & 4.1  & \num{2.3}{-3} & 4.0  & \num{1.7}{-3} & 4.1  & \num{1.8}{-3} & 4.1 \\ \hline
  1/40  & \num{8.1}{-4} & 4.0  & \num{4.5}{-4} & 3.8  & \num{5.6}{-4} & 4.0  & \num{4.2}{-4} & 4.0  & \num{4.5}{-4} & 4.0 \\ \hline
  rate        &    2.04       &      &    1.97       &      &    2.01       &      &    2.02       &      &    2.02       &   \\  \hline
\end{tabular}
}
% -------------------------------------------------------------------------------------------------------
\newcommand{\pistonTableII}{%
\begin{tabular}{|c|c|c|c|c|c|c|c|c|c|c|} \hline 
  \multicolumn{11}{|c|}{Piston motion, $\rhob=1$} \\ \hline 
\strutt$\hj$&     \errp     &  r   &     \erruv     &  r   &     \errxA     &  r   &     \errvA     &  r   &     \erraA     &  r    \\[3pt] \hline 
  1/10  & \num{5.9}{-2} &      & \num{1.2}{-2} &      & \num{9.6}{-3} &      & \num{1.2}{-2} &      & \num{3.4}{-2} &     \\ \hline
  1/20  & \num{1.4}{-2} & 4.2  & \num{2.8}{-3} & 4.2  & \num{2.4}{-3} & 4.1  & \num{2.8}{-3} & 4.2  & \num{8.1}{-3} & 4.2 \\ \hline
  1/40  & \num{3.5}{-3} & 4.1  & \num{6.9}{-4} & 4.0  & \num{5.8}{-4} & 4.0  & \num{6.9}{-4} & 4.1  & \num{2.0}{-3} & 4.1 \\ \hline
  rate         &    2.04       &      &    2.04       &      &    2.02       &      &    2.04       &      &    2.04       &  \\    \hline
\end{tabular}
}
% -------------------------------------------------------------------------------------------------------
\newcommand{\pistonTableIII}{%
\begin{tabular}{|c|c|c|c|c|c|c|c|c|c|c|} \hline 
  \multicolumn{11}{|c|}{Piston motion, $\rhob=0.001$} \\ \hline 
\strutt$\hj$&     \errp     &  r   &     \erruv     &  r   &     \errxA     &  r   &     \errvA     &  r   &     \erraA     &  r    \\[3pt] \hline 
  1/10  & \num{9.9}{-2} &      & \num{1.7}{-2} &      & \num{1.0}{-2} &      & \num{1.6}{-2} &      & \num{5.8}{-2} &     \\ \hline
  1/20  & \num{2.4}{-2} & 4.1  & \num{4.0}{-3} & 4.1  & \num{2.6}{-3} & 4.0  & \num{4.0}{-3} & 4.1  & \num{1.4}{-2} & 4.1 \\ \hline
  1/40  & \num{6.0}{-3} & 4.0  & \num{9.8}{-4} & 4.0  & \num{6.4}{-4} & 4.0  & \num{9.8}{-4} & 4.0  & \num{3.5}{-3} & 4.0 \\ \hline
  rate  &    2.02       &      &    2.03       &      &    2.00       &      &    2.02       &      &    2.03       &     \\ \hline
\end{tabular}
}
{
\begin{table}[hbt]\tableFont % you should set \tableFont to \footnotesize or other size
\begin{center}
  \pistonTableI \\
\bigskip
  \pistonTableII \\
\bigskip
  \pistonTableIII
\caption{Piston motion in three dimensions. Maximum errors and estimated convergence rates at $t=0.8$ computed using the~\ampRB~scheme
   for a heavy, $\rhob=10$, medium, $\rhob=1$, and very light, $\rhob=0.001$, moving piston. 
    The column labeled "r" provides the ratio of the errors at the current grid spacing to that on the next coarser grid.
 }
\label{table:piston3D}
\end{center}
\end{table}
}
} % =================================== END TABLE ===============================================

A refinement study is conducted on a sequence of grids of increasing resolution to check the accuracy of the implementation of the~\ampRB~scheme for three dimensions.
The time-step is taken as $\dt_j=0.4/(10 j)$ for the composite grid $\Gcp^{(j)}$, and the equations are integrated to $t_{{\rm final}}=0.8$. 
Note the time step is smaller than the refinement study conducted in~\cite{rbins2017} due to the stability constraint from the additional cylindrical grid which has cells with smaller volumes.
However, the numerical solutions still exhibit superconvergence at $t=n/2$, where $n$ is a nonnegative integer, 
corresponding to the half-period of the oscillatory motion of the body specified in~\eqref{eq:prescribedPistonMotion} when $x_b=a_b=0$, which is the case found in~\cite{rbins2017} at $t_{{\rm final}} = 1$.
Therefore, to avoid this superconvergence, the equations in the current study are integrated to $t_{{\rm final}} = 0.8$ instead.
Table~\ref{table:piston3D} presents the max-norm errors and estimated convergence rates
of the refinement study for a heavy, medium and light body. 
The max-norm error of a quantity $q$ is denoted by $E_j^{q}$ for the results on grid $\Gcp^{(j)}$, and if $q$ is a vector, then the maximum is taken over all of its components.  The results of the table show approximately second-order accuracy for all solution quantities and for all three densities of the body.  As an added check, the calculations were also performed without the extra cylindrical body-fitted grid and the errors were found to be nearly identical to the errors given in the table.  This verifies  the accuracy of the discrete surface integrals when computing the fluid force on the rigid body.

\subsection{One particle settling in a box} 
\label{sec:settlingParticle}

\newcommand{\Gcsp}{\Gc_{\rm sp}} % changed "f" (falling) to "r" (rising)

We consider the common benchmark problem of a moderately light 
spherical particle settling due to gravity in a small container. 
This problem was first studied experimentally~\cite{Cate2002} by considering a spherical particle settling in silicon oil towards the bottom of a box, 
and it has later been widely used as a test problem to validate numerical FSI algorithms, including~\cite{gibou2012efficient, kempe2012improved, yang2015non, Koblitz2016}.
The setup of the problem is as follows.
A spherical particle of density $\tilde \rhos=1120\;{\rm kg}/{\rm m}\sp3$ and diameter $D_b =2\,\diskRadius = 15\;{\rm mm}$ is immersed in a fluid channel of a rectangular box.
The fluid density is $\tilde \rho = 970 \;{\rm kg}/{\rm m}\sp3$ and its viscosity is $\tilde\mu=0.373\;{\rm kg}/({\rm m}\,{\rm s})$.
Initially the fluid and particle are quiescent,
and at $t = 0$ the particle is allowed to move under the influence of gravity as specified by
the external body force,
\begin{equation}
\label{eq:particleForce}
\fvbe(t)=\frac 4 3 \pi \diskRadius\sp3(\rho_b-\rho) \, \gv,
\end{equation}
where the acceleration is $\gv=[0,\,0,\,-\tilde g]$ with $\tilde g=9.81\;{\rm m}/{\rm s}\sp2$.
In the previous studies, the problem is typically presented in terms of dimensional units. 
Here we prefer to nondimensionalise the problem by introducing suitable reference scales. 
The reference scale for length is the diameter of the particle $\tilde d = D_b$.
The reference velocity is taken to be $\tilde v_{\rm term}= 0.038\;{\rm m}/{\rm s}$ given in the original work~\cite{Cate2002},
which is the terminal velocity of the particle in an infinitely long channel with the corrections due to the side boundary effects.
The densities are all scaled by the fluid density of $\tilde \rho$.
The dimensionless fluid density is therefore $\rho=1$ and the body density is $\rhob=1120/970\approx1.1546$.
The dimensionless fluid viscosity is hence $\mu=\tilde\mu/(\tilde\rho\, \tilde d\, \tilde v_\text{term})=0.67462$, corresponding to the Reynolds number ${\rm Re} = 1.5$.
The dimensionless gravity is  $g=\tilde g\,\tilde d/(\tilde v_\text{term})\sp2=101.90$,
which is turned on instantaneously at $t = 0$ in the simulation.
The resulting fluid domain is $\Omega(t) = [-x_c, x_c] \times [-y_c, y_c] \times [z_0, z_1]$, with $x_c = y_c = 10/3$, $z_0=0$ and  $z_1 = 32/3$.
The particle of radius $0.5$  is initially  located at $(x,y,z) = (0,0,8.5)$.
The top boundary condition of the fluid domain is a zero pressure condition, and the remaining boundary conditions are no-slip walls.

Numerical solutions are computed using a {\em sphere-in-a-box} composite grid, denoted by $\Gcsp^{(j)}$, which consists of a background Cartesian grid and two body-fitted grids attached to the surface of the rigid body.
The target background grid spacing $h^{(j)}$ is set to be $2/(15j)$, corresponding to the grid size of $(50 j)\times(50j)\times(80j)$ covering the entire domain.
The boundary-fitted grids are two curvilinear grids surrounding the particle, 
each of which covers more than half of its surface.
The thin body-fitted grids have a target grid spacing approximately equal to $0.71 h^{(j)}$ with a radial width of $5 h^{(j)}$.
The time step $\dt^{(j)}$ is taken as $1/(40j)$.
The interested reader is referred to~\cite{pog2008a} for more details on the construction of the sphere-in-a-box grid.

The added-damping tensors for this problem are computed based on a discrete surface integral on the two body-fitted grids at the beginning of the simulation.
%\dws{Qi: Lets discuss this bit.  I have made some minor changes but I'm not sure I have this quite right.}  
Formulas for the discrete added-damping tensors corresponding to the exact forms {in~\eqref{eq:ADvv} through~\eqref{eq:ADww}} for a spherical particle of radius $r$ are given in the appendix of~\cite{rbins2017} as
\begin{align*}
 & \Dvva = \frac{\mu}{\dn}\, \frac{8}{3} \pi r^2 \, 
      \begin{bmatrix} 1 & 0 & 0 \\ 0 & 1 & 0 \\ 0 & 0 & 1\end{bmatrix},  \quad 
   \Dvwa = \zerov,  \quad 
   \Dwwa =   \frac{\mu}{\dn}\, \frac{8}{3} \pi r^4 \, 
         \begin{bmatrix} 1 & 0 & 0 \\ 0 & 1 & 0 \\ 0 & 0 & 1 \end{bmatrix}. 
\end{align*}
where
\begin{equation}
\dn \eqdef \frac{\Delta s_{n}}{1-e\sp{-\delta}},\qquad \delta \eqdef \frac{\Delta s_{n}}{\sqrt{\nu\dt/2}},
\label{eq:dn}
\end{equation}
with $\Delta s_{n}$ being the mesh spacing in the normal direction, i.e.~$\Delta s_n\approx 0.71h^{(j)}$
for the present calculations. 
%See~\cite{rbins2017} for further details on the added-damping length scale $\dn$.
%\dws{I have added this last bit to define $\dn$.  Check if this is correct.}  
Here, $r=0.5$ so that the diagonal elements of the added-damping tensors are
\[
    \Dvva_{ii} = \frac8{3} \pi r^2\,\frac{\mu}{\dn} \approx 2.0944\,\frac{\mu}{\dn} , 
    \qquad \Dwwa_{ii} = \frac{8}{3}\pi r^4\,\frac{\mu}{\dn} \approx 0.52360\,\frac{\mu}{\dn} , \qquad i=1, 2, 3,
\]
which is in excellent agreement with the values computed by the discrete surface integral.  
For instance, the diagonal elements of the tensors for the grid $\Gcsp^{(4)}$ are computed in the beginning of the simulation as 
\[
\Dvva_{h,ii}\approx 2.0886\,\frac{\mu}{\dn} , 
   \qquad \Dwwa_{h,ii}\approx 0.52215\,\frac{\mu}{\dn} , \qquad i=1, 2, 3.
\]
The off-diagonal elements are zero for the exact tensors, and this is in excellent agreement with the computed tensors for which the maximum off-diagonal element is approximately $7.3\times10^{-6}$ in absolute value.

{% ------ MATLAB CURVES LIGHT BODY ------
\newcommand{\figWidth}{7.75cm}
\newcommand{\trimfig}[2]{\trimFig{#1}{#2}{.0}{.0}{.0}{.0}}
\newcommand{\figWidtha}{8.5cm}
\begin{figure}[htb]
\begin{center}
\resizebox{14cm}{!}{% START resize box
\begin{tikzpicture}[scale=1]
  \useasboundingbox (0.0,0) rectangle (16,6.3);  % set the bounding box (so we have less surrounding white space)
  \draw(-.5,-.75) node[anchor=south west,xshift=0pt,yshift=+0pt] {\trimfig{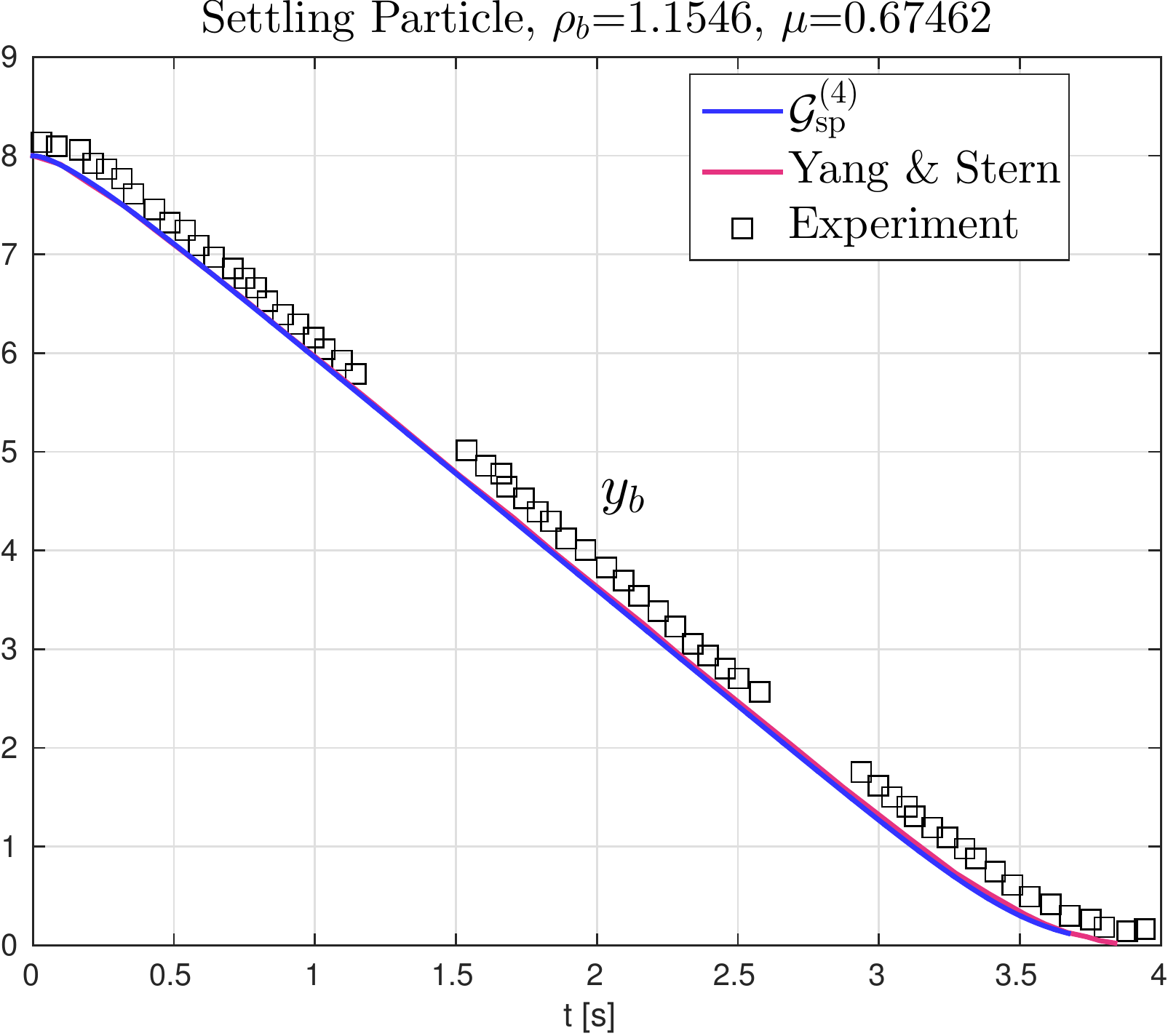}{\figWidth}};
  \draw(8.0,-.75) node[anchor=south west,xshift=0pt,yshift=+0pt] {\trimfig{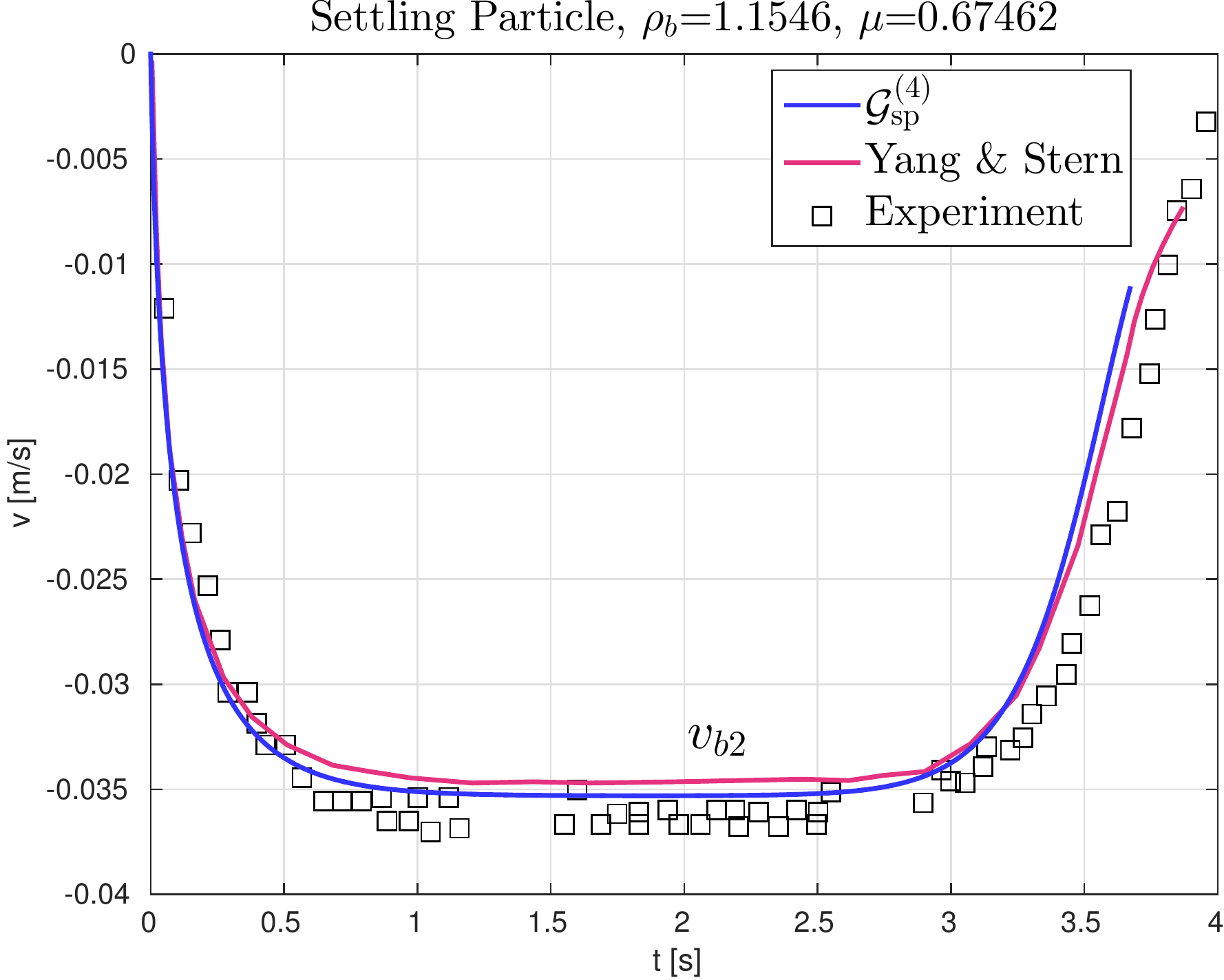}{\figWidtha}};
%
% grid:
%\draw[step=1cm,gray] (0,0) grid (16,6.3);
\end{tikzpicture}
} % END resize box
\end{center}
  \caption{Settling particle. Time history of the position (left) and velocity (right).
      The solutions are computed on grid $\Gc^{(4)}_{\rm sp}$.
      The black square is the experimental data taken from Fig.~5 in~\cite{Cate2002}
      and the red curve is the numerical solutions taken from Fig.~13 in~\cite{yang2015non}.
  }
  \label{fig:settlingParticleCompareCurves}
\end{figure}
}

{% ------ RISING RECTANGLE CONTOUR PLOTS  ------
% --------- START drawContour -----------
\newcommand{\drawContour}[7]{%
\begin{scope}[#1]
\draw(0.0,0) node[anchor=south west,xshift=-4pt,yshift=+0pt] {\trimfig{sphereDrop/#2}{\figWidth}};
  \draw(.5,5.4) node[draw,fill=white,anchor=west,xshift=2pt,yshift=1pt] {\scriptsize #3};
  % \draw(0,0.0) node[draw,fill=white,anchor=south west,xshift=12pt,yshift=12pt] {\scriptsize #4};
  \draw(1.2,5.4) node[draw,fill=white,anchor=west,xshift=2pt,yshift=1pt] {\scriptsize #5};
 % colour bar:
\begin{scope}[xshift=9pt,yshift=-1pt]
  \draw (\xcb,\ycb) node[anchor=south west,xshift=0.cm,yshift=.5cm,rotate=-90] {\trimfigcb{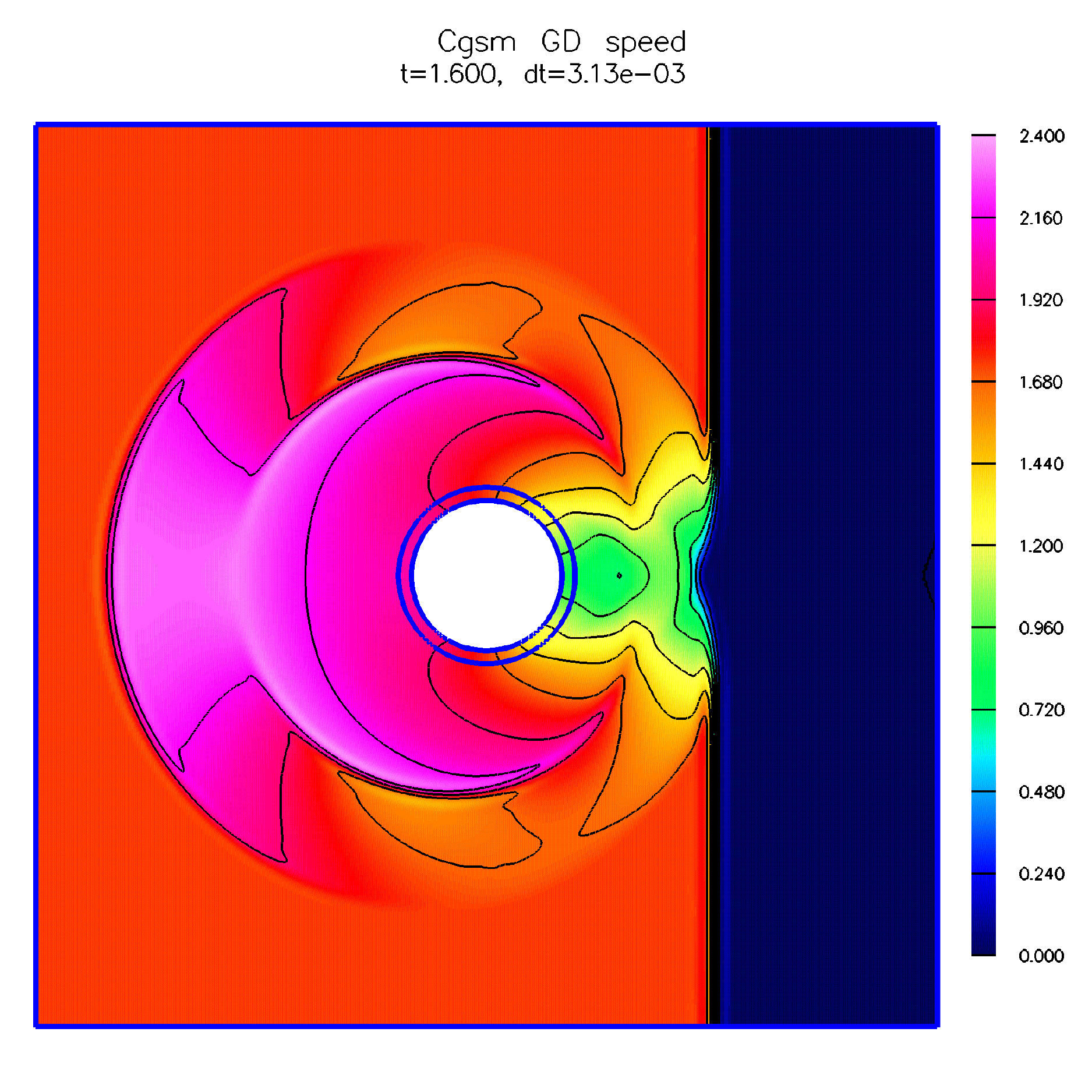}{\cbWidth}{\cbHeight}};
  \draw (.8,0) node[anchor=north,xshift=+3pt,yshift=+2pt] {\scriptsize $#6$};
  \draw (4.4,0) node[anchor=north,xshift=+0pt,yshift=+2pt] {\scriptsize $#7$};
\end{scope}
% \draw (\xcb,\ycb) node[anchor=south west,xshift= +8pt,yshift=+1pt] {\scriptsize $#6$};
% \draw (\xcb,\ycbTop) node[anchor=south west,xshift= +8pt,yshift=-6pt] {\scriptsize $#7$};
\end{scope}
}
% --------- END drawContour -----------
% -- for colour bar ----
\newcommand{\cbWidth}{.2cm}% colour bar width
\newcommand{\cbHeight}{4cm}% colour bar height
\newcommand{\xcb}{.5cm}% colour bar lower left corner
\newcommand{\ycb}{-.2cm}% colour bar lower left corner
\setlength{\ycbTop}{\ycb+\cbHeight}% colour bar top label position
\setlength{\ycbMid}{\ycb+\cbHeight*\real{.5}}% colour bar top label position
\newcommand{\trimfigcb}[3]{\includegraphics[width=#2, height=#3, clip, trim=17cm 2.35cm 1.65cm 2.35cm]{#1}}
%
% ========================= DRAW ================
%
\newcommand{\figWidth}{5cm}
\newcommand{\trimfig}[2]{\trimw{#1}{#2}{.12}{.12}{.0}{.14}}
\begin{figure}[htb]
\begin{center}
\resizebox{14cm}{!}{% START resize box
\begin{tikzpicture}[scale=1]
  \useasboundingbox (0.4,1) rectangle (16.,13);  % set the bounding box (so we have less surrounding white space)
%  top row: 
 \drawContour{xshift= 0.0cm,yshift=7cm}{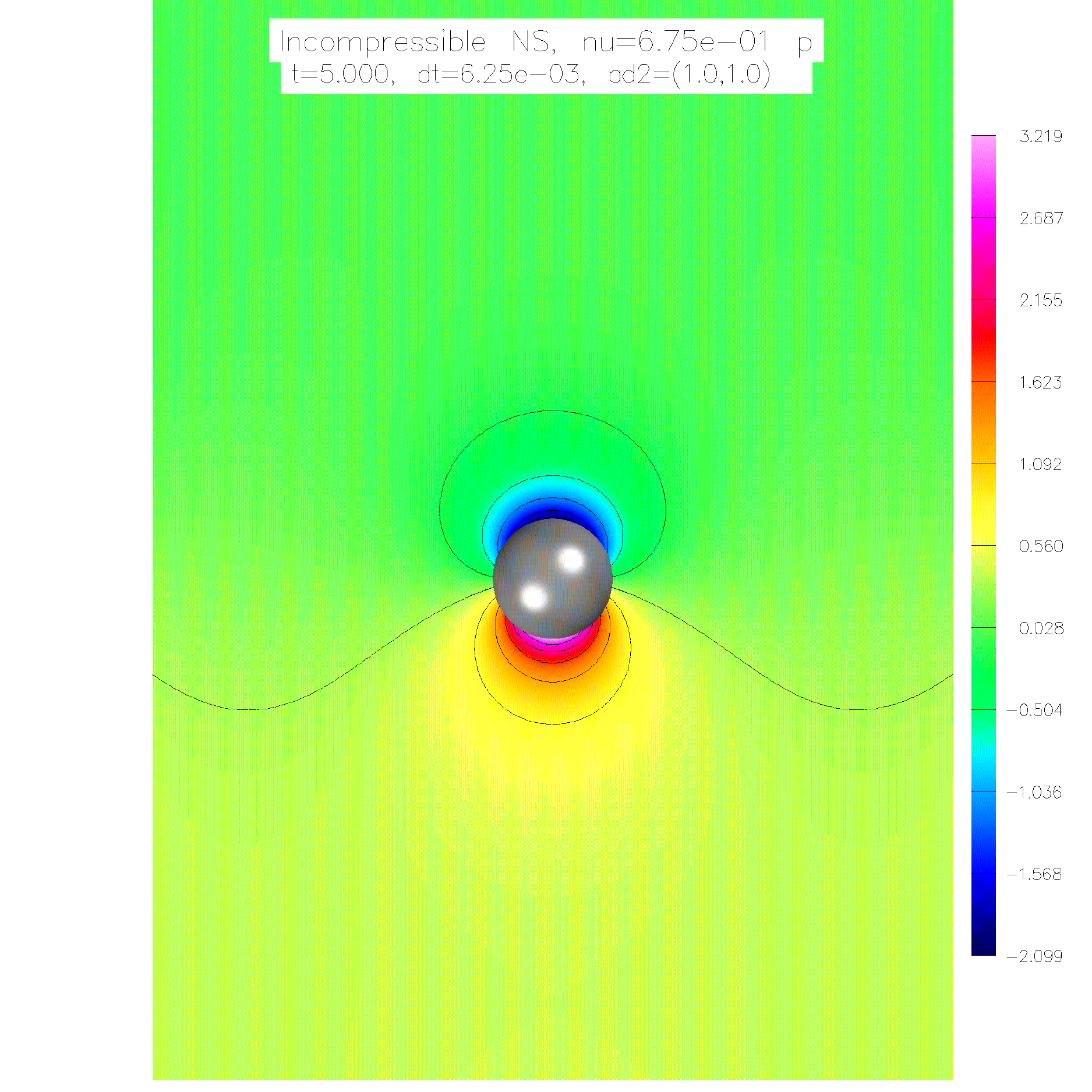}{$p$}{$p$}{$t=5.0$}{-2.01}{3.22};
 \drawContour{xshift=5.25cm,yshift=7cm}{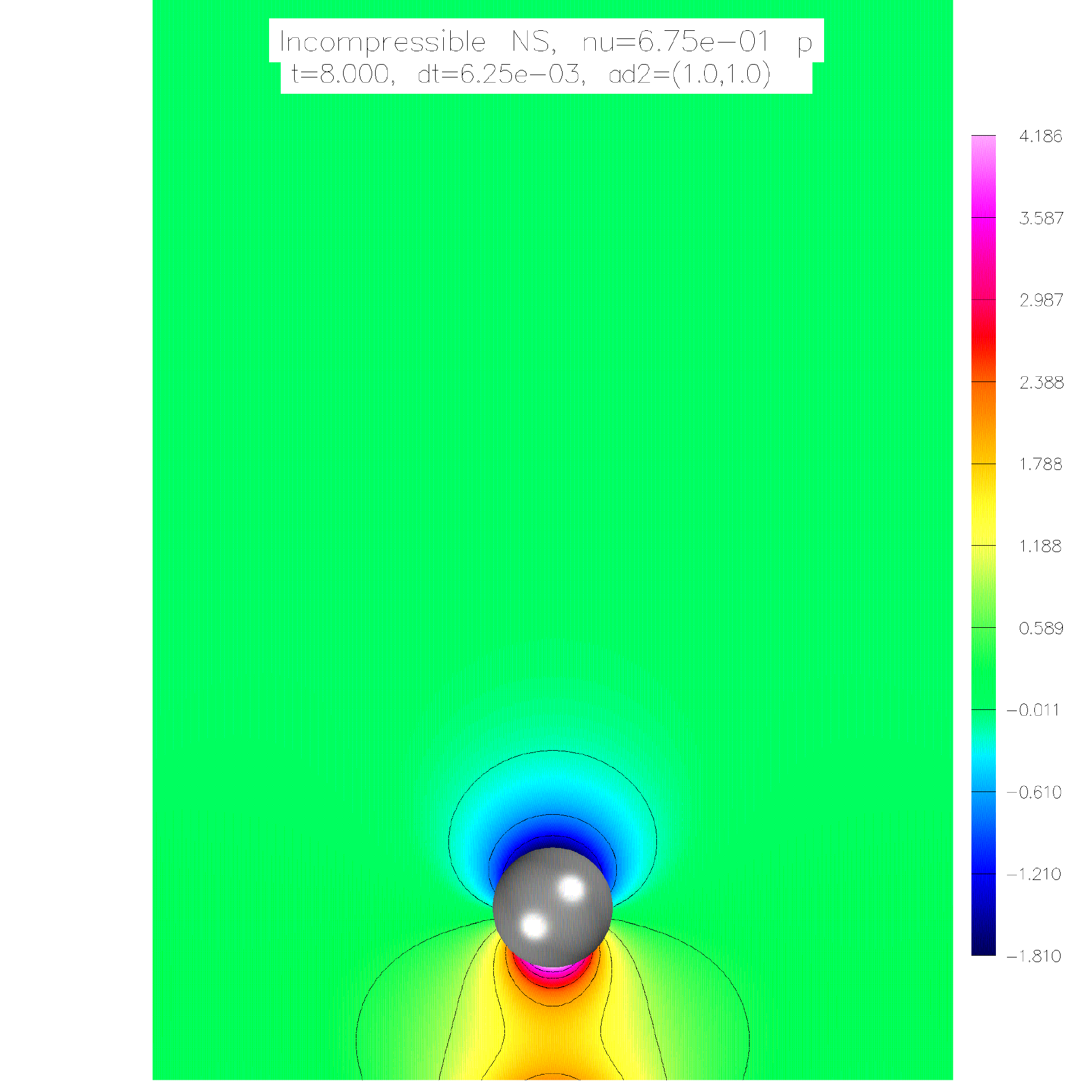}{$p$}{$p$}{$t=8.0$}{-1.81}{4.19};
 \drawContour{xshift=10.5cm,yshift=7cm}{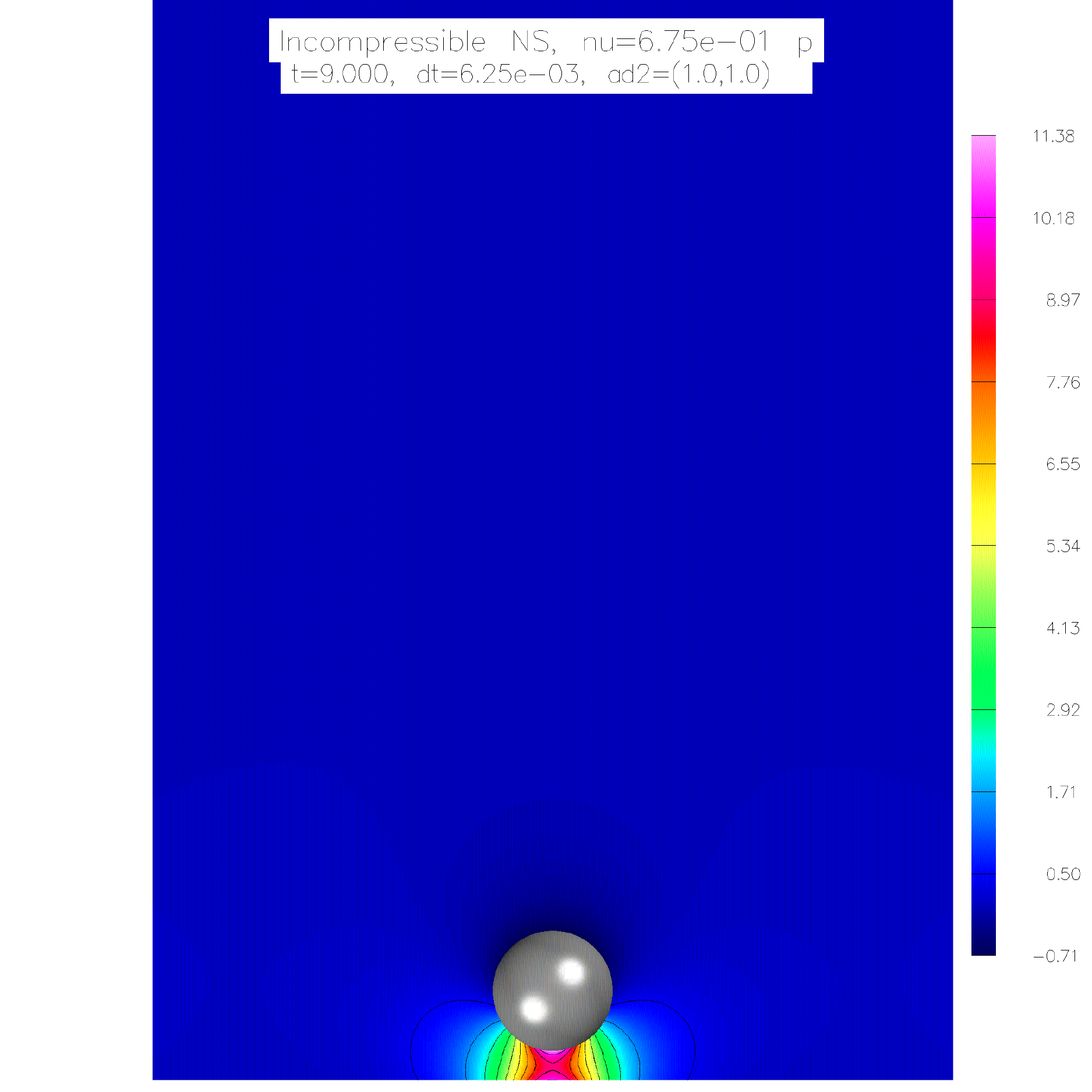}{$p$}{$p$}{$t=9.0$}{-0.71}{11.38};
%
% bottom row: 
 \drawContour{xshift= 0.0cm,yshift=0.4cm}{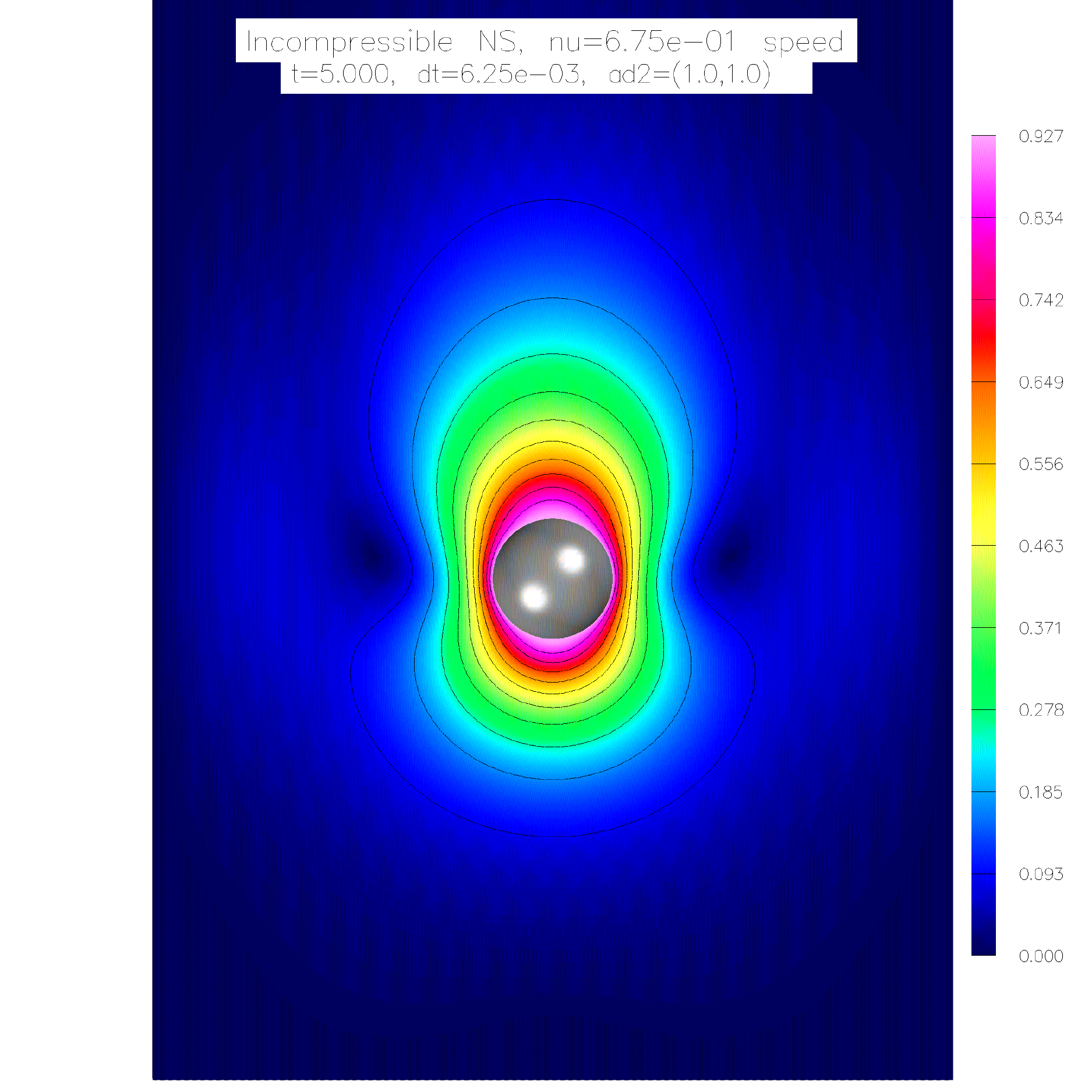}{$|\vv|$}{$|\vv|$}{$t=5.0$}{0}{.927};
 \drawContour{xshift=5.25cm,yshift=0.4cm}{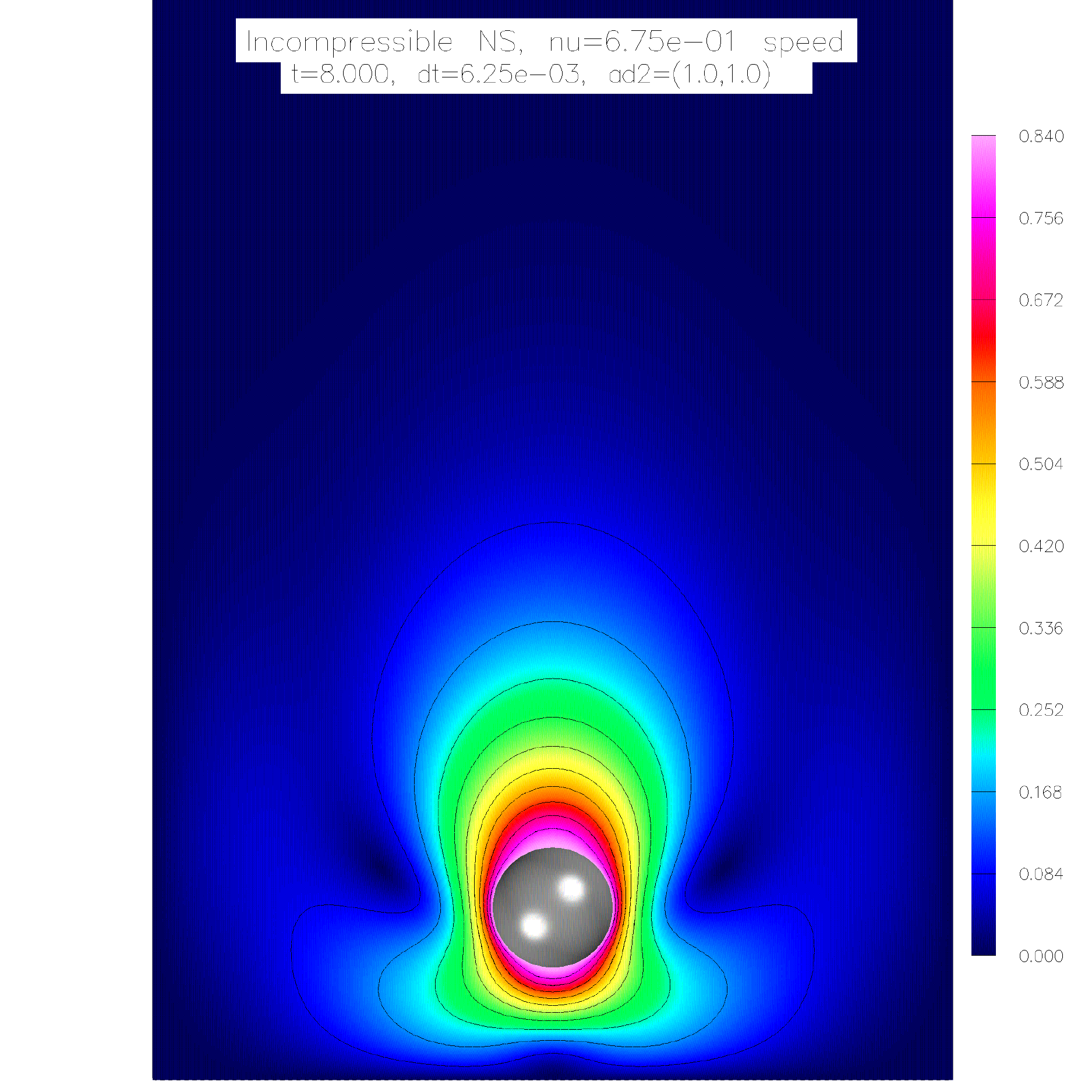}{$|\vv|$}{$|\vv|$}{$t=8.0$}{0}{.840};
 \drawContour{xshift=10.5cm,yshift=0.4cm}{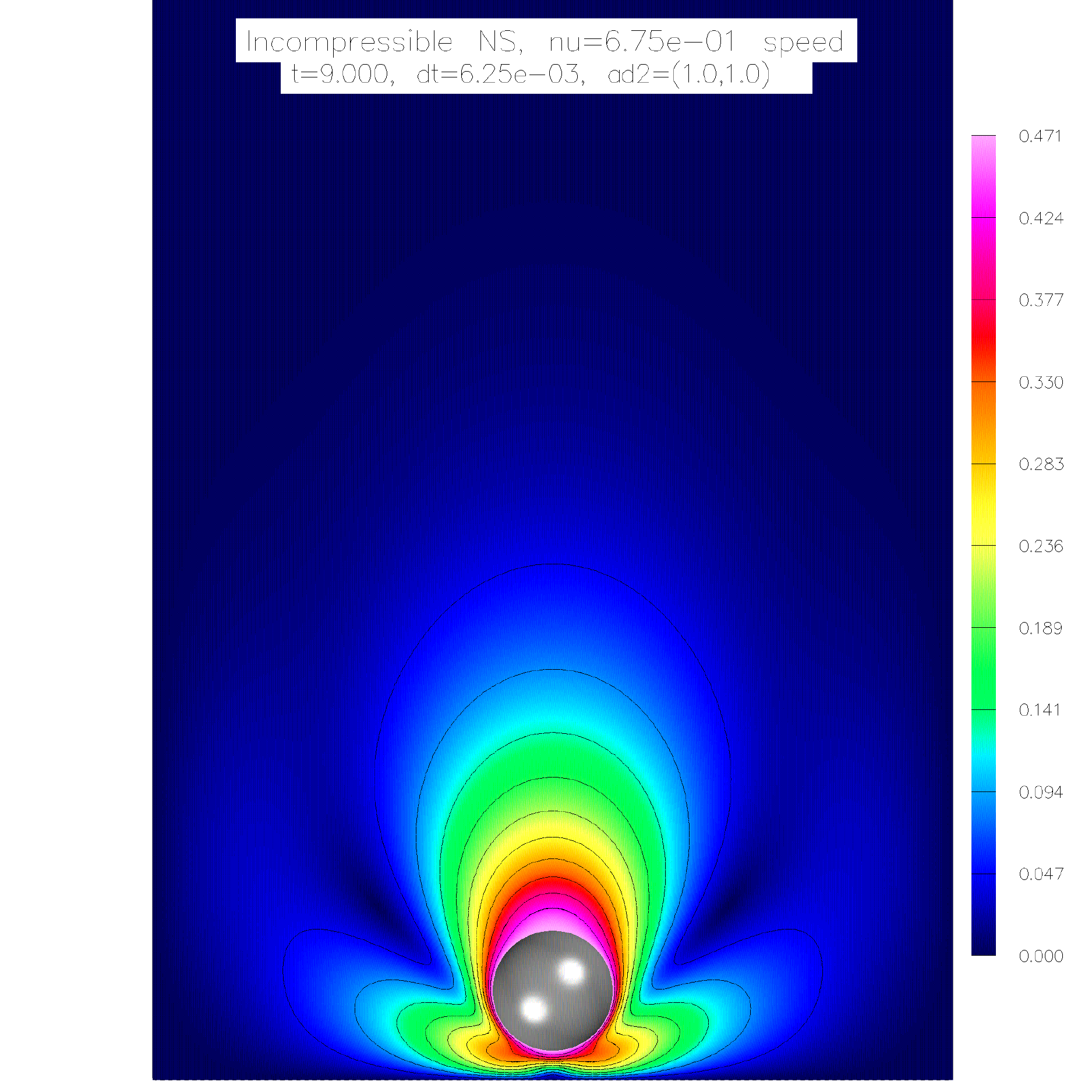}{$|\vv|$}{$|\vv|$}{$t=9.0$}{0}{.471};
%
% grid:
%\draw[step=1cm,gray] (0,0) grid (16,13);
\end{tikzpicture}
}% end resize box
\end{center}
  \caption{Settling particle. Contours of the pressure (top row) and speed (bottom row) at times $t=5$, $8$ and~$9$ computed
      using the composite grid $\Gc_{\rm sp}^{(4)}$. }
  \label{fig:settlingParticleContours}
\end{figure}
}

%This problem has been commonly used as a validation exercise and numerical solutions are compared to the experimental data from~\cite{Cate2002}.
Figure~\ref{fig:settlingParticleCompareCurves} presents the time history of the vertical position, $y_b(t)$, and vertical velocity, $v_{b2}(t)$, of the body computed using the~\ampRB~scheme with the composite grid $\Gcsp^{(4)}$.
To compare with the available data, the time and velocity in Figure~\ref{fig:settlingParticleCompareCurves} have been scaled properly so that their units are ${\rm s}$ and ${\rm m/s}$, respectively.
The present results are in good agreement with the experimental data in~\cite{Cate2002} and the numerical results in~\cite{gibou2012efficient, kempe2012improved, yang2015non, Koblitz2016}.
The numerical solutions from~\cite{yang2015non}, for example, are also plotted in Figure~\ref{fig:settlingParticleCompareCurves} for comparison.  
The numerical solutions are taken directly from Fig.~5 in~\cite{yang2015non}, although the velocity therein was 
plotted in a larger scale, which may lead to potential errors in reproducing the data.
Despite these potential errors, the numerical solution of the~\ampRB~scheme appears to be in good agreement with the results in~\cite{yang2015non}.
%\dws{shouldn't something be said about the comparison with~\cite{yang2015non} in the figure before mentioning a comparison in the next sentence that is not shown?}
%In particular, 
It is also found that the results of the~\ampRB~scheme are almost indistinguishable with the results given by the TP-RB scheme in~\cite{Koblitz2016}.
The contours of the pressure and speed are presented in Figure~\ref{fig:settlingParticleContours} for three different times.
The body accelerates initially due to its negative buoyancy which is turned on impulsively at $t=0$.  It then approaches a steady velocity before the body eventually slows down as it approaches the bottom wall.
As in the previous work~\cite{yang2015non, Koblitz2016}, collisions between the particle and the wall are not considered.
Therefore, the simulation is stopped when there are insufficient grid lines in the gap between the particle and the bottom wall.
Note that the other three cases from~\cite{Koblitz2016} with difference Reynolds numbers have also been computed using the~\ampRB~scheme, and it has been found that
 the results are almost identical to the TP-RB scheme if the same grids are used.
 %\dws{Again there is mention of the results in~\cite{Koblitz2016} and again no comparisons are presented.}

{% ------ MATLAB CURVES : TWO BODIES MOTION  ------
\newcommand{\figWidth}{8cm}
\newcommand{\figWidthb}{7.8cm}
\newcommand{\figWidtha}{6cm}
\newcommand{\figWidthc}{8.1cm}
\newcommand{\trimfig}[2]{\trimFig{#1}{#2}{.0}{.0}{.0}{.0}}
\begin{figure}[htb]
\begin{center}
\resizebox{14cm}{!}{% START resize box
\begin{tikzpicture}[scale=1]
  \useasboundingbox (0.0,0) rectangle (16,13.6);  % set the bounding box (so we have less surrounding white space)

  \draw(-.3,6.4) node[anchor=south west,xshift=0pt,yshift=+0pt] {\trimfig{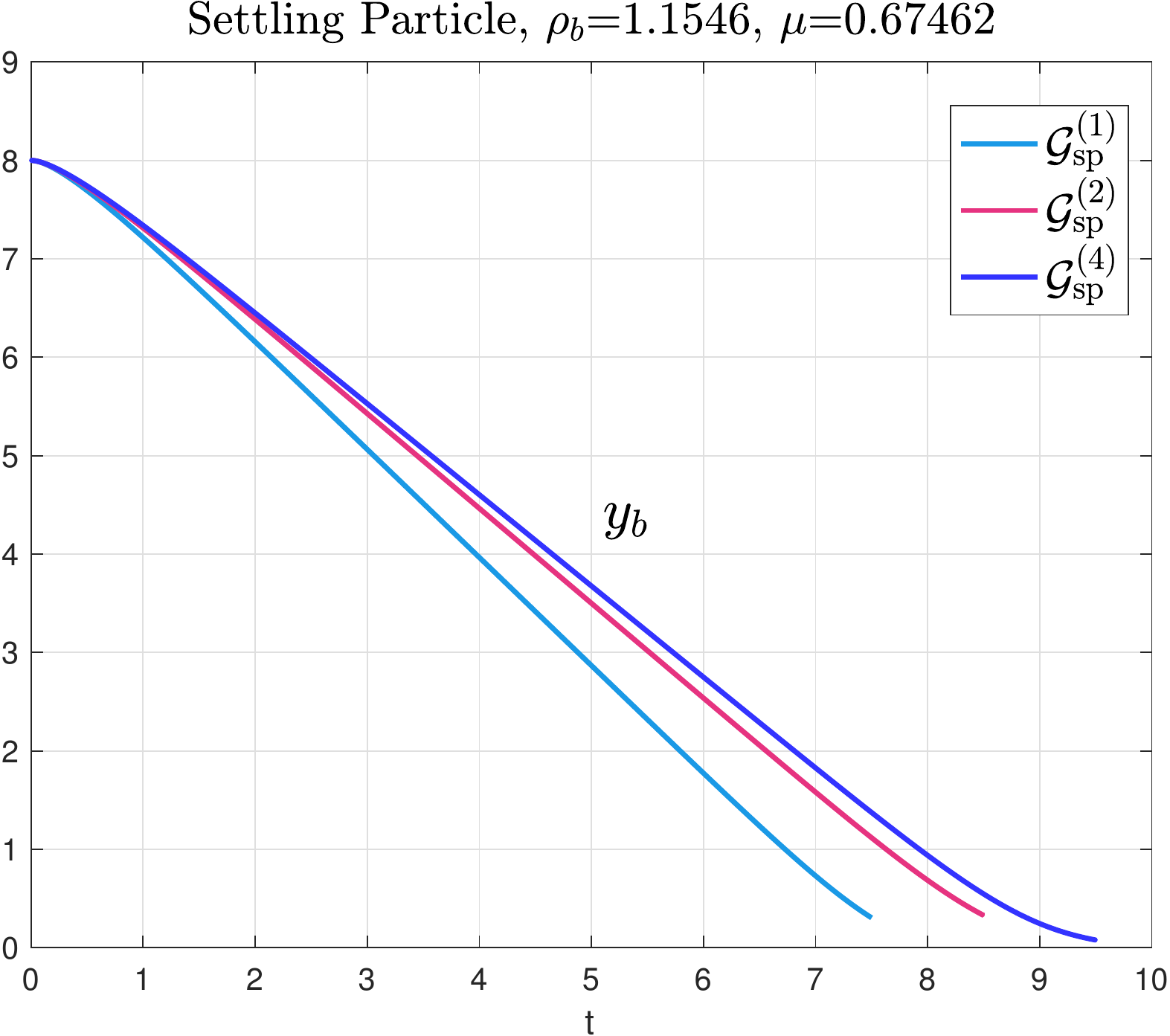}{\figWidthb}};
  \draw(8.0,6.4) node[anchor=south west,xshift=0pt,yshift=+0pt]   {\trimfig{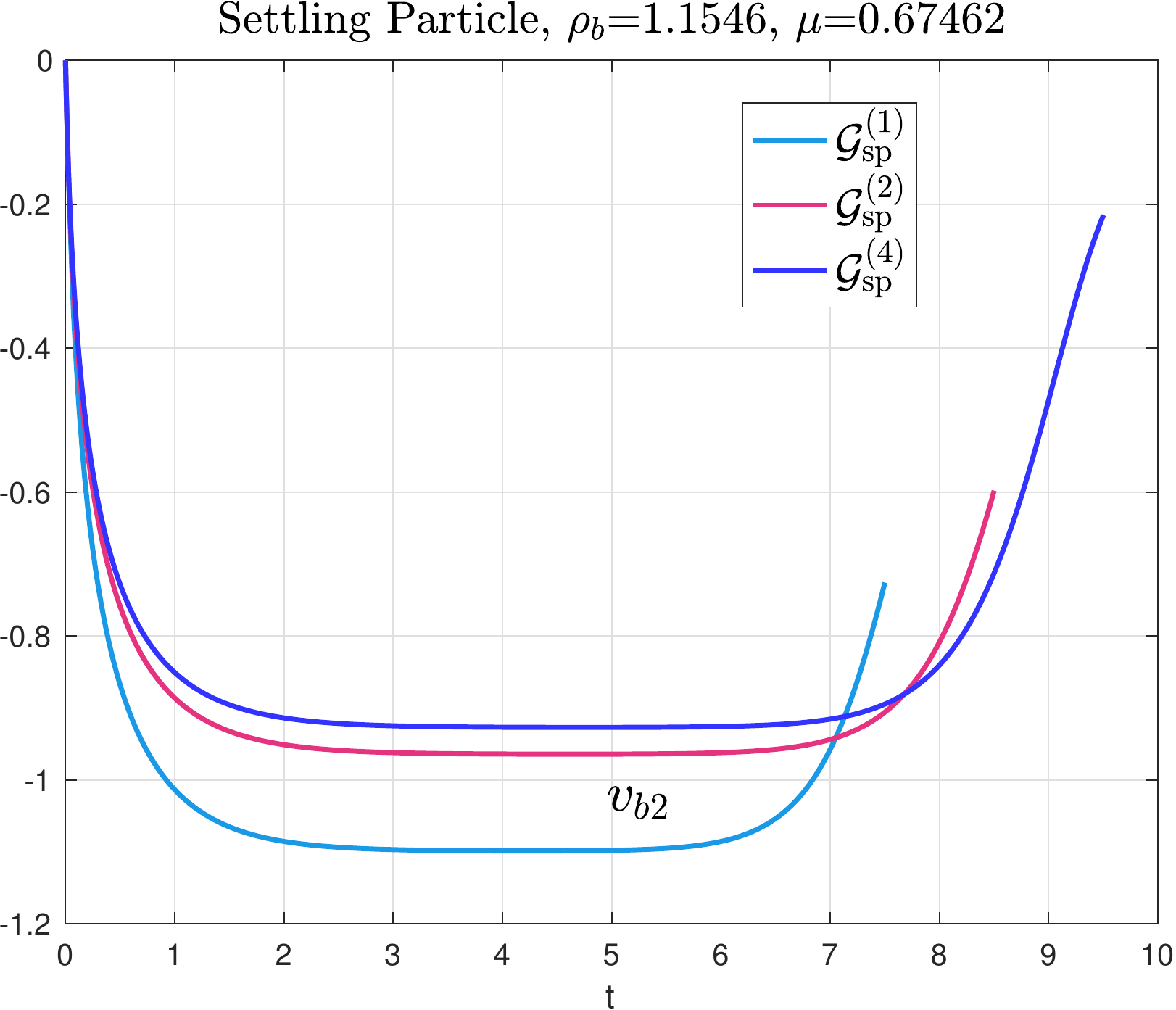}{\figWidthc}};
  \draw(-.5,-.75) node[anchor=south west,xshift=0pt,yshift=+0pt] {\trimfig{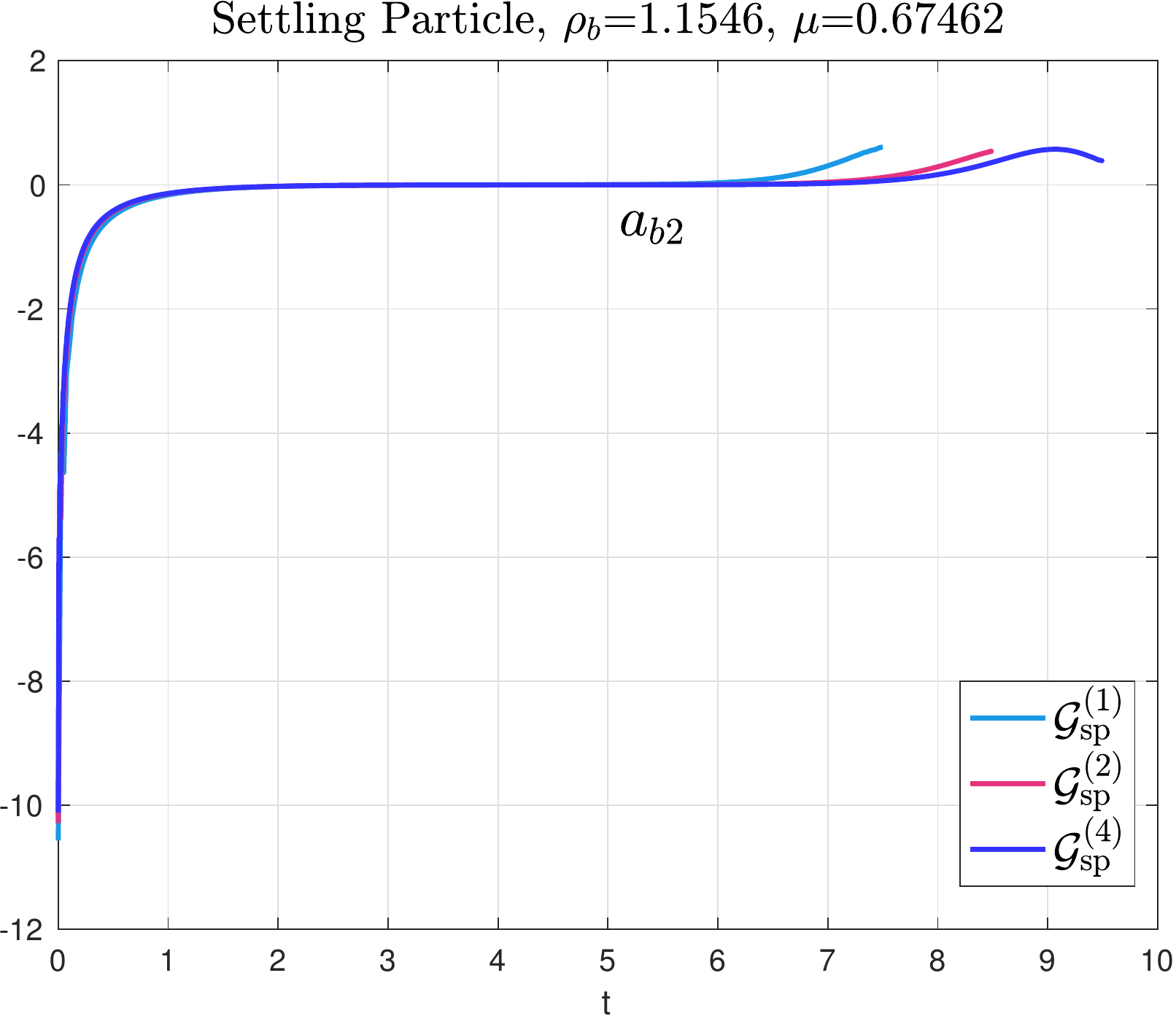}{\figWidth}};
  \draw(9.0,.25) node[anchor=south west,xshift=0pt,yshift=+0pt]{\trimfig{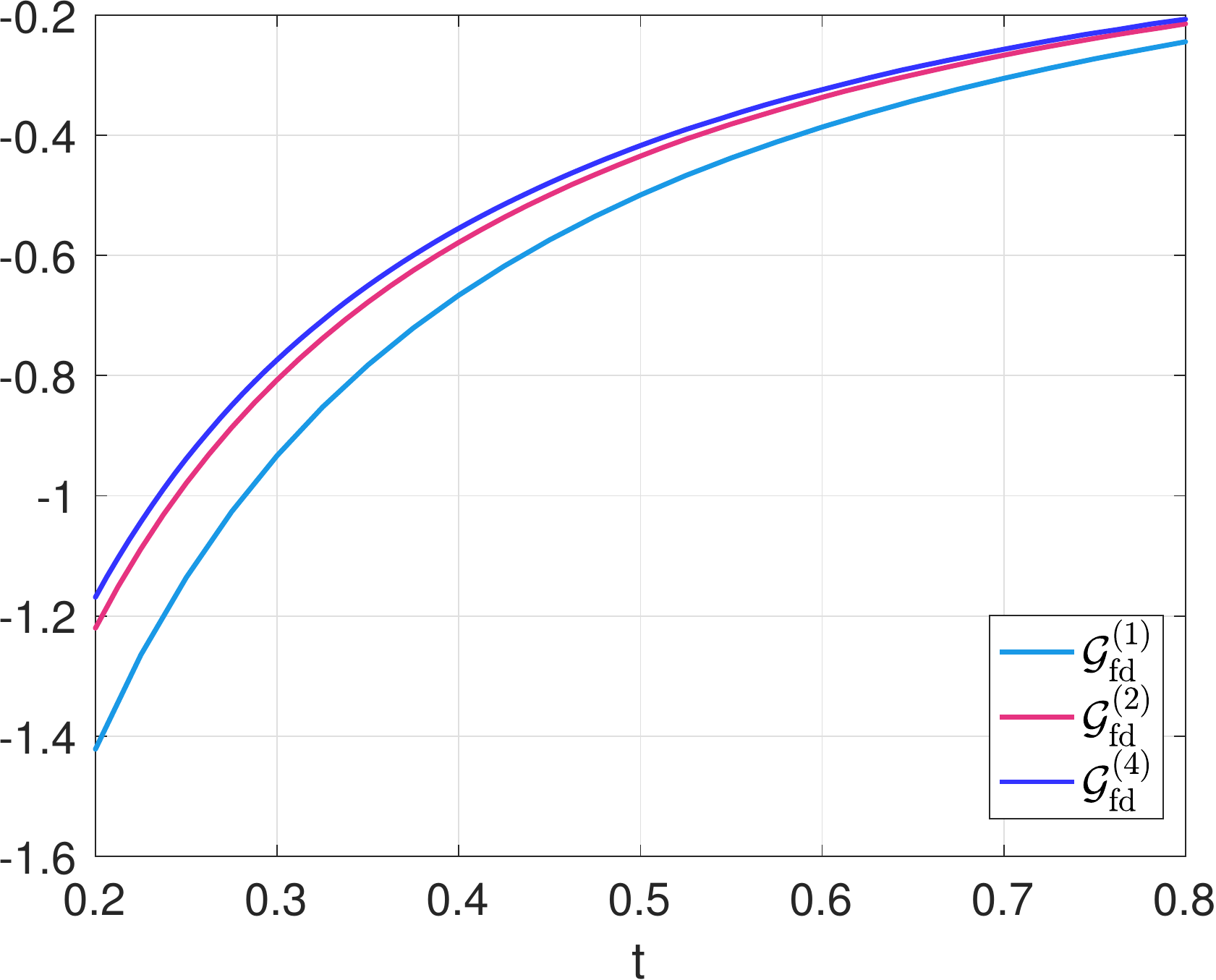}{\figWidtha}};
% grid:
%\draw[step=1cm,gray] (0,0) grid (16,13.6);
\end{tikzpicture}
}
\end{center}
  \caption{Settling particle. Time history of the position (top left), velocity (top right), acceleration (bottom left) and its zoomed view of the particle.  
     }
  \label{fig:settlingParticleCurves}
\end{figure}
}

{% ------ RISING RECTANGLE CONTOUR PLOTS  ------
% --------- START drawContour -----------
\newcommand{\drawContour}[8]{%
\begin{scope}[#1]
\draw(0.0,0) node[anchor=south west,xshift=-4pt,yshift=+0pt] {\trimfig{sphereDrop/#2}{\figWidth}};
  \draw(.5,5.4) node[draw,fill=white,anchor=west,xshift=2pt,yshift=1pt] {\scriptsize #3};
  % \draw(0,0.0) node[draw,fill=white,anchor=south west,xshift=12pt,yshift=12pt] {\scriptsize #4};
  \draw(1.2,5.4) node[draw,fill=white,anchor=west,xshift=2pt,yshift=1pt] {\scriptsize #5};
  \draw(4.4,5.4) node[draw,fill=white,anchor=west,xshift=2pt,yshift=1pt] {\scriptsize #8};
 % colour bar:
\begin{scope}[xshift=9pt,yshift=-1pt]
  \draw (\xcb,\ycb) node[anchor=south west,xshift=0.cm,yshift=.5cm,rotate=-90] {\trimfigcb{fig/colourBarLines}{\cbWidth}{\cbHeight}};
  \draw (.8,0) node[anchor=north,xshift=+3pt,yshift=+2pt] {\scriptsize $#6$};
  \draw (4.4,0) node[anchor=north,xshift=+0pt,yshift=+2pt] {\scriptsize $#7$};
\end{scope}
% \draw (\xcb,\ycb) node[anchor=south west,xshift= +8pt,yshift=+1pt] {\scriptsize $#6$};
% \draw (\xcb,\ycbTop) node[anchor=south west,xshift= +8pt,yshift=-6pt] {\scriptsize $#7$};
\end{scope}
}
% --------- END drawContour -----------
% -- for colour bar ----
\newcommand{\cbWidth}{.2cm}% colour bar width
\newcommand{\cbHeight}{4cm}% colour bar height
\newcommand{\xcb}{.5cm}% colour bar lower left corner
\newcommand{\ycb}{-.2cm}% colour bar lower left corner
\setlength{\ycbTop}{\ycb+\cbHeight}% colour bar top label position
\setlength{\ycbMid}{\ycb+\cbHeight*\real{.5}}% colour bar top label position
\newcommand{\trimfigcb}[3]{\includegraphics[width=#2, height=#3, clip, trim=17cm 2.35cm 1.65cm 2.35cm]{#1}}
%
% ========================= DRAW ================
%
\newcommand{\figWidth}{5cm}
\newcommand{\trimfig}[2]{\trimw{#1}{#2}{.12}{.12}{.0}{.14}}
\begin{figure}[htb]
\begin{center}
\resizebox{14cm}{!}{% START resize box
\begin{tikzpicture}[scale=1]
  \useasboundingbox (0.4,1) rectangle (16.,13);  % set the bounding box (so we have less surrounding white space)
%  top row: 
  \drawContour{xshift= 0.0cm,yshift=7cm}{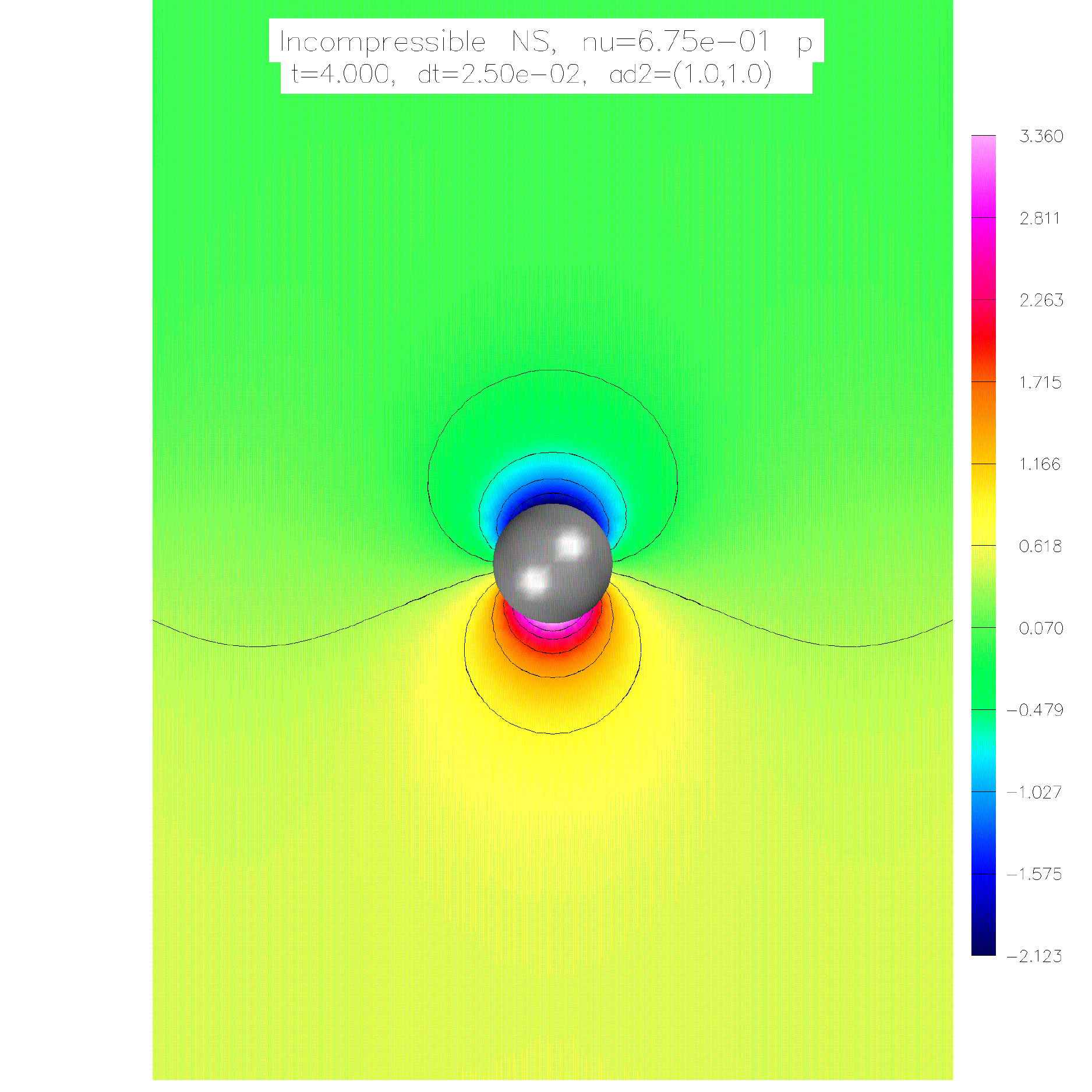}{$p$}{$p$}{$t=4$}{-2.12}{3.36}{$\Gc_{\rm sp}^{(1)}$};
  \drawContour{xshift=5.25cm,yshift=7cm}{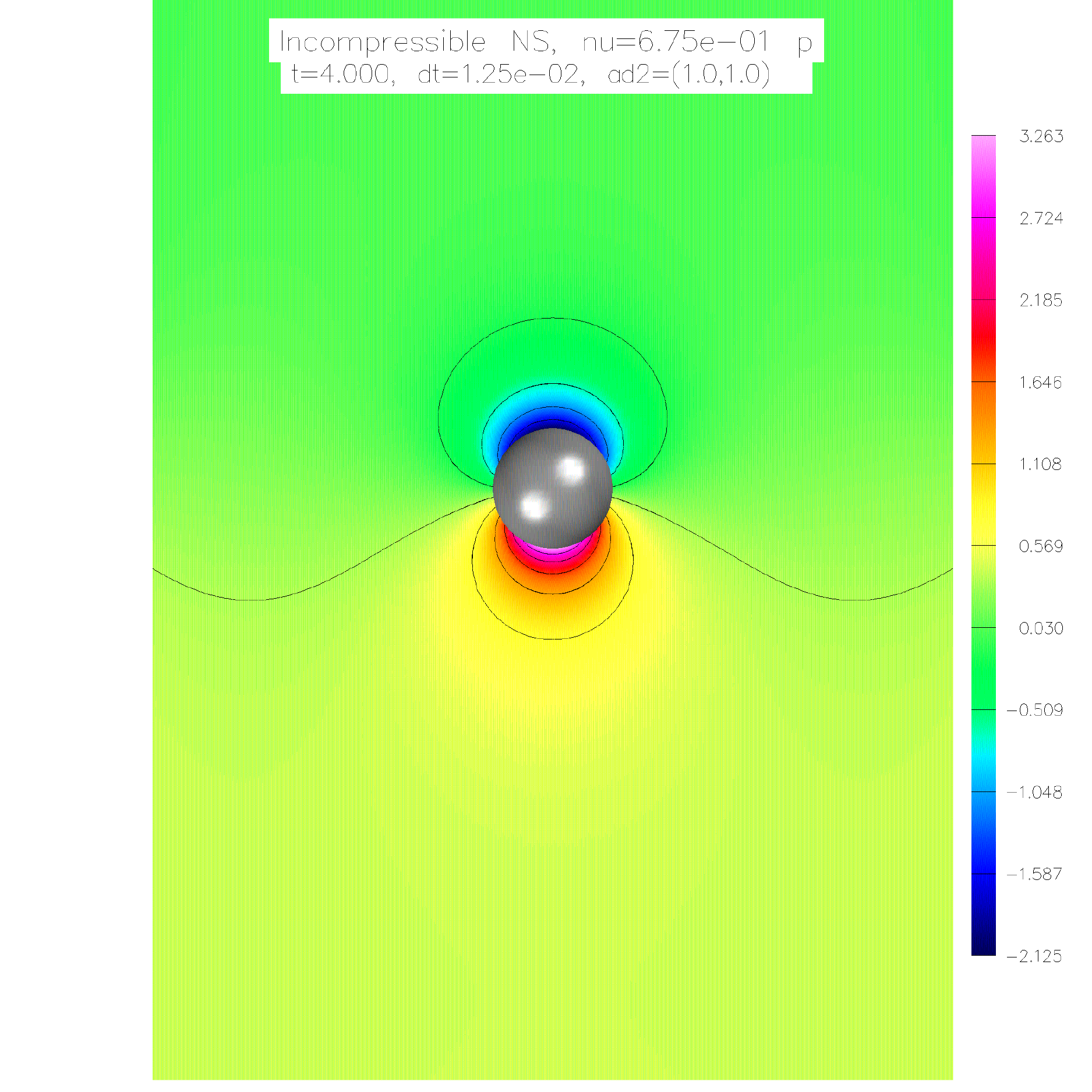}{$p$}{$p$}{$t=4$}{-2.12}{3.26}{$\Gc_{\rm sp}^{(2)}$};
  \drawContour{xshift=10.5cm,yshift=7cm}{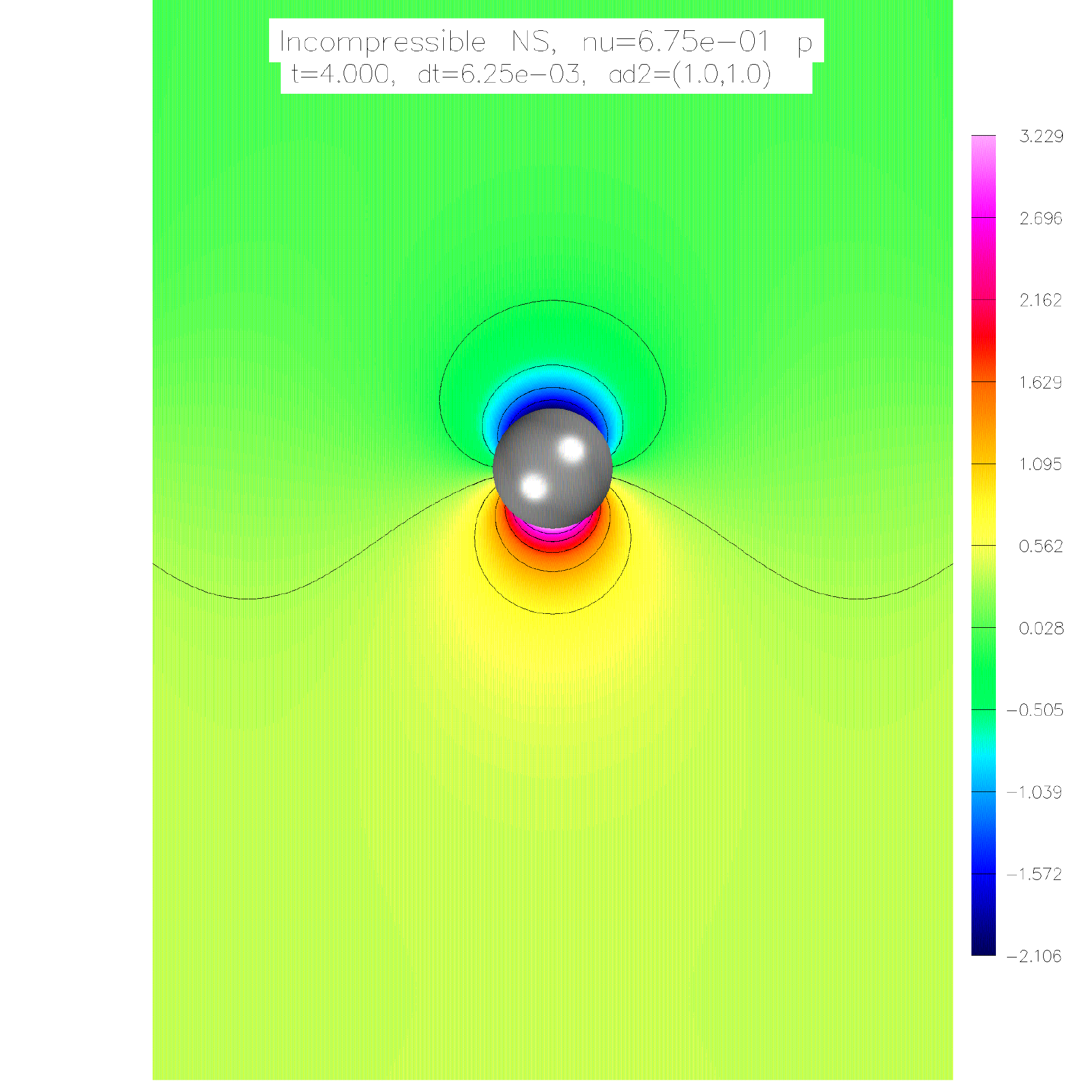}{$p$}{$p$}{$t=4$}{-2.11}{3.23}{$\Gc_{\rm sp}^{(4)}$};
%
% bottom row: 
 \drawContour{xshift= 0.0cm,yshift=0.4cm}{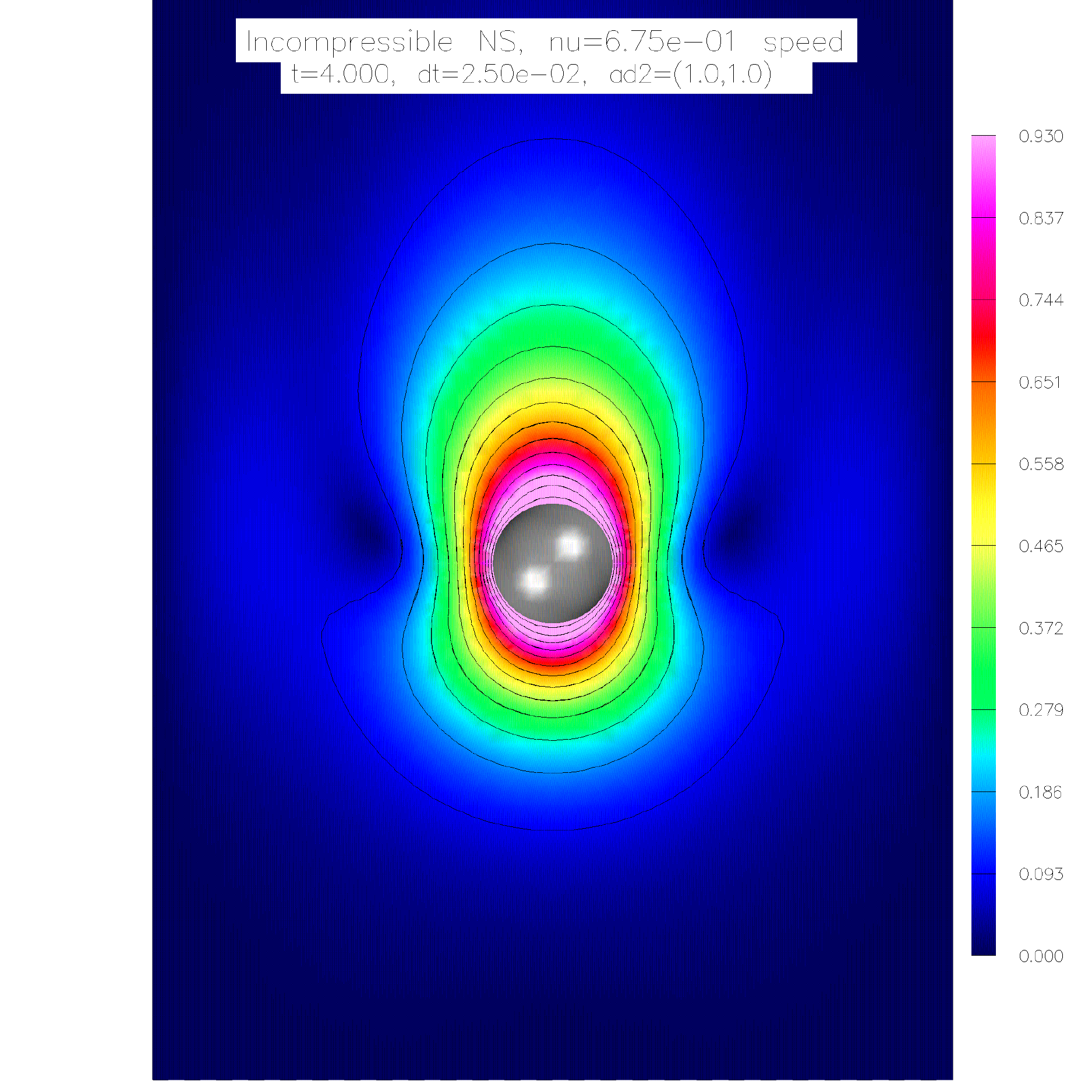}{$|\vv|$}{$|\vv|$}{$t=4$}{0}{.93}{$\Gc_{\rm sp}^{(1)}$};
 \drawContour{xshift=5.25cm,yshift=0.4cm}{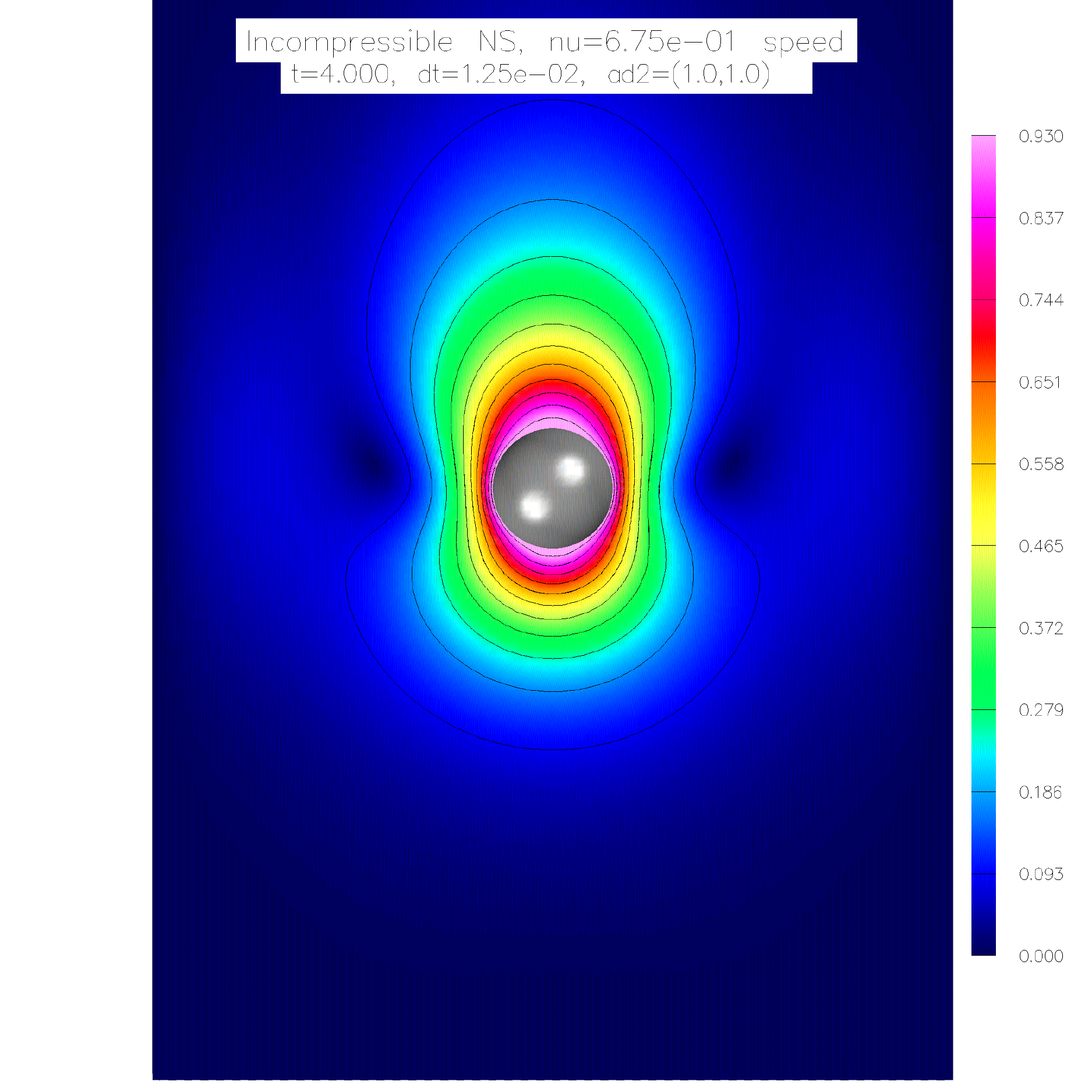}{$|\vv|$}{$|\vv|$}{$t=4$}{0}{.93}{$\Gc_{\rm sp}^{(2)}$};
 \drawContour{xshift=10.5cm,yshift=0.4cm}{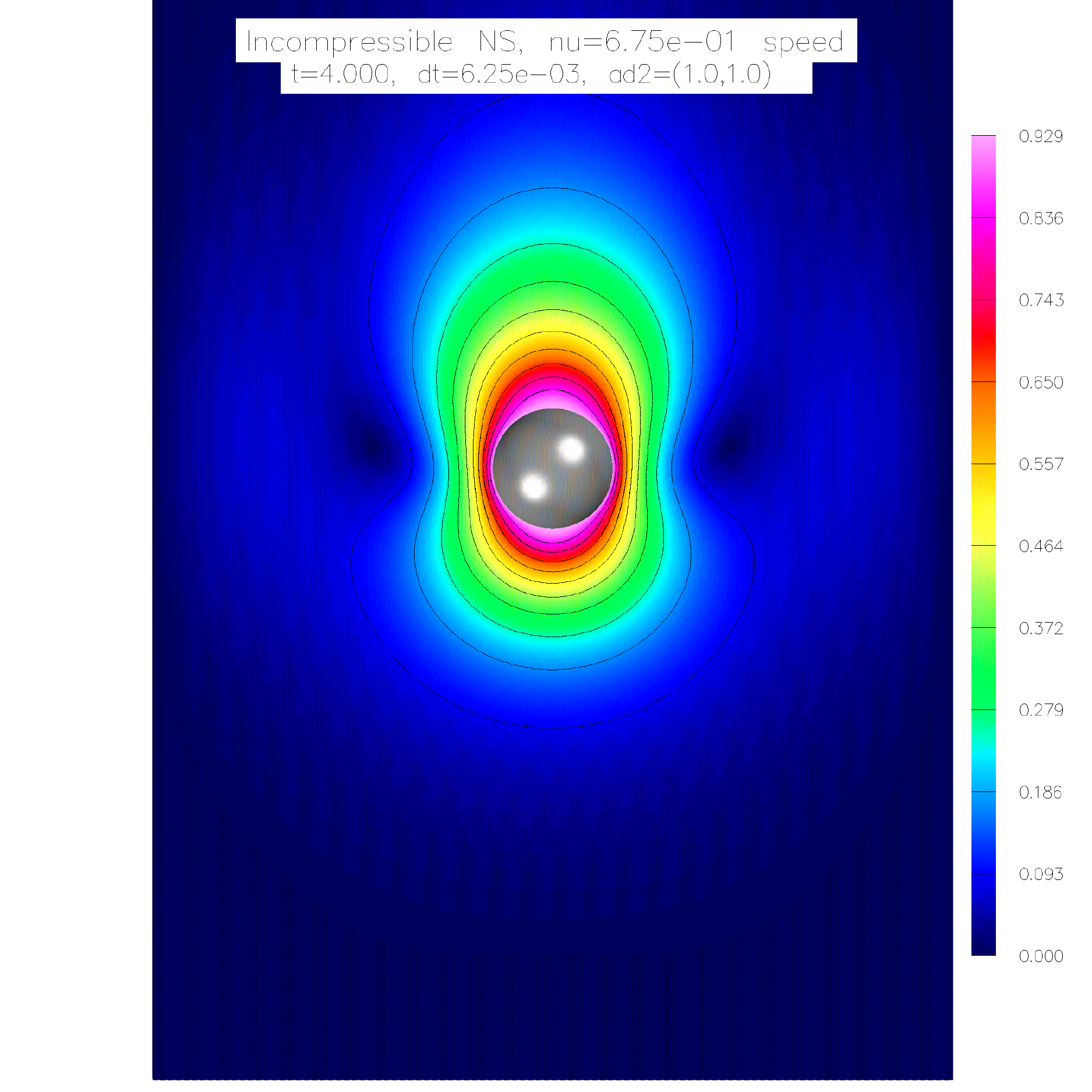}{$|\vv|$}{$|\vv|$}{$t=4$}{0}{.93}{$\Gc_{\rm sp}^{(4)}$};
%
% grid:
%\draw[step=1cm,gray] (0,0) grid (16,13);
\end{tikzpicture}
}% end resize box
\end{center}
  \caption{Settling particle. Contours of the pressure (top row) and speed (bottom row) at times $t=4$ computed
      using the composite grids  $\Gc_{\rm sp}^{(1)}$,  $\Gc_{\rm sp}^{(2)}$ and $\Gc_{\rm sp}^{(4)}$. }
  \label{fig:settlingParticleRefine}
\end{figure}
}

Figure~\ref{fig:settlingParticleCurves} presents the behaviour of the vertical position, velocity and acceleration of the particle versus time for different grid resolutions.
This figure is complemented by Figure~\ref{fig:settlingParticleRefine} that shows the fluid quantities and body positions at $t=4$ for different grid resolutions.
The numerical solutions in Figure~\ref{fig:settlingParticleCurves} are very smooth and show good convergence as the grid resolution increases.
In fact, Richardson extrapolation estimates for the convergence rate of the position and velocity are both around 1.90 at the times $t= 2$, 4 and 6. 
However, it is found that the resulting particle tends to drop faster at lower resolutions,
and so the simulations are halted earlier (e.g., the calculation for grid $\Gcsp^{(1)}$ is stopped at $t \approx 7.5$). 
%Therefore, the simulations are truncated even earlier since there are fewer resolutions between the particle and wall when the particle is settling. 
This makes the self-convergence studies challenging at the later times. 
Self-convergence studies of the accelerations are not as clean and  yield a convergence rate of approximately 1.52 at $t=2$ and 4, while the rate improves slightly to 1.65 at $t=6$.
There are at least  two  plausible reasons leading to this reduced convergence rate. 
One reason is that a very large gravity is turned on impulsively in this test.
Due to a sudden change of the body acceleration, the first time step may involve a first-order error, and in fact we observe that the accelerations oscillate slighlty during the first few time steps.
%which indicates this first order error in the first predictor is relatively large thanks to the impulsive start.
The second reason is that the problem becomes steady during most of the simulation. 
This makes it more challenging for the iterative Krylov solvers used here, since the solvers may take less Krylov iterations than required to obtain the designed order of accuracy.
Note that  all of the tests we considered in the previous work~\cite{rbins2017} were solved by a direct sparse matrix solver,
which eliminates this potential issue.

{% ------ MATLAB CURVES UNSTABLE LIGHT BODY ------
\newcommand{\figWidth}{8.cm}
\newcommand{\trimfig}[2]{\trimFig{#1}{#2}{.0}{.0}{.0}{.0}}
\newcommand{\figWidthz}{6cm}
\newcommand{\figWidtha}{6cm}
\newcommand{\trimfiga}[2]{\trimFig{#1}{#2}{.25}{.25}{.25}{.25}}
\begin{figure}[htb]
\begin{center}
\resizebox{16.5cm}{!}{% START resize box
\begin{tikzpicture}[scale=1]
  \useasboundingbox (0.0,1) rectangle (21.,7);  % set the bounding box (so we have less surrounding white space)
  \draw(0.0,0) node[anchor=south west,xshift=-4pt,yshift=+0pt] {\trimfig{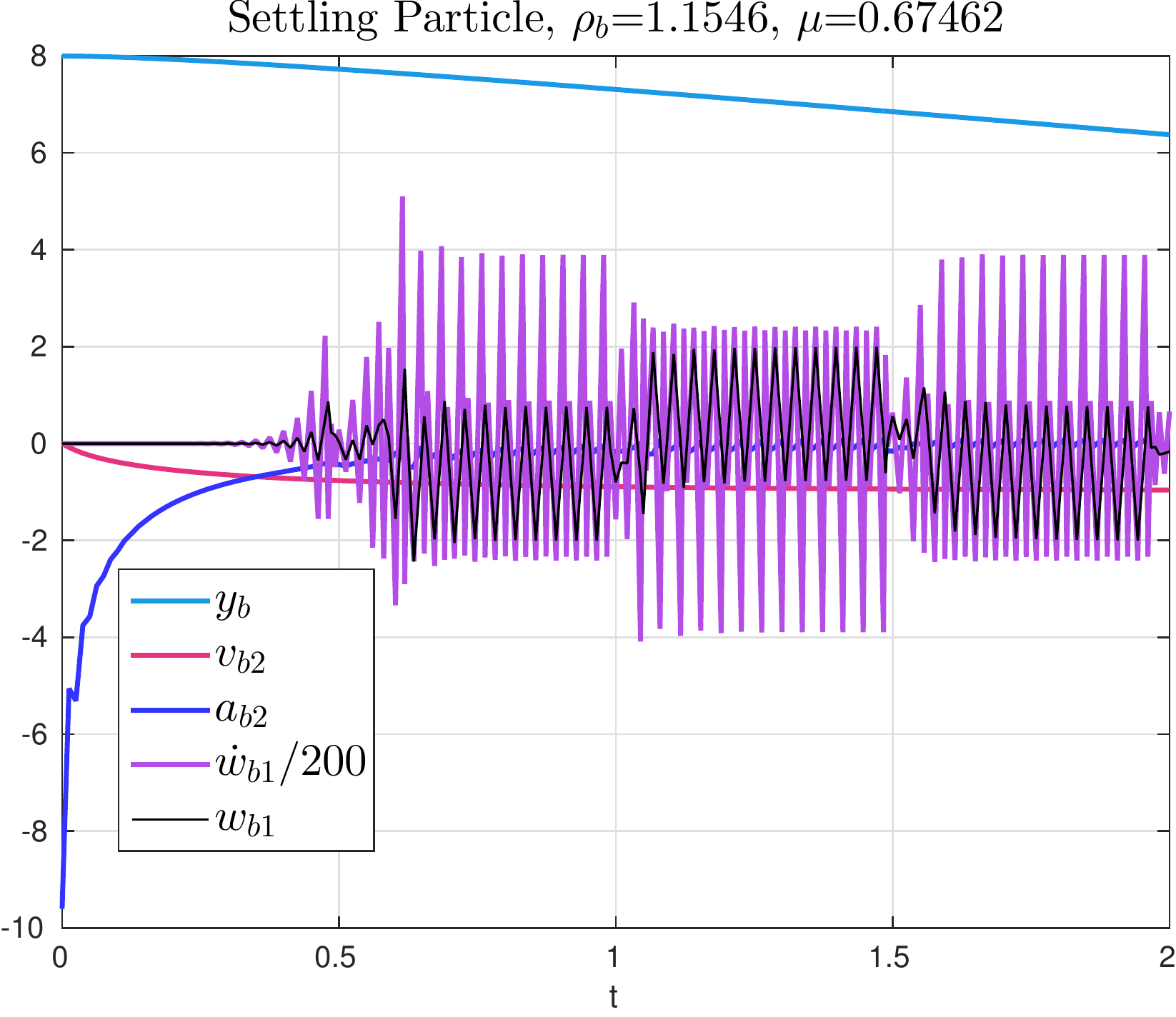}{\figWidth}};
  \draw(8.2,.5) node[anchor=south west,xshift=-4pt,yshift=+0pt] {\trimfig{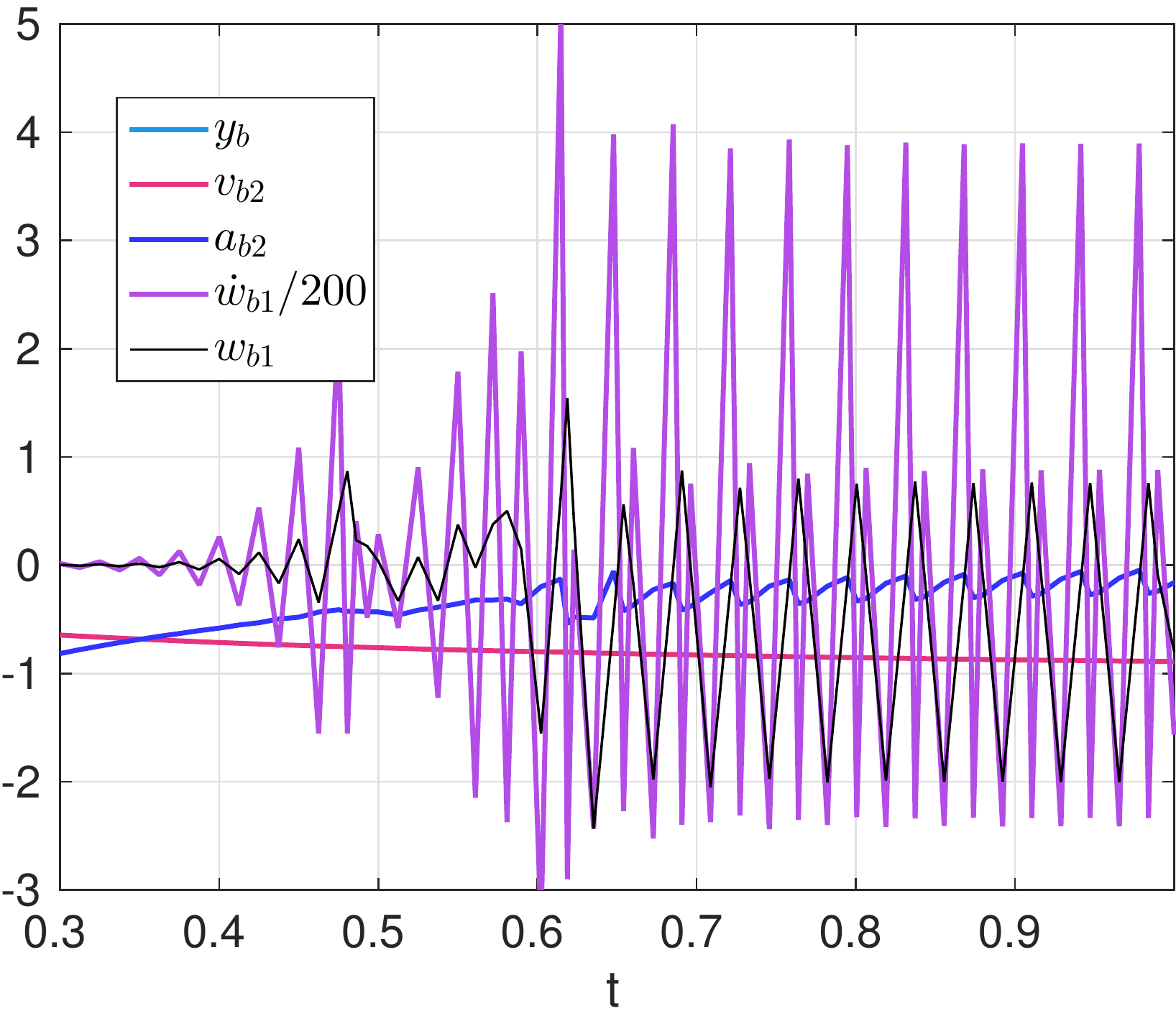}{\figWidthz}};
  \draw(15,.5) node[anchor=south west,xshift=-4pt,yshift=+0pt] {\trimfiga{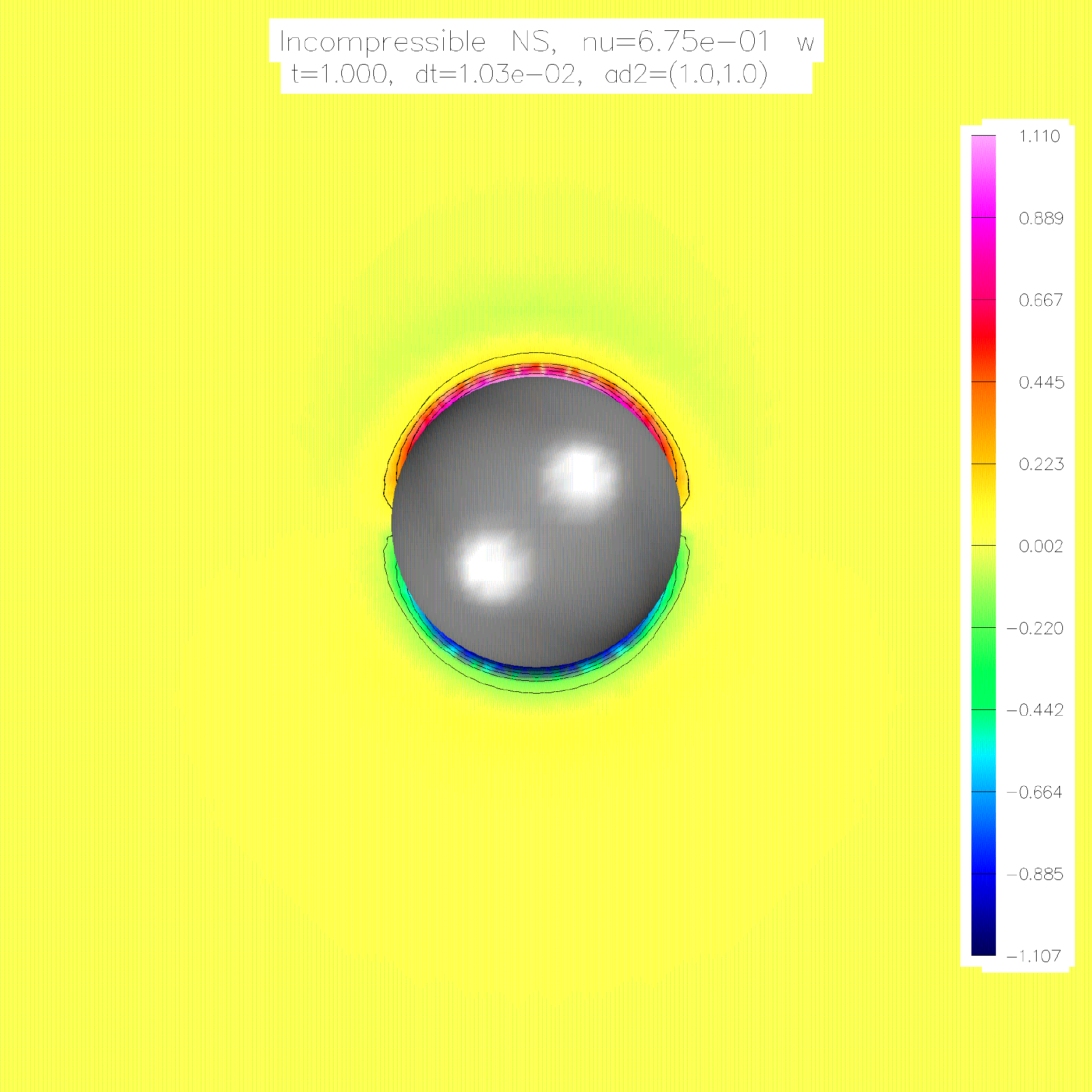}{\figWidtha}};
  \draw(15,6.2) node[draw,fill=white,anchor=west,xshift=2pt,yshift=1pt] {$v_{3}$};
  \draw(16 ,6.2) node[draw,fill=white,anchor=west,xshift=2pt,yshift=1pt] {$t=1$};
%
% grid:
%\draw[step=1cm,gray] (0,0) grid (21,7);
\end{tikzpicture}
}% END resize box
\end{center}
  \caption{
      Illustration of the added-damping instability for the settling particle test when
      the added-damping tensor is turned off.
   The particle starts to develop an unphysical rotation around $t=0.4$.
These oscillations in the angular acceleration starts to affect the translational acceleration at $t=0.6$, 
and eventually dominate the overall motion of the particle.
The right plot is the contour plot of $v_{3}$ at $t=1$, which shows the fluid velocity has been affected by the particle rotation.
    Results are shown for grid $\Gc_{\rm sp}^{(2)}$ solved by the~\ampRB~scheme with the added-damping term turned off intentionally.
    %while the exact motion should be dominated by the vertical translation
 }
  \label{fig:settlingParticleInstability}
\end{figure}
}

This test is further used to demonstrate the importance of including the added-damping tensors in the~\ampRB~scheme.
It is found that instabilities due to added-damping effects are still important for this problem 
even though the body density is heavier than that of the fluid and the motion of the body is primarily translational.
To illustrate the importance of these tensors, Figure~\ref{fig:settlingParticleInstability} presents the results from the~\ampRB~scheme but with 
the added-damping tensors intentionally turned off.
It is found that the particle starts to develop an unphysical rotation around $t=0.4$,
and the unphysical rotation eventually pollutes the overall motion of the particle.
Different from the rising body case in~\cite{rbins2017}, 
the counter-acting effects of the pressure is not very strong in the current test due to its  geometry and the oscillations do not saturate as time evolves.

\subsection{Performance comparisons between~\ampRB~and TP-RB schemes }
\label{sec:compareSchemes}

We now use the benchmark problem described in the previous section to compare several aspects of the performance of the~\ampRB~and TP-RB schemes.
%All the following simulations were performed on a workstation with a 2.4 GHz Xeon processor.
To eliminate the complexity due to the parallel implementation, the comparisons are performed using the serial implementations of the algorithms.
%while the performance of the current parallel version of the code will be discussed later in Section~\ref{sec:parallel}.
The main modification for the~\ampRB~scheme is the application of the AMP interface conditions for the pressure Poisson equation, and 
so we first compare the corresponding linear system from the~\ampRB~scheme to that for the TP-RB scheme. Without loss of generality, we consider the two linear systems at $t=0$ 
when there is an impulsive gravity applied to the particle in the quiescent fluid. 
For this special case when $\vv=\zerov$, the discrete Poisson problem for the~\ampRB~scheme becomes
\begin{align}
    &    \Delta_h p_\iv = 0 , \qquad \iv \in \Omega_h\cup\Gamma_h, \label{eq:discretePressure}
\end{align}
with
\begin{align}
   \nv_\iv^T\grad_h p_\iv + \rho{\nv_\iv^T} \Big( \avb + \bvb\times(\rv_\iv^{0}-\xvb^0 )   \Big)  
               &= 0, \qquad\iv\in\Gamma_h, \label{eq:discretePressureInterfaceEquation} 
               \\
 \left\{ 
\begin{bmatrix} \mrb\, I_{3\times 3} & 0 \\  0 & \Ib  \end{bmatrix}
+ \dt\,\ADtensor^0
 \right\} 
\begin{bmatrix} \avb \\ \bvb  \end{bmatrix}
+
\Fv(p_\iv) 
& =  \dt\, 
     \ADtensor^0 \begin{bmatrix} \avb^* \\ \bvb^*  \end{bmatrix} + \frac 4 3 \pi \diskRadius\sp3(\rho_b-\rho) \, \gv,
    \label{eq:discreteRBaccelerationInterface}
\end{align}
together with the remaining boundary conditions at physical boundaries of the fluid domain.
Here $\rv_\iv^{0}$ and $\xvb^0$ are the initial locations of the interface and rigid body, 
and $\avb^*$ and $\bvb^*$ denote initial guesses for the accelerations of the body.  These latter quantities are taken to be
%The initial guess of the body acceleration in the~\ampRB~scheme can be simply taken as
\begin{align*}
    \avb^* = \left(\frac{\rhos-\rho}{\rhos+0.5\,\rho}\right)\gv, \qquad \bvb^*= \zerov,
\end{align*}
which are the exact initial accelerations of a spherical particle in a viscous quiescent fluid of infinite extent.
%therefore it serves as a very good initial guess in the current geometry. 
%In practice, \QT{the body accelerations at $t = 0$ can be computed without prediction from the given initial condition} \old{the initial acceleration can be directly solved without prediction} since the added-damping tensors are not involved at $t = 0$.
The computed added-damping tensor, denoted by $\ADtensor^0$, is included in~\eqref{eq:discreteRBaccelerationInterface} even though it is not required at $t=0$ since added-damping effects are not instantaneous
so that the discrete problem in~\eqref{eq:discretePressure}--\eqref{eq:discreteRBaccelerationInterface} without $\ADtensor^0$ is well-posed.
However, we have included the tensor intentionally to test the full linear system in a general case when $t>0$.
The corresponding discrete Poisson problem used in the TP-RB scheme  involves  the pressure equation in~\eqref{eq:discretePressure} with the interface condition
\begin{align}
   \nv_\iv^T\grad_h p_\iv =- \rho{\nv_\iv^T} \Big( \avb^* + \bvb^*\times(\rv_\iv^{0}-\xvb^0 )   \Big)  
               , \qquad\iv\in\Gamma_h, %\label{eq:pressureInterfaceEquation} 
\end{align}
and the remaining physical boundary conditions. 
%The readers are referred to our previous work~\cite{rbinsmp2017} for more discussions on a simplified version of this system.
%since the nonslip boundary condition is satisfied exactly on the interface of the sphere.
We note that the initial guesses for the accelerations in the interface condition for the TP-RB scheme is not as important since sub-iterations are needed at $t=0$ to obtain a converged initial acceleration for this scheme.

{% ------ MATLAB CURVES UNSTABLE LIGHT BODY ------
\newcommand{\figWidth}{7cm}
\newcommand{\trimfig}[2]{\trimFig{#1}{#2}{.25}{.26}{.0}{.0}}
\newcommand{\figWidthz}{5cm}
\newcommand{\figWidtha}{6cm}
\newcommand{\trimfiga}[2]{\trimFig{#1}{#2}{.28}{.26}{.0}{.0}}
\def\xma{-.15}
\def\xmb{5.3}
\def\yma{1.3}
\def\ymb{6.75}
\tikzstyle{ann} = [font=\footnotesize,inner sep=1pt]
\begin{figure}[htb]
\begin{center}
\resizebox{16cm}{!}{% START resize box
\begin{tikzpicture}[scale=1]
  \useasboundingbox (0.0,1) rectangle (16.,7.7);  % set the bounding box (so we have less surrounding white space)
  \draw(-1,0) node[anchor=south west,xshift=-4pt,yshift=+0pt] {\trimfig{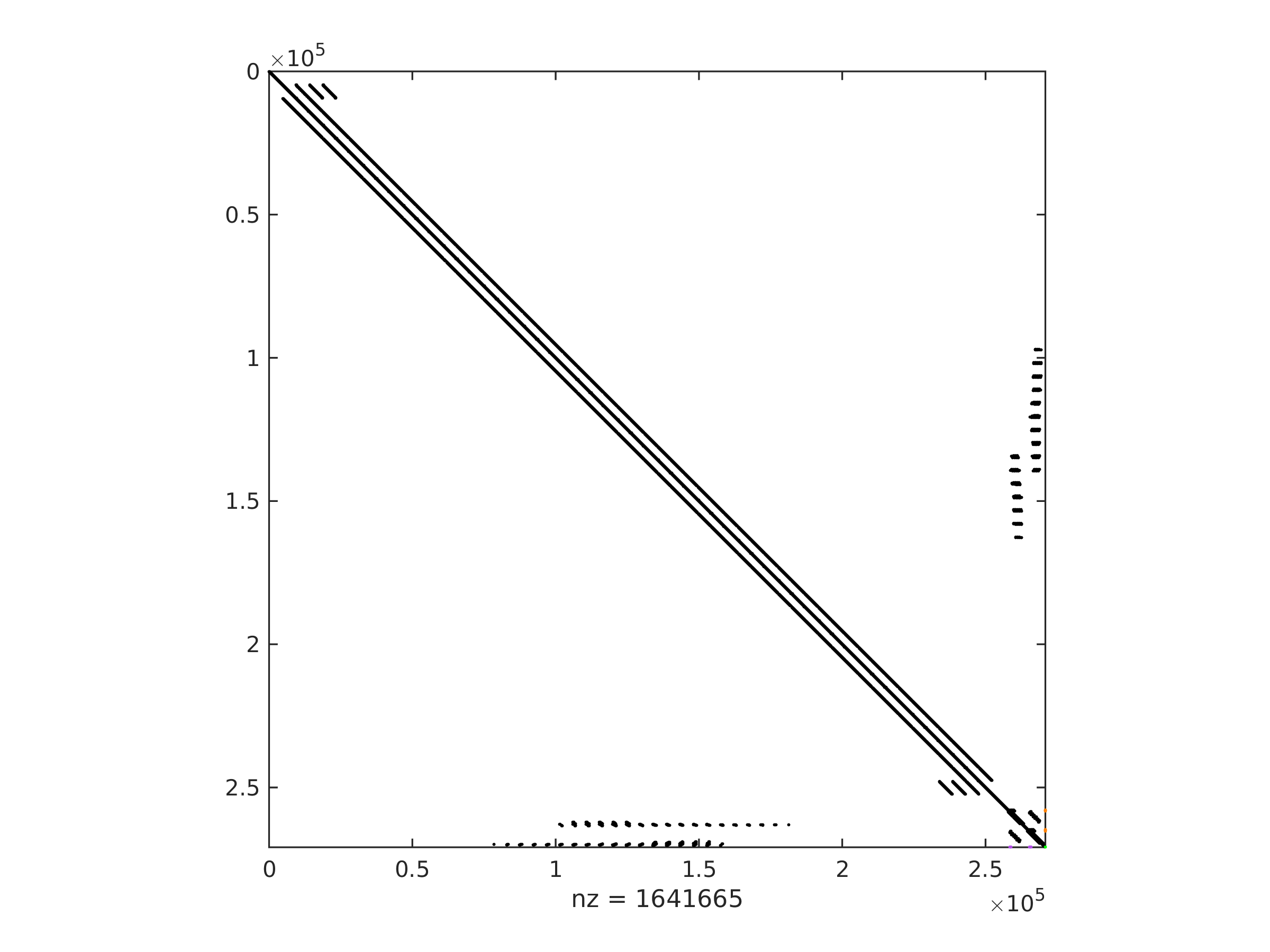}{\figWidth}};
  \draw(5.8,1) node[anchor=south west,xshift=-4pt,yshift=+0pt] {\trimfiga{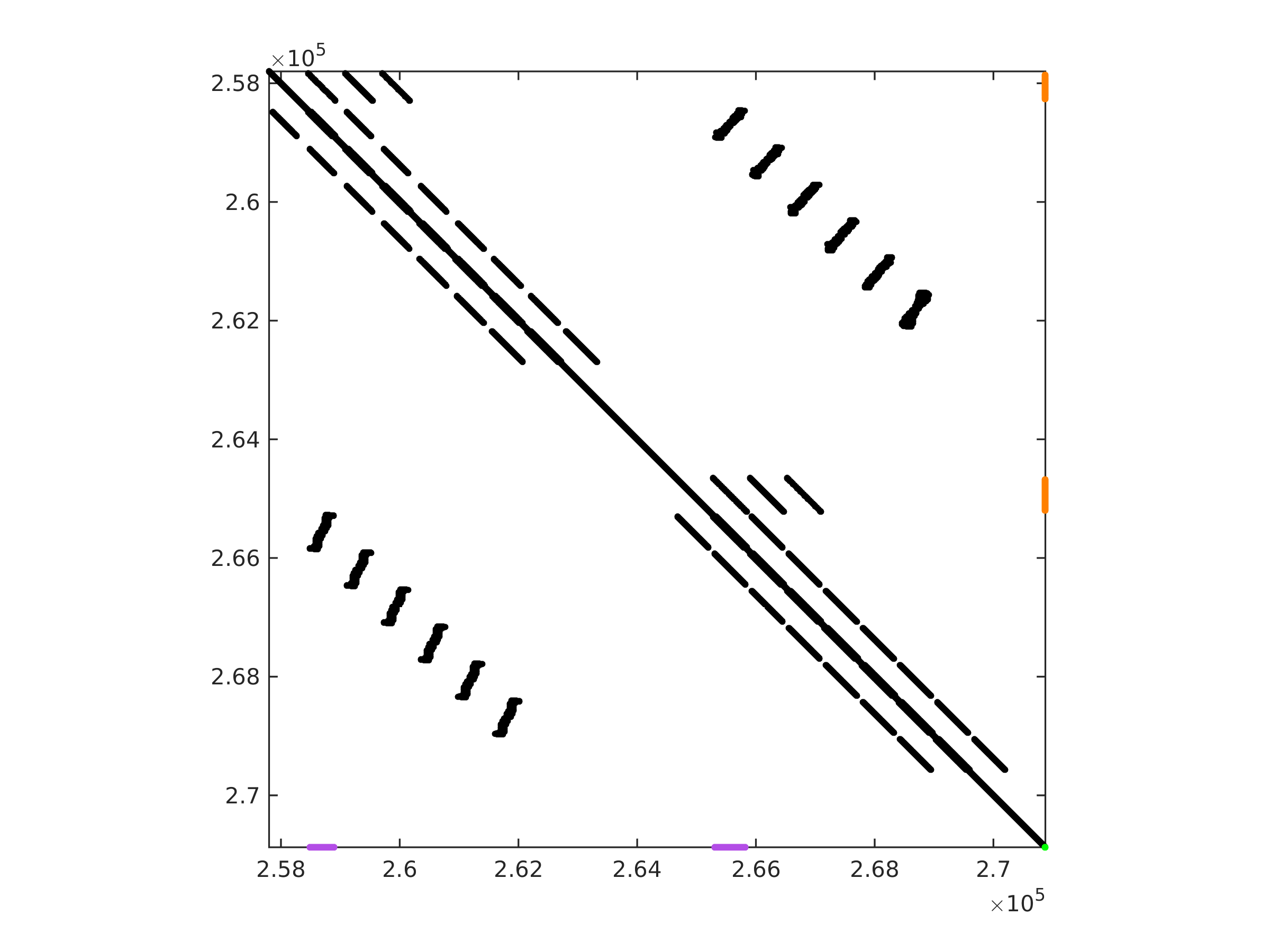}{\figWidtha}};
  \draw(11.5,2) node[anchor=south west,xshift=-4pt,yshift=+0pt] {\trimfig{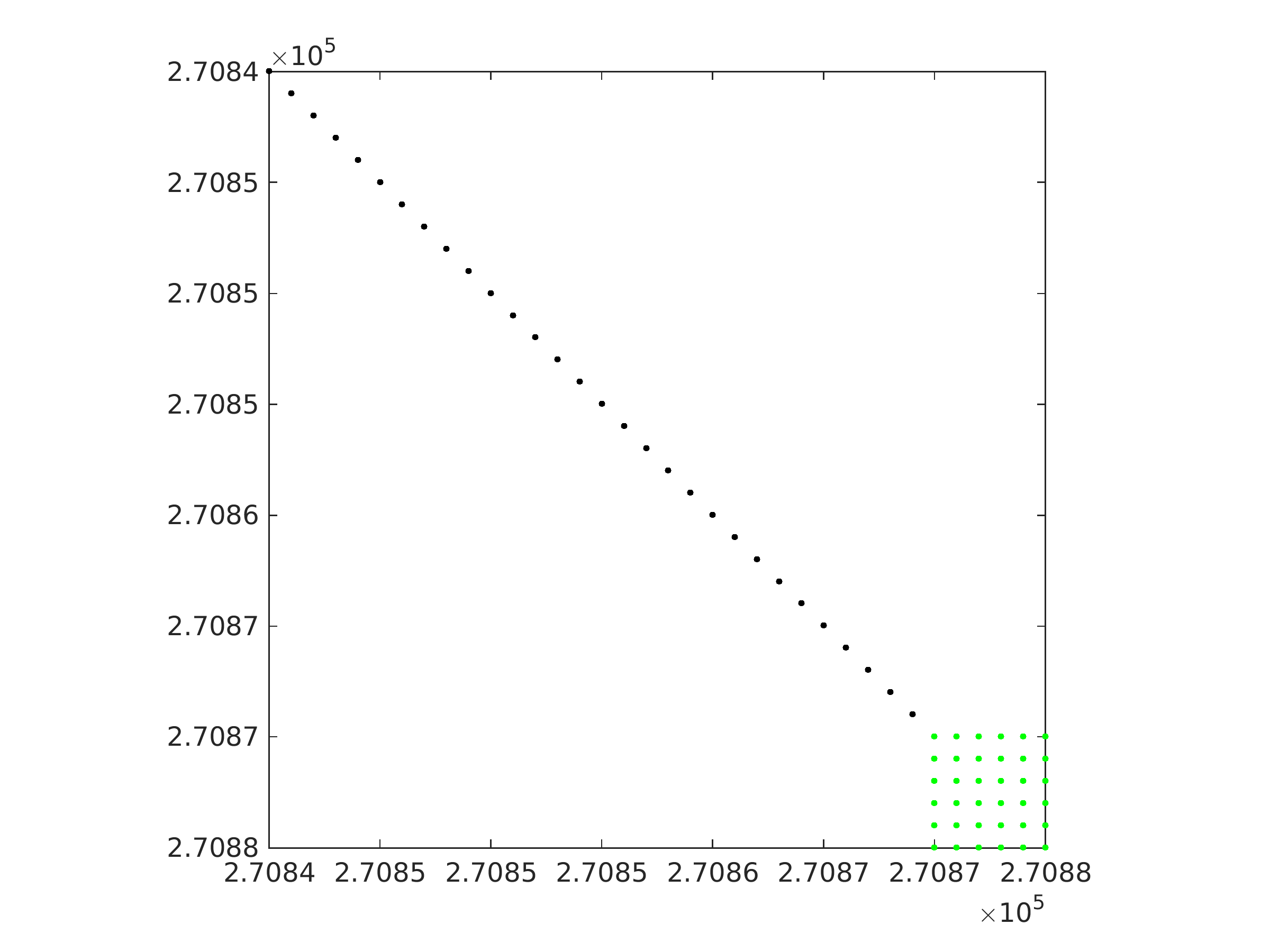}{\figWidthz}};
%  \draw(16 ,6.2) node[draw,fill=white,anchor=west,xshift=2pt,yshift=1pt] {$t=1$};
  %\draw[red,thick,dashed] (\xma,\yma) -- (\xma,\ymb) -- (\xmb,\ymb) -- (\xmb,\yma) -- (\xma,\yma);
  \draw[red,thick,dashed](\xma,\yma)rectangle(\xmb,\ymb);
  \draw[blue,thick,dashed](6.5,4.85) rectangle ++(1.85,1.85);
  \draw[blue,thick,dashed](9.05,2.3) rectangle ++(1.85,1.85);
   %\node[ann] at (13,.8) {$\begin{bmatrix} \mrb\, I_{3\times 3} & 0 \\  0 & \Ib  \end{bmatrix}
   \draw(12.4,1.2) node[draw,fill=white,anchor=west,xshift=0pt,yshift=0pt]{\scriptsize$\begin{bmatrix} \mrb\, I_{3\times 3} & 0 \\  0 & \Ib  \end{bmatrix} +\dt\,\ADtensor$};
   \draw[arrows=->,line width=1pt](14.2,1.65)--(15.4,2.7);
   %
   %\node[ann] at (7.1,.8) {$\Fv(p_\iv)$};
   \draw(7.5,.6) node[draw,fill=white,anchor=west,xshift=0pt,yshift=0pt]{\scriptsize$\Fv(p_\iv)$};
   \draw[arrows=->,line width=1pt](8.,.9)--(6.9,1.75);
   \draw[arrows=->,line width=1pt](8.,.9)--(9.3,1.78);
   %
   %\node[ann] at (10,6.8) {$\rho{\nv_\iv^T} \Big( \avb + \bvb\times(\rv_\iv^{0}-\xvb^0 ) \Big)$};
   \draw(7.0,7.25) node[draw,fill=white,anchor=west,xshift=0pt,yshift=0pt]{\scriptsize$\rho{\nv_\iv^T} \Big( \avb + \bvb\times(\rv_\iv-\xvb ) \Big)$};
   \draw[arrows=->,line width=1pt](10.7,7.2)--(11.3,6.6);
   \draw[arrows=->,line width=1pt](10.7,7.2)--(11.3,4.2);
  %\draw[red,thick,dashed] (\xza,\yza) -- (\xza,\yzb) -- (\xzb,\yzb) -- (\xzb,\yza) -- (\xza,\yza);
%
% grid:
%\draw[step=1cm,gray] (0,0) grid (16,7.7);
\end{tikzpicture}
}% END resize box
\end{center}
  \caption{
   The sparsity pattern of the pressure matrix of $270875\times270875$ in the~\ampRB~scheme.
   Left: the sparsity pattern from the grid $\Gcsp^{(1)}$.
   Middle and right: two zoomed views of the pattern.
   %This matrix has a size of $270875\times270875$.
   The entries corresponding to the background Cartesian grid are marked by the red square and the entries of the two boundary fitted grids are marked by the blue squares.
   The important terms in the AMP interface conditions are marked in different colors and pointed out in the sparsity pattern.
 }
  \label{fig:sparsePattern}
\end{figure}
}

Figure~\ref{fig:sparsePattern} shows the sparse matrix pattern of the linear system of the~\ampRB~scheme for the grid $\Gcsp^{(1)}$,
which consists of one background grid and two boundary-fitted grids.
The pressure matrix has a size of $270875\times270875$ and can be divided into several blocks.
The entries corresponding to the background Cartesian grid are marked by the red box in the figure, and 
the entries of the two boundary-fitted grids are marked by the blue boxes.
The majority of the diagonal blocks correspond to the discrete Poisson operator, while the off-diagonal blocks correspond primarily to the interpolations between the different grids.
The six AMP constraints in~\eqref{eq:discreteRBaccelerationInterface} related to the accelerations are positioned in the last rows of the matrix,
and various terms in the constraints~\eqref{eq:discretePressureInterfaceEquation}--\eqref{eq:discreteRBaccelerationInterface} have been marked by different colors 
and pointed out explicitly in the sparsity pattern.
Apart from the six extra rows corresponding to the equations in~\eqref{eq:discreteRBaccelerationInterface} and the six extra columns for the accelerations of the body, the matrix of the new system is identical to that used in the TP-RB scheme.

{
\newcommand{\figWidth}{7.75cm}
\newcommand{\trimfig}[2]{\trimFig{#1}{#2}{.0}{.0}{.0}{.0}}
\begin{figure}[hbt]
\begin{center}
\begin{minipage}[b]{0.44\linewidth}
\begin{tikzpicture}[scale=1]
  \useasboundingbox (0,0.) rectangle (8.,6);  % set the bounding box (so we have less surrounding white space)
  \draw(0,-0.5) node[anchor=south west,xshift=0pt,yshift=+0pt] {\trimfig{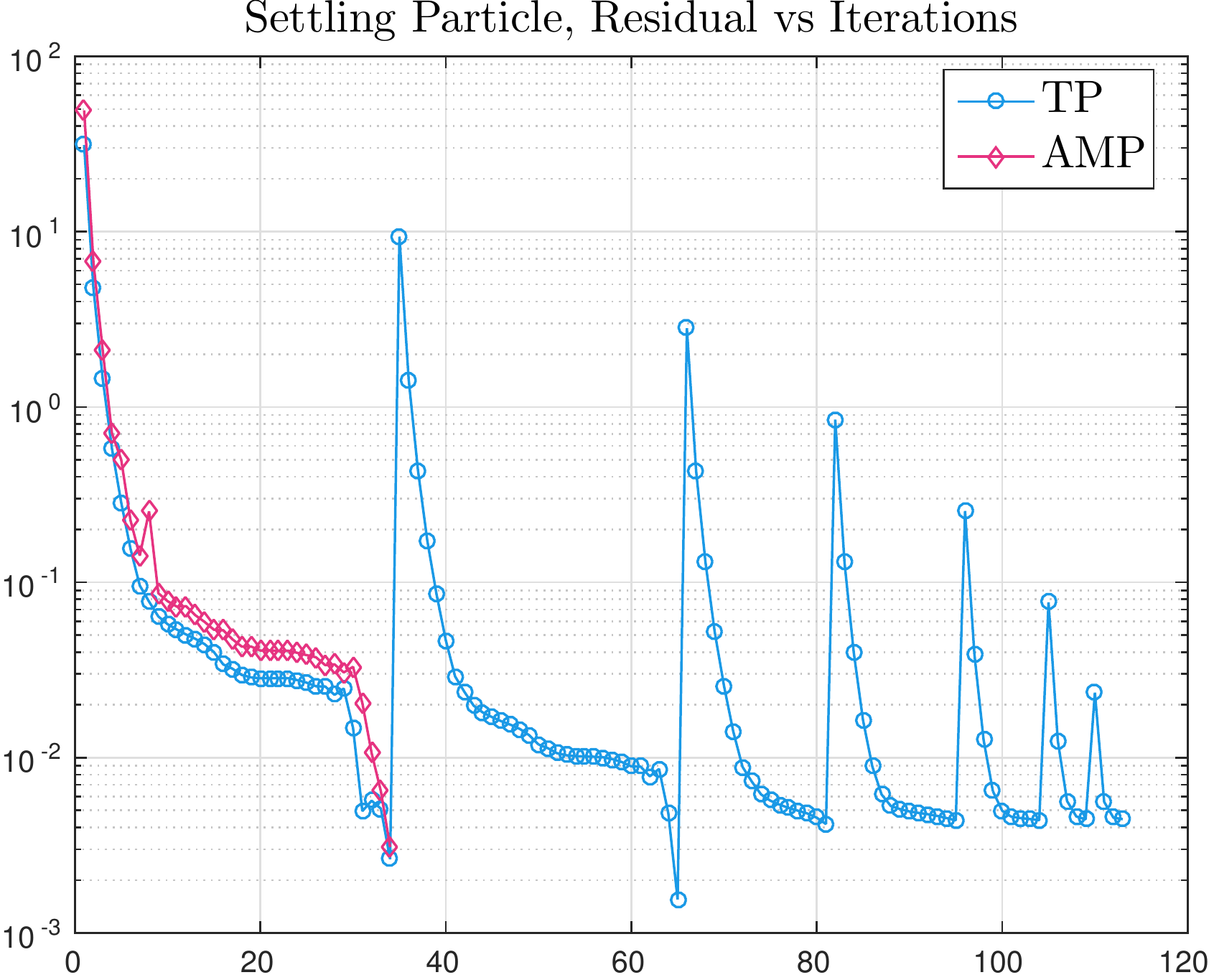}{\figWidth}};
  % grid:
% \draw[step=1cm,gray] (0,0) grid (8,6);
 \end{tikzpicture}
\end{minipage}
\hfill
\begin{minipage}[b]{0.52\linewidth}
    \centering
\tableFont
\begin{tabular}{| c | c | c | c|}
\hline
\multicolumn{4}{|c|}{Pressure matrices of \ampRB~and TP schemes } \\ \hline
Scheme and Grid & $\sigma_{\rm max}$ & $\sigma_{\rm min}$ & Condition number  \\ \hline
\hline
AMP on $\Gcsp^{(1)}$ & 1.71 &\num{3.19}{-5} & \num{5.34}{4}  \\ \hline
TP  on $\Gcsp^{(1)}$ & 1.00 &\num{3.19}{-5} & \num{3.13}{4}  \\ \hline
AMP on $\Gcsp^{(2)}$ & 1.49 &\num{8.04}{-6} & \num{1.86}{5}  \\ \hline
TP  on $\Gcsp^{(2)}$ & 1.01 &\num{8.02}{-6} & \num{1.26}{5}  \\ \hline
\end{tabular}
\vspace{1.5cm}
\end{minipage}
\end{center}
  \caption{
Left: residual verse iterations in the~\ampRB~and TP schemes. Right: estimated condition numbers and other quantities of the pressure matrices in the two schemes. 
The sudden jump in the residual of the TP scheme on the left plot is due to the restart of the pressure solver in each sub-iteration.
$\sigma_{\rm max}$ and $\sigma_{\rm min}$ are the estimated largest and smallest singular values.
All the quantities in the right table are estimated through GMRES.
}
  \label{fig:ampTpCompare}
\end{figure}
}

The linear systems for the two schemes are solved by a serial Bi-CGSTAB solver with an ILU(1) preconditioner.
The Bi-CGSTAB solver is chosen because the matrix for a composite grid is typically not symmetric, and the Bi-CGSTAB solver often shows better performance than other Krylov solvers in this case.
The ILU(1) preconditioner is chosen here because the system is relatively easily inverted and only a few Krylov iterations are needed to converge.
%that other preconditioners such as ILU(3) show slower performance. 
%This is because the linear equations of the current test in each step are considered relatively easy and do not need many Krylov iterations even with a simple ILU(1) preconditioner.
%For some hard problem we considered later, ILU(3) are used to reduce the number of Krylov iterations.
%This is probably because the overall matrices in three dimensions are more banded than those in two dimensions.
Figure~\ref{fig:ampTpCompare} presents the convergence of the residual for each sub-iteration.
The~\ampRB~scheme uses only one Krylov solve to obtain the solution, while the TP-RB scheme uses a Krylov solve for each of the seven sub-iterations needed for convergence of the whole system (corresponding to the sudden jumps of the residual in the figure).
Note the performance of the Krylov iterations in the~\ampRB~and TP-RB schemes is similar in the first solve, both
using 33 Krylov iterations to obtain a convergent solution. This indicates the both systems are similarly conditioned.
The subsequent sub-iterations for the TP-RB scheme require fewer and fewer Krylov iterations as the previous computed pressure serves as a better initial guess for the next sub-iteration.

To further understand the nature of the linear system, the condition numbers are estimated by the extremal singular values.
The extremal singular values of the matrices are estimated by the singular values in the Hessenberg matrix of the GMRES iterations after about 1000 iterations.
The table in Figure~\ref{fig:ampTpCompare} shows that the largest singular value of the \ampRB~matrix is slightly larger than the TP-RB scheme
while both schemes have almost identical smallest singular values which  scale as $O(h^2)$.
This confirms that the two systems have very similar conditioning.
Recall that the condition number of the system is the ratio between the  maximal and minimal singular values if the $2$-norm is used.
For the systems at a later time when $t > 0$,
we find that the performance of the linear solver applied to the two systems is very similar to the case at $t = 0$,
since the difference between the linear systems at two different times lies mainly in the right-hand side.
%Those similar results are hence not presented to save space.

{ % =================================== START TABLE ===============================================
% \setlength\extrarowheight{5pt}
% -------------------------------------------------------------------------------------------------------
\begin{table}[hbt]\tableFont % you should set \tableFont to \footnotesize or other size
%\begin{center}
%
\hspace{-.4cm}
\begin{tabular}{|c|c|c||c|c|c|c|c|c||c|c|c|c|}
\hline
\multicolumn{13}{|c|}{\strutt Run-time performance of AMP versus TP } \\[3pt] \hline
\multirow{2}{*}{Grid} &
\multirow{2}{*}{rtol} & 
\multirow{2}{*}{atol} & 
\multicolumn{6}{c||}{TP} &
\multicolumn{4}{c|}{AMP}
%%\multicolumn{1}{c|}{Domain 0} &
\\ 
\cline{4-13} 
& & & rtolc & atolc & p/step & v/step & its/step & sec/step  & p/step & v/step & sec/step & speed-up \\
\hline
\hline
%new
 $\Gcsp^{(1)}$  & \num{1}{-5} & \num{1}{-7} & \num{1}{-4} & \num{1}{-6} & 25 & 12 & 7 & 5.55 & 18 & 5 & 3.36 &1.65\\ 
\hline                                                                                                        
 $\Gcsp^{(2)}$  & \num{2.5}{-6} & \num{2.5}{-8} & \num{2.5}{-5} & \num{2.5}{-7} & 53 & 16 & 8 & 40.8 & 39 & 6 & 23.6 & 1.73\\ 
\hline
\end{tabular}
%
%old
% $\Gcsp^{(1)}$  & \num{1}{-4} & \num{1}{-6} & \num{1}{-4} & \num{1}{-6} & 25 & 12 & 7 & 4.79 & 17 & 5 & 2.88 &1.66\\ 
%$\Gcsp^{(2)}$ &  \num{1}{-5} & \num{1}{-7} & \num{1}{-5} & \num{1}{-7} & 46 & 51 & 26 & 74.6& 20 & 6 & 19.8 &3.77\\ 
%\hline                                                                                                        
%$\Gcsp^{(2)}$  & \num{1}{-6} & \num{1}{-8} & \num{1}{-5} & \num{1}{-7} & 72 & 14 & 6 & 40.9 & 58 & 7 & 28.7 &1.43\\ 
 \caption{Settling particle. Comparison of the run-time performance of the AMP scheme versus the TP scheme.
The relative and absolute tolerances of the velocity and pressure solvers are denoted by ``rtol''and ``atol'',
while the tolerances of the sub-iterations in the TP scheme are denoted by ``rtolc'' and ``atolc''.
``P/step'' and ``v/step'' stand for the averaged number of Krylov iteration in the pressure and velocity solvers per time-step, while
``sec/step'' for the averaged CPU run-time per time-step.
The averaged number of iteration in the TP scheme are also given in the column of ``its/step''.
 }
\label{table:settlingParticleCompareSchemes}
%\end{center}
\end{table}
}
%\QT{JCP will fix the table width. So no need to worry about it for now. Or we can remove the column of ``atol''.}

We further compare the total CPU time to solve the full problems given by the~\ampRB~and TP-RB schemes. 
Table~\ref{table:settlingParticleCompareSchemes} presents detailed information of the two schemes on the grids $\Gcsp^{(1)}$ and $\Gcsp^{(2)}$.
The grid $\Gcsp^{(1)}$ has about $2.2 \times 10^5$ total grid points, while $\Gcsp^{(2)}$ has about $1.7 \times 10^6$ points.
%We leave out the result of $\Gcsp^{(4)}$ (about $1.3\times10^7$ grid points) due to the very long run-time  (the \ampRB~scheme needs about 100 hours of CPU time and the TP-RB scheme needs even longer).
Other information such as the tolerances in the Bi-CGSTAB solvers and sub-iterations in the TP-RB scheme are also provided in Table~\ref{table:settlingParticleCompareSchemes}.
The under-relaxation parameter in the TP-RB scheme (see~\cite{Koblitz2016}) is chosen to be 0.5, which we found to give the best performance for this problem.
The results for $\Gcsp^{(1)}$ show that the CPU time of the~\ampRB~scheme is about 3.36 seconds per time step, while the time of the TP-RB scheme is about 5.55 seconds per time step,
which corresponds to a speed-up factor of 1.65 the~\ampRB~scheme.
Note the TP-RB scheme for this problem requires about 7 sub-iterations for each time step, but the performance of the~\ampRB~scheme only has increased by a factor of 1.65.
This is expected since, as demonstrated in Figure~\ref{fig:ampTpCompare}, the later correction steps in the TP-RB scheme require fewer Krylov iterations than first two steps.
In a typical time step, the predictor and the first corrector in the TP-RB scheme account for most of the Krylov iterations in the pressure solver (``p/step'' in Table~\ref{table:settlingParticleCompareSchemes})
and the further correction steps only require one or two Krylov iterations to obtain a converged local solution.
Note that for both schemes, each Krylov iteration takes almost the same CPU time on average, although the linear system for the~\ampRB~scheme is slightly larger.
The $\Gcsp^{(2)}$ results show a similar speed-up factor of 1.73.  The CPU time of the~\ampRB~scheme is about 23.6 seconds per time step for this grid, while the time of TP-RB scheme is about 40.8 seconds per time step.
Note also that the tolerances of each sub-iteration in the TP-RB scheme are chosen larger than the tolerances of the Krylov solver to avoid too many sub-iterations.
We have found that this choice is sufficient to yield satisfactory numerical solutions.
%and we found it is enough to produce numerical solutions of the good quality for this test.
However, if the tolerances for the sub-iterations and Krylov solvers are chosen to be the same, the~\ampRB~scheme is about 4 times faster than the TP-RB scheme.
Finally, note that the current overhead in various other parts of the schemes, such as computing information related to the overlapping grids and evaluating the advection terms, are roughly the same in both schemes.
For instance, the overhead accounts for about 15\% of the ``sec/step'' in the~\ampRB~scheme and it accounts for about 10\% for the TP-RB scheme.

%The speed-up will be even larger if all those facts are considered.
%It also confirms that the AMP scheme shows much better performance than the TP-RB scheme even on a relatively simple case when the density of the body is heavier than the fluid.
%All the performance results clearly indicate that the sub-iterations in the TP-RB scheme should be avoided as much as possible. 

In summary, the above tests show the new linear system in the~\ampRB~scheme to solve the discrete Poisson problem for the pressure requires only minor modifications to the corresponding linear system for pressure in the TP-RB scheme.  Both systems have similar conditioning and can be solved effectively using preconditioned iterative Krylov solvers.  The~\ampRB~scheme requires no sub-iterations per time step unlike the TP-RB scheme so the cost per time step for the~\ampRB~scheme is lower.  As the ratio of the density of the rigid body to that the fluid becomes smaller, the TB-RB scheme requires more and more sub-iterations (and ultimately fails).  In contrast, the~\ampRB~scheme does not require sub-iterations and so the comparative performance of the~\ampRB~scheme is more pronounced for lighter bodies, such as in the problem discussed in the next section.

\subsection{One particle falling or rising in a long container} 
\label{sec:risingParticle}

The next validation problem we consider is a spherical particle falling or rising 
due to a buoyant force in a very long container.
The original problem was discussed in an experimental study~\cite{mordant2000velocity}, but then used later in~\cite{uhlmann2005immersed} to validate an immersed boundary method.
In many subsequent studies of FSI algorithms, such as those in~\cite{kempe2012improved, yang2014sharp, yang2015non, tschisgale2017non}, Case~2 in~\cite{uhlmann2005immersed} has been adopted 
as a standard benchmark problem to demonstrate numerical stability for problems involving low-density particles.
The problem shares a very similar geometry as the particle dropping test discussed in the previous two sections, but
the size of the particle here is much smaller compared to the computational domain.
%its nondimensionalization in different work can be slightly different. 
We follow the nondimensional setup from the work in~\cite{uhlmann2005immersed, kempe2012improved}.
An incompressible fluid with viscosity $\mu=0.00104238$ and density $\rho=1$ occupies the domain $[-x_c, x_c] \times [-y_c, y_c] \times [z_0, z_1]$, where $x_c = y_c = 0.64$, $z_0=0$ and  $z_1 = 10$.
The particle has a density of $\rhos$ and a radius of $\diskRadius=1/12$.
The fluid and the particle are initially quiescent and an impulsive external body force satisfying~\eqref{eq:particleForce} 
with $\tilde g=9.81$ is applied at $t = 0$.
In~\cite{uhlmann2005immersed} the particle was initially located at ${\bf x}_0=(0, 0, 9.5)$ and the density of particle is $\rhos = 2.56$. 
However, in more recent work the particle density has varied, and the velocity profile of the particle for different values of $\rhos$ are presented. 
In our simulations, the particle is initially located at ${\bf x}_0$ for cases when the particle falls ($\rhos>1$), and it is located at $(0, 0, 0.5)$ initially for cases when it rises ($\rhos<1$).
%Therefore, the initial location of the particle is not important as long as the particle achieves a steady state inside the fluid domain.
Periodic boundary conditions are applied in all three directions as in the previous work~\cite{uhlmann2005immersed, kempe2012improved}.

The original problem is designed to model the situation when the particle is immersed in an infinite fluid domain.
%and therefore the flow structure in this test all surrounds the particle.
%The computational cost can be significantly reduced if a coarser grid is used as the background grid.
%and therefore there are no boundary layers in the other part of the fluid domain.
The composite grid approach allows us to use a coarser background grid than the previous work~\cite{uhlmann2005immersed, kempe2012improved} and 
at the same time attach much finer boundary-fitted grids to the particle so that its boundary layer and the particle wake are well resolved.
The sphere-in-a-box composite grid similar to the previous case is used with slight modifications.
A target mesh spacing $h=0.02$ is used for the background Cartesian grid, corresponding to a grid of $64\times64\times500$ covering the entire domain.
Much finer boundary-fitted grids are attached to the particle and cover the fluid region of length $\Delta r = 0.3$ in the radial (normal) direction.
The target mesh spacing for the boundary-fitted grids is $h/4 = 0.005$. 
Here we use a three-patch sphere to avoid excessively small time steps resulting from small cells in the two-patch version.
The resulting composite grid\footnote{The interested reader is again referred to~\cite{pog2008a} for more details on the sphere-in-a-box grid.}, consisting of the background grid and three boundary-fitted grids, is denoted by $\Gc_{\rm sr}$.
A fixed time step, $\dt = 0.001$, is used in all the simulations. 

%will be sufficient to obtain a numerical solution of similar quality.

{% ------ MATLAB CURVES LIGHT BODY ------
\newcommand{\figWidth}{8.5cm}
\newcommand{\trimfig}[2]{\trimFig{#1}{#2}{.0}{.0}{.0}{.0}}
\newcommand{\figWidtha}{8.5cm}
\begin{figure}[htb]
\begin{center}
\resizebox{16cm}{!}{% START resize box
\begin{tikzpicture}[scale=1]
  \useasboundingbox (0.0,0) rectangle (18,6.5);  % set the bounding box (so we have less surrounding white space)
  \draw(-.1,-.75) node[anchor=south west,xshift=0pt,yshift=+0pt] {\trimfig{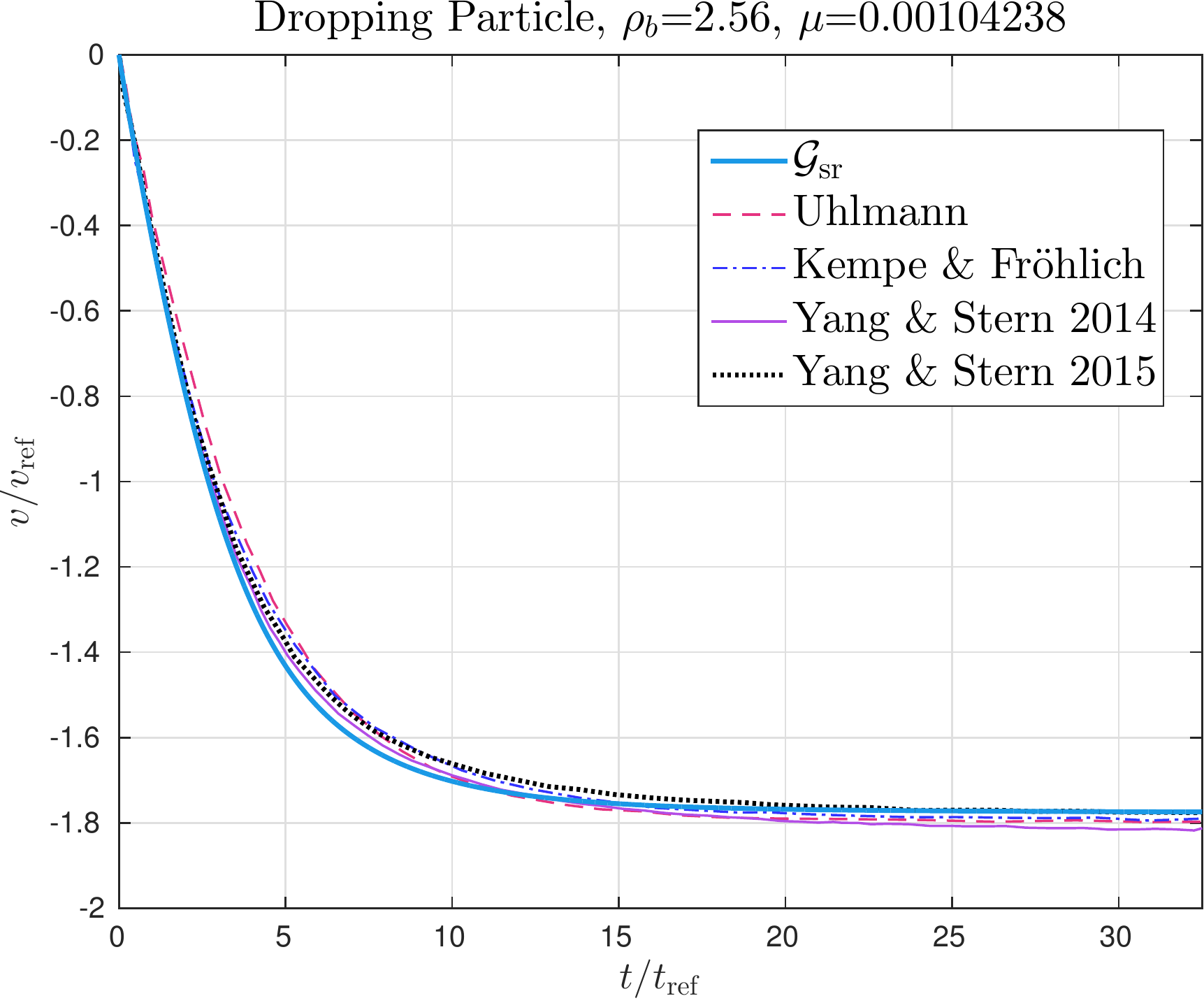}{\figWidth}};
  \draw(9.1,-.75) node[anchor=south west,xshift=0pt,yshift=+0pt] {\trimfig{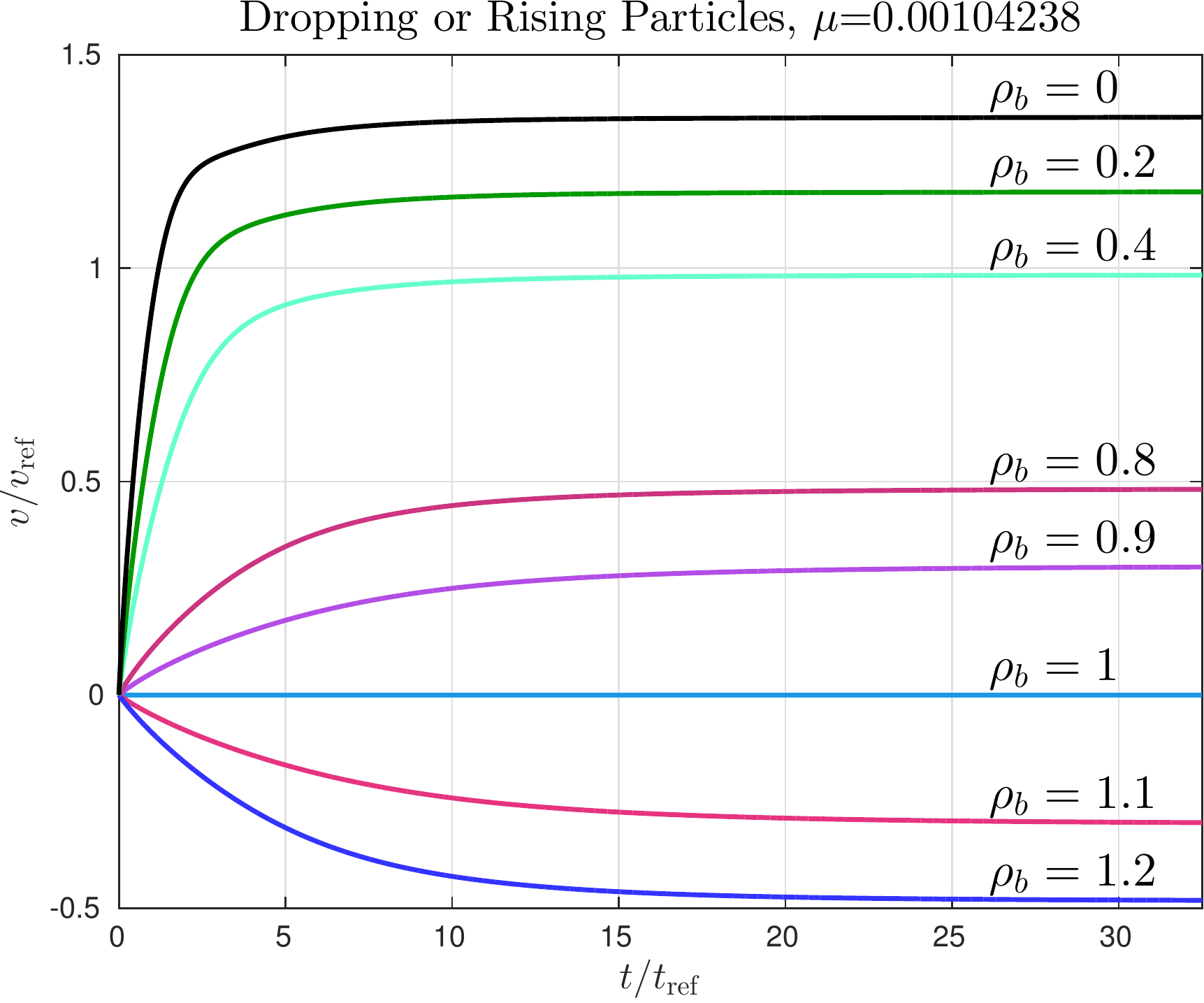}{\figWidtha}};
%
% grid:
%\draw[step=1cm,gray] (0,0) grid (18,6.5);
\end{tikzpicture}
} % END resize box
\end{center}
  \caption{Dropping or rising particle. Time history of the velocity (left) of the case $\rhos=2.56$ and velocity (right) of the varying density cases.  The reference velocity and time are $v_{\rm ref} = \sqrt{2\diskRadius\tilde g}$ and $t_{\rm ref}=\sqrt{2\diskRadius/\tilde g}$, where $\diskRadius=1/12$ and $\tilde g=9.81$.
      The solutions are computed on grid $\Gc_{\rm rs}$.
      In the left plot, the red dashed curve is the numerical solutions taken from Fig.~16 in~\cite{uhlmann2005immersed},
      the blue dash-dotted curve is from Fig.~11 in~\cite{kempe2012improved},
      the purple curve is from Fig.~5 in~\cite{yang2014sharp}
      and the black dotted curve is from Fig.~10 in~\cite{yang2015non}.
  }
  \label{fig:risingParticleCompareCurves}
\end{figure}
}

The left plot in Figure~\ref{fig:risingParticleCompareCurves} shows the time history of the particle velocity for the case $\rhos=2.56$.
%The velocity and time presented in Figure~\ref{fig:risingParticleCompareCurves} have been scaled by the reference values of $v_{\rm ref} = \sqrt{D |\gv|}$ and $t_{\rm ref}=\sqrt{D/|\gv|}$ with the diameter of the particle $D=1/6$.
The numerical solution given by the~\ampRB~scheme is compared with the results found in the literature~\cite{uhlmann2005immersed, kempe2012improved, yang2014sharp, yang2015non} using different immersed boundary methods.
Our results are in good overall agreement with these results.
The computed terminal speed is found to be $v_{\rm term} = 1.774 \, v_{\rm ref}$, which is in good agreement with the experimental terminal speed of $1.797$ reported in~\cite{mordant2000velocity}. 
%Refining the grid $\Gc_{\rm sr}$ by a factor of two results into a better terminal speed of $v_{\rm term}=???$.
%However, the simulations become significantly more expensive and therefore the finer grid is not used in the simulations of varying density tests..
We observe that the particle drops slightly faster at early times in our simulation as compared to the results from the immersed boundary methods, 
although all of the simulations produce similar terminal speeds.
We also observe that the particle wake does not break in our simulations and remains axisymmetric, even when a finer grid is used.
Whether the particle wave remains axisymmetric or not can effect the profile of the particle velocity, and this differs in the various simulations.  For example, the numerical results in~\cite{uhlmann2005immersed} show an asymmetric wake, which triggers much larger horizontal velocities than what we observe,
while the wake remains axisymmetric for the simulations discussed in~\cite{yang2014sharp}.
Whether the wake remains axisymmetric or not is likely due to a hydrodynamic instability and whether this instability is triggered by perturbations in the flow due to the treatment of the moving grids in the vicinity of the particle.  We observe, for example, that a similar simulation using the TP-RB scheme (with sub-time-step iterations) with the same composite grid $\Gc_{\rm sr}$ shows an axisymmetric wake and a nearly identical velocity profile to that given by the~\ampRB~scheme.  This suggests that the issue is not related to the treatment of the interface conditions in the~\ampRB~scheme.

The results from varying density tests are also presented in Figure~\ref{fig:risingParticleCompareCurves},
where the time history of velocities of selected particle densities are presented.
The~\ampRB~scheme remains stable for any density ratio, including the case $\rhos=0$, while the previous methods
appear to be unable to compute solutions for low density ratios.
% work of the immersed boundary method all have \QT{??} \old{a lower stability bound} in terms of density ratios.
%In particular, the results of a zero-density case are presented along with the other cases, 
%which shows a significant improvement compared to previous work.
For example, to our knowledge, the lowest stability bound of this particular benchmark test was previously presented in~\cite{yang2015non}, 
in which their scheme is stable for $\rhos/\rho \ge 0.1$.
On the other hand, the TP-RB scheme in our implementation also struggles to remain stable for small density ratios.
For instance, when $\rhos/\rho = 0.1$, the TP-RB scheme needs about 9 sub-iterations to converge on average.
When the body density is $0.01$, the TP-RB scheme struggles and needs about 85 sub-iterations to converge on average, which makes the scheme significantly more expensive.  It is difficult to stabilize the TP-RB scheme for even small density ratios.
%It becomes very difficult to stabilize the TP scheme when the density is even smaller.

%\input texFiles/twoParticles.tex

\subsection{Numerical simulations of a mechanical heart valve}
\label{sec:heartValveIntro}

%the basic setup and why we want to do that: complex geometry, multiple bodies, more mechanism is involved such as repulsive torque
The last example we consider is a model of a bi-leaflet mechanical heart valve, which has been studied previously in~\cite{borazjani2008curvilinear, de2009direct}, for example, using other numerical approaches.
As shown in Figure~\ref{fig:heartValveDemo}, the heart valve consists of a rigid cylindrical ring supporting two solid leaflets, which can rotate about fixed hinge axes.
The heart valve is placed in a fluid channel which models the aorta.
The structure of each leaflet is designed to be very thin so that its inertia is small compared to that of the fluid flowing through the device.
Previous studies have found this system to be very challenging to 
simulate mainly due to large added-mass effects~\cite{borazjani2008curvilinear}.
Therefore, this problem serves as a good test to demonstrate the applicability of the~\ampRB~scheme together with the
moving composite grid approach for a challenging engineering problem of significant interest.
The problem also can be used as a good benchmark problem for other FSI algorithms.
Unfortunately, a well-specified benchmark problem based on this application is not currently available
in the literature to the best of our knowledge.
The results of previous simulations are available, but the input to these simulations is based on
experimental data, and there is incomplete information given on the initial conditions or boundary conditions in many cases, so that it is difficult to
reproduce the results for the purpose of quantitative comparisons.
Thus, another focus of the present study is to fully specify benchmark problems based on the bi-leaflet
geometry for both two and three-dimensional FSI simulations.
From the point of view of FSI algorithms and added-mass effects, the difficult interval of the full cardiac cycle involves the
rapid rotation of the leaflets as the value opens and closes.  Hence, the designed benchmark problems and our simulations focus on this aspect of the
application.

%- 
{
\newcommand{\xcorner}{-0}% x of lower left corner of top leaflet
\newcommand{\ycorner}{-2.8}% y lower left corner of top leaflet
\newcommand{\xbcorner}{-.2}% x of lower left corner of bottom leaflet
\newcommand{\ybcorner}{-3.12}% y of lower left corner of bottom leaflet
\newcommand{\figWidth}{7.cm}
\newcommand{\figWidtha}{7.5cm}
\newcommand{\trimfig}[2]{\trimh{#1}{#2}{.0}{.0}{.0}{.0}}
\newcommand{\trimfiga}[2]{\trimh{#1}{#2}{.1}{.1}{.1}{.1}}
\newcommand{\labelsize}{\small}
\begin{figure}[htb]
\begin{center}
\resizebox{16cm}{!}{% START resize box
\begin{tikzpicture}[scale=1]
  \useasboundingbox (-.5,1) rectangle (18,7.75);  % set the bounding box (so we have less surrounding white space)
  \begin{scope}[xshift=14.5cm,yshift=6.7cm]
      \fill[fill=blue!15] (-1.5,-6.15) -- (3.6,-6.15) -- (3.6,.3) -- (-1.5,.3) -- cycle; 
      \draw[->,very thick,blue] (-2,-2) -- (-1,-2); 
      \draw[->,very thick,blue] (-2,-3.8) -- (-1,-3.8); 
      \draw[->,thick] (-2,-2.95) -- (4.,-2.95) node [anchor=north] {$x$}; % x-axis
      \draw[->,thick] (0.06,-6.5) -- (0.06,.8) node [anchor=east] {$y$}; % y-axis
      \draw[-,very thick] (-1.5,-6.15) -- (3.6,-6.15);
      \draw[-,very thick] (-1.5,.3) -- (3.6,.3);
      \draw[fill=green!20,draw=green,line width=2pt,rounded corners=2pt,rotate around={70:(\xcorner,\ycorner)}]
        (\xcorner,\ycorner) rectangle ++(2.9,.2);
        \draw[fill=red!20,draw=red,line width=2pt,rounded corners=2pt,rotate around={-70:(\xbcorner,\ybcorner)}]
        (\xbcorner,\ybcorner) rectangle ++(2.9,.2);
    \fill[black] (0.06,-2.95) circle (2 pt); % origin
    \fill[black] (0.06,-2.4) circle (2 pt); % center of top
    \fill[black] (0.06,-3.5) circle (2 pt); % center of bottom
    \coordinate (origo) at (0.06,-2.4);
    \coordinate (pivot) at (1,-2.4);
    %draw arrows
    \draw[thick,dashed,->] (origo) -- ++(2.5,0) node (horizontal) {};
    %\draw[thick,dashed,->] (origo) -- ++(2,0) node (horizontal) {};
    \draw[thick] (origo) -- ++(10:2.5) coordinate (angle1);
    \draw[thick] (origo) -- ++(75:2.5) coordinate (angle2);
    \pic[draw,thick, ->, "\labelsize$\theta_{\rm min}$",angle radius = 1.4cm ,angle eccentricity=1.5] {angle = horizontal--origo--angle1};
    \pic[draw,thick, ->, "\labelsize$\theta_{\rm max}$",angle eccentricity=1.8] {angle = horizontal--origo--angle2};
    %\draw [color=black](pivot)+(10:1.2) node[rotate=0] {$\beta$};
    %\fill[black] (origo) circle (0.05);
        %\draw[thick,gray,->] (origo) -- ++(4,0) node[black,right] {$x$};
        %\draw[thick,gray,->] (origo) -- ++(0,-4) node (mary) [black,below] {$y$};
    %label:
    \draw(-.4,-2.) node[anchor=north] {\small$0.15$}; 
    \draw(-1.5,-6.2) node[anchor=north] {\labelsize$-0.5$};
    \draw(3.6,-6.2) node[anchor=north] {\labelsize$1.25$};
  \end{scope}
  % demo ---
  \draw(-1.8,0) node[anchor=south west,xshift=-4pt,yshift=+0pt] {\trimfig{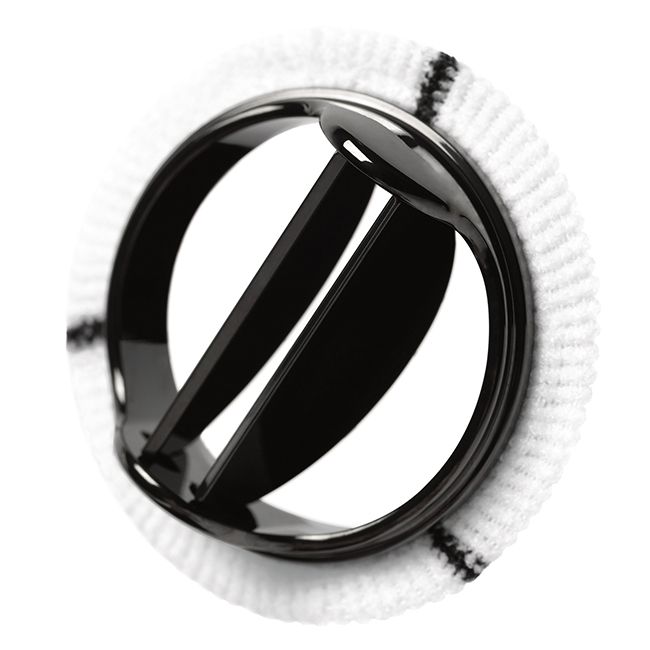}{\figWidtha}};
  \draw(5.2,0) node[anchor=south west,xshift=-4pt,yshift=+0pt] {\trimfig{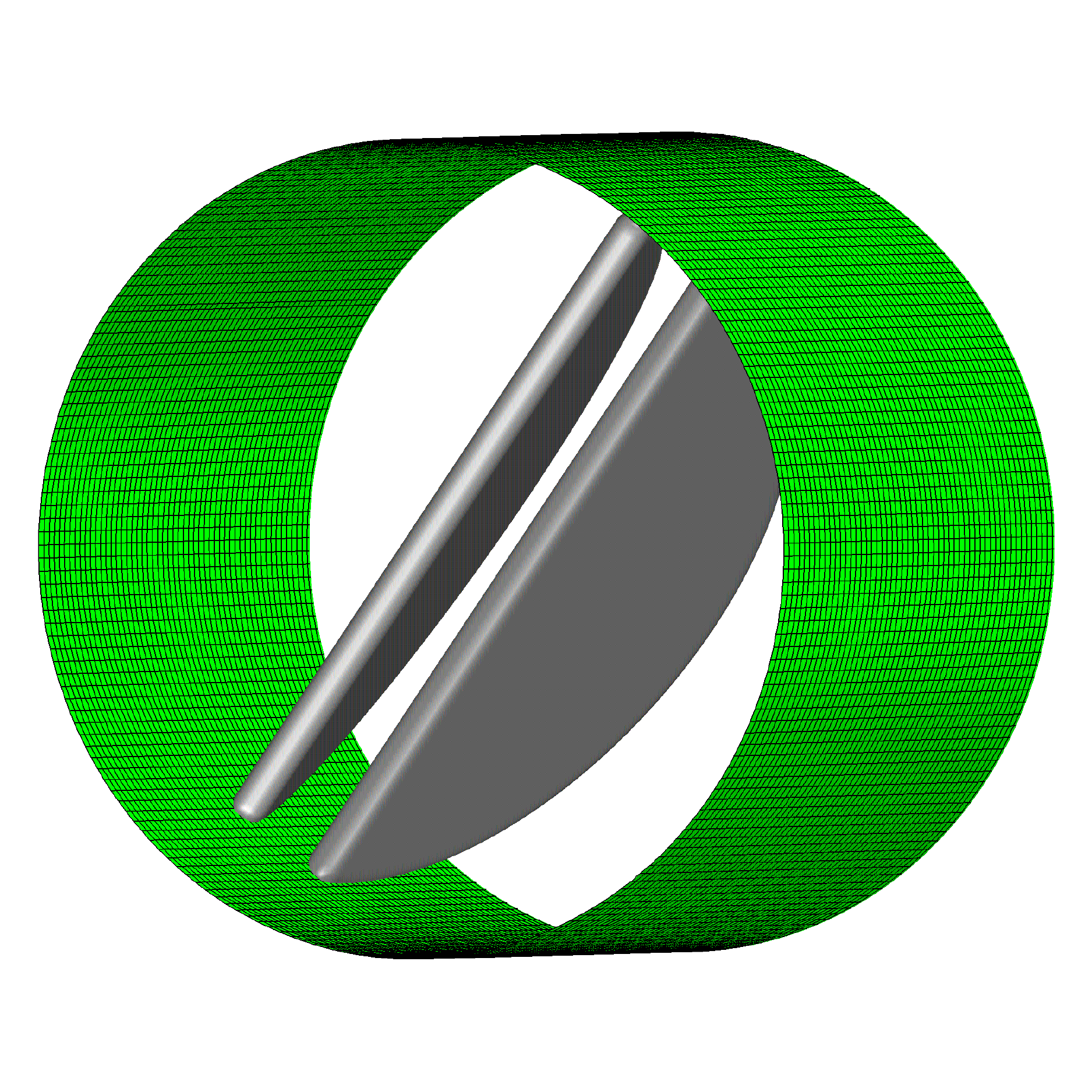}{\figWidth}};
{
\def\thetax{2}
\def\thetay{145}
\def\thetaz{55}
\def\radx{.6}
\def\rady{1.2}
\def\radz{1.2}
%plot coordiantes
\begin{scope}[xshift=10.4cm,yshift=.2cm]
\draw[thick,->] (0,0) -- ({\radx*cos(\thetax)},{\radx*sin(\thetax)}) node[anchor=west]{$x$};
\draw[thick,->] (0,0) -- ({\rady*cos(\thetay)},{\rady*sin(\thetay)}) node[anchor=north]{$y$}; 
\draw[thick,->] (0,0) -- ({\radz*cos(\thetaz)},{\radz*sin(\thetaz)}) node[anchor=west]{$z$};
\end{scope}
}
  %
% grid:
%\draw[step=1cm,gray] (0,0) grid (18,7.75);
\end{tikzpicture}
}% end resize box
\end{center}
  \caption{Mechanical heart valve. Left: a typical bi-leaflet mechanical heat valve (SJM Regent).
  Middle: leaflets and the computational domain. 
  Right: geometrical configuration at the symmetry plane of the domain ($t=0$).
}
  \label{fig:heartValveDemo}
\end{figure}
}

%However, due to the lack of the experimental data, our attention is not to simulate the practical problem exactly and match with experimental observation. Instead, particular attentions are focused on the FSI portion in a full cardiac cycle and  the resulting numerical stabilities of different schemes.

As illustrated in Figure~\ref{fig:heartValveDemo}, the fluid channel covers a cylindrical domain with length $L_c = 1.75$ and radius $R_{c} = 1.1$.
The left boundary is located at $x_0=-0.5$ and the right boundary is $x_1=1.25$.
The valve is fully immersed in an incompressible  fluid with density $\rho=1$ and viscosity $\nu=0.1$.
Two leaflets are allowed to rotate along the hinge axes, which are located at $(x_b, y_b) = (0,\pm0.15)$ and parallel to $z$-axis,
and the leaflets are attached to the hinge at the position 0.1 away from the bottom tips.
%\QT{The length of leaflets are slightly different in 2D and 3D due to dimensional effects, so I prefer to talk about the size of leaflets later when the leaflets and grids are introduced.}
The range of angles for the leaflets are $[\theta_{\rm min}, \theta_{\rm max}]$ with $\theta_{\rm min}=10^{\circ}$ and $\theta_{\rm max}=75^{\circ}$,
which corresponds to the fully open and closed heart valve.
While the precise shape of the leaflets are given later for the specific cases of two and three-dimensional flow,
they are assumed to have the same moment of inertia  given by  $I_b = 0.001$ as used in~\cite{borazjani2008curvilinear}.
%, in which the leaflet inertia has been normalized with respect to the diameter of the aorta and the fluid density.
The leaflets are considered to be ``light'' as added-mass effects corresponding to their rotational motion are large for this test.

%collisions and governing equations 
%To successfully simulate this problem, several numerical or modeling treatments have been implemented to address various issues.
%To successfully simulate this problem, a repulsive torque is applied to the leaflets to model the collisions between the leaflets with the hinges or housing and therefore the moving angles of leaflets can be restricted to a certain range.
%A damping term is added to the system during the collisions to model the loss of energy of the system.

During the operation of the actual mechanical heart valve, the limits of the rotation angles of the leaflets are reached at points of solid-solid contact.  In our numerical simulations, such contacts are avoided by imposing the tighter limits on the rotation angles as given above.  These chosen limits could be imposed in the numerical solution using a post-processing step immediately following any time step for which the rotation limits have been exceeded.  While this approach appears to be a simple choice to implement, we have found that the resulting abrupt perturbation to the solution can lead to numerical instabilities.  Thus, we have adopted an alternate approach in which a continuous repulsive force is applied to the equations of motion of the leaflets to smoothly restrict their rotation to the chosen range.  In addition, we apply a damping term to the equations of motion to model the loss of energy that would occur during a collision.  The equation of motion for the angular acceleration of each leaflet is thus taken to be
\begin{align}
    I_b \,{\dot\omega}_b = \int_{\GammaB} \ev_z^T \Big[ (\rv-\xvb)\times \sigmav \nv \Big] \, \dArea+ g_{\rm rt}(\theta_b) - B(\theta_b)\, \omega_b,
\label{eq:leafletTorque}
\end{align}
where $g_{\rm rt}(\theta_b)$ is the repulsive torque and $B(\theta_b) \, \omega_b$ is the damping term.  Formulas for the functions, $g_{\rm rt}(\theta_b)$ and $B(\theta_b)$, are given in~\eqref{eq:repulsiveTorque} and~\eqref{eq:repulsiveDamp}, respectively, of~\ref{sec:collision} where a full description of the contact model is given.  Here, we note that these terms are set to zero except for a small angle $\delta$ near the limits of the range of rotation.  This angle in~\eqref{eq:repulsiveTorque} and~\eqref{eq:repulsiveDamp}  is taken to be $\delta=3^\circ$ for all of our simulations.
The stiffness parameters in~\eqref{eq:repulsiveTorque} are taken to be $\epsilon_1 = 0.05$ and $\epsilon_2=0.01$,
for the open and closed states, respectively,
and the damping parameter in~\eqref{eq:repulsiveDamp} is $B_0 = 20$.
The stiffness parameter near the closed state, $\epsilon_2$, is chosen to be smaller than that near the open state, $\epsilon_1$, because the collision between the leaflets and hinges is more violent as the valve closes.
The parameters given above specify our collision model, and they are fixed for all of our choices of the grid resolution in our FSI simulations of this problem.

The fluid flow through the heart valve is driven by a pressure difference between the left and right boundaries, which models the changing of blood pressure in a cardiac cycle.
On the left boundary, $x_0=-0.5$, the transverse velocity is set to be zero, while the pressure is taken to be
\begin{align}
    p(x_0,y,z,t) =  \begin{cases}
        p_{\rm max}\, \sin(\pi t) & \text{if $0\le t\le 1$} ,\\
        p_{\rm min}\, \sin(2\pi t) & \text{if $1<t\le 1.25$} , \\
        p_{\rm min} & \text{if $t>1.25$},
\end{cases}
    \label{eq:pressureProfile}
\end{align}
with $p_{\rm max} = 20$ and $p_{\rm min}=-40$.
The magnitude of the minimum pressure $p_{\rm min}$ is chosen to be larger than the maximum pressure $p_{\rm max}$ to mimic a typical cardiac cycle, 
in which the magnitude of the difference between ventricular and aortic pressure is larger in diastole than systole.
The boundary condition at $x_1=1.25$ is a zero-pressure outflow/inflow condition, which is implemented by setting
\begin{align*}
    %\label{eq:inflowOutlfow}
    p(x_1,y,z,t) = 0, \qquad \frac{\partial \vv}{\partial n}(x_1,y,z,t) = \zerov.
\end{align*}
This boundary condition allows the fluid to flow from left to right during the opening process, and then flow backwards during the closing process.
%This boundary condition allows both inflow and outflow at the right boundary.
%The right boundary condition is a zero-pressure outflow/inflow condition, implemented by setting
%\begin{align*}
%    %\label{eq:inflowOutlfow}
%    p(x_1,y,z,t) = 0, \qquad \frac{\partial \vv}{\partial n}(x_1,y,z,t) = \zerov.
%\end{align*}
%This boundary condition allows both inflow and outflow at the right boundary.
The remaining boundary conditions on the surface of the leaflets and on the walls of the cylinder are taken as no-slip walls.
The fluid is at rest initially and the leaflets are also at rest in a closed state at the rotation angles $\pm70^{\circ}$.
%and the fluid and leaflets are all at rest initially.  
%As illustrated in Figure~\ref{fig:heartValve2DGrid},
%the leaflets start from a closed state, rotated by $\pm70^{\circ}$ from the $x$-direction.
The time step is determined using a CFL number of 0.9, 
which is based on the advection terms in the fluid momentum equations.

The~\ampRB~sheme used here avoids sub-time-step-iterations between the fluid and solid equations by handling added-mass and added-damping effects in the AMP interface conditions, whereas the previous work for this problem employed various TP-RB schemes~\cite{borazjani2008curvilinear, de2009direct}.
Some discussion of added-mass effects for this heart-valve problem can be found in~\cite{borazjani2008curvilinear}.
In this paper, both a loosely-coupled FSI algorithm (i.e., a TP-RB scheme without sub-iterations) and the strongly-coupled FSI algorithm (i.e., a TP-RB scheme with sub-iterations) are considered.
It was concluded that the TP-RB scheme without sub-iterations is unstable regardless of the size of the time step for the problems involving bodies of small inertia,
and the TP-RB scheme with under-relaxed sub-iterations and the Aitken's acceleration technique is required for certain difficult problems with strong added-mass effects.
Our finding for this problem, and the theoretical results in our previous work~\cite{rbinsmp2017}, are consistent with the results described in~\cite{borazjani2008curvilinear}.

%basic setup
To better isolate the different issues and to clearly illustrate the results, we first consider simulations of the heart-value problem in two dimensions.  This reduction can be viewed as an approximation of the flow along the symmetry plane of the three-dimensional domain.  We note that the two-dimensional problem shares many of the important numerical difficulties as the full three-dimensional problem, such as strong added-mass effects, and this problem is reasonable as a starting point and forms a good benchmark problem.
We then conclude our numerical investigations with simulations of the heart-valve problem in three dimensions.
The three-dimensional benchmark problem is closer to a realistic geometry in which
two half-disk leaflets are used to model the heart valve. 

%and the blood flow in the aorta.
%and serves as a good benchmark test for the FSI algorithms in three dimensions.

\subsubsection{Benchmark problem in two dimensions}
\label{sec:heartValve2D}

%domain and grid
For the two-dimensional version of the heart-valve problem, we consider a fluid channel with horizontal and vertical extent corresponding to the symmetry plane of the three-dimensional domain described previously.  This implies a channel covering the domain $[x_0, x_1]\times[-R_c, R_c]$ with $x_0=-0.5$, $x_1=1.25$ and  $R_c = 1.1$.
Here, each leaflet is a rectangle (with rounded corners) of length and width equal to $1.0$ and $0.1$, respectively.
%The hinges are located at $(x_b, y_b) = (0,\pm0.15)$ and the leaflets are allowed to rotate around the hinges.
The composite grid for the domain, denoted by $\Gcts^{(j)}$, consists of one background Cartesian grid and two boundary-fitted grids as shown in Figure~\ref{fig:heartValve2DGrid}.
The background grid has a grid spacing given by $h_j = 1/(40j)$, where $j$ is a refinement factor.
The boundary-fitted grids are stretched in the normal direction to the surface of the leaflets.
The resulting grid spacing is close to $h_j$ for the grid lines away from the surface and overlapping with the background grid, while the grid spacing in the normal direction is reduced by a factor of approximately~4 near the surface.
The equation of motion in~\eqref{eq:leafletTorque} for each leaflet reduces to
\begin{align}
    I_b \,{\dot\omega}_b = \int_{\Gamma_b} \bigl(\tilde \rv\sp{T}\sigmav\nv\bigr)\,ds + g_{\rm rt}(\theta_b) - B(\theta_b)\, \omega_b, 
\label{eq:leafletTorqueTwoD} 
\end{align}
where $\tilde \rv=(-r_2+y_b, r_1-x_b)^T$ for the line integral around the body $\Gamma_b$ in~\eqref{eq:leafletTorqueTwoD}.  The center of mass of the two leaflets given by $(x_b,y_b)=(0,\pm0.15)$ initially is nearly fixed during the simulations due to very large value for $m_b$ taken in~\eqref{eq:linearAccelerationEquation}.
%where the two-dimensional vector $\tilde \rv$ for the line integral around the body $\Gamma_b$ in~\eqref{eq:leafletTorqueTwoD} is defined as $\tilde \rv=(-r_2+y_b, r_1-x_b)^T$.

%- 
{% --------- cylDrop -----
% --------- START drawContour -----------
\newcommand{\drawContour}[7]{%
\begin{scope}[#1]
\draw(0.0,0) node[anchor=south west,xshift=-4pt,yshift=+0pt] {\trimfig{heartValve/2d/#2}{\figWidth}};
  \draw(1.,7.1) node[draw,fill=white,anchor=west,xshift=4pt,yshift=0pt] {\scriptsize #3};
  \draw(1.7,7.1) node[draw,fill=white,anchor=west,xshift=2pt,yshift=0pt] {\scriptsize #5};
 % colour bar:
\begin{scope}[xshift=28pt,yshift=18pt]
  \draw (\xcb,\ycb) node[anchor=south west,xshift=0.cm,yshift=.5cm,rotate=-90] {\trimfigcb{fig/colourBarLines}{\cbWidth}{\cbHeight}};
  \draw (.8,0) node[anchor=north,xshift=+3pt,yshift=+2pt] {\scriptsize $#6$};
  \draw (4.4,0) node[anchor=north,xshift=+0pt,yshift=+2pt] {\scriptsize $#7$};
\end{scope}
\end{scope}
}
% --------- END drawContour -----------
% -- for colour bar ----
\newcommand{\cbWidth}{.2cm}% colour bar width
\newcommand{\cbHeight}{4cm}% colour bar height
\newcommand{\xcb}{.5cm}% colour bar lower left corner
\newcommand{\ycb}{-.2cm}% colour bar lower left corner
\newcommand{\xcorner}{-0}% x of lower left corner of top leaflet
\newcommand{\ycorner}{-2.8}% y lower left corner of top leaflet
\newcommand{\xbcorner}{-.2}% x of lower left corner of bottom leaflet
\newcommand{\ybcorner}{-3.12}% y of lower left corner of bottom leaflet
\setlength{\ycbTop}{\ycb+\cbHeight}% colour bar top label position
\setlength{\ycbMid}{\ycb+\cbHeight*\real{.5}}% colour bar top label position
\newcommand{\trimfigcb}[3]{\includegraphics[width=#2, height=#3, clip, trim=17cm 2.35cm 1.65cm 2.35cm]{#1}}
\newcommand{\figWidth}{7.5cm}
\newcommand{\trimfig}[2]{\trimh{#1}{#2}{.14}{.19}{.025}{.1}}
\newcommand{\trimfiga}[2]{\trimh{#1}{#2}{.27}{.27}{.025}{.1}}
\newcommand{\labelsize}{\small}
\begin{figure}[htb]
\begin{center}
\resizebox{16cm}{!}{% START resize box
\begin{tikzpicture}[scale=1]
  \useasboundingbox (-.5,1) rectangle (18,7.75);  % set the bounding box (so we have less surrounding white space)
  \draw(2,0) node[anchor=south west,xshift=-4pt,yshift=+0pt] {\trimfig{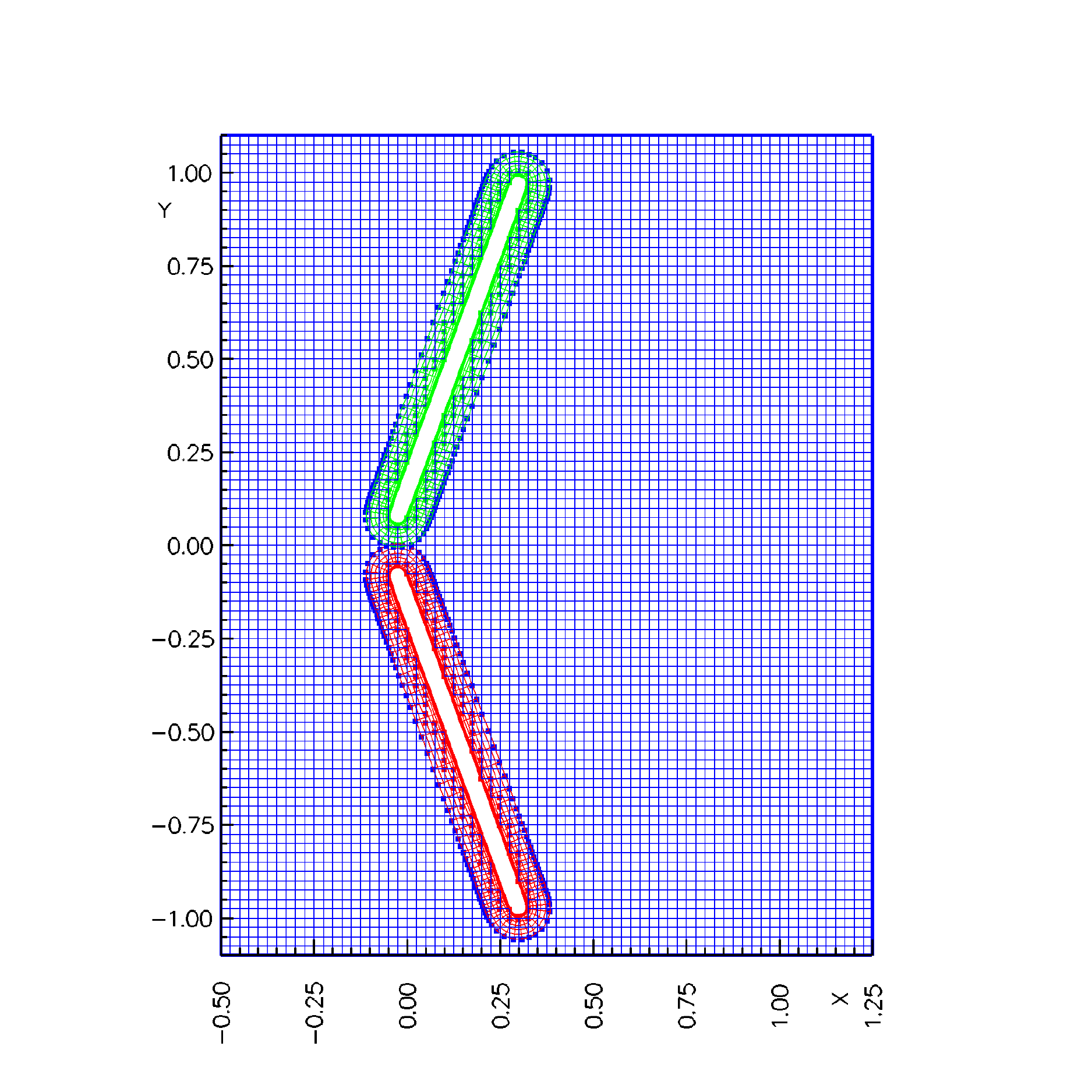}{\figWidth}};
  \drawContour{xshift=9.2cm,yshift=0cm}{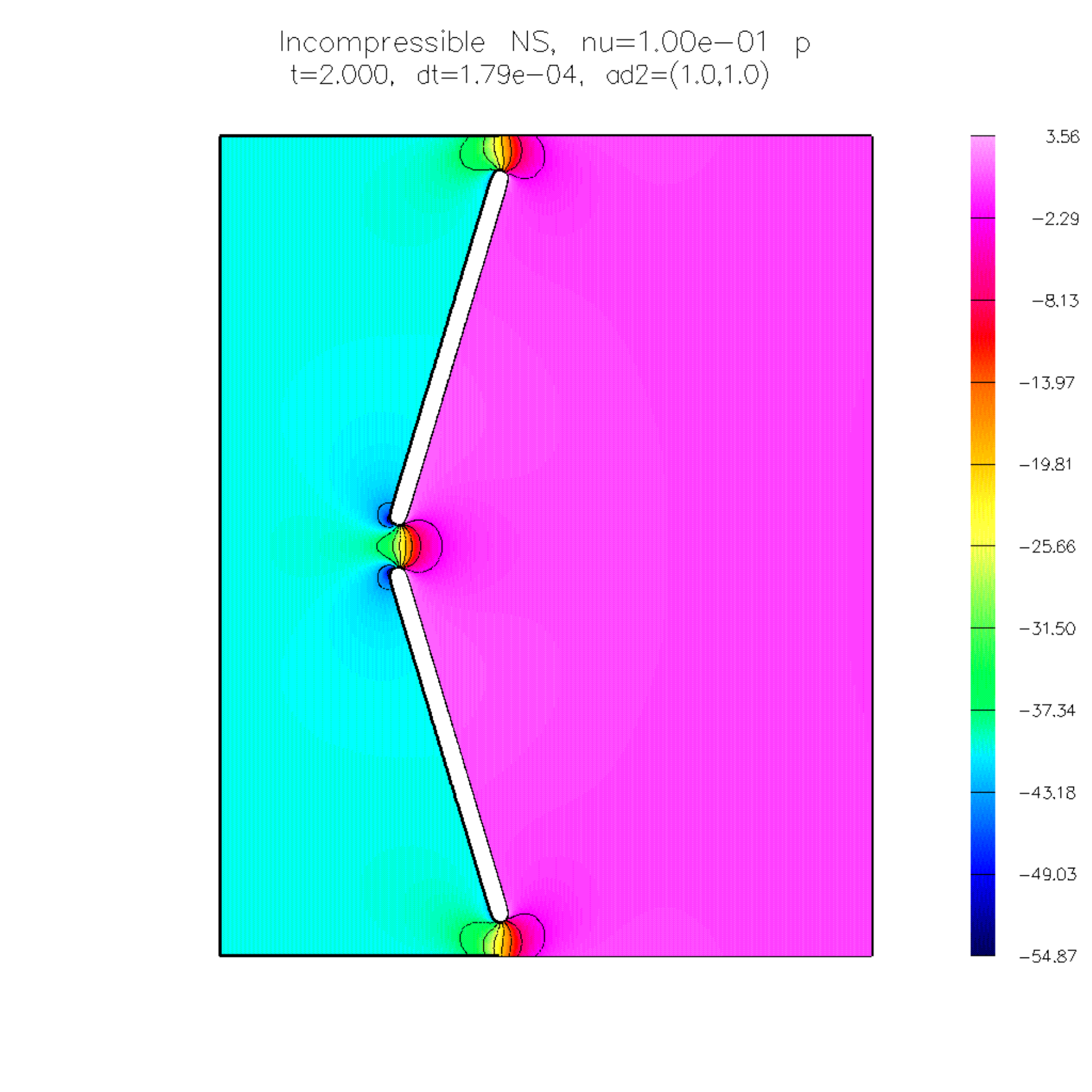}{$p$}{$p$}{$t=2$}{-53.9}{3.0};
% grid:
%\draw[step=1cm,gray] (0,0) grid (18,7.75);
\end{tikzpicture}
}% end resize box
\end{center}
  \caption{Mechanical heart valve in 2D. %Left: Geometrical configuration at $t=0$.
  Left: Composite grid $\Gcts^{(1)}$ at $t=0$ (coarse grid). 
  Right: contours of the pressure at $t=2$ using grid $\Gcts^{(8)}$.
}
  \label{fig:heartValve2DGrid}
\end{figure}
}

The instantaneous streamlines of the numerical solution at selected times are presented in Figure~\ref{fig:heartValve2Dstreamlines}, and the corresponding rotation angle, angular velocity and acceleration of the
upper leaflet are given in Figure~\ref{fig:heartValve2DCurves}.  (The latter figure also indicates grid convergence which is discussed below.)
%complemented by the pressure at $t = 2$ in  Figure~\ref{fig:heartValve2DGrid}
%as well as the time history of the motion of the top leaflet in Figure~\ref{fig:heartValve2DCurves}.
At early times, the pressure at the left boundary increases according to the prescibed profile in~\eqref{eq:pressureProfile}, and this drives
the fluid to the left causing the leaflets to open (see the streamline plots at $t = 0.2$ and $t=0.5$).
%The speed of the flow also increases during this stage.
At approximately $t = 0.6$, the leaflets approach the lower limit of the rotation angle given by $\pm\theta_{\rm min}$,
and they exhibit a small oscillation due to combined influences of the fluid, repulsive force and damping term as shown in the time history of the leaflet motion in Figure~\ref{fig:heartValve2DCurves}.
By about $t = 0.8$, the fluid force on the leaflets balances the repulsive force, and hence the leaflets appear to reach a steady state. 
%At the same time all the flow passes two leaflets and move from the left boundary  to the right quickly.
Through this time, starting at $t=0.5$, the applied pressure on the fluid at the left boundary has been decreasing, and after $t=1$ the pressure is negative.
%At $t = 1$, the left and right pressure are both zero but the fluid continue to move from left to right.
%After $t = 1$, the pressure gradient is reversed and the flow quickly turns around.
The effect of this negative pressure is first seen in the plot of the instantaneous streamlines at $t=1.1$, where
a relatively small region of reverse flow is observed at the right-hand side of the channel near the centerline, $y=0$.
%At $t = 1.1$, some of the flow between the leaflets has a reversed direction,
%while most of the flow still moves from left to right.
By $t=1.2$, the flow is in the midst of reversal and several vortices have formed in the channel as a result. 
Meanwhile, the leaflets have started to swing back towards their closed position as a result of the reversed flow.
The leaflets approach their closed position at approximately $t=1.6$, and again an oscillation in the rotation
angle about the limit state is observed. 
This motion is similar to that during the opening stage, but the oscillation is more severe since the leaflets close faster due to a larger magnitude of the pressure difference in the prescribed pressure profile during the closing process.
After $t = 1.6$ the oscillations are significantly reduced by the damping term,
and the leaflets reach another steady state corresponding to a closed valve.

{% ------ STREAMLINES ------
\newcommand{\drawContour}[7]{%
\begin{scope}[#1]
\draw(0.0,0) node[anchor=south west,xshift=-4pt,yshift=+0pt] {\trimfig{heartValve/2d/#2}{\figWidth}};
\draw(.7,4.75) node[draw,fill=white,anchor=west,xshift=2pt,yshift=-1pt] {\scriptsize #5};
 % colour bar:
\begin{scope}[xshift=5pt,yshift=1pt]
  \draw (\xcb,\ycb) node[anchor=south west,xshift=0.cm,yshift=.5cm,rotate=-90] {\trimfigcb{fig/colourBarLines}{\cbWidth}{\cbHeight}};
  \draw (.8,0) node[anchor=north,xshift=+3pt,yshift=+2pt] {\scriptsize $#6$};
  \draw (3.9,0) node[anchor=north,xshift=+0pt,yshift=+2pt] {\scriptsize $#7$};
\end{scope}
\end{scope}
}
% -- for colour bar ----
\newcommand{\cbWidth}{.2cm}% colour bar width
\newcommand{\cbHeight}{3.4cm}% colour bar height
\newcommand{\xcb}{.5cm}% colour bar lower left corner
\newcommand{\ycb}{-.2cm}% colour bar lower left corner
\setlength{\ycbTop}{\ycb+\cbHeight}% colour bar top label position
\setlength{\ycbMid}{\ycb+\cbHeight*\real{.5}}% colour bar top label position
\newcommand{\trimfigcb}[3]{\includegraphics[width=#2, height=#3, clip, trim=17cm 2.35cm 1.65cm 2.35cm]{#1}}
%
% ============================ DRAW ===================================
\newcommand{\figWidth}{4.2cm}
\newcommand{\trimfig}[2]{\trimw{#1}{#2}{.16}{.16}{.1}{.1}}
\begin{figure}[htb]
\begin{center}
%\resizebox{16cm}{!}{% START resize box
\begin{tikzpicture}[scale=1]
  \useasboundingbox (0.0,.25) rectangle (16.,11);  % set the bounding box (so we have less surrounding white space)
% 
%  top row: 
\begin{scope}[yshift=5.8cm]
 \drawContour{xshift=-0.5cm,yshift=0.0cm}{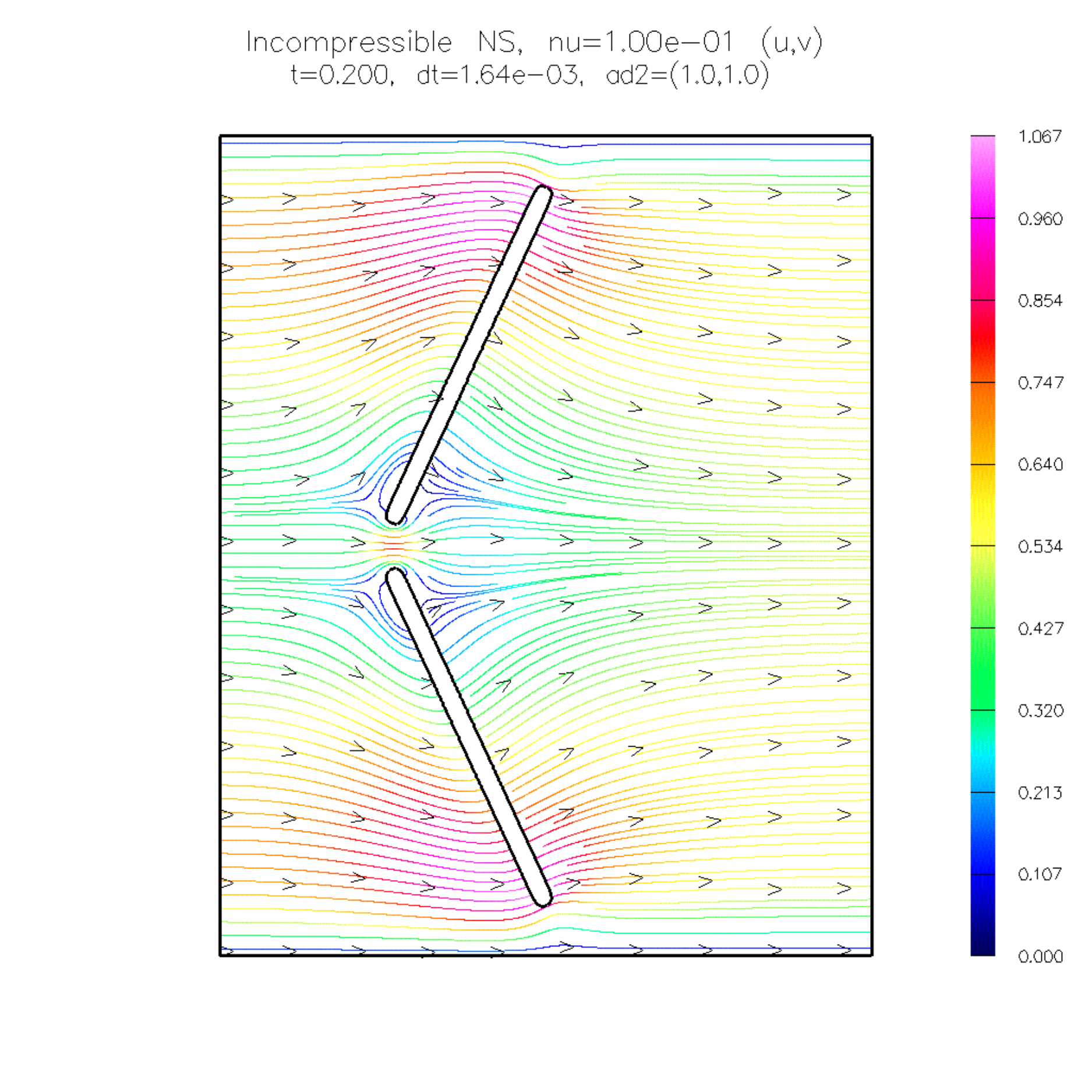}{$p$}{$p$}{$t=0.2$}{0.0}{1.06};
 \drawContour{xshift= 3.5cm,yshift=0.0cm}{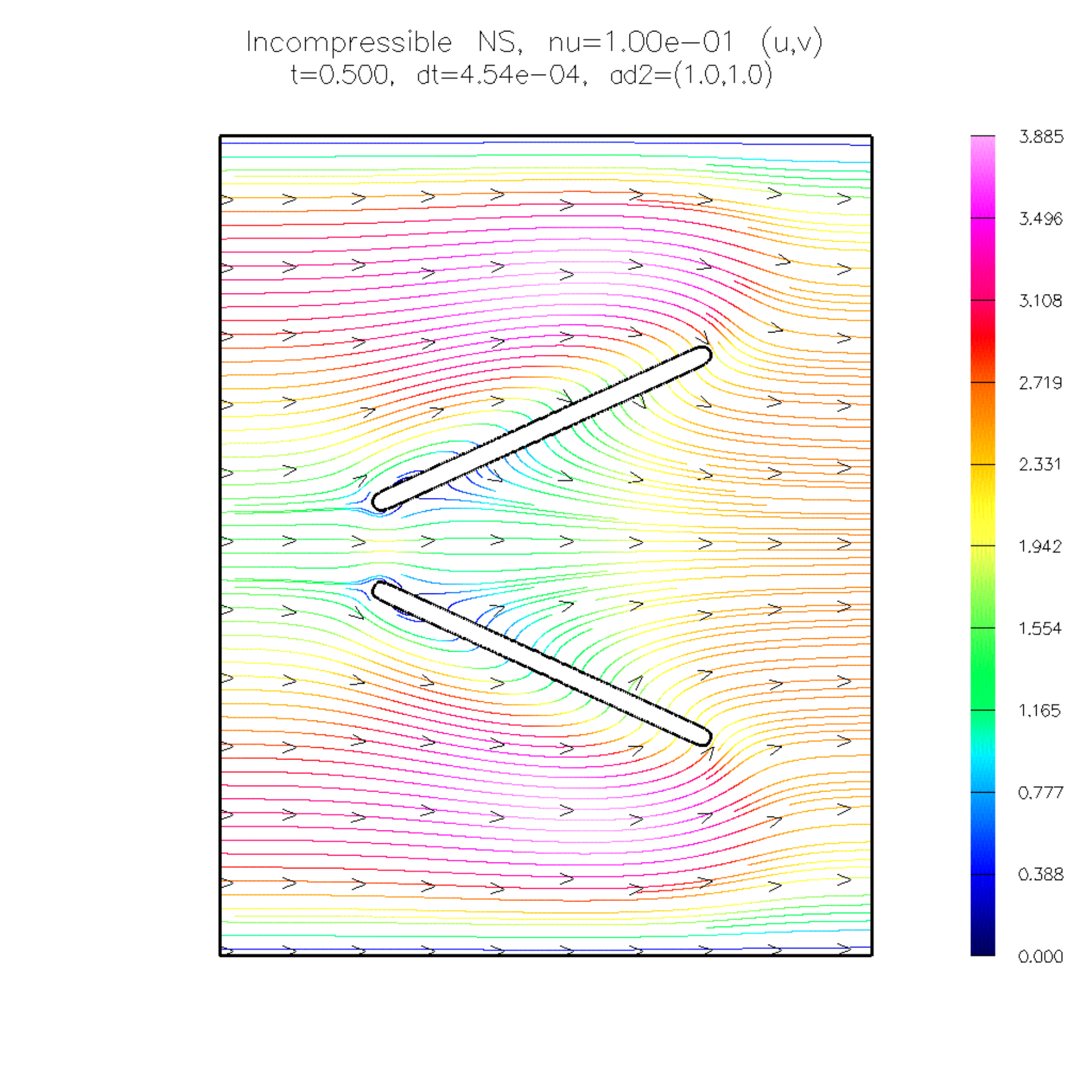}{$p$}{$p$}{$t=0.5$}{0.0}{3.88};
 \drawContour{xshift= 7.5cm,yshift=0.0cm}{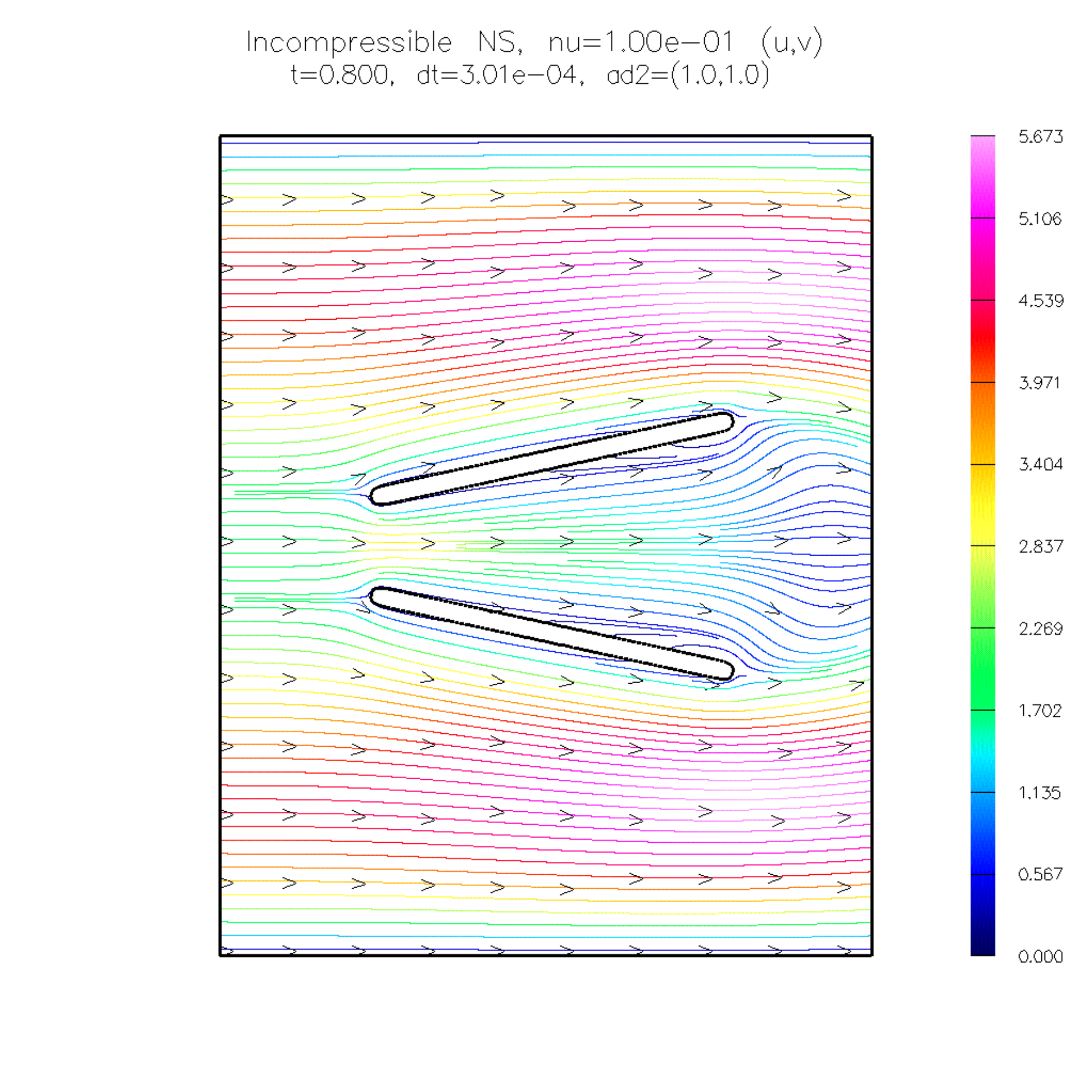}{$p$}{$p$}{$t=0.8$}{0.0}{5.66};
 \drawContour{xshift=11.5cm,yshift=0.0cm}{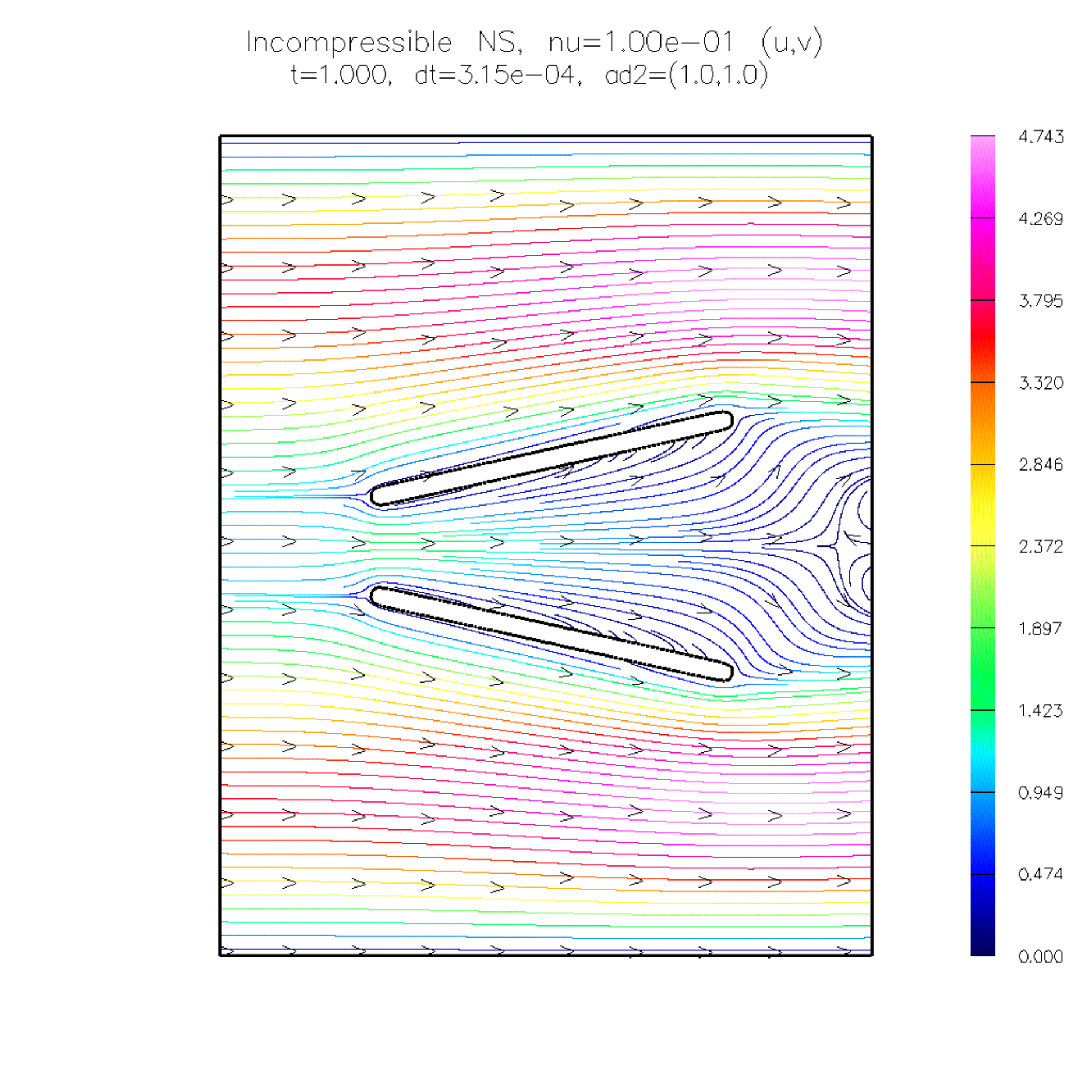}{$p$}{$p$}{$t=1.0$}{0.0}{4.74};
\end{scope}
%  bottom row: 
 \drawContour{xshift=-0.5cm,yshift=0.0cm}{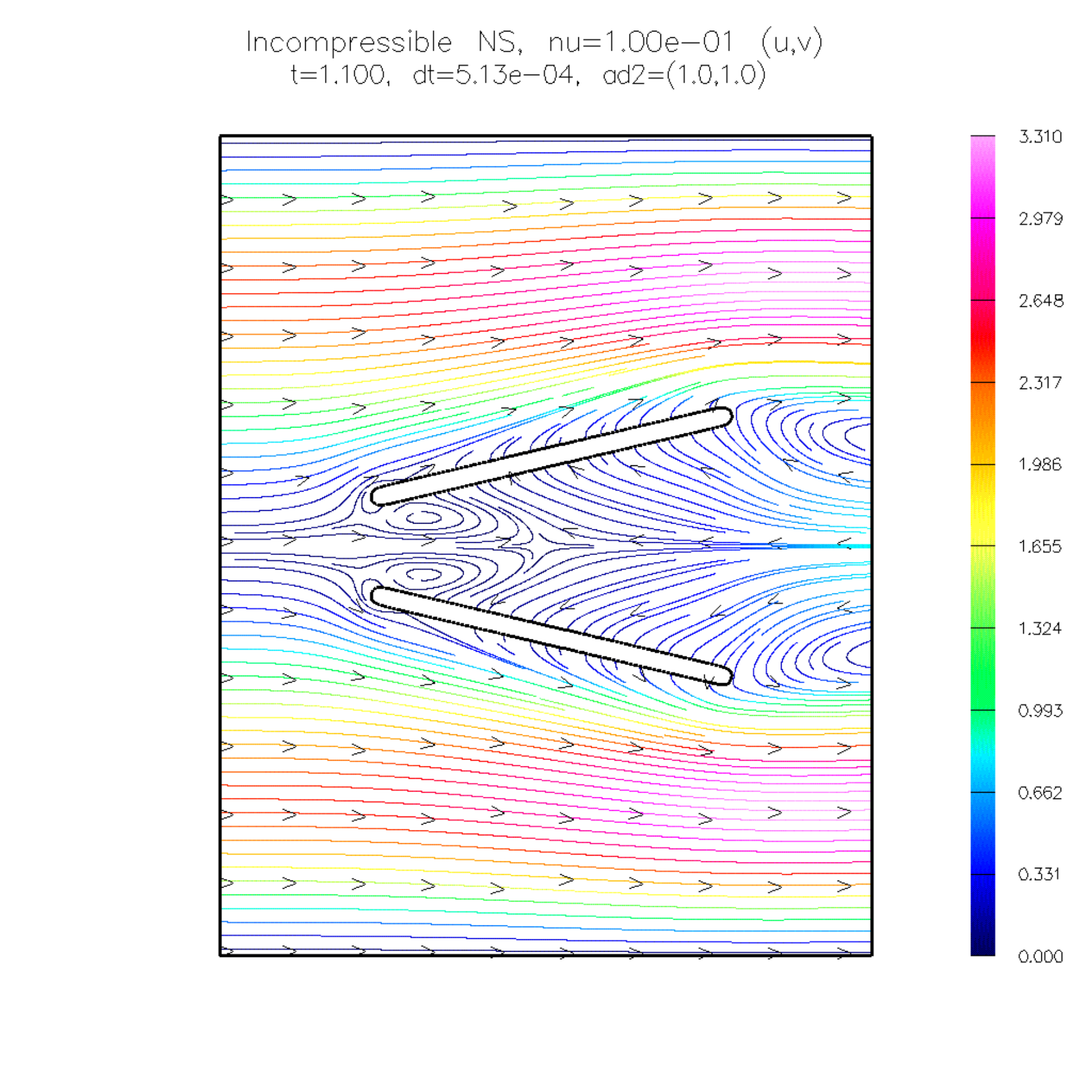}{$p$}{$p$}{$t=1.1$}{0.0}{3.30};
 \drawContour{xshift= 3.5cm,yshift=0.0cm}{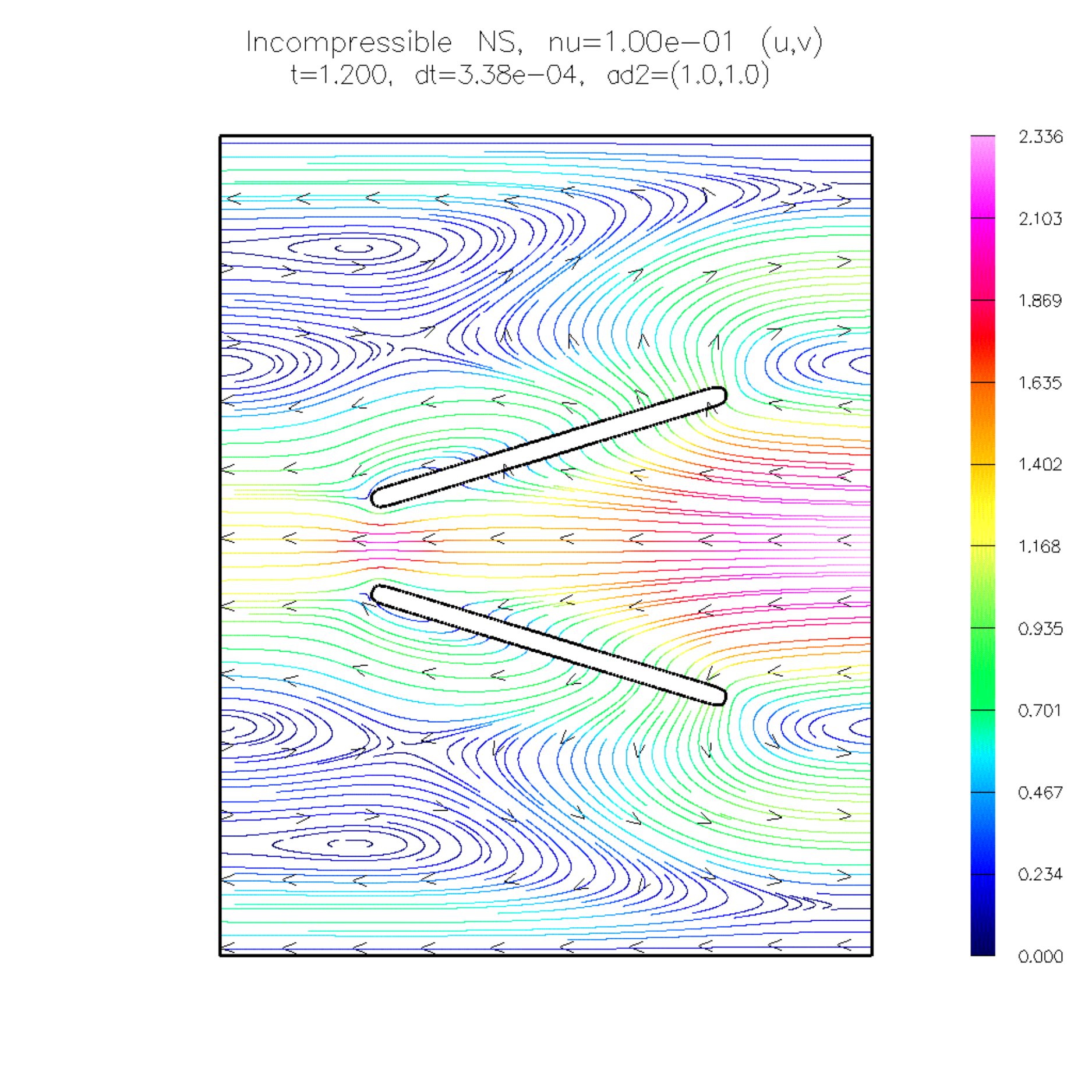}{$p$}{$p$}{$t=1.2$}{0.0}{2.33};
 \drawContour{xshift= 7.5cm,yshift=0.0cm}{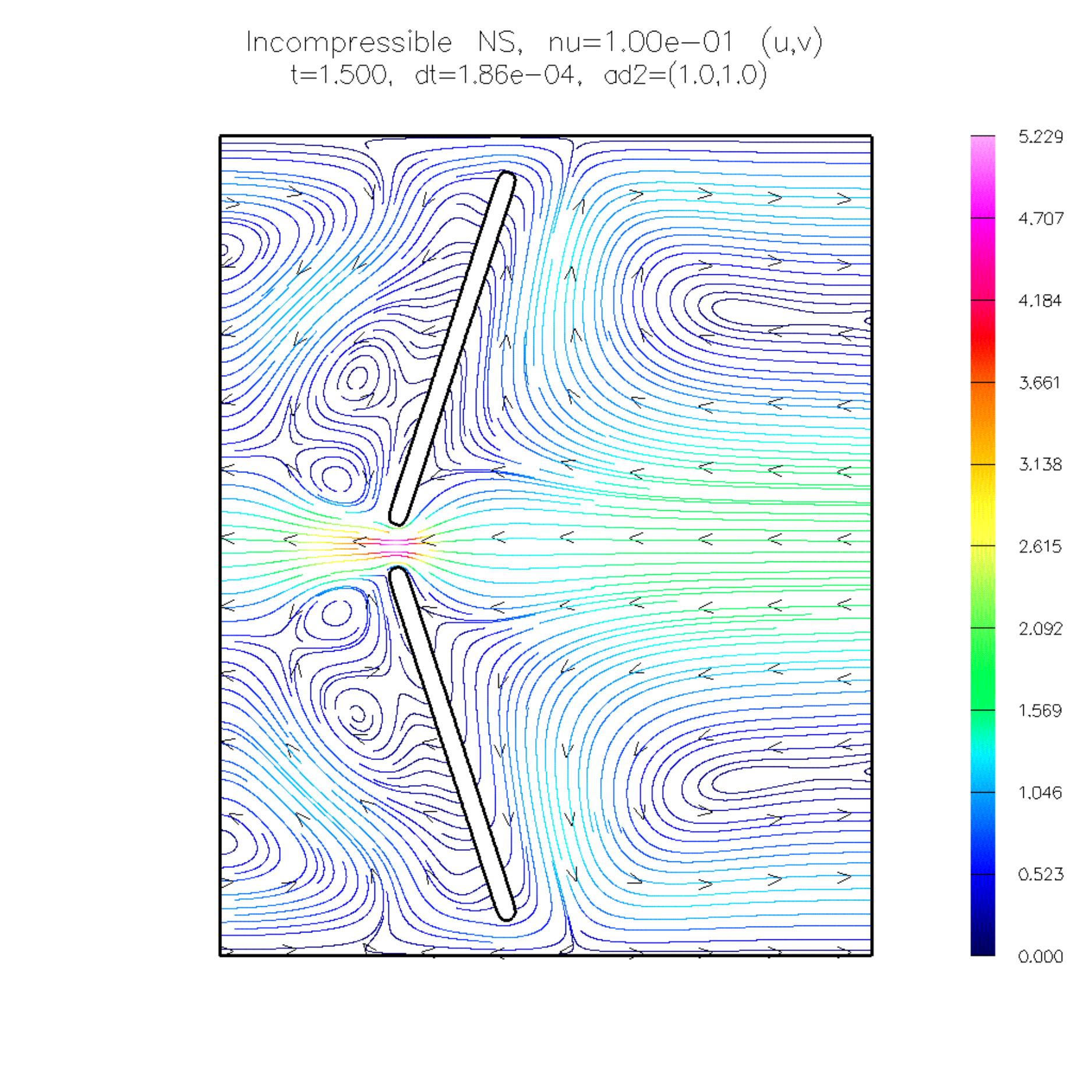}{$p$}{$p$}{$t=1.5$}{0.0}{5.15};
 \drawContour{xshift=11.5cm,yshift=0.0cm}{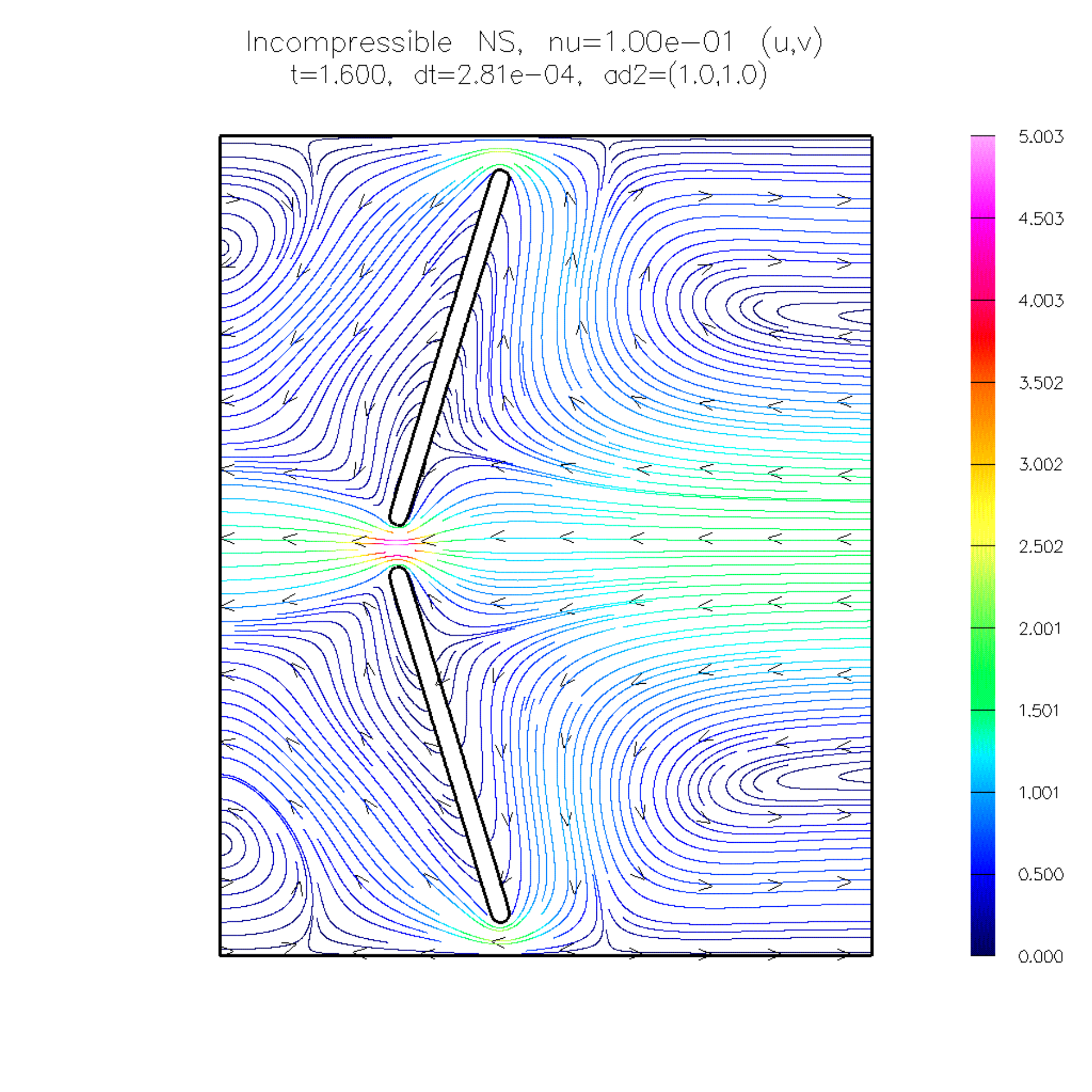}{$p$}{$p$}{$t=1.6$}{0.0}{4.93};
% 
% ---------------
% grid:
%\draw[step=1cm,gray] (0,0) grid (16,11);
\end{tikzpicture}
%}% end resize box
\end{center}
  \caption{Mechanical heart valve in 2D. Computed instantaneous streamlines (colored by the flow speed) at selected times using the composite grid $\Gcts^{(8)}$.  
}
  \label{fig:heartValve2Dstreamlines}
\end{figure}
}

{% ------ MATLAB CURVES LIGHT BODY ------
\newcommand{\figWidth}{8.cm}
\newcommand{\trimfig}[2]{\trimFig{#1}{#2}{.0}{.0}{.0}{.0}}
\newcommand{\figWidthz}{4.cm}
\newcommand{\figWidthb}{3.5cm}
\newcommand{\trimfigz}[2]{\trimFig{#1}{#2}{.0}{.0}{.24}{.0}}
\begin{figure}[htb]
\begin{center}
\resizebox{14cm}{!}{% START resize box
\begin{tikzpicture}[scale=1]
  \useasboundingbox (0.0,0) rectangle (16,13.8);  % set the bounding box (so we have less surrounding white space)
  \draw(-.5,6.7) node[anchor=south west,xshift=0pt,yshift=+0pt] {\trimfig{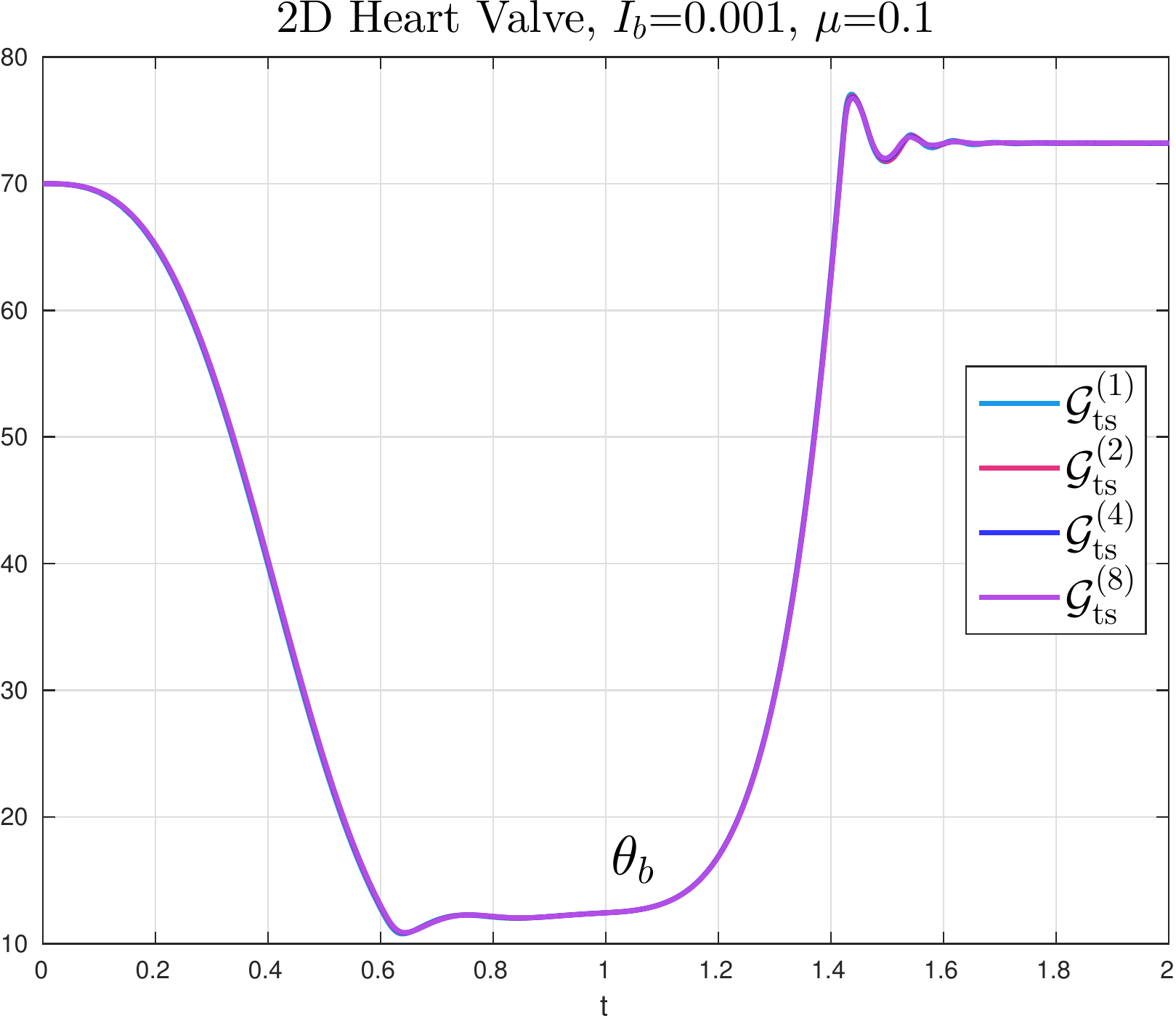}{\figWidth}};
  %zoom
  \draw(1.4, 10.4) node[anchor=south west,xshift=0pt,yshift=+0pt] {\trimfigz{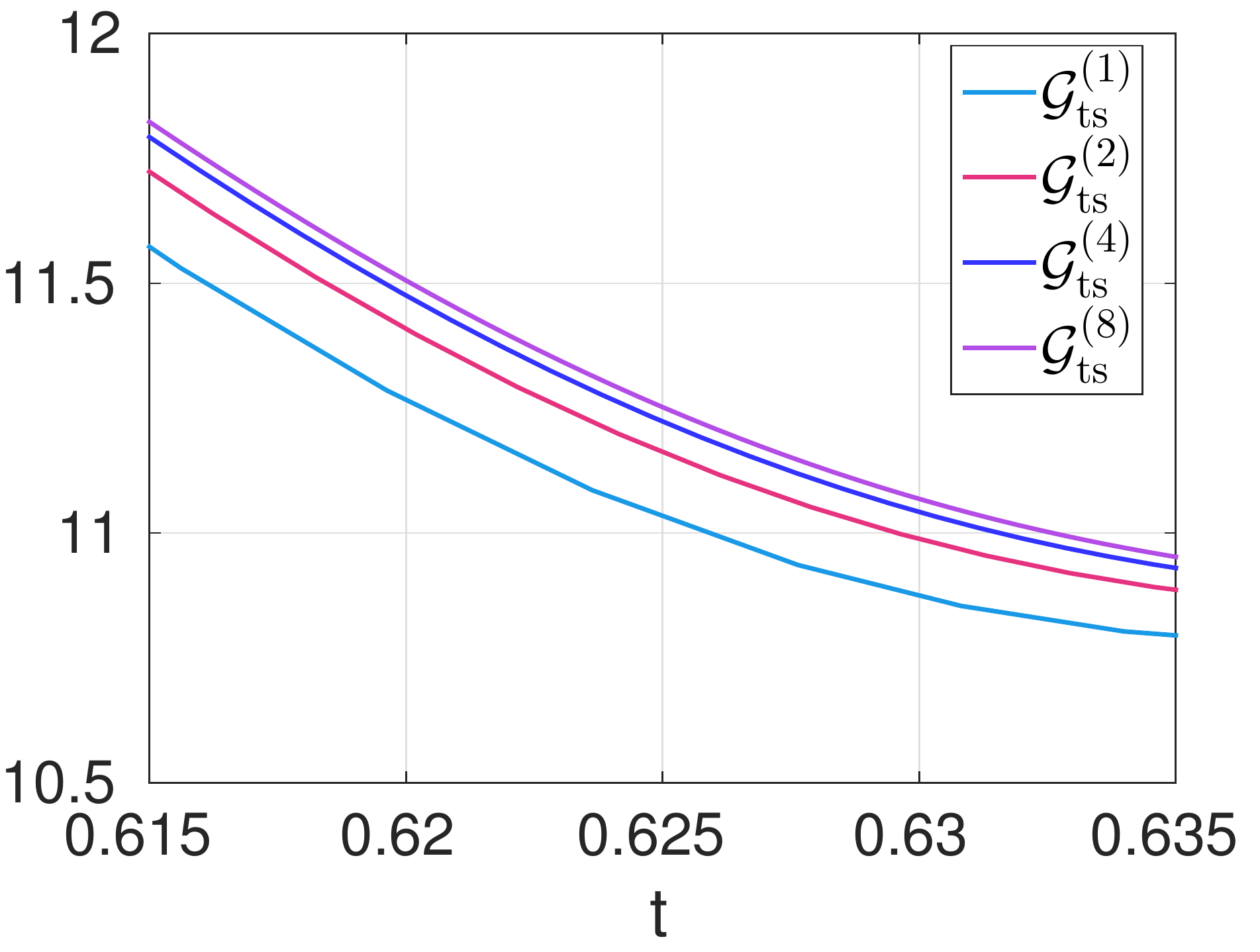}{\figWidthb}};
  \draw(8.0,6.7) node[anchor=south west,xshift=0pt,yshift=+0pt] {\trimfig{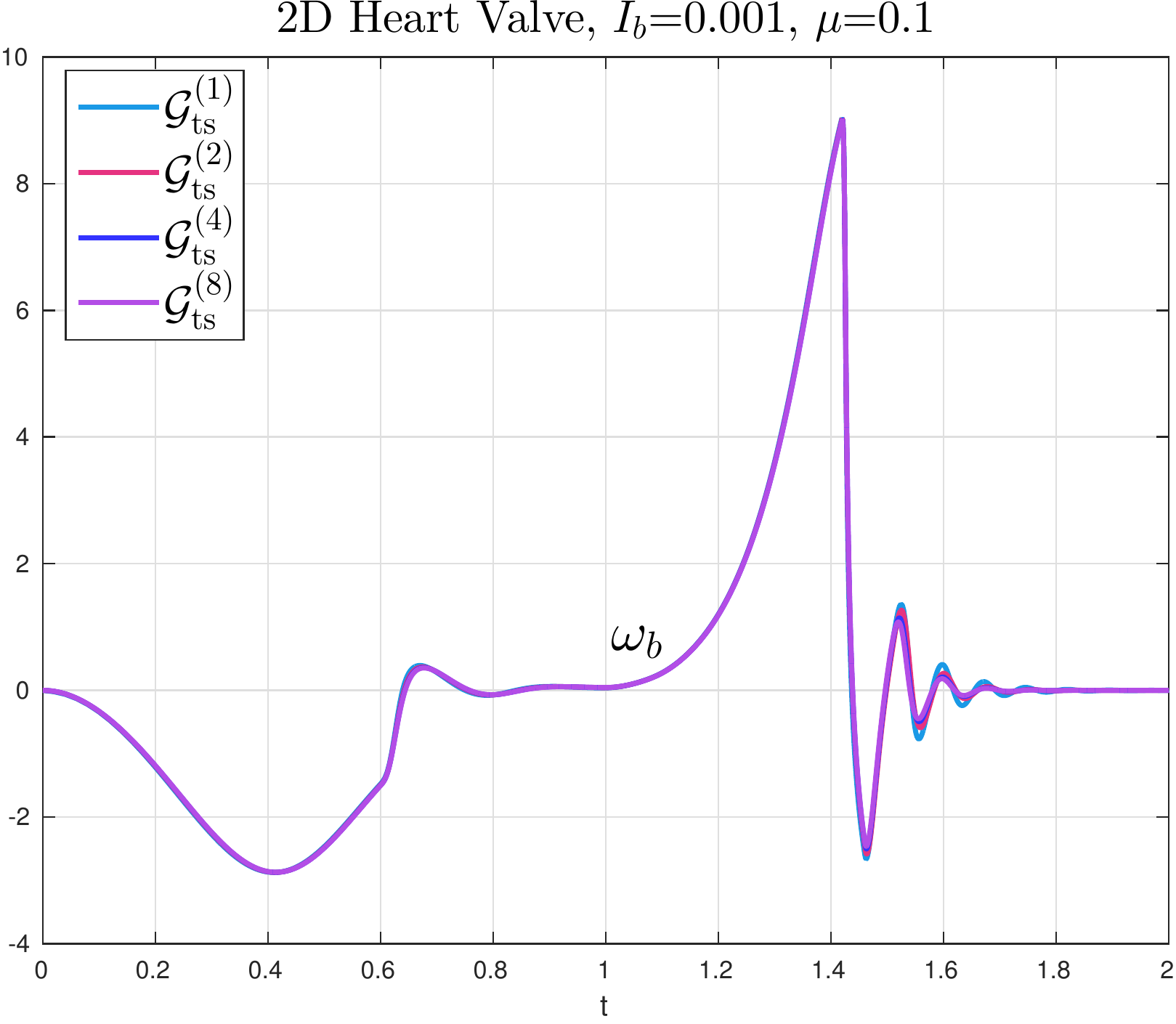}{\figWidth}};
  %zoom
  \draw(9.7, 10.3) node[anchor=south west,xshift=0pt,yshift=+0pt] {\trimfigz{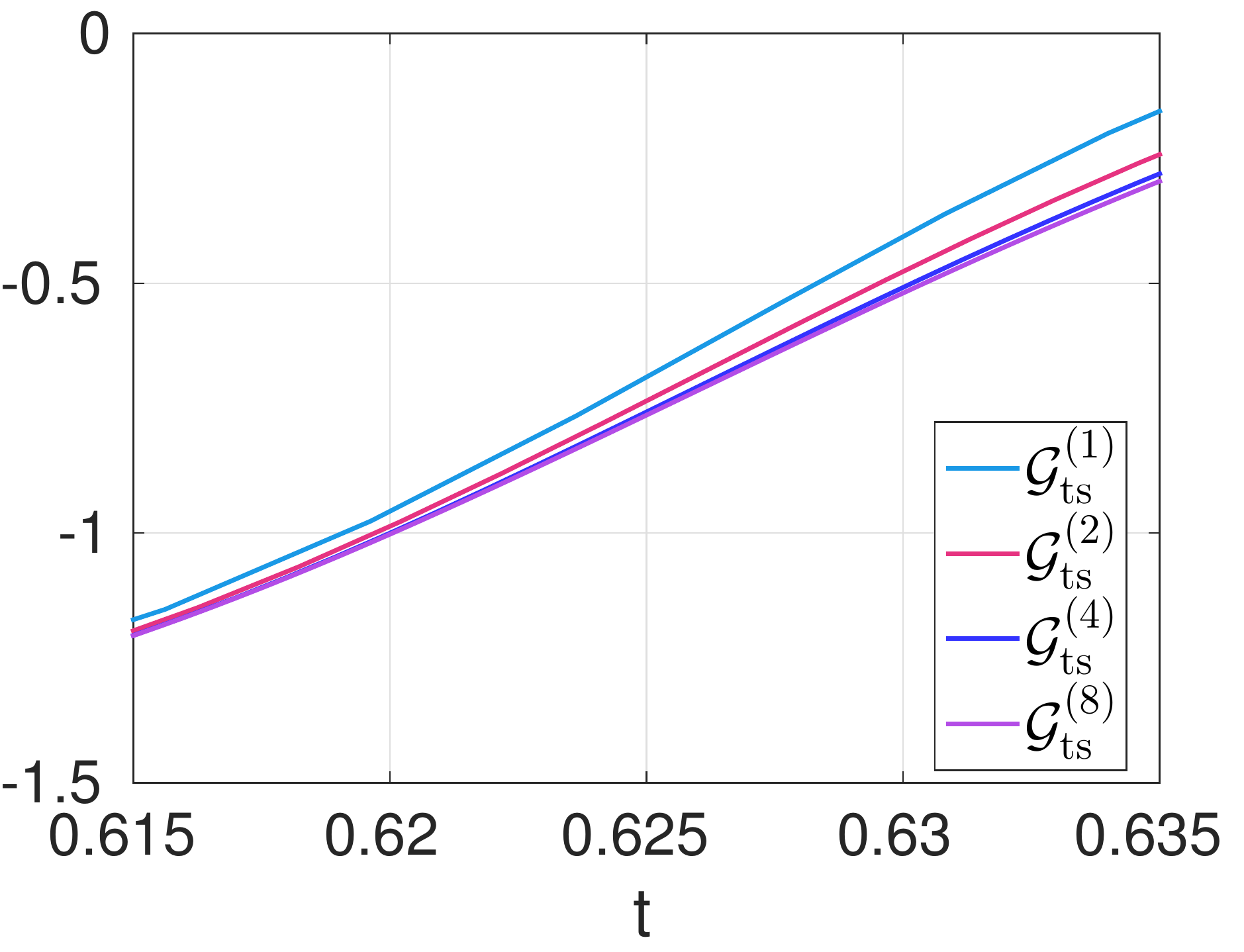}{\figWidthb}};
  %---
  \draw(-.5,-.75) node[anchor=south west,xshift=0pt,yshift=+0pt] {\trimfig{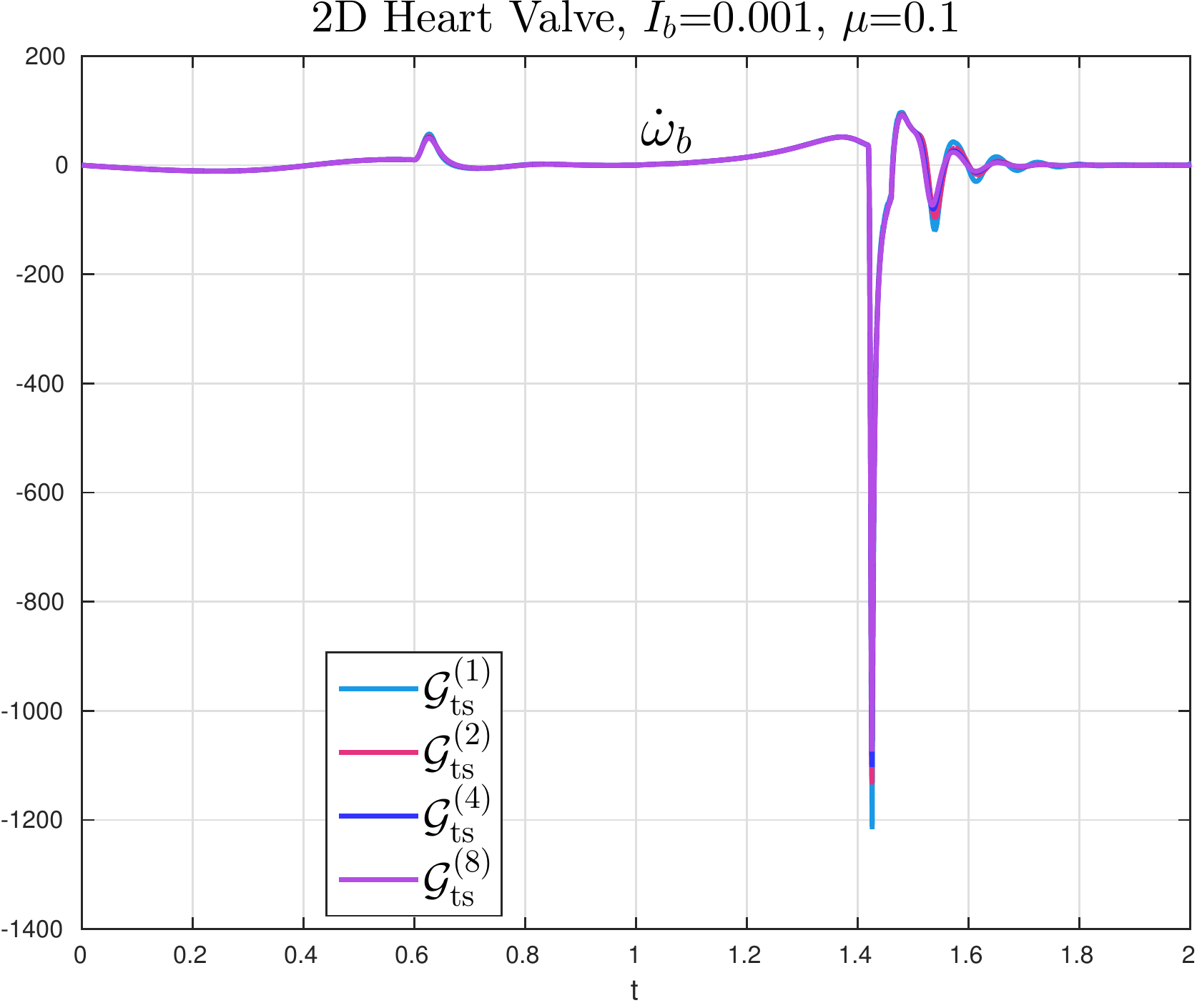}{\figWidth}};
  %zoom
  \draw(.5, 1.7) node[anchor=south west,xshift=0pt,yshift=+0pt] {\trimfigz{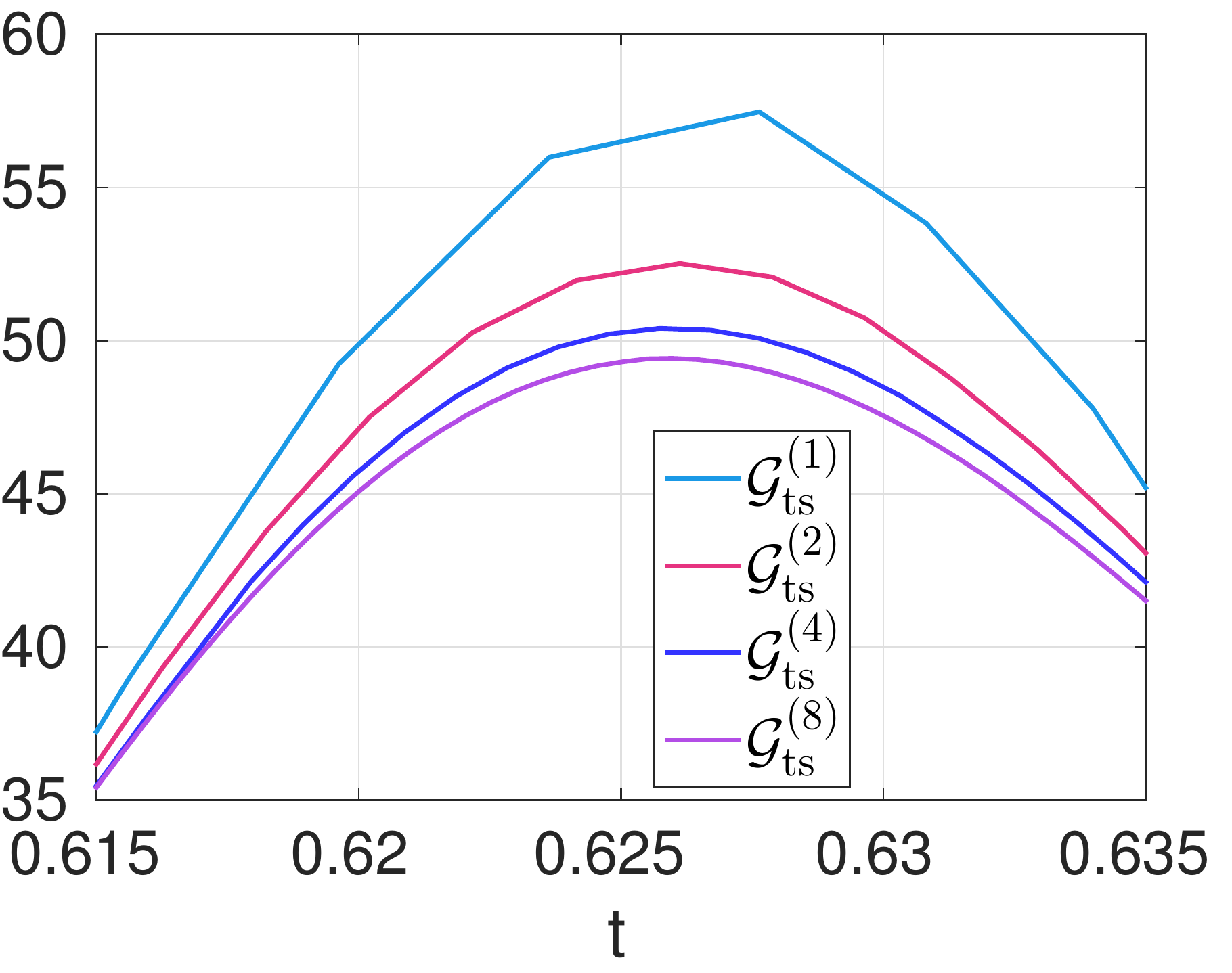}{\figWidthz}};
  \draw(8,-.75) node[anchor=south west,xshift=0pt,yshift=+0pt] {\trimfig{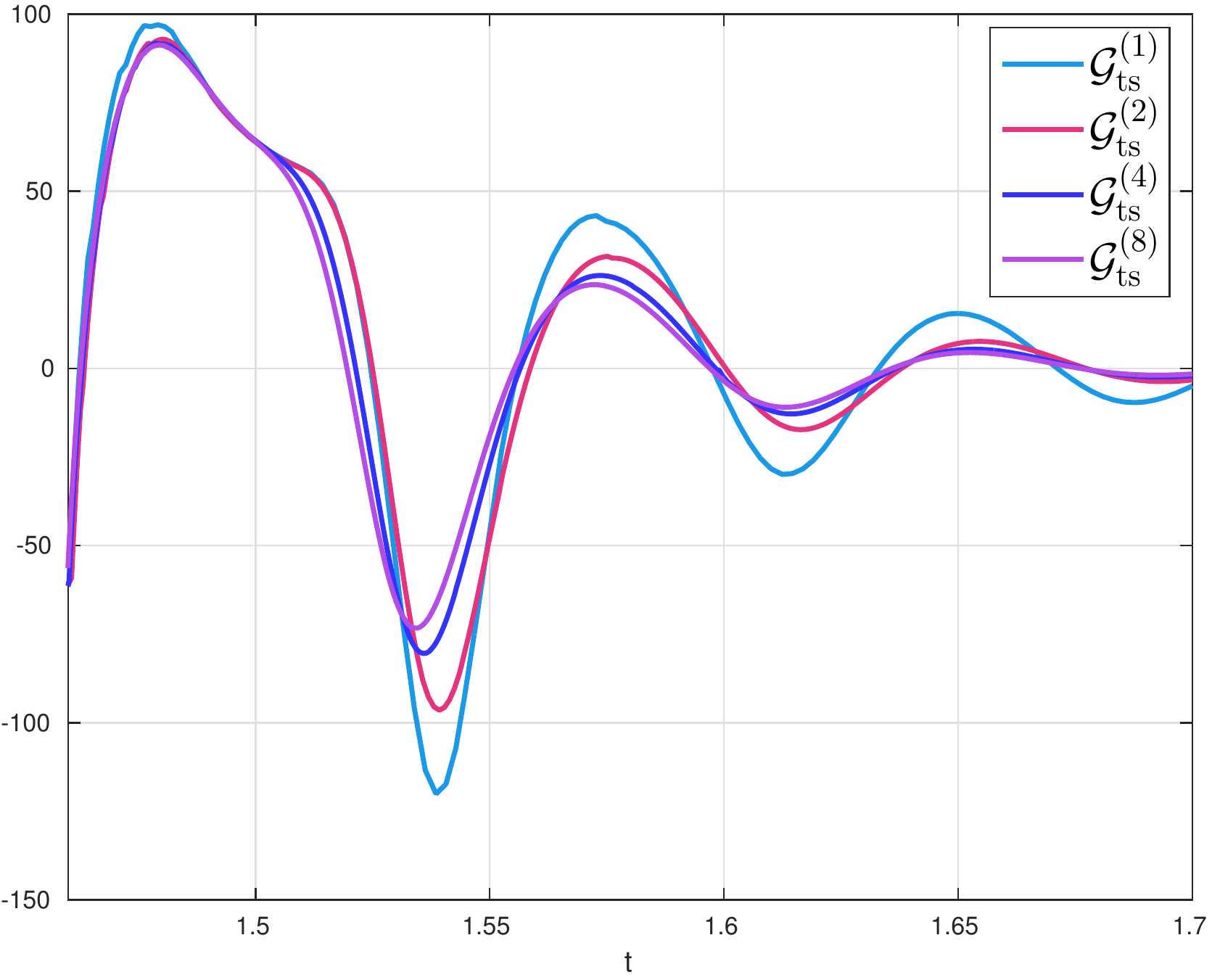}{\figWidth}};
% grid:
%\draw[step=1cm,gray] (0,0) grid (16,13.8);
\end{tikzpicture}
} % END resize box
\end{center}
  \caption{Mechanical heart valve in 2D. Time history of the rotation angle (top left),  angular velocity (top right) and acceleration (bottom left) of the top leaflet.
      Several zoomed views are presented including the zoomed view of the accelerations at the later time (bottom right).
      Results are shown from calculations using the composite grid, $\Gcts^{(j)}$, $j = 1, 2, 4, 8$.
  }
  \label{fig:heartValve2DCurves}
\end{figure}
}

{% ------ STREAMLINES ------
\newcommand{\drawContour}[7]{%
\begin{scope}[#1]
\draw(0.0,0) node[anchor=south west,xshift=-4pt,yshift=+0pt] {\trimfig{heartValve/2d/#2}{\figWidth}};
\draw(.7,4.75) node[draw,fill=white,anchor=west,xshift=2pt,yshift=-1pt] {\scriptsize #5};
\draw(3.4,4.65) node[draw,fill=white,anchor=west,xshift=2pt,yshift=-1pt] {\scriptsize #4};
 % colour bar:
\begin{scope}[xshift=5pt,yshift=1pt]
  \draw (\xcb,\ycb) node[anchor=south west,xshift=0.cm,yshift=.5cm,rotate=-90] {\trimfigcb{fig/colourBarLines}{\cbWidth}{\cbHeight}};
  \draw (.8,0) node[anchor=north,xshift=+3pt,yshift=+2pt] {\scriptsize $#6$};
  \draw (3.9,0) node[anchor=north,xshift=+0pt,yshift=+2pt] {\scriptsize $#7$};
\end{scope}
\end{scope}
}
% -- for colour bar ----
\newcommand{\cbWidth}{.2cm}% colour bar width
\newcommand{\cbHeight}{3.4cm}% colour bar height
\newcommand{\xcb}{.5cm}% colour bar lower left corner
\newcommand{\ycb}{-.2cm}% colour bar lower left corner
\setlength{\ycbTop}{\ycb+\cbHeight}% colour bar top label position
\setlength{\ycbMid}{\ycb+\cbHeight*\real{.5}}% colour bar top label position
\newcommand{\trimfigcb}[3]{\includegraphics[width=#2, height=#3, clip, trim=17cm 2.35cm 1.65cm 2.35cm]{#1}}
%
% ============================ DRAW ===================================
\newcommand{\figWidth}{4.2cm}
\newcommand{\trimfig}[2]{\trimw{#1}{#2}{.16}{.16}{.1}{.1}}
\begin{figure}[htb]
\begin{center}
%\resizebox{16cm}{!}{% START resize box
\begin{tikzpicture}[scale=1]
  \useasboundingbox (0.0,.25) rectangle (16.,5.1);  % set the bounding box (so we have less surrounding white space)
% 
%  bottom row: 
  \drawContour{xshift=-0.5cm,yshift=0.0cm}{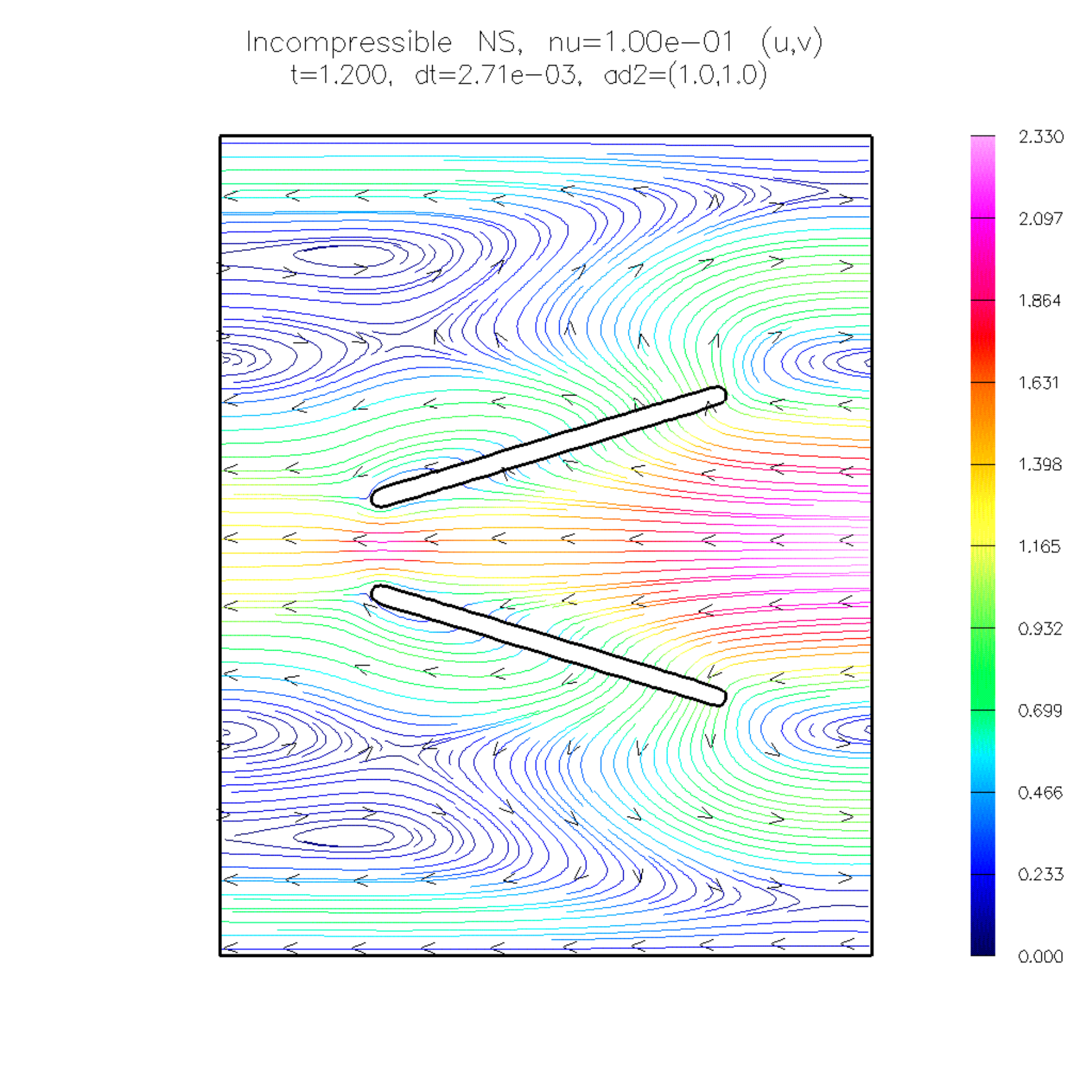}{$p$}{$\Gcts^{(1)}$}{$t=1.2$}{0.0}{2.33};
  \drawContour{xshift= 3.5cm,yshift=0.0cm}{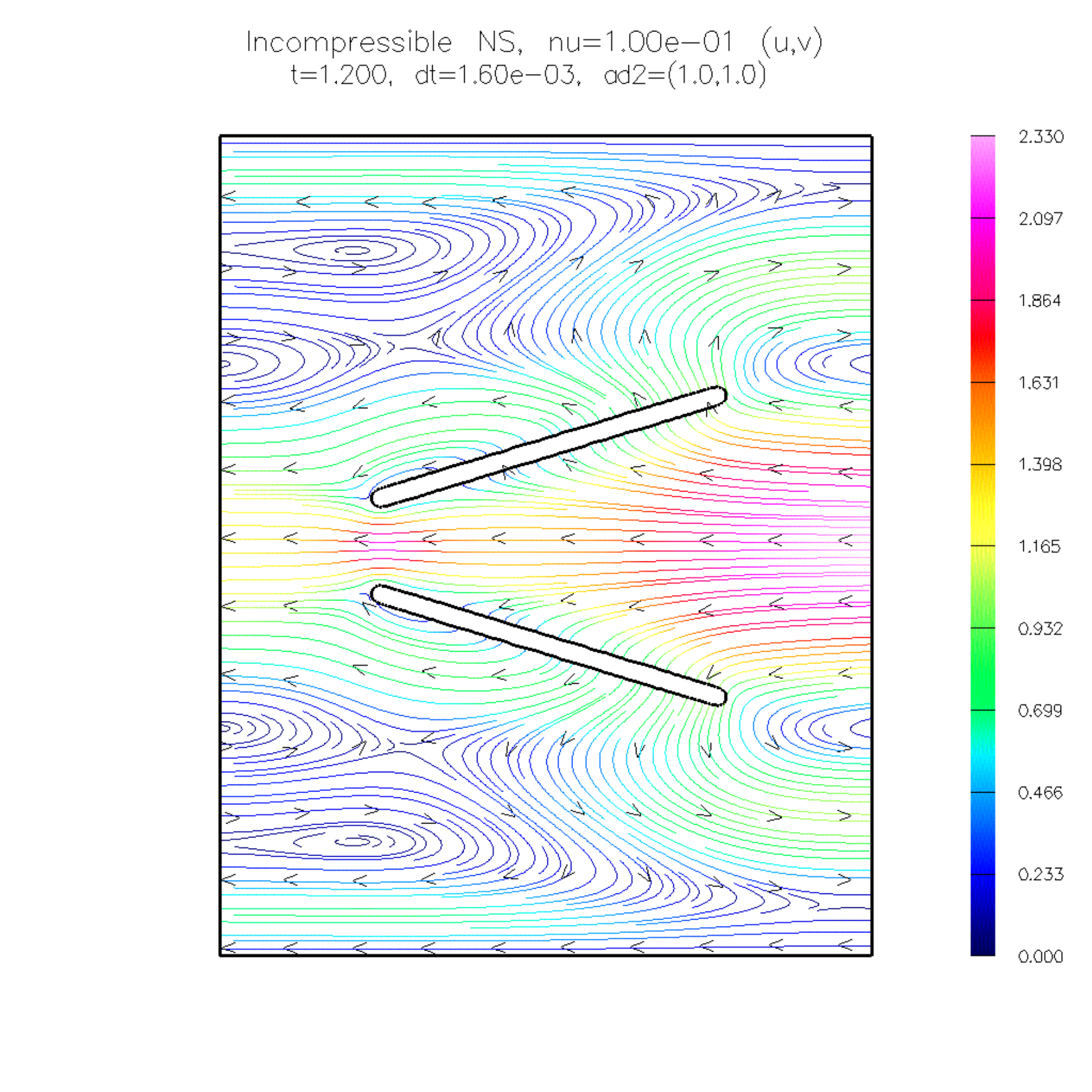}{$p$}{$\Gcts^{(2)}$}{$t=1.2$}{0.0}{2.33};
 \drawContour{xshift= 7.5cm,yshift=0.0cm}{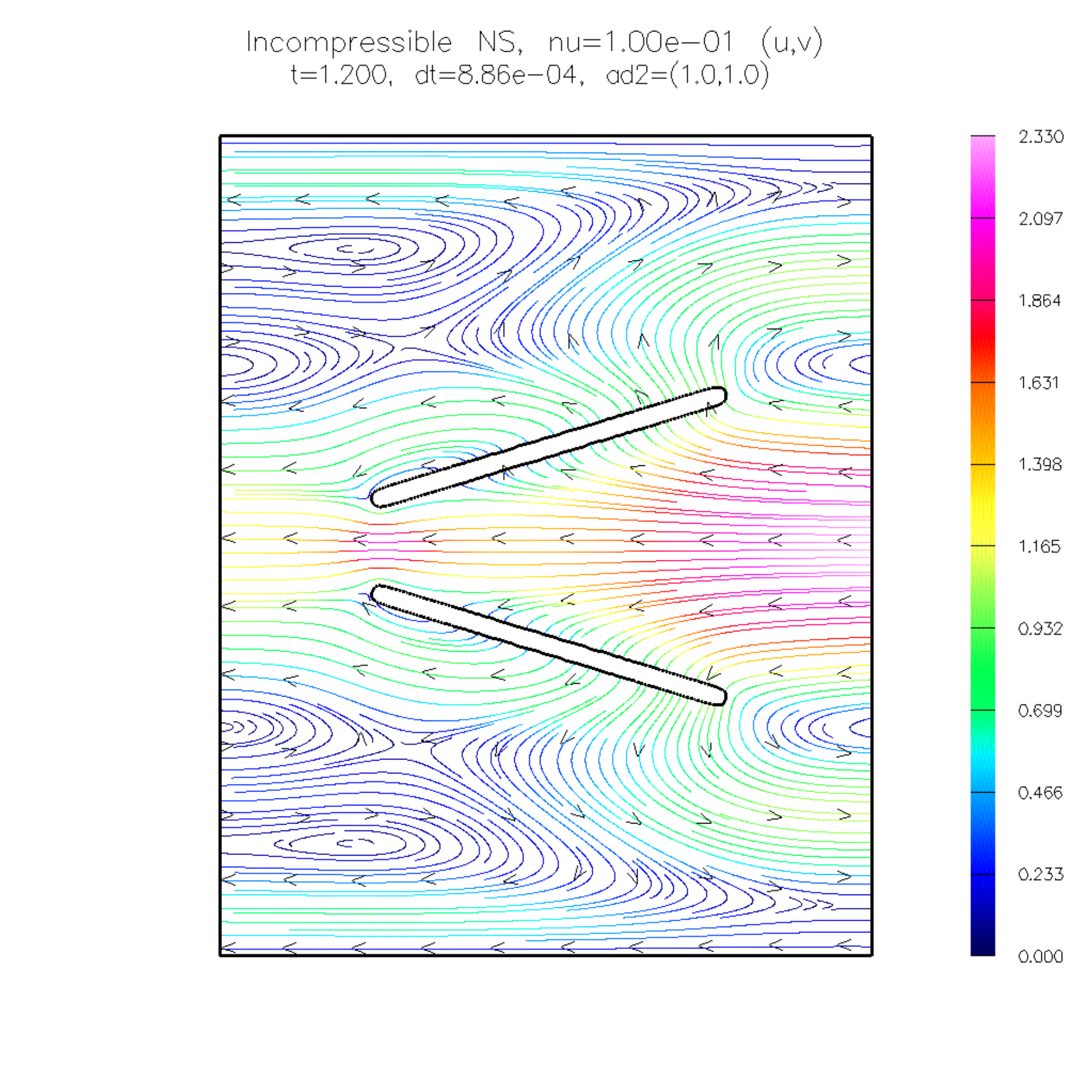}{$p$}{$\Gcts^{(4)}$}{$t=1.2$}{0.0}{2.33};
 \drawContour{xshift=11.5cm,yshift=0.0cm}{heartValve2dt1p2streamline}{$p$}{$\Gcts^{(8)}$}{$t=1.2$}{0.0}{2.33};
% 
% ---------------
% grid:
%\draw[step=1cm,gray] (0,0) grid (16,5.1);
\end{tikzpicture}
%}% end resize box
\end{center}
  \caption{Mechanical heart valve in 2D. Computed instantaneous streamlines (colored by the flow speed) at $t=1.2$ using the composite grids $\Gcts^{(j)}$, $j=1,2,4,8$.}
  \label{fig:heartValve2DRefine}
\end{figure}
}

The time history of the behavior of the top leaflet, presented in Figure~\ref{fig:heartValve2DCurves}, is used
to illustrate the grid convergence of solutions of the~\ampRB~scheme for this problem.  In each plot, four curves are shown corresponding to the time history of the rotation angle, angular velocity and acceleration from calculations using
composite grids of increasing resolution.  We note that the separation between curves decreases rapidly as the grid
resolution increases.
However, the zoomed view in Figure~\ref{fig:heartValve2DCurves} of the acceleration at later times of the simulation, shows that
the convergence of this component of the solution is not as clean.
To give a quantitative sense of the convergence rate,  Richardson extrapolation is used to estimate the time-averaged convergence rate of the three finest grids.
At an early time, $t=0.5$, the estimated convergence rates of the rotation angle, angular velocity and acceleration are 1.51, 1.40 and 1.44, respectively.
The rates drop to 1.36,  1.32 and   1.24 at $t = 1$,
and then to 1.19,  1.15 and 0.89 at the final time $t = 2$.
The reduced convergence rates at later times are largely due to the collision model.
For instance, the convergence of the acceleration is poor during the closing stage, as noted earlier and
shown in the zoomed view of the acceleration.
Generally, we find that the grid convergence is sensitive to the choices of the parameters in the collision model, such as $\epsilon$ and $B_0$.  Realistic simulations require a stiff collision model, so that very fine grids with correspondingly small time steps would be needed to get better convergence results, especially during the closing stage of the simulation.
Figure~\ref{fig:heartValve2DCurves} is again complemented by Figure~\ref{fig:heartValve2DRefine}, which presents the fluid quantities and leaflets positions at $t=1.2$ for different grid resolutions.
It shows the solutions on different grids match well.
%and it is necessary to fix all those parameters during the refinement study to obtain a reasonable convergent result.

We note that because our FSI model simulation has avoided solid-solid contacts, the channel is not complete blocked by the two leaflets in their closed position.  The result of this choice is that there is some small leakage of the fluid between the two leaflets and between the leaflets and the walls of the channel when the valve is in its closed position.
%The issue is a challenging one and our subsequent work will consider approaches to remedy the situation.
We also note that the
Reynolds number of flow for the current simulation is ${\rm Re} = U_{\rm max}D/\nu \approx 100$, where $D=2R_c$ is the channel width and $U_{\rm max}$ is the peak speed of the flow at the inlet.  The Reynolds number for the simulations in~\cite{borazjani2008curvilinear, de2009direct} was reported to be larger, approximately 6000 or more.  We have chosen
to compute the flow in the heart value using a larger viscosity, and thus a smaller Reynolds number, so that sufficient grid resolution of the flow can be achieved, especially for the subsequent benchmark calculations in three dimensions.
Finally, we note that our version of the TP-RB scheme is not able to compute solutions for this problem since we are not able to find a fixed value for the relaxation parameter to stabilize the scheme through the under-relaxed iterations.
As suggested in~\cite{borazjani2008curvilinear}, Aiken's convergence acceleration may be necessary to stabilize the TP-RB scheme for this case.
In contrast, the~\ampRB~scheme remains stable throughout the whole simulations.

\subsubsection{Benchmark problem in three dimensions}
\label{sec:heartValve3D}

%\QT{The grid of 2D and 3D use the same name now, $\Gchv^{(j)}$.}

The geometry of the two leaflets and their position in the cylindrical channel for the three-dimensional heart-value problem are illustrated in Figure~\ref{fig:heartValve3Dgrid}. %the fluid channel covers a cylindrical domain with length $L_c = 1.75$ and radius $R_{c} = 1.1$.
A thin disk of width $L_\ell = 0.1$ and radius $R_\ell = 1$ is split into two halves to obtain two leaflets.
%Two leaflets are allowed to rotate along a hinge axis parallel to $z$-axis.
To allow the leaflets to rotate without colliding into each other, a thin portion of the half-disk along its equator is removed so that the height of the resulting solid is 0.95, which is shorter than its radius $R_\ell$.
The edges of the two leaflets are rounded to avoid sharp corners (similar to that done in the two-dimensional problem).
The above process to define the geometry is done using a computer-aided design (CAD) program.
Figure~\ref{fig:heartValve3Dgrid} also shows the surface grids on the leaflets and the composite grid.
The composite grid, denoted by $\Gchv^{(j)}$ with resolution factor $j$, consists of two background grids and the collection of boundary-fitted grids that represent the leaflets as shown in the figure.
The first background grid is a Cartesian grid, with target  grid spacing $h_{j} = 1/(40j)$, covering the bulk of the fluid channel, while the second grid is a thin boundary-fitted grid, with grid spacing $2h_{j}/3$, attached to the cylindrical boundary of the channel.  (The Cartesian grid is not shown in the figure so that the position of the leaflets can be seen, and only the boundary of the cylindrical boundary-fitted grid is shown.)
The surface of each leaflet is divided into three grids, two surface grids of grid spacing $h_j$ on the two faces of the half-disks and one edge grid of grid spacing $2h_j/3$ forming the perimeter of the solid.
The two surface grids and the edge grid are then extended to the volume, with target grid spacing $2h_j/3$, to form the three-dimensional boundary-fitted grids around the leaflet.
The composite grid is generated using the {\bf Ogen} grid generator.  The IGES file generated by the CAD program provides the input to form the surface and edge grids for the leaflets, and then Ogen generates the body-fitted volume grids and creates the composite grid from these grids along with the two background grids as described above.  For more details on generating grids from a CAD geometry, see~\cite{cadProjection02}.

{% ------ BODY ACCELERATION CURVES ------
\newcommand{\figWidth}{7.cm}
\newcommand{\figWidthb}{6.5cm}
\newcommand{\trimfig}[2]{\trimw{#1}{#2}{.0}{.0}{.0}{.0}}
% zoom: 
\newcommand{\figWidthz}{4.5cm}
\newcommand{\trimfigz}[2]{\trimwb{#1}{#2}{.1}{.1}{.05}{.3}}
%coordinates
\def\thetax{-25}
\def\thetay{187}
\def\thetaz{90}
\def\radx{1}
\def\rady{1}
\def\radz{1}
\begin{figure}[htb]
\begin{center}
\resizebox{14cm}{!}{% START resize box
\begin{tikzpicture}[scale=1]
  \useasboundingbox (0.,1) rectangle (16.,12.2);  % set the bounding box (so we have less surrounding white space)
%  \draw(-1.5,0.0) node[anchor=south west,xshift=-4pt,yshift=+0pt] {\trimfig{heartValve/3d/heartValve3dGrid1}{\figWidth}};
%  \draw(4.5,-0.2) node[anchor=south west,xshift=-4pt,yshift=+0pt] {\trimfig{heartValve/3d/heartValve3dGridLeaflet}{\figWidth}};
%  \draw(11.5,1.0) node[anchor=south west,xshift=-4pt,yshift=+0pt] {\trimfigz{heartValve/3d/heartValve3dGridZoom}{\figWidthz}};
%old:
  \draw(0.2,5.5) node[anchor=south west,xshift=-4pt,yshift=+0pt] {\trimfig{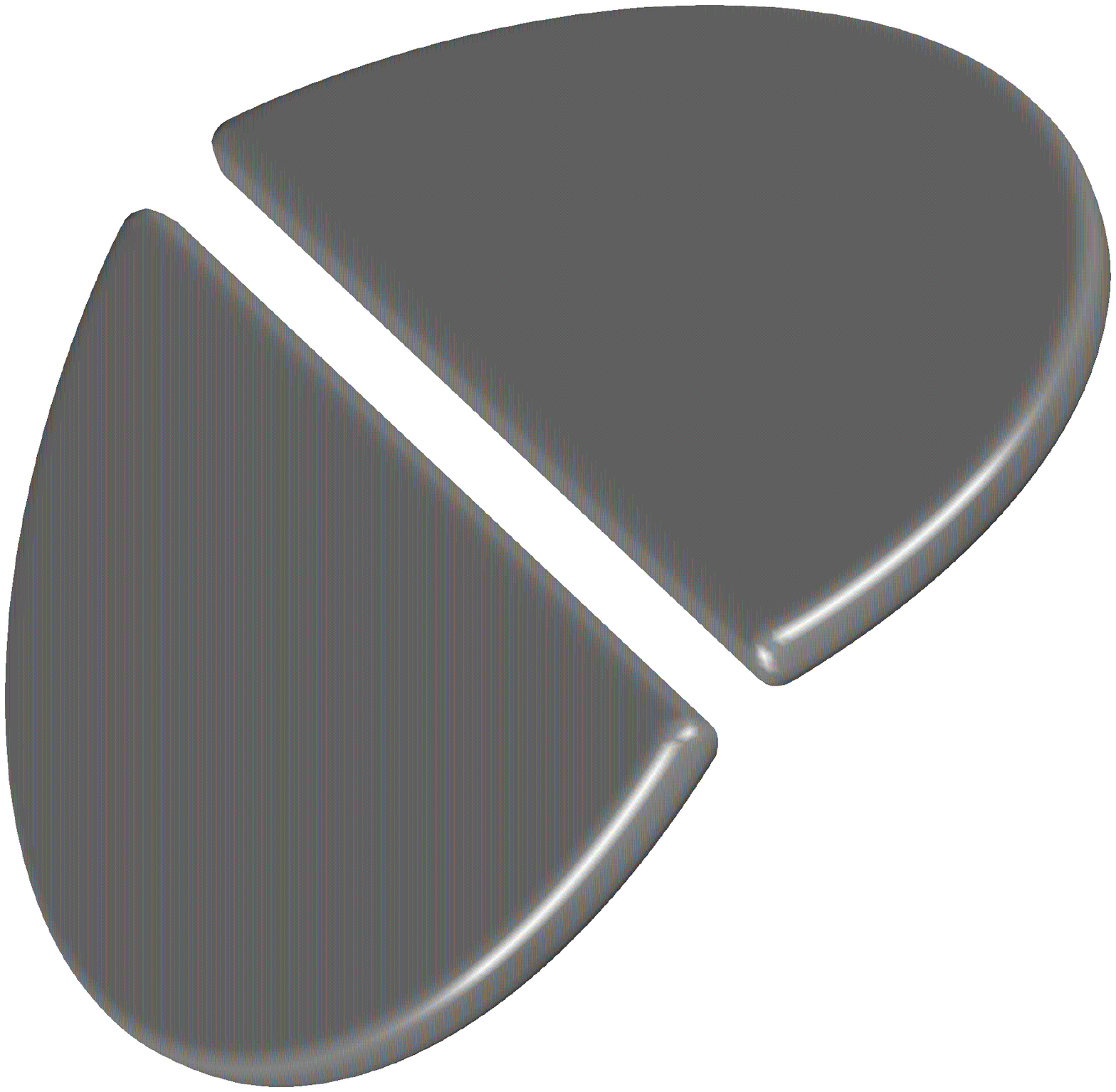}{\figWidth}};
  \draw(8.2,5.5) node[anchor=south west,xshift=-4pt,yshift=+0pt] {\trimfig{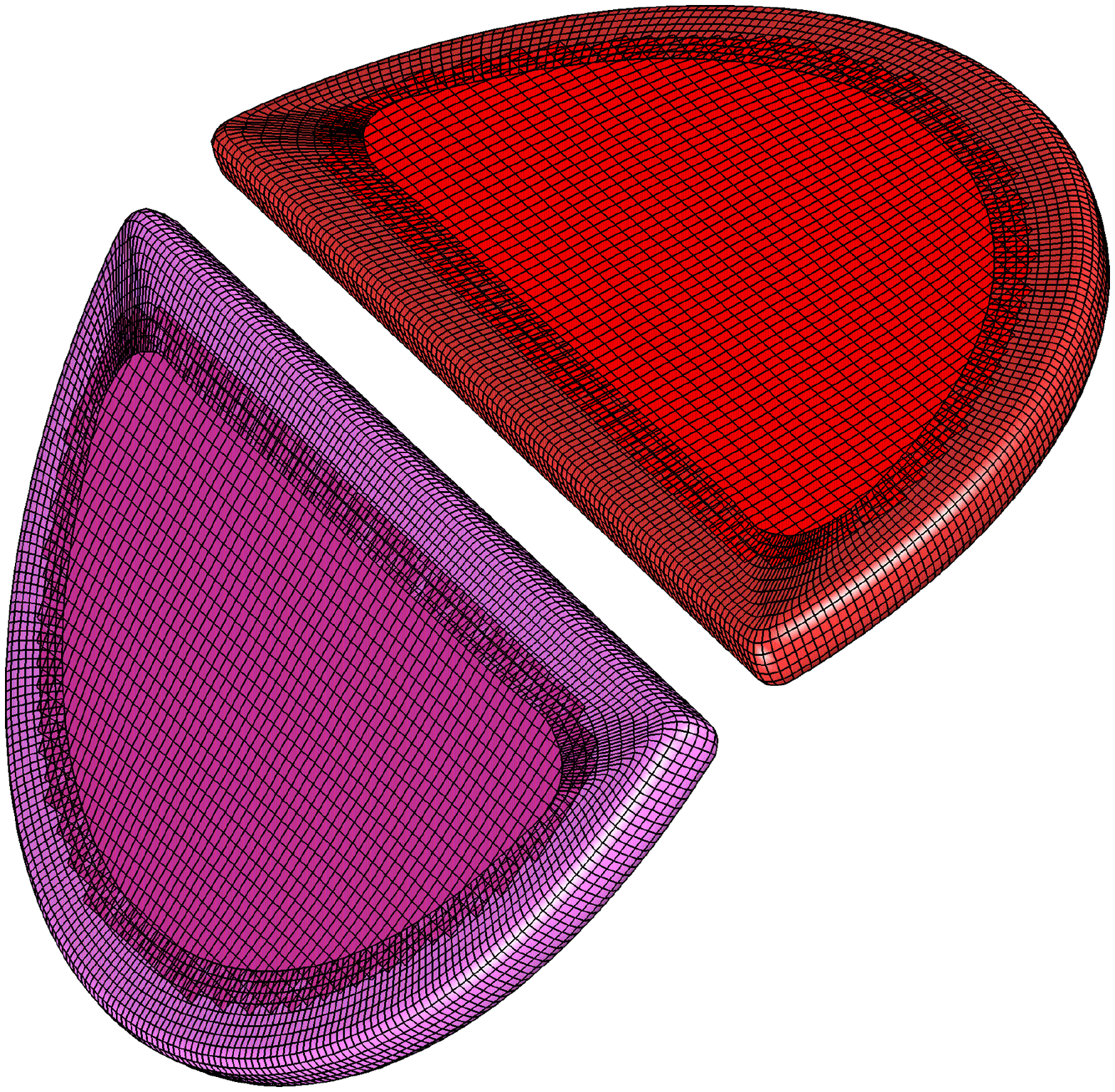}{\figWidth}};
  \draw(0.5,-0.4) node[anchor=south west,xshift=-4pt,yshift=+0pt] {\trimfig{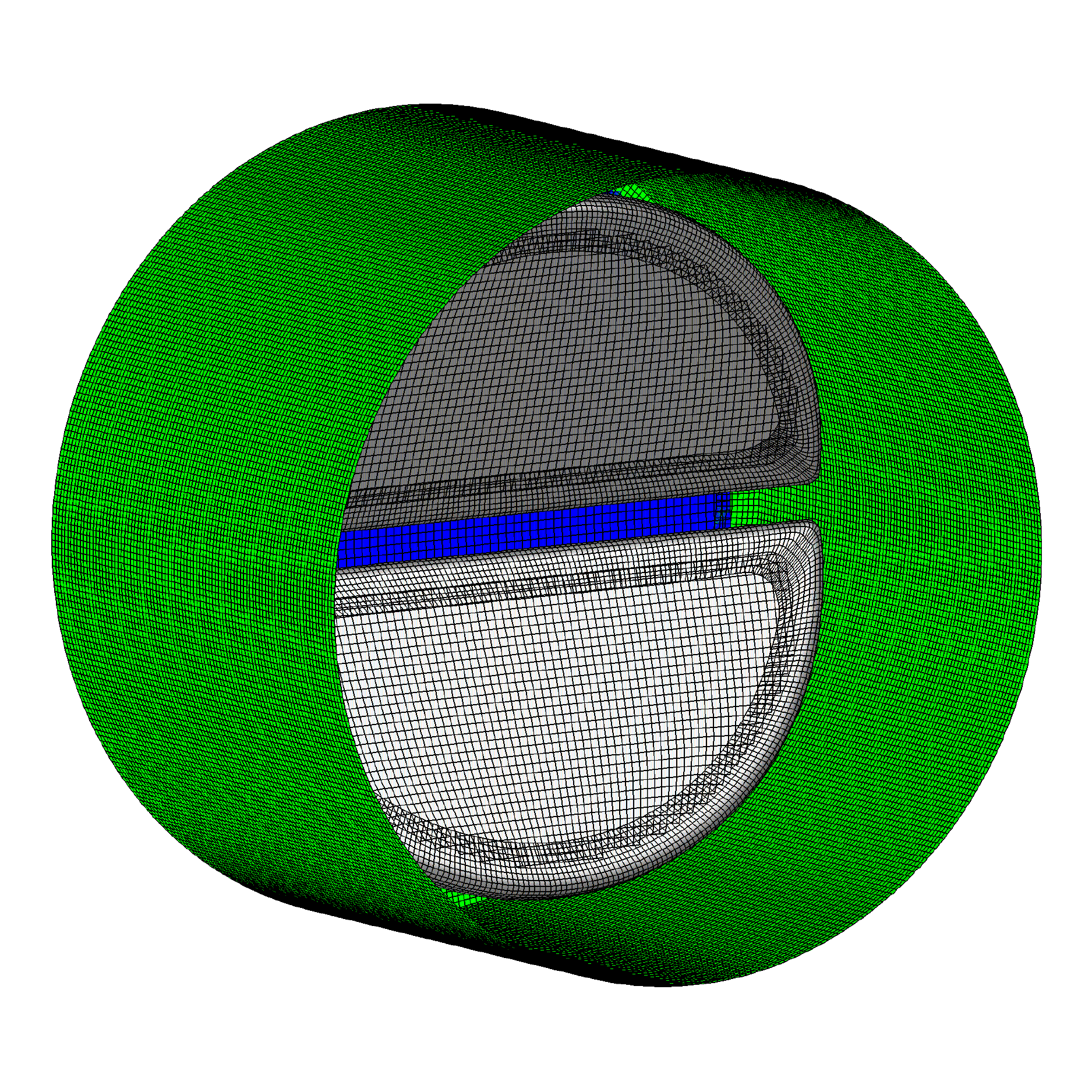}{\figWidthb}};
  \draw(9.3,.5) node[anchor=south west,xshift=-4pt,yshift=+0pt] {\trimfigz{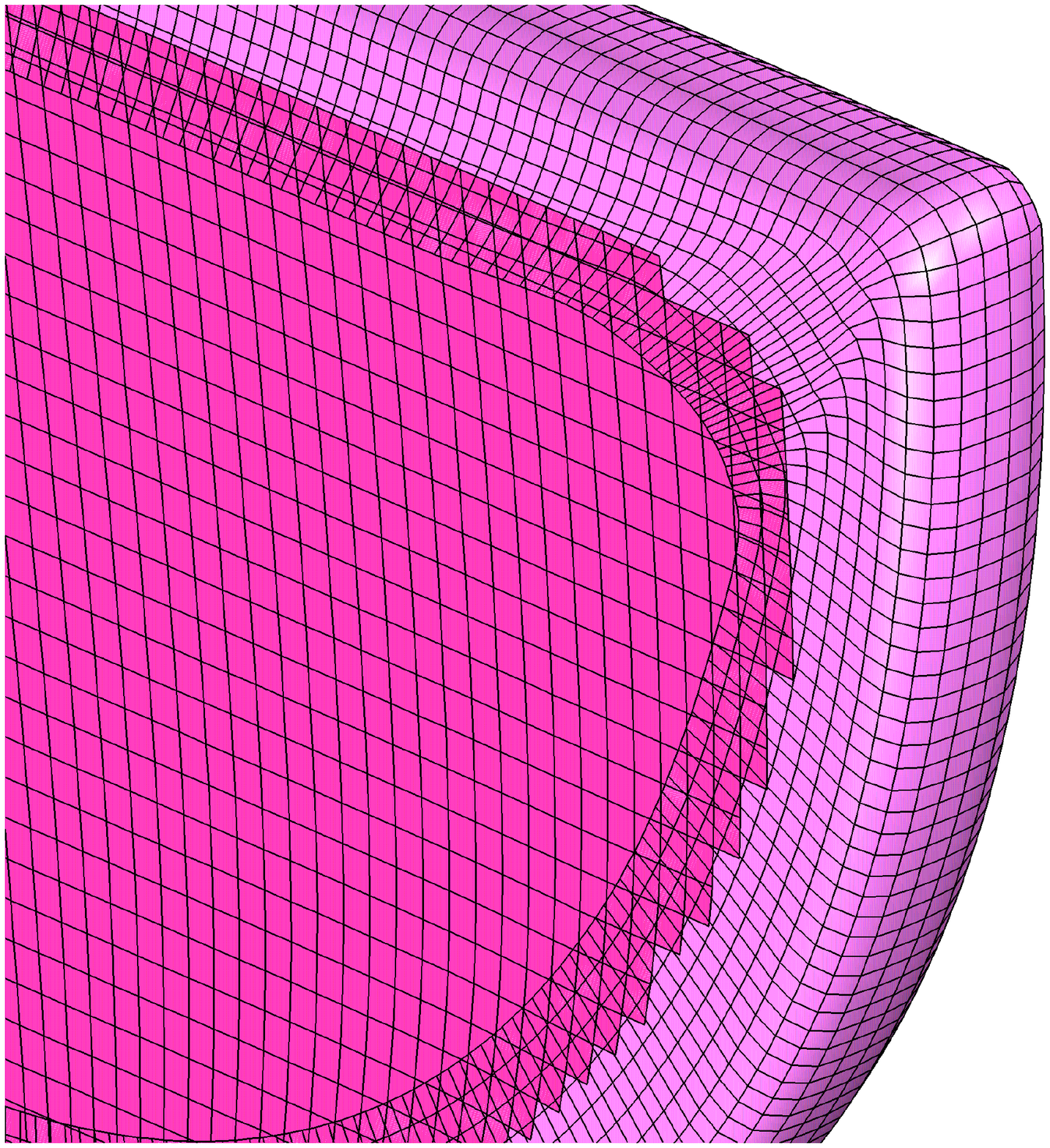}{\figWidthz}};
%plot coordiantes
\begin{scope}[xshift=1.cm,yshift=1.cm]
\draw[thick,->] (0,0) -- ({\radx*cos(\thetax)},{\radx*sin(\thetax)}) node[anchor=north west]{$x$};
\draw[thick,->] (0,0) -- ({\rady*cos(\thetay)},{\rady*sin(\thetay)}) node[anchor=south east]{$z$}; 
\draw[thick,->] (0,0) -- ({\radz*cos(\thetaz)},{\radz*sin(\thetaz)}) node[anchor=south]{$y$};
\end{scope}
%
% grid:
%\draw[step=1cm,gray] (0,0) grid (16,12.2);
\end{tikzpicture}
}% end resize box
\end{center}
\caption{Mechanical heart valve grid in 3D.  
    Top left:  the CAD geometry represented as the surface of two leaflets. 
    Top right: surface grids generated with the hyperbolic grid generator.
    Different surface grids have been colored by different colors.
    Bottom left: the composite grid at $t = 0$. 
    Bottom right: zoomed view of the corner of the bottom leaflet.
    The presented grids are from $\Gchv^{(1)}$ (coarse grid).
     }
  \label{fig:heartValve3Dgrid}
\end{figure}
}

{% ------ STREAMLINES ------
\newcommand{\drawContour}[7]{%
\begin{scope}[#1]
\draw(-.3,0) node[anchor=south west,xshift=-4pt,yshift=+0pt] {\trimfig{heartValve/3d/#2}{\figWidth}};
\draw(.2,5.4) node[draw,fill=white,anchor=west,xshift=2pt,yshift=-1pt] {\scriptsize #5};
 % colour bar:
\begin{scope}[xshift=5pt,yshift=-2pt]
  \draw (\xcb,\ycb) node[anchor=south west,xshift=0.cm,yshift=.5cm,rotate=-90] {\trimfigcb{fig/colourBarLines}{\cbWidth}{\cbHeight}};
  \draw (.8,0) node[anchor=north,xshift=+3pt,yshift=+2pt] {\scriptsize $#6$};
  \draw (3.7,0) node[anchor=north,xshift=+0pt,yshift=+2pt] {\scriptsize $#7$};
\end{scope}
\end{scope}
}
% -- for colour bar ----
\newcommand{\cbWidth}{.2cm}% colour bar width
\newcommand{\cbHeight}{3.3cm}% colour bar height
\newcommand{\xcb}{.5cm}% colour bar lower left corner
\newcommand{\ycb}{-.2cm}% colour bar lower left corner
\setlength{\ycbTop}{\ycb+\cbHeight}% colour bar top label position
\setlength{\ycbMid}{\ycb+\cbHeight*\real{.5}}% colour bar top label position
\newcommand{\trimfigcb}[3]{\includegraphics[width=#2, height=#3, clip, trim=17cm 2.35cm 1.65cm 2.35cm]{#1}}
%
% ============================ DRAW ===================================
\newcommand{\figWidth}{5cm}
\newcommand{\trimfig}[2]{\trimw{#1}{#2}{.12}{.12}{.1}{.1}}
\begin{figure}[htb]
\begin{center}
%\resizebox{16cm}{!}{% START resize box
\begin{tikzpicture}[scale=1]
  \useasboundingbox (0.0,.25) rectangle (15.,12.1);  % set the bounding box (so we have less surrounding white space)
% 
%  top row: 
\begin{scope}[yshift=6.5cm]
 \drawContour{xshift= 0cm,yshift=0.0cm}{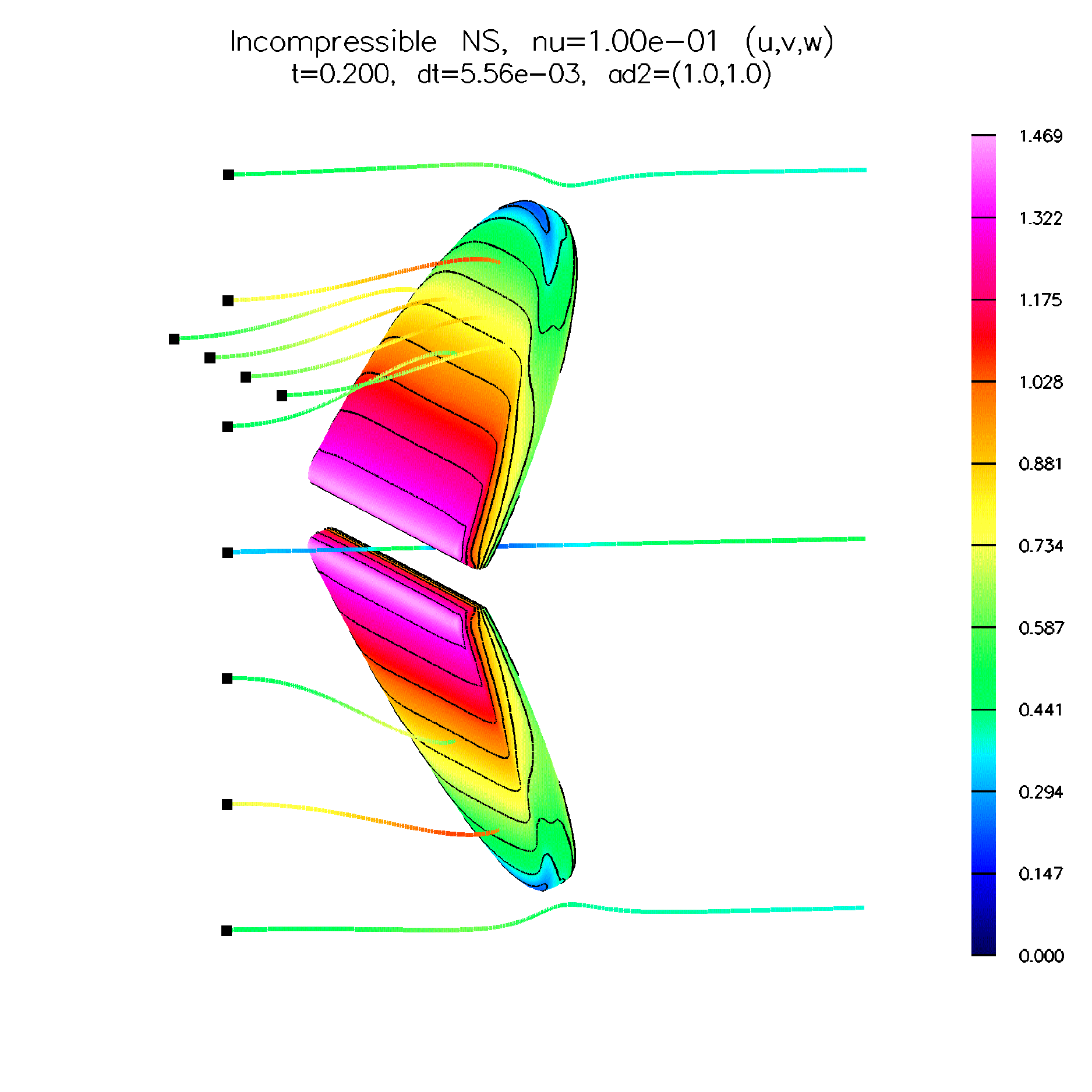}{$p$}{$p$}{$t=0.2$}{0.0}{1.47};
 \drawContour{xshift= 5cm,yshift=0.0cm}{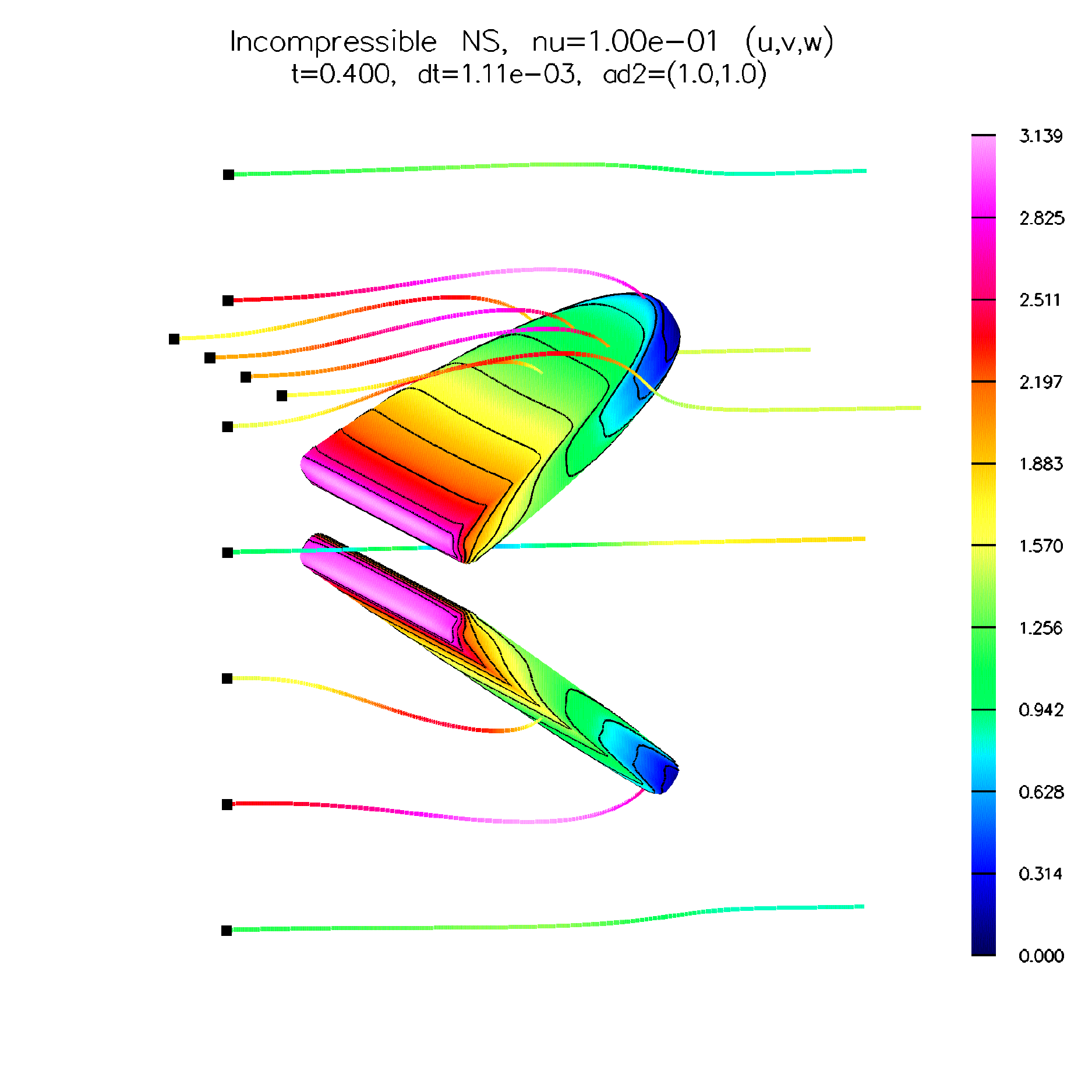}{$p$}{$p$}{$t=0.4$}{0.0}{3.14};
 \drawContour{xshift= 10cm,yshift=0.0cm}{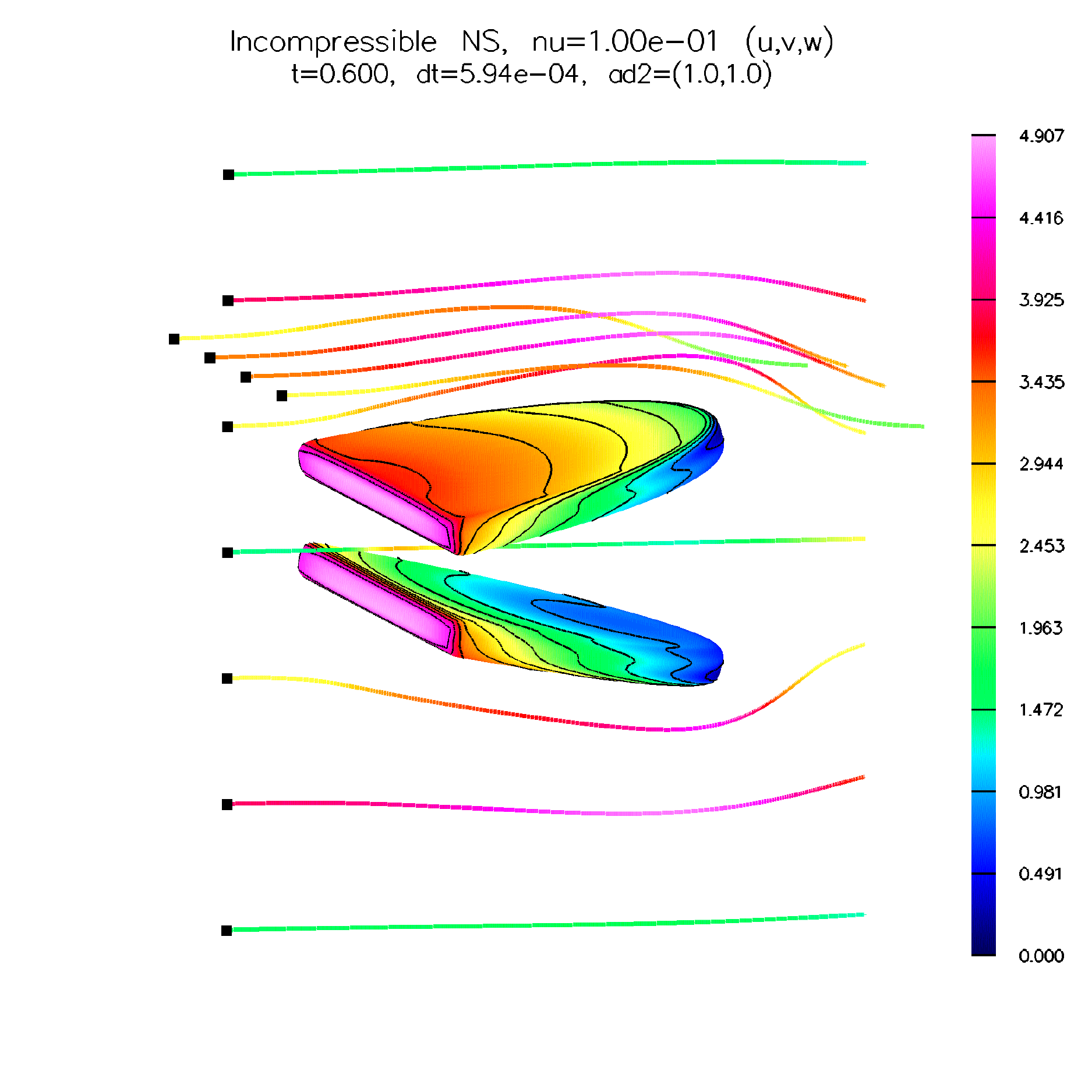}{$p$}{$p$}{$t=0.6$}{0.0}{4.91};
\end{scope}
%  bottom row: 
 \drawContour{xshift= 0cm,yshift=0.0cm}{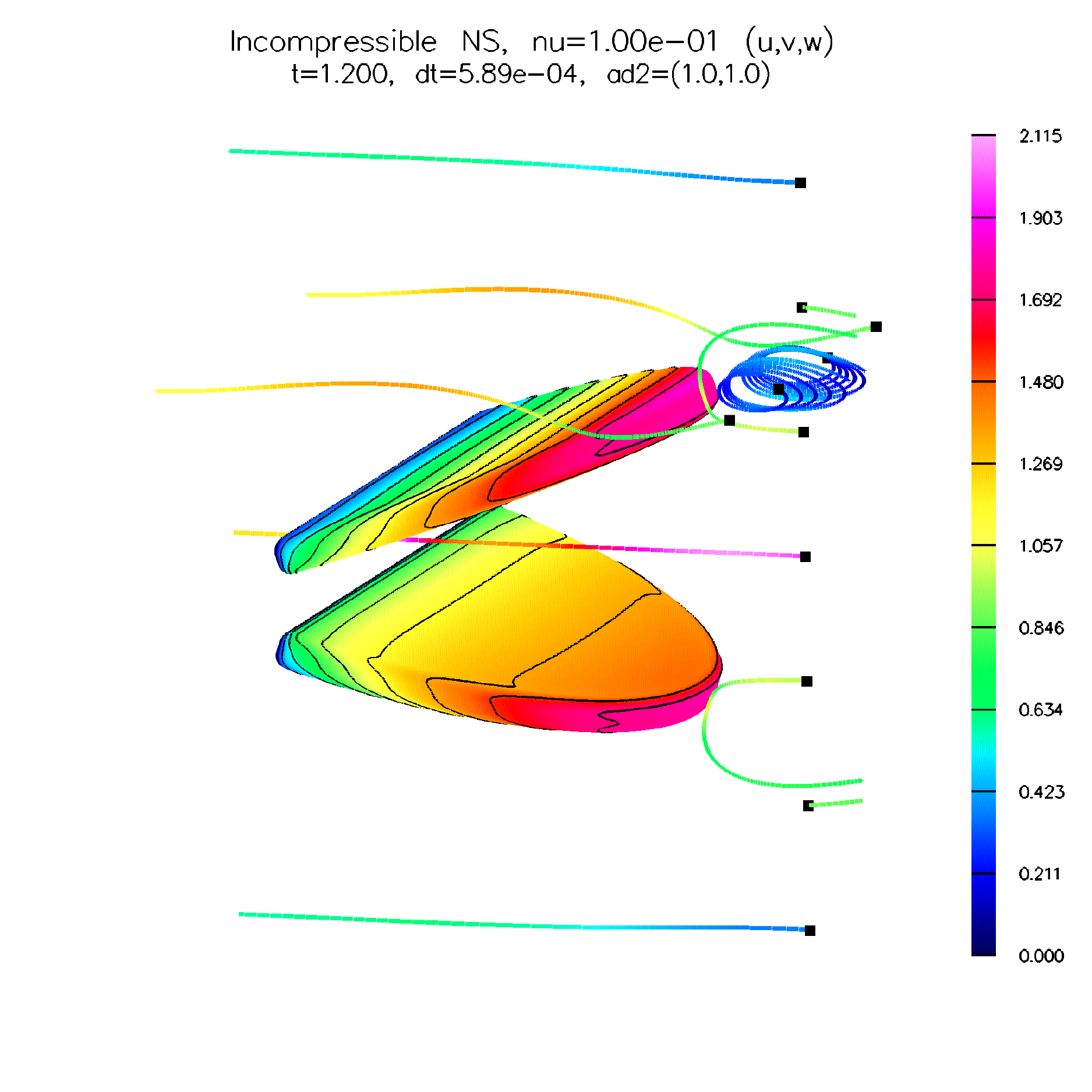}{$p$}{$p$}{$t=1.2$}{0.0}{2.12};
 \drawContour{xshift= 5cm,yshift=0.0cm}{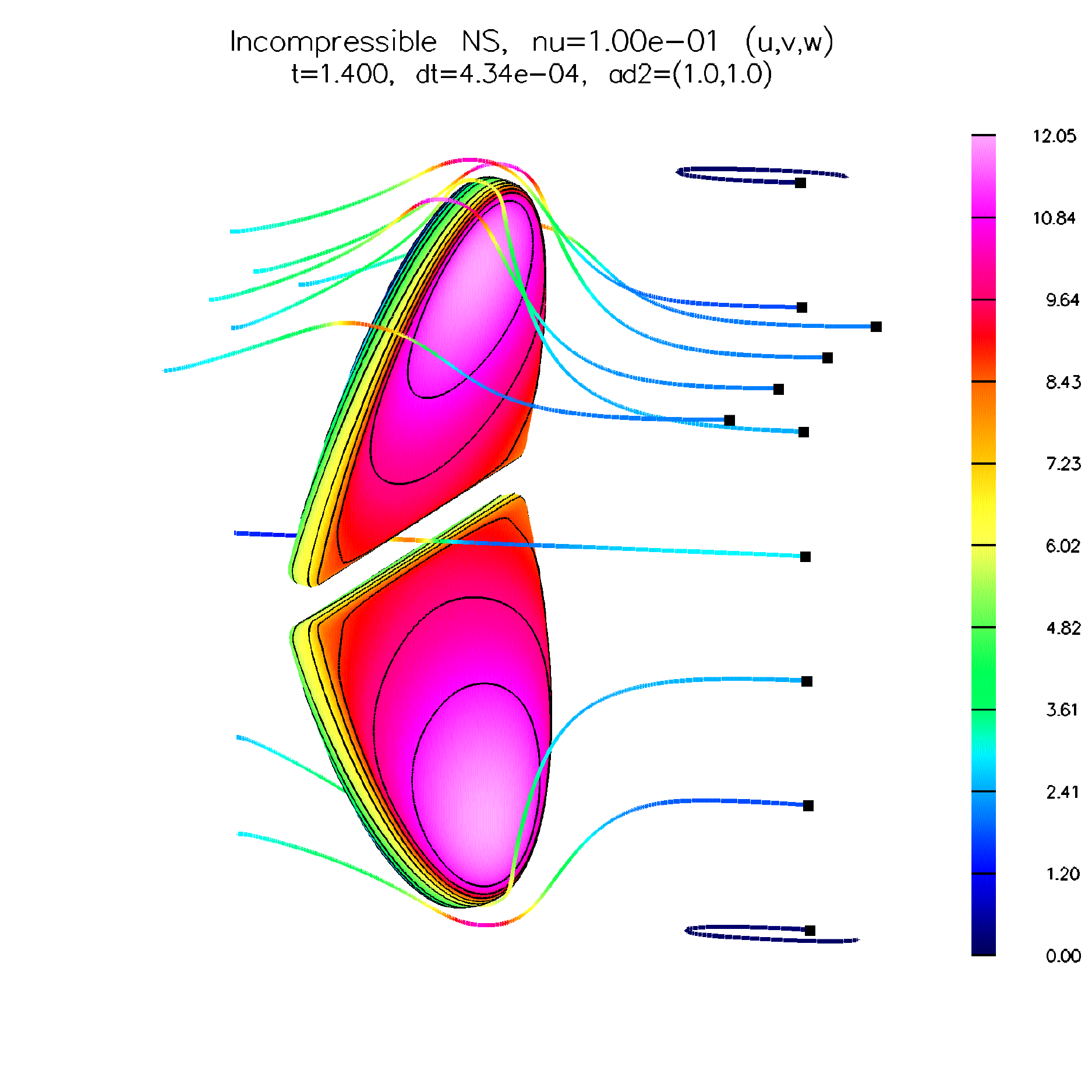}{$p$}{$p$}{$t=1.4$}{0.0}{12.05};
 \drawContour{xshift= 10cm,yshift=0.0cm}{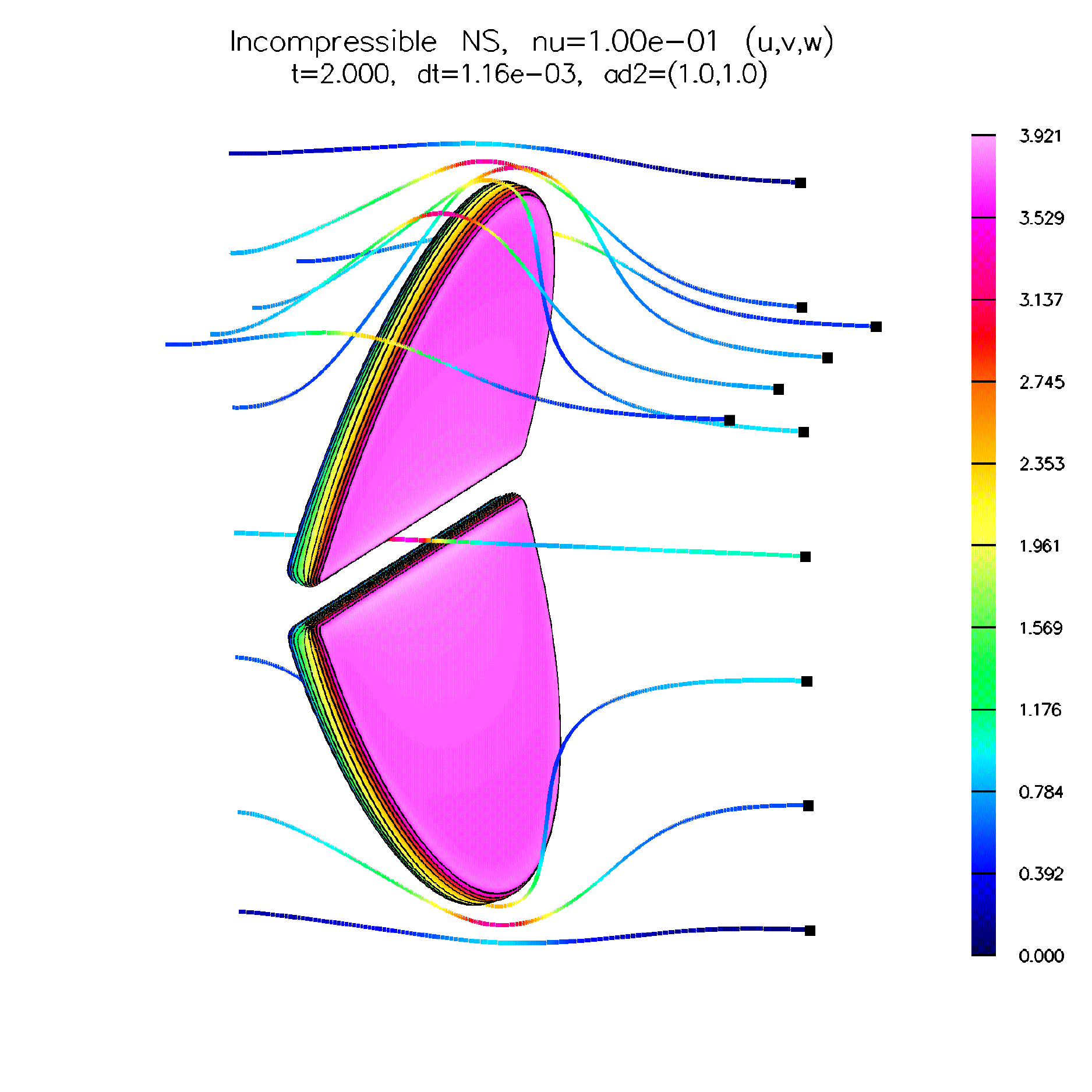}{$p$}{$p$}{$t=2.0$}{0.0}{3.92};
% 
% ---------------
% grid:
%\draw[step=1cm,gray] (0,0) grid (15,12.1);
\end{tikzpicture}
%}% end resize box
\end{center}
  \caption{Mechanical heart valve in 3D. Computed instantaneous streamlines (colored by the flow speed) at selected times using the composite grid $\Gchv^{(2)}$.  
  The colorbars are provided for the flow speed. The pressure is plotted on surface of the leaflets. 
}
  \label{fig:heartValve3Dstreamlines}
\end{figure}
}

%what is going on in the solution, what happens at different times
The behaviour of the numerical solution at selected times during the simulation of the three-dimensional heart-value problem is shown in Figure~\ref{fig:heartValve3Dstreamlines}.  Colour contours of the fluid pressure on the surface of the leaflets is shown along with the instantaneous streamlines from a collection of tracer particles which illustrate the flow around the leaflets.
The images along the top row of the figure at times $t=0.2$, $0.4$ and $0.6$ are chosen to illustrate the behaviour of the solution as the valve opens.
The starting points for the tracer particles in these images lie on the plane $x=x_0=-0.5$ at the left boundary of the fluid channel, and are positioned in the shape of a cross.  There are seven tracers along the vertical line segment of the cross with $x=-0.5$ and $z=0$, and there are four tracers along its horizontal line segment with $x=-0.5$ and $y=0.5$.  At the earliest time, $t=0.2$, we observe that the instantaneous streamlines from some of the tracers stagnate on the surface of the leaflets indicating that the flow is blocked, while at later times the streamlines from these tracers show that the flow passes over and around the leaflets as the valve opens.
%The figure features the opening and closing processes of the leaflets.
%During the opening process (the top row in Figure~\ref{fig:heartValve3Dstreamlines}), 11 starting points  are selected in the region close to the left boundary.
%Seven points are located along the line $(-0.5,y,0)$ and the other four points are located along the line $(-0.5,0.5,z)$.
The bottom row of images in Figure~\ref{fig:heartValve3Dstreamlines} at the times $t=1.2$, $1.4$ and $2.0$ are shown to illustrate the flow while the valve closes.  For these images, the starting points of the tracers are positions on the
plane $x=1.1$ near the right boundary of the channel.  Seven tracers are located along the verticle line with $x=1.1$ and $z=0$, while four more tracers are located along the horizontal line with $x=1.1$ and $y=0.5$, all forming a cross shape similar to before.  At $t=1.2$, we note a recirculation zone near the tip of the top leaflet as the flow reverses due to a change in the sign of the prescribed pressure profile in~\eqref{eq:pressureProfile}.  At the two later times, the flow again becomes blocked as the valve closes, but with some leakage as noted earlier.  However, there is a significant region of near stagnant flow on the face of leaflet as indicated by the near uniform pressure there.

The time history of the motion of the top leaflet is presented in Figure~\ref{fig:heartValve3DCurves}.
Overall, the motion of leaflet is similar to that observed in the two-dimensional case, and shown previously
in Figure~\ref{fig:heartValve2DCurves}.  For example, there is a slight oscillation seen as the leaflet approaches
its lower limit, $\theta_{\rm min}$, at the open position.  There is a brief interval of time where the leaflet is
in equilibrium at its open position, prior to a rapid rotation of the leaflet towards the upper limit, $\theta_{\rm max}$, at the closed positon.  Near the upper limit there is a stronger oscillation, also observed in the two-dimensional case, as the stress on the leaflet from the fluid achieves a balance with the stiff repulsive force of the contact model.  The oscillation dampens in time due to the damping term in the model, and an equilibrium is again reached but now with the valve in the closed position.
In Figure~\ref{fig:heartValve3DCurves} for the three-dimensional case, there are three curves drawn in each plot corresponding to the results computed using the composite grid, $\Gchv^{(j)}$, with increasing values of the resolution factor $j=1$, 2 and~4.  This is done to indicate the grid convergence of the solutions.  In the plots of the rotation angle $\theta_b$, the angular velocity $\omega_b$ and the acceleration ${\dot\omega}_b$, we see that the curves from the solutions on the different grids lie nearly on top of one another, as was seen for the two-dimensional benchmark problem.
The differences in the curves can be seen in the zoomed views, and in these views we observe a convergence in the behaviours as the grid is refined.  Estimates of the convergence rate for the angular velocity and acceleration are $0.92$ and $1.46$, respectively, while the convergence rate for the rotation angle is $0.30$.  The convergence rates found in the three-dimensional case are generally lower than the ones found in the two-dimensional case for this difficult problem.  This is likely due to the coarser grids used to compute the convergence estimates for the three-dimensional case.
Figure~\ref{fig:heartValve3DCurves} is again complemented by Figure~\ref{fig:heartValve3DRefine}, which shows good agreement of the solutions at $t=1.2$ for different grid resolutions.

{% ------ MATLAB CURVES LIGHT BODY ------
\newcommand{\figWidth}{8.cm}
\newcommand{\trimfig}[2]{\trimFig{#1}{#2}{.0}{.0}{.0}{.0}}
\newcommand{\figWidthz}{3.5cm}
\newcommand{\figWidthb}{3.5cm}
\newcommand{\trimfigz}[2]{\trimFig{#1}{#2}{.0}{.0}{.05}{.0}}
\begin{figure}[htb]
\begin{center}
\resizebox{14cm}{!}{% START resize box
\begin{tikzpicture}[scale=1]
  \useasboundingbox (0.0,0) rectangle (16,13.8);  % set the bounding box (so we have less surrounding white space)
  \draw(-.5,6.7) node[anchor=south west,xshift=0pt,yshift=+0pt] {\trimfig{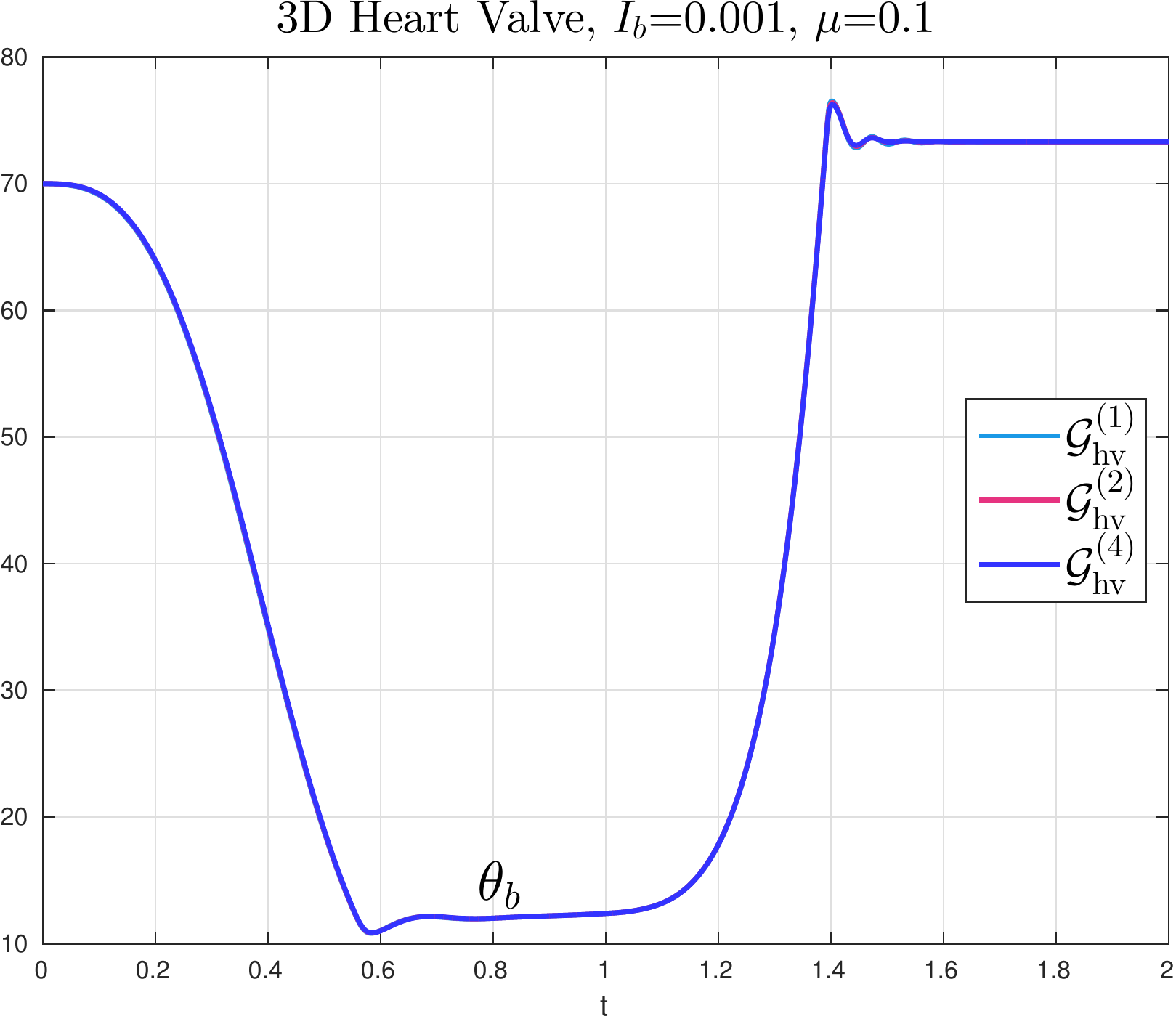}{\figWidth}};
  %zoom
  \draw(1.4, 10.4) node[anchor=south west,xshift=0pt,yshift=+0pt] {\trimfigz{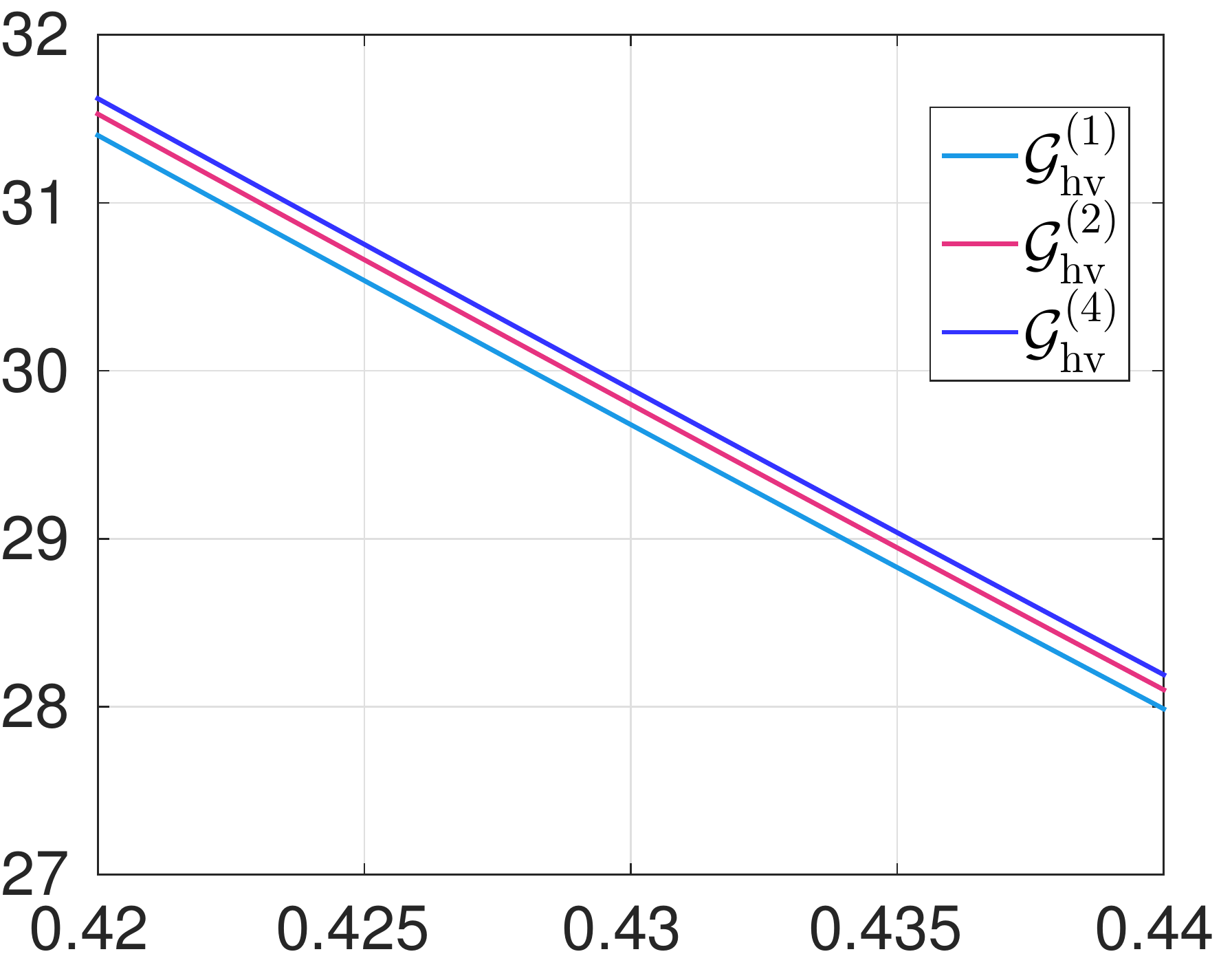}{\figWidthb}};
  \draw(8.0,6.7) node[anchor=south west,xshift=0pt,yshift=+0pt] {\trimfig{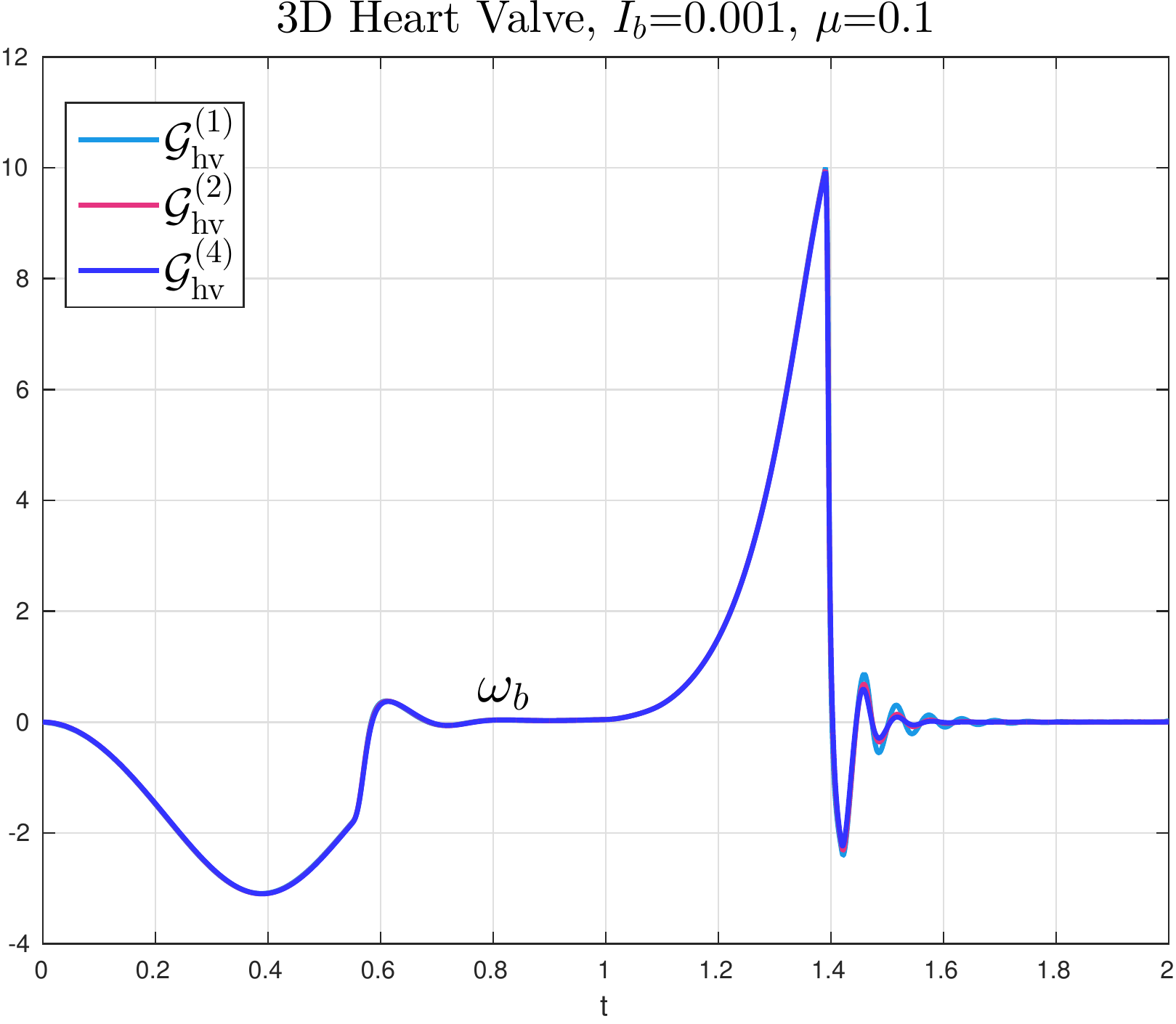}{\figWidth}};
  %zoom
  \draw(9.7, 10.3) node[anchor=south west,xshift=0pt,yshift=+0pt] {\trimfigz{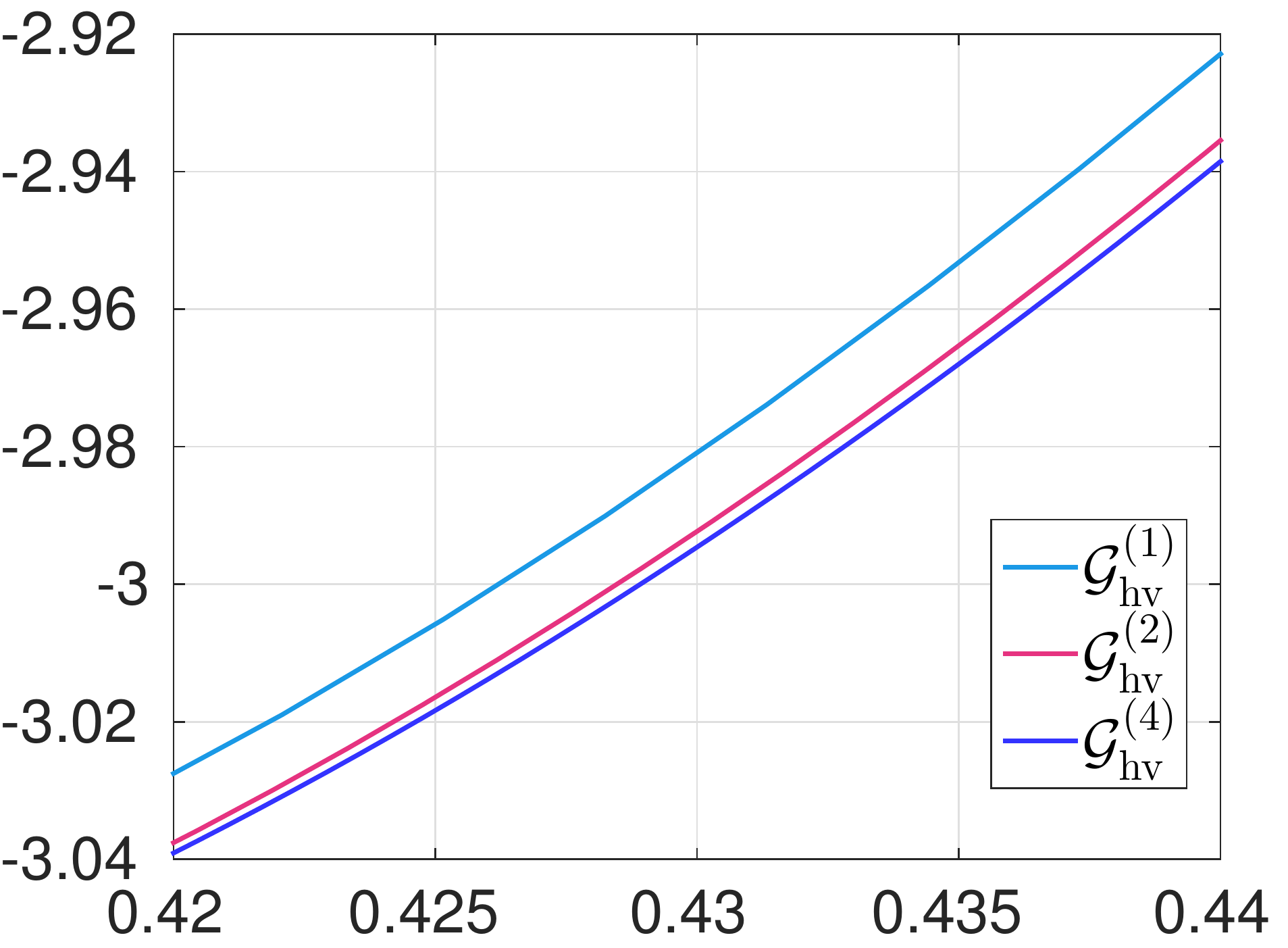}{\figWidthb}};
  %---
  \draw(-.5,-.75) node[anchor=south west,xshift=0pt,yshift=+0pt] {\trimfig{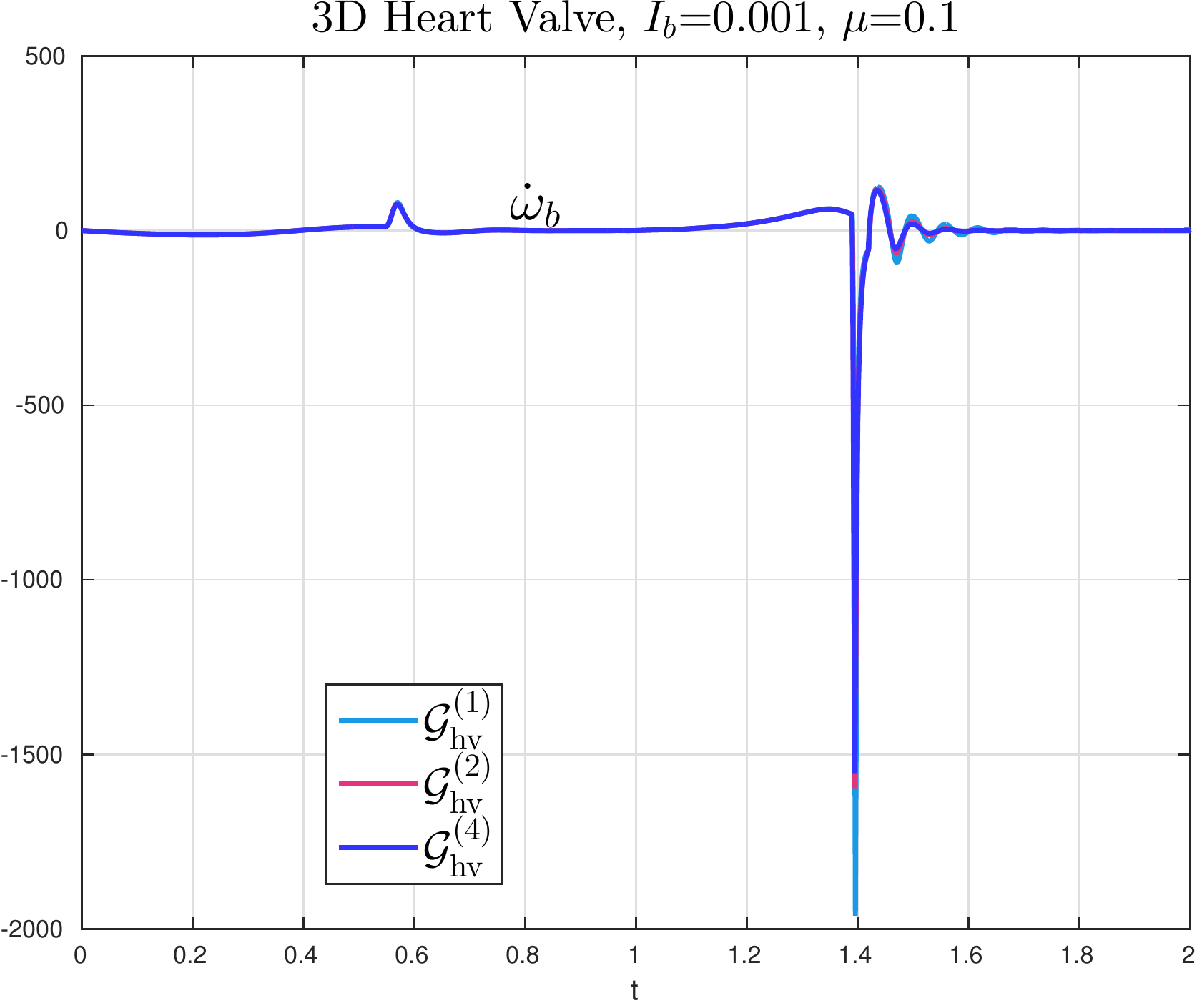}{\figWidth}};
  %zoom
  \draw(.5, 1.7) node[anchor=south west,xshift=0pt,yshift=+0pt] {\trimfigz{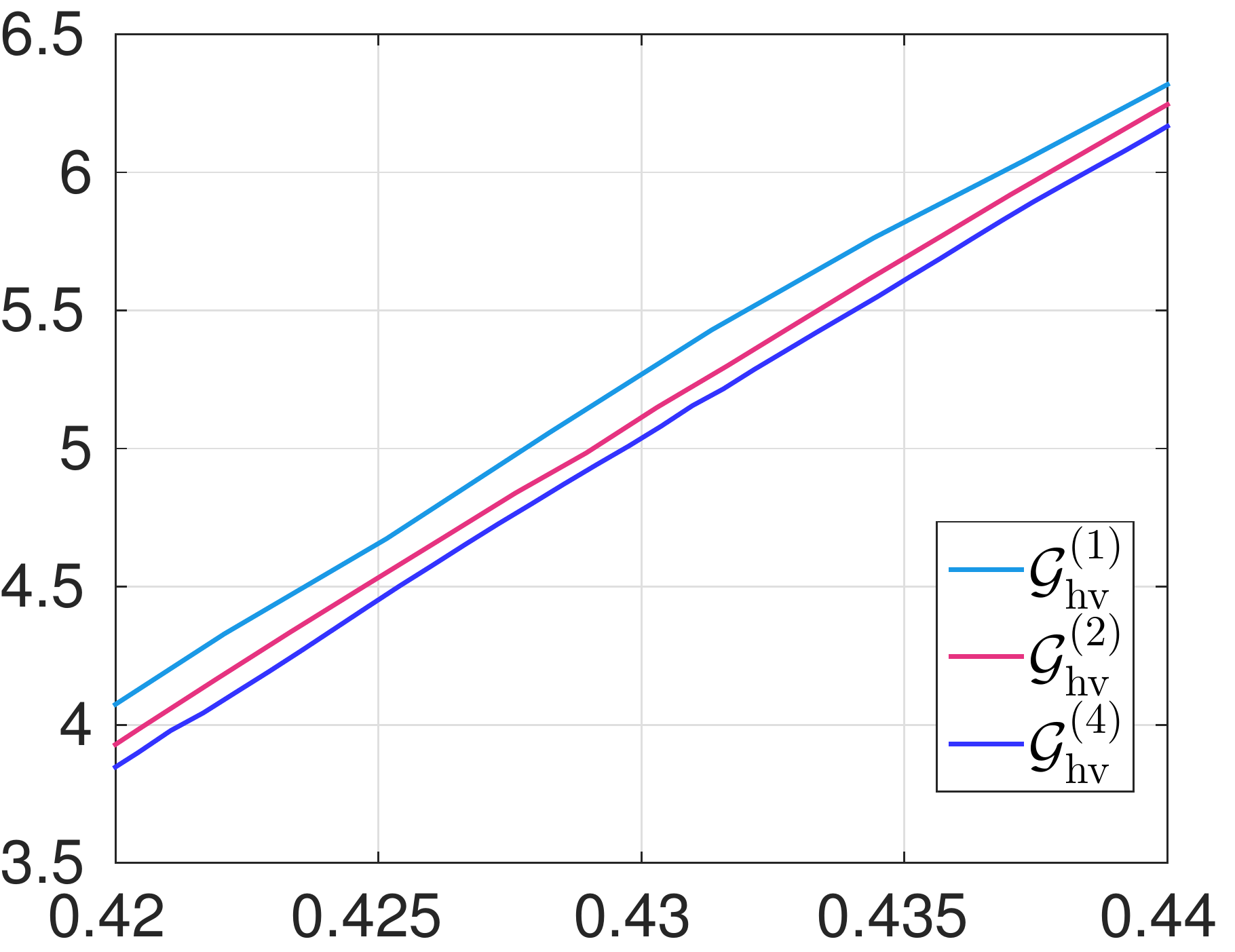}{\figWidthz}};
  \draw(8,-.75) node[anchor=south west,xshift=0pt,yshift=+0pt] {\trimfig{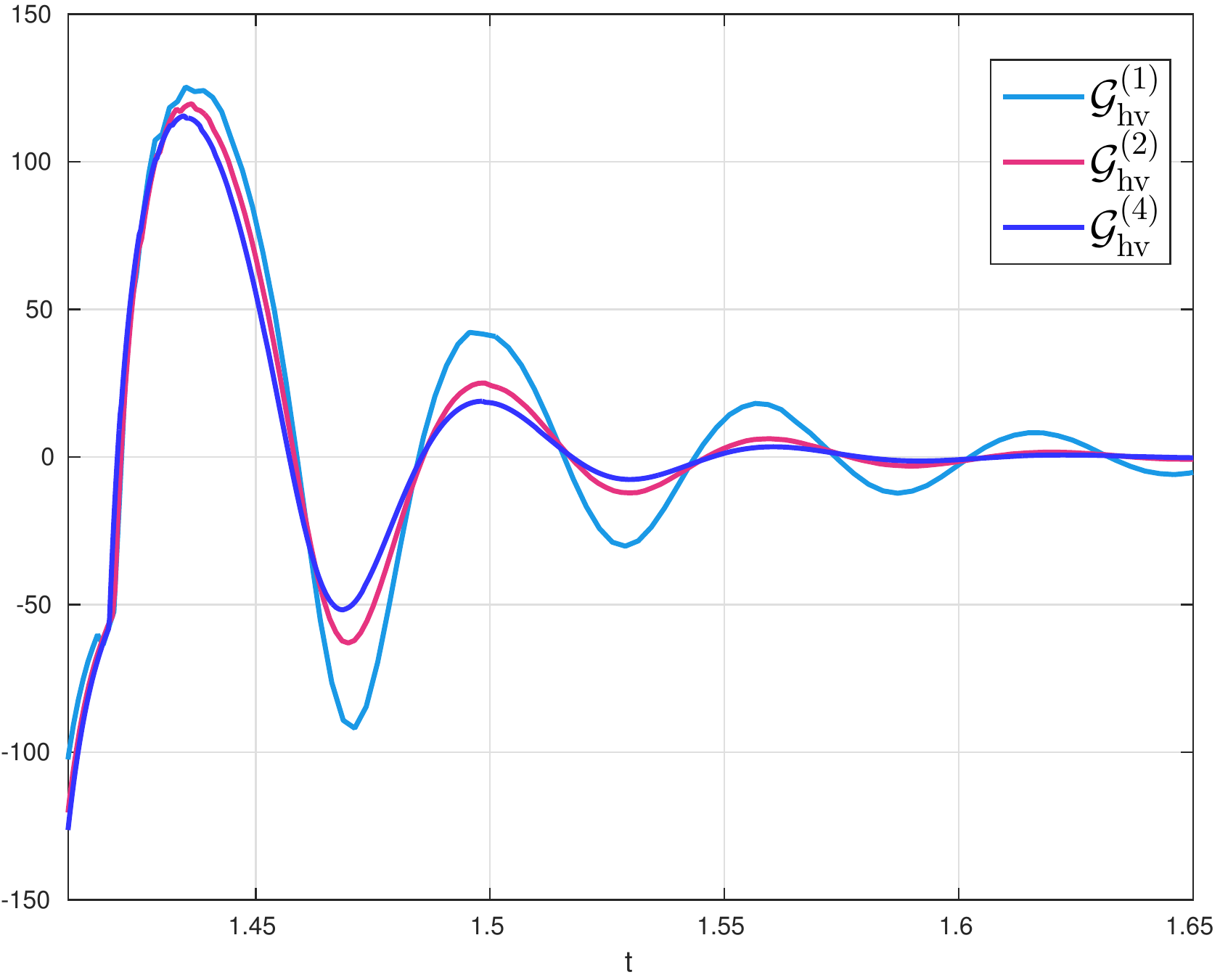}{\figWidth}};
%
% grid:
%\draw[step=1cm,gray] (0,0) grid (16,13.8);
\end{tikzpicture}
} % END resize box
\end{center}
  \caption{Mechanical heart valve in 3D. Time history of the rotation angle (top left),  angular velocity (top right) and acceleration (bottom left) of the top leaflet.
      Several zoomed views are presented including the zoomed view of the accelerations at the later time (bottom right).
      Results are shown from calculations using the composite grid, $\Gchv^{(j)}$, $j = 1, 2, 4$.
  }
 
  \label{fig:heartValve3DCurves}
\end{figure}
}

{% ------ STREAMLINES ------
\newcommand{\drawContour}[7]{%
\begin{scope}[#1]
\draw(-.3,0) node[anchor=south west,xshift=-4pt,yshift=+0pt] {\trimfig{heartValve/3d/#2}{\figWidth}};
\draw(.2,5.4) node[draw,fill=white,anchor=west,xshift=2pt,yshift=-1pt] {\scriptsize #5};
\draw(3.2,5.3) node[draw,fill=white,anchor=west,xshift=2pt,yshift=-1pt] {\scriptsize #4};
 % colour bar:
\begin{scope}[xshift=5pt,yshift=-2pt]
  \draw (\xcb,\ycb) node[anchor=south west,xshift=0.cm,yshift=.5cm,rotate=-90] {\trimfigcb{fig/colourBarLines}{\cbWidth}{\cbHeight}};
  \draw (.8,0) node[anchor=north,xshift=+3pt,yshift=+2pt] {\scriptsize $#6$};
  \draw (3.7,0) node[anchor=north,xshift=+0pt,yshift=+2pt] {\scriptsize $#7$};
\end{scope}
\end{scope}
}
% -- for colour bar ----
\newcommand{\cbWidth}{.2cm}% colour bar width
\newcommand{\cbHeight}{3.3cm}% colour bar height
\newcommand{\xcb}{.5cm}% colour bar lower left corner
\newcommand{\ycb}{-.2cm}% colour bar lower left corner
\setlength{\ycbTop}{\ycb+\cbHeight}% colour bar top label position
\setlength{\ycbMid}{\ycb+\cbHeight*\real{.5}}% colour bar top label position
\newcommand{\trimfigcb}[3]{\includegraphics[width=#2, height=#3, clip, trim=17cm 2.35cm 1.65cm 2.35cm]{#1}}
%
% ============================ DRAW ===================================
\newcommand{\figWidth}{5cm}
\newcommand{\trimfig}[2]{\trimw{#1}{#2}{.12}{.12}{.1}{.1}}
\begin{figure}[htb]
\begin{center}
%\resizebox{16cm}{!}{% START resize box
\begin{tikzpicture}[scale=1]
  \useasboundingbox (0.0,.25) rectangle (15.,5.6);  % set the bounding box (so we have less surrounding white space)
% 
%  bottom row: 
 \drawContour{xshift= 0cm,yshift=0.0cm}{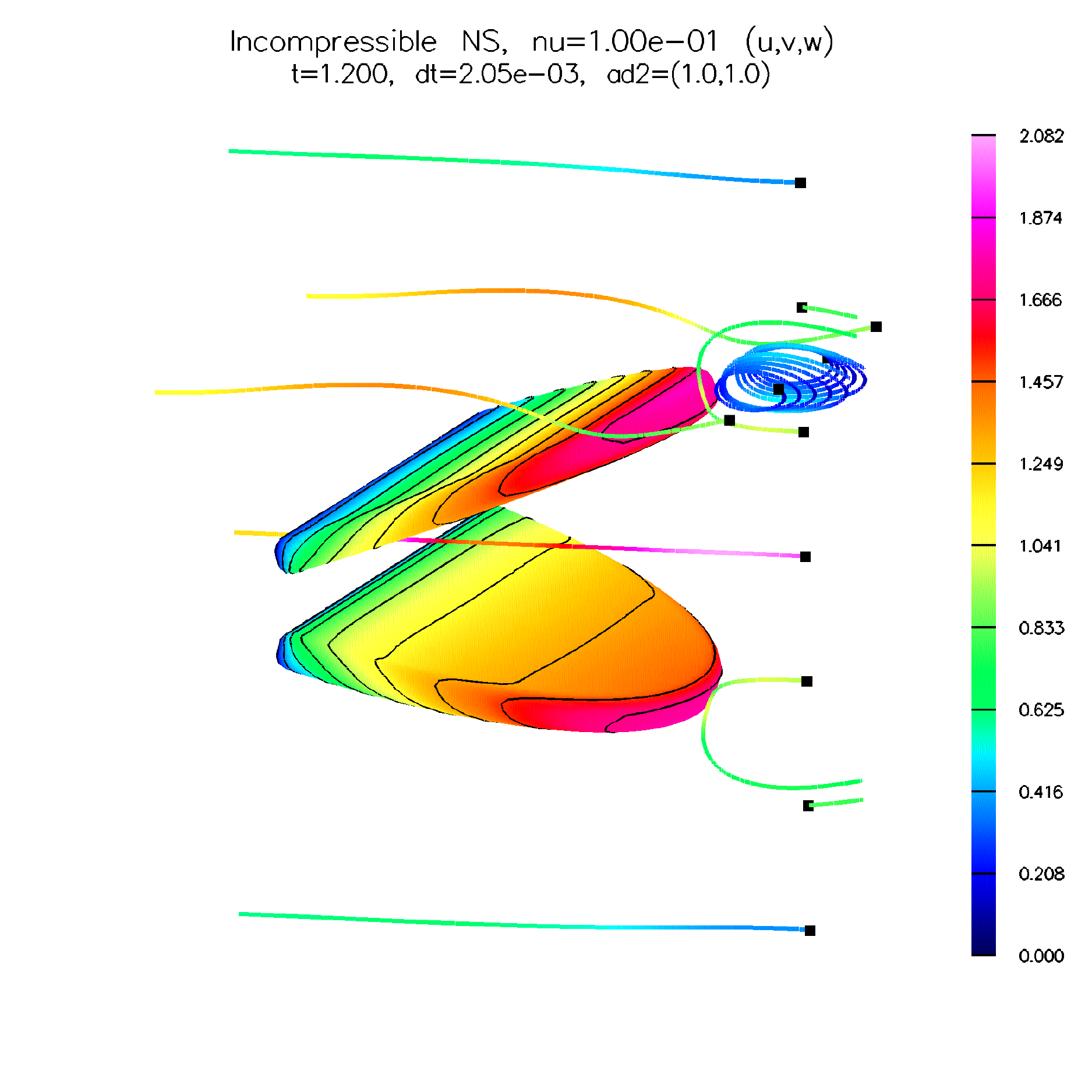}{$p$}{$\Gchv^{(1)}$}{$t=1.2$}{0.0}{2.08};
 \drawContour{xshift= 5cm,yshift=0.0cm}{heartValve3dt1p2streamline}{$p$}{$\Gchv^{(2)}$}{$t=1.2$}{0.0}{2.12};
 \drawContour{xshift= 10cm,yshift=0.0cm}{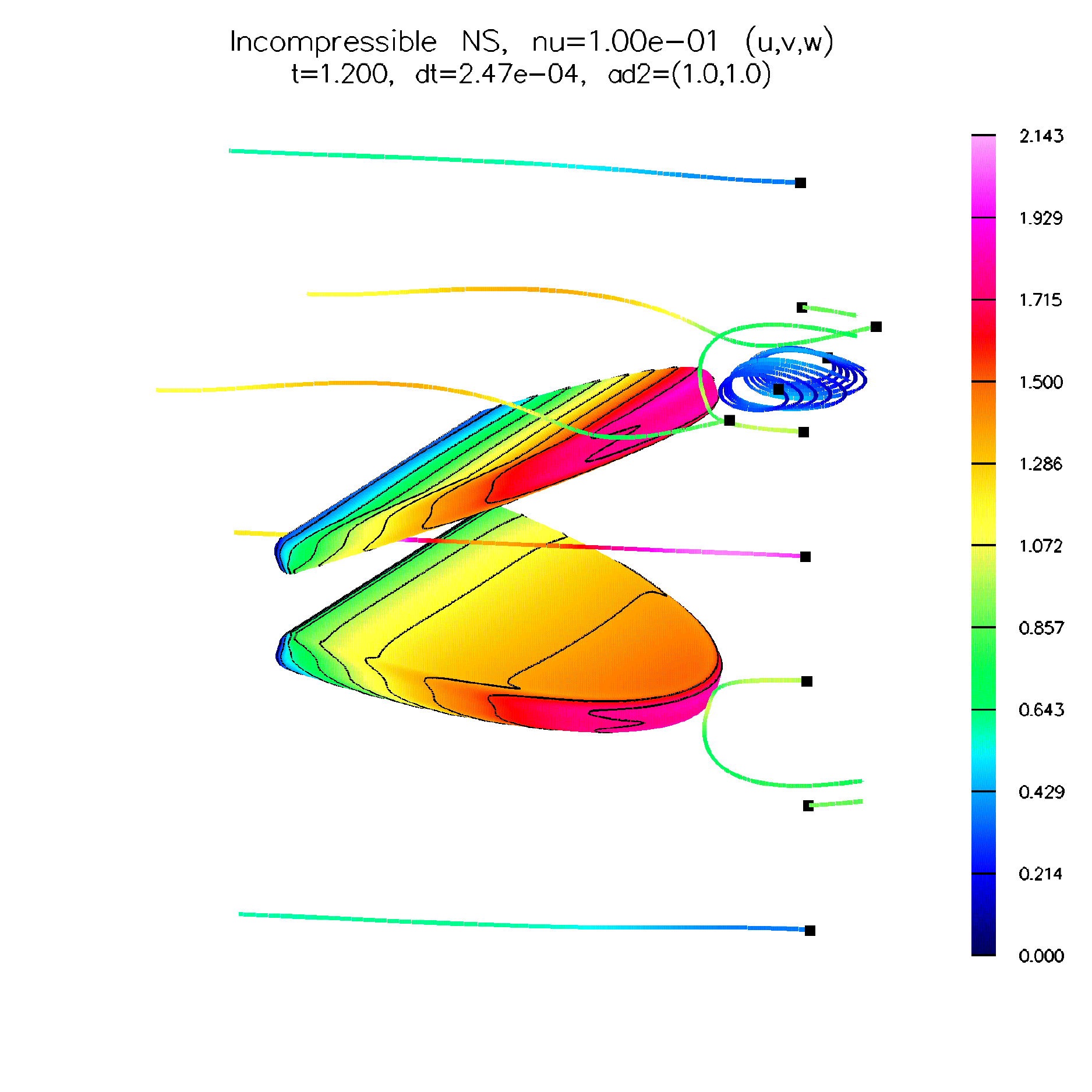}{$p$}{$\Gchv^{(4)}$}{$t=1.2$}{0.0}{2.14};
% 
% ---------------
% grid:
%\draw[step=1cm,gray] (0,0) grid (15,5.6);
\end{tikzpicture}
%}% end resize box
\end{center}
  \caption{Mechanical heart valve in 3D. Computed instantaneous streamlines (colored by the flow speed) at $t=1.2$ using the composite grids $\Gchv^{(j)}$, $j=1,2,3$.  
  The colorbars are provided for the flow speed. The pressure is plotted on surface of the leaflets.}
  \label{fig:heartValve3DRefine}
\end{figure}
}

As in the two-dimensional benchmark problem, we are not able to find a stable under-relaxed parameter for our version of the TP-RB scheme when simulating the problem in three dimensions. 
The TP-RB scheme typically becomes unstable after a couple of time steps even with under-relaxed sub-iterations.
In contrast, the~\ampRB~scheme remains stable throughout the whole simulation.

\section{Conclusions}
\label{sec:conclusions}

In this work we have extended the~\ampRB~scheme derived in~\cite{rbinsmp2017, rbins2017} to three dimensions 
for FSI problems with rigid bodies moving in an incompressible fluid.
The extension relies on the moving overlapping grid approach, 
parallel implementation based on MPI and various building blocks handled properly in three dimensions.
Several details of the full three-dimensional algorithm have been described,
including the AMP interface condition handling the added-mass and added-damping effects,
a strategy to find the surface quadrature on composite grids,
and parallel implementation such as parallel partitioning and sparse linear solvers.
In addition, an analysis of two model problems for a spherical rigid-body showed the stability of the
three-dimensional scheme with respect to added-mass and added-damping effects.
The analysis and careful numerical results confirmed that the~\ampRB~scheme remains stable without sub-iterations 
for light and even zero-mass bodies of general geometries in three dimensions.

%summary of numerical examples
The resultant algorithm has been verified and validated through several benchmark problems.
It has been shown that the algorithm achieves a second-order accuracy and remains stable for some challenging problems.
In particular, the benchmark problems involving particles have confirmed the algorithm is stable for very light bodies including a zero-mass body
and the added-damping effects have to been handled properly to successfully simulate the problems.
The linear system derived from the AMP interface condition has been shown to have similar conditioning of a pressure Poisson system in a TP scheme.
Benchmark problems involving mechanical heart valves show the applicability of the~\ampRB~scheme to practical engineering problems
and also demonstrate the advantageous of the~\ampRB~scheme.
All the benchmark problems have shown that the~\ampRB~scheme is significantly more efficient than the TP scheme that requires sub-iterations to remain stable.
The algorithm has been implemented in both serial and parallel, and the parallel performance test shows reasonable scaling.

A variety of future directions exist for this work.
From a computational perspective, improving the parallel performance of the algorithm would be quite useful, and an important first step would
be to address the poor scaling of the linear system solves.
The multigrid solver in~\cite{CGMG,automg} has been shown to be much faster and more memory efficient than
Krylov solvers, but this multigrid solver would need to be extended to handle the AMP interface conditions and optimized
in parallel. Optimizations of the grid-generator for moving grids and for parallel may also yield significant performance gains. 
There are also many interesting applications of the present work, for example the heart valve simulations, and steps toward
realistic Reynolds numbers, more accurate representation of leaflet geometry, and inlet boundary conditions generalized from
experimental data are all interesting topics for future work.

\appendix
\section{Collision model}
\def\Dmn{D}
\label{sec:collision}

Glowinski et al.~\cite{glowinski2000distributed} proposed a simple collision model based on a continuous repulsive force that depends on  the distance between two rigid particles
\begin{align*}
    \Fv(\Dmn) = \begin{cases} 
        \zerov & \text{if}\quad \Dmn > R_1 +R_2+\delta, \\
        \frac{1}{\epsilon}  (\frac{\Dmn- R_1-R_2-\delta}\delta)^2 \, \frac{\xv_{b,1}-\xv_{b,2}}\Dmn & \text{if} \quad \Dmn \le R_1+R_2+\delta.
    \end{cases}
    %\label{eq:repulsiveForce}
\end{align*}
Here $\Dmn$ stands for the distance between the centers of two particles, of which the radii are $R_1$ and $R_2$ respectively. 
%Based on the results derived from a simple 1D model problem, 
It is further suggested that the small distance is taken as $\delta \simeq \ds$ and the stiffness parameter is $\epsilon \simeq \ds^2$.
This approach is widely used in  FSI due to its simple formulation and straightforward implementation. 
In the implementation the force is only applied to the right hand side of the Newton-Euler equations as an external body force, while the FSI algorithms remain the same.
In many simulations using this idea, $\delta$ is taken as $3\ds$ and $\epsilon$ is taken between $10^{-4}$ and $10^{-8}$.

In the current work,  this collision model has been extended to a repulsive torque applied to the leaflet in the mechanical heart valve in Section~\ref{sec:heartValveIntro}.
In the literature, it appears that a post-processing approach is often used to restrict the rotation angle in the designed range, for instance, simply setting a hard limit on the range of the rotation angle.
However, we find that such a post-processing approach leads to a significant mismatch between the leaflet speed and fluid velocity, and eventually the fluid solver becomes unstable. 
Note the mismatch in the~\ampRB~scheme  will not be diminished by sub-iterations while the sub-iterations may help to stabilize the TP-RB scheme in such a case.

\newcommand{\rfplanar}{{MP-RF}}
\subsection{A model problem for repulsive forces}
\label{sec:collisionModel}

{
\newcommand{\lbfont}{\small}
\def\ysb{-.24} % bottom of rigid body 
\def\ysa{-2}   % top of rigid body 
\def\ya{0} % bottom of fluid domain 
\def\yb{2} % top of fluid domain
\def\xL{8}
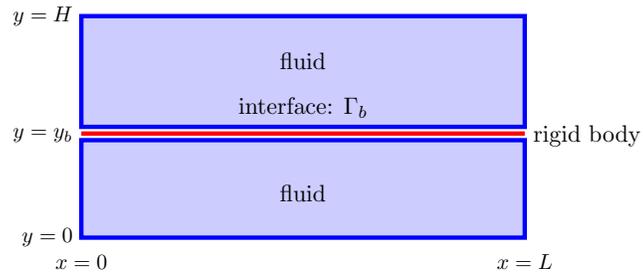
\begin{figure}[hbt]
	\newcommand{\textFont}{\normalss}
	\begin{center}
           \resizebox{6cm}{!}{% START resize box
		\begin{tikzpicture}[scale=.9]
		\useasboundingbox (0,\ysa) rectangle (\xL,\yb+.2);  % set the bounding box (so we have less surrounding white space)
		\draw[thick,fill=blue!20,draw=blue,line width=2pt] (0,\ysa) rectangle (\xL,\ysb);
		% \draw[thick,black] (0,\ysa) rectangle (\xL,\ysb);
                \draw[thick,black] (4,-1.2) node {fluid};
		\draw[thick,fill=blue!20,draw=blue,line width=2pt] (0,\ya) rectangle (\xL,\yb);
		\draw[thick,blue] (0,\ya) rectangle (\xL,\yb);
		\draw[thick,black] (4,1.2) node {fluid};
		\draw[thick,black] (4,\ya) node[anchor=south] {interface: $\GammaB$};

		% -- solid
		\draw[-,thick,blue] (0  ,\ysa) -- (\xL,\ysa);
		\draw[-,thick,red,line width=2pt] (0  ,\ysb+.12) -- (\xL,\ysb+.12);
		\draw[-,thick,blue] (0  ,\ya) -- (0  ,\ya) node[anchor=east,black,yshift=-4pt] {\lbfont$y=y_b$};
		\draw[-,thick,blue] (\xL,\ya) -- (\xL,\ya) node[anchor=west,black,yshift=-4pt] {rigid body};
		\draw[-,thick,blue] (\xL,\ysa) -- (\xL,\ysb);
		\draw[-,thick,blue] (0,\yb) -- (\xL,\yb);
		\draw[-,thick,blue] (\xL,0) -- (\xL,\yb);
		\draw (0  ,\ysa) node[anchor=north,black,yshift=-4pt] {\lbfont$x=0$};
		\draw (\xL,\ysa) node[anchor=north,black,yshift=-4pt] {\lbfont$x=L$};
		\draw (0,  \ysa) node[anchor=east,black] {\lbfont$y=0$};
		\draw (0,   \yb) node[anchor=east,black] {\lbfont$y=H$};
		%
		%\draw[blue] (4,\yH) node[anchor=south,yshift=5pt] {\textFont $y=\eta(x,t)$};
		% grid:
		%\draw[step=1cm,gray] (0,\ysa) grid (8,\yb);
		%
		\end{tikzpicture}
            }% end resize box
	\end{center}
	\caption{The geometry for the model problem for repulsive forces.} \label{fig:collisionModelProblem}
\end{figure}
}

In this section, a model problem is used to investigate the effects of the repulsive force and the damping term applied to the rigid body.
Consider a rectangular rigid body of mass $m_b$ immersed in an incompressible fluid. 
The geometry of the model problem is given in  Figure~\ref{fig:collisionModelProblem}, in which
a rigid body is located at $y_b$.
The height of this rigid body has no effects in this case, 
and so the rigid body is assumed to be infinitely thin and denoted as a red line.
The motion of the full system is driven by the pressure gradient between the top and bottom boundaries of fluid. 
We assume the pressure at $y=H$ is $p_H(t)$, the pressure at $y=0$ is $p_0(t)$, and further assume $p_H(t)>p_0(t)$.
A repulsive force $g_\text{rf}(y)$ is applied to the rigid body with
\begin{align}
    g_\text{rf}(y) = \begin{cases} 
        0 & \text{if}\quad y > y_0+\delta, \\
        \frac{1}{\epsilon}  (\frac{y- y_0-\delta}\delta)^2 & \text{if} \quad y_0 \le y \le y_0+\delta, \\
        \frac{1}{\epsilon}   & \text{if} \quad  y < y_0.
    \end{cases}
    \label{eq:MP-RF-grf}
\end{align}
Initially, $y_b(0)>y_0+\delta$.
This repulsive force can be interpreted as a penalty that restricts the rigid body motion.
%For instance, the leaflet of BMHV can only rotate in certain range due to the special hinge system, in which case this repulsive force is modeling the force between the leaflets and the hinge system when they collide.
Following similar arguments for the model problems in~\cite{rbinsmp2017}, 
this model problem can be assumed to be only $y$-dependent.
Therefore, the model with repulsive forces included (\rfplanar) is given by
\begin{equation}
\text{{\rfplanar}:}\;  
\left\{ 
  \begin{alignedat}{3}
 & \rho \frac{\partial v}{\partial t} +\frac{\partial p}{\partial y} = 0 , \quad&&  y\in(0,y_b)\cup(y_b,H),  \\
 &  \frac{\partial v}{\partial y} = 0   , \quad&& y\in(0,y_b)\cup(y_b,H), \\
 &   m_b \frac{d v_b}{dt} = L\,(p(y_b^-,t)- p(y_b^+,t)) + g_\text{rf}(y_b) -B(y_b)\,v_b, \\
 &  v(y_b^+,t)=v(y_b^-,t)=v_b, \quad p(H,t)=p_H(t), \quad p(0,t) = p_0(t).
  \end{alignedat}  \right. 
 \label{eq:MP-RF}
\end{equation}
Note that there is an extra damping term $B(y_b)\,v_b$ involved in the governing equation of the rigid body.
Through numerical experiments, we find that it is important to involve a damping term to model the energy loss during the collision. The damping term is applied to the rigid body when the repulsive force is turned on.
The damping coefficient $B(y)$ follows a similar fashion to $g_\text{rf}(y)$ and is given by
\begin{align}
    B(y) = \begin{cases} 
        0 & \text{if}\quad y > y_0+\delta, \\
        B_0  (\frac{y- y_0-\delta}\delta)^2 & \text{if} \quad y_0 \le y \le y_0+\delta, \\
        B_0   & \text{if} \quad  y < y_0,
    \end{cases}
    \label{eq:dampingTerm}
\end{align}
with the coefficient $B_0>0$.

It is easy to show that the solutions of this model problem satisfy
\begin{align}
    v(y,t) & = v_b(t),\\
    (m_b+M_a) \frac{d v_b}{dt} & = L\,(p_0(t)- p_H(t)) + g_\text{rf}(y_b) -B(y_b)\,v_b,
    \label{eq:MP-RFsolution}
\end{align}
with the added mass being $M_a = \rho L H$.
Therefore, the effect of the whole fluid in this model problem is fully described 
by its added mass in the isolated equation of the rigid body motion~\eqref{eq:MP-RFsolution}.
The isolated equation~\eqref{eq:MP-RFsolution} can be normalized as 
\begin{align}
    \frac{d v_b}{dt} & = \tilde f(t) + \tilde g_\text{rf}(y_b) - \tilde B(y_b)\,v_b.
\end{align}
Here $\tilde f(t) = L(p_t(t) - p_H(t))/(m_b + M_a)$, $\tilde g_\text{rf}(y_b)$ and $\tilde B(y_b)$ follow the same formulations given in~\eqref{eq:MP-RF-grf} and~\eqref{eq:dampingTerm}, respectively, except $\tilde \epsilon = \epsilon (m_b + M_a)$ and $\tilde B_0 = B_0/(m_b + M_a)$.

{% ------ MATLAB CURVES  ------
\newcommand{\labelFont}{\scriptsize}
\newcommand{\trimfig}[2]{\trimw{#1}{#2}{.0}{.0}{.0}{.0}}
\newcommand{\trimfigz}[2]{\trimw{#1}{#2}{.05}{.05}{.03}{.05}}
\newcommand{\figWidth}{3.9cm}
\begin{figure}[htb]
\begin{center}
\begin{tikzpicture}[scale=1]
  \useasboundingbox (0.0,1) rectangle (16,3.6);  % set the bounding box (so we have less surrounding white space)
  \begin{scope}[xshift=-.5cm,yshift=0cm]
    \draw(0.0,0.0) node[anchor=south west,xshift=-4pt,yshift=+0pt] {\trimfig{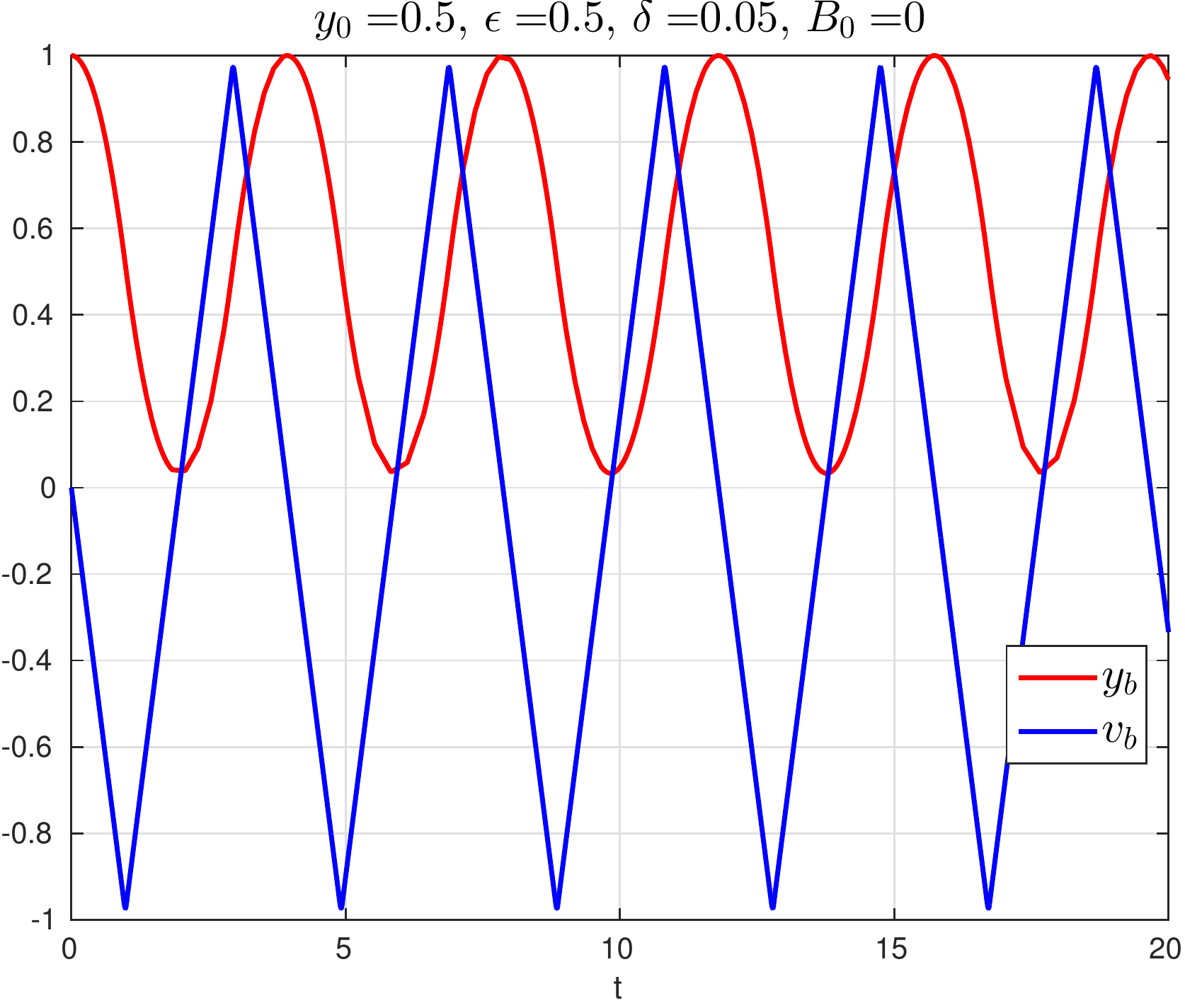}{\figWidth}};
    \draw(2.4,0.8) node[draw,fill=white] {\labelFont I};
    \draw(4.0,0.0) node[anchor=south west,xshift=-4pt,yshift=+0pt] {\trimfig{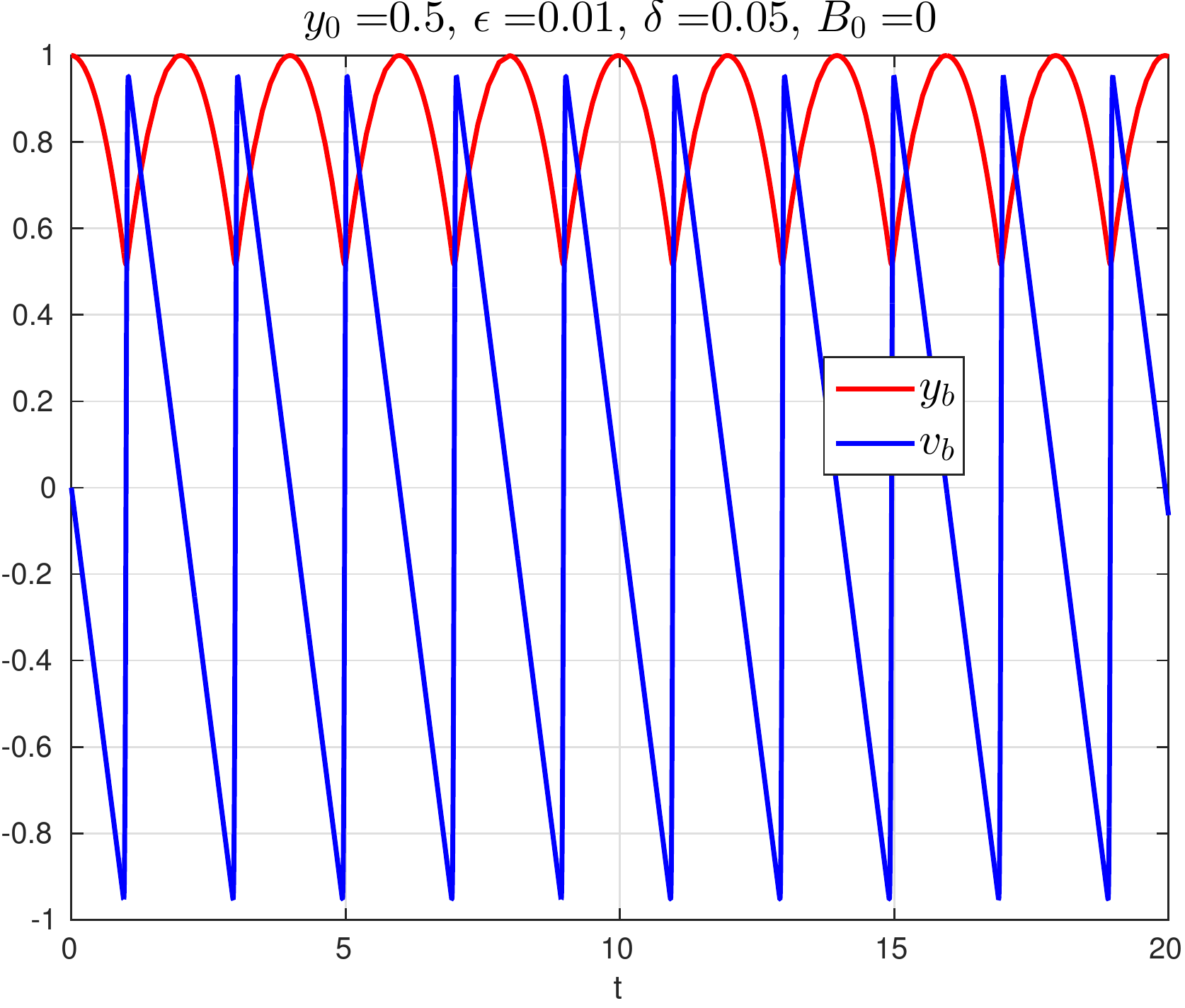}{\figWidth}};
    \draw(6.4,0.8) node[draw,fill=white] {\labelFont II};
    \draw(8.0,0.0) node[anchor=south west,xshift=-4pt,yshift=+0pt] {\trimfig{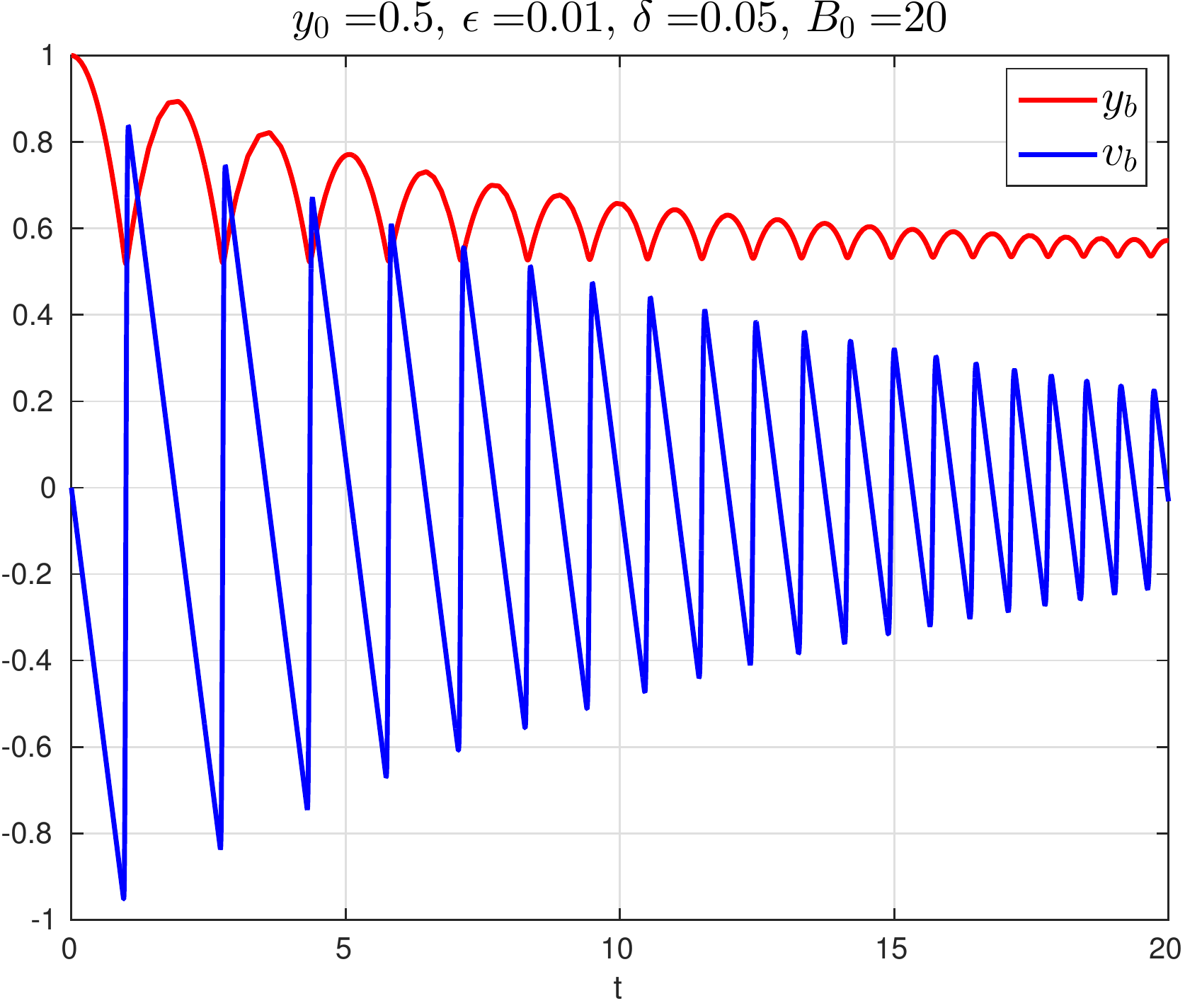}{\figWidth}};
    \draw(10.4,0.8) node[draw,fill=white] {\labelFont III};
    \draw(12.0,0.0) node[anchor=south west,xshift=-4pt,yshift=+0pt] {\trimfig{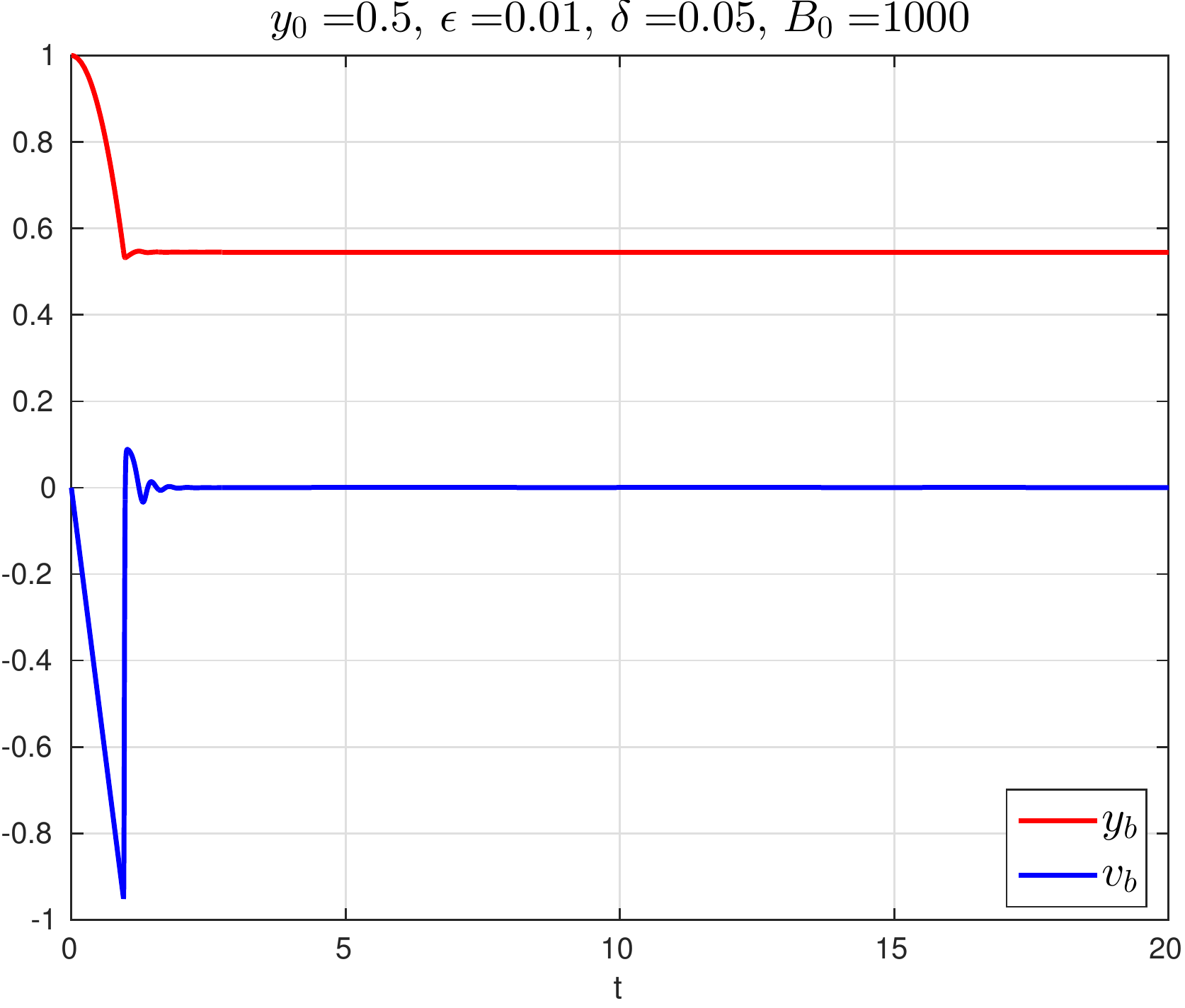}{\figWidth}};
    \draw(14.4,0.8) node[draw,fill=white] {\labelFont IV};
  \end{scope}
%
% grid:
%\draw[step=1cm,gray] (0,0) grid (16,3.6);
\end{tikzpicture}
\end{center}
  \caption{MP-RF.  Velocity and location of the rigid body.
   (I)~$\epsilon$ is too large and there is no damping;
   (II)~$\epsilon$ is small and there is no damping;
   (III)~$\epsilon$ is small and the solution is under-damped;
   (IV)~$\epsilon$ is small and the solution is over-damped.
    }
  \label{fig:collisionModel}
\end{figure}
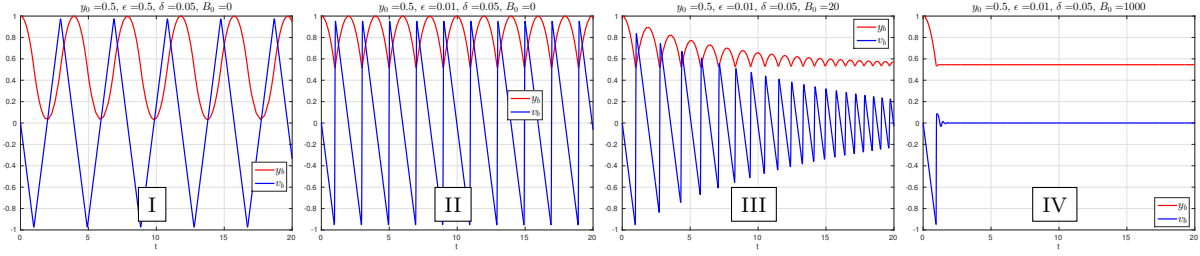
}

For simplicity, we consider the case when $f(t) = -1$, the initial location of the rigid body is $y_b(0)=1$ with the speed $v_b(0) = 0$,
and the rigid body is assumed to collide at $y_0=0.5$.
%with some parts of the system that are not explicitly described.
The results of some typical choices of the parameters $\tilde \epsilon$, $\delta$ and $\tilde B_0$ are presented in Figure~\ref{fig:collisionModel}.
In Case I, $\tilde \epsilon$ is too big so that the resulting repulsive force is too small to slow down the rigid body above $y_0$. In practice, this may lead to the interpenetration of multiple grids in the composite grid framework, which typically leads to the failure of simulation.
In this case, the energy is not conserved due to the formulation of the repulsive force.
In Case II, $\tilde \epsilon$ is small enough so that the repulsive force is able to slow down the rigid body before it hits the bound. It is easy to show that the repulsive force conserves the total energy of the system in this case.
Assume at the time $\tau$, the height of the rigid body is $y_0+\delta$ where the repulsive force starts to affect. Then the system satisfies,
\begin{align}
    \ddot{y}_b & = -1 + \frac{1}{\tilde \epsilon} (\frac{y_b-y_0-\delta}{\delta})^2 ,
    \label{eq:simpleMR-RF}
\end{align}
with $\dot y_b(\tau) = v_\tau$ and $y_b(\tau) = y_0+\delta$.
Multiplying $\dot y_b$ at both side of~\eqref{eq:simpleMR-RF} and integrating leads to
\begin{align}
    \frac{\dot y_b^2}2 + y_b- \frac{1}{3 \tilde \epsilon \delta^2} ({y_b-y_0-\delta})^3=C ,
    \label{eq:simpleMR-RFeng}
\end{align}
which is true when the rigid body is in $[y_0, y_0+\delta]$.
The left hand side of~\eqref{eq:simpleMR-RFeng} can be interpreted as the total energy of the system
that consists of the kinetic energy, the potential energy and the energy related to the repulsive force.
Therefore, the energy has been shown to be conserved without the damping term.

%Although it is useful to conserve the energy in the collision, 
However, in some cases a significant amount of  total energy can be lost through other complicated mechanisms during the collision. The damping term plays an important role in our model to model the energy loss. In Case (III), a damping term is turned on and the total energy of the system decays as the collision happens.
The solution is under-damped and slowly converges to a constant.
In Case (IV), a much larger damping term is used so that the motion of the rigid body is over-damped, in which case the solution quickly converges to a constant.
This is  useful to model the situation where significant amount of energy is lost during the collisions.

\subsection{Applications to the heart valve simulations}
In the mechanical heart valve simulations, the leaflets are designed to rotate only in a certain range of angles. 
Therefore, an extra mechanism is needed to force the leaflets stop before they cross the designed bounds.
In the current simulations, this collision model has been used as the constraint on the range of the leaflet angle. 
%Consider two dimensions for simplicity. due to the special hinge system and the collisions between the housing and leaflets
Assume the range of the rotational angle of leaflets $\theta$ be $[\theta_{\rm min}, \theta_{\rm max}]$.
The following external torque and damping term are used
\begin{align}
    g_{\rm rt}(\theta) &= \begin{cases} 
        \frac{1}{\epsilon_1}    & \text{if} \quad \theta \le \theta_{\rm min}, \\
        \frac{1}{\epsilon_1}  (\frac{\theta - \theta_{\rm min}-\delta}\delta)^2  & \text{if} \quad \theta_{\rm min} < \theta \le \theta_{\rm min}+\delta, \\
        0 & \text{if}\quad \theta_{\rm min}+\delta <\theta < \theta_{\rm max}-\delta, \\
        - \frac{1}{\epsilon_2}  (\frac{\theta- \theta_{\rm max}+\delta}\delta)^2  & \text{if} \quad \theta_{\rm max}-\delta \le \theta \le \theta_{\rm max}, \\
       - \frac{1}{\epsilon_2}    & \text{if} \quad \theta > \theta_{\max}. \\
    \end{cases}
    \label{eq:repulsiveTorque}
    \\
    B(\theta) &= \begin{cases} 
        B_0    & \text{if} \quad \theta < \theta_{\min} \text{ or }  \theta > \theta_{\max}. \\
        B_0  (\frac{\theta - \theta_{\rm min}-\delta}\delta)^2  & \text{if} \quad \theta_{\rm min} < \theta \le \theta_{\rm min}+\delta, \\
        0 & \text{if}\quad \theta_{\rm min}+\delta <\theta < \theta_{\rm max}-\delta, \\
        B_0  (\frac{\theta- \theta_{\rm max}+\delta}\delta)^2  & \text{if} \quad \theta_{\rm max}-\delta \le \theta \le \theta_{\rm max}, \\
    \end{cases}
    \label{eq:repulsiveDamp}
\end{align}
%which is used to reduce the angular speed of the leaflet when it gets close to the bounds of the range.
This simple model can be interpreted as the results of the collision between the leaflets  and  the hinge, although the hinge is not explicitly described in the current simulation. 
Note the repulsive torque and damping term are a direct extension of the above collision model.
This model shows much better performance in our implementation than
simply setting a hard limit of the angles.
The damping term is particularly important during the process of closing, where the significant amount of energy are lost through the collision between the leaflet and the hinge. 
See Section~\ref{sec:heartValveIntro}.

\bibliographystyle{../../elsart-num}
\bibliography{../../journal-ISI,../../jwb,../../henshaw,../../henshawPapers,../../fsi}

\begin{thebibliography}{10}
\expandafter\ifx\csname url\endcsname\relax
  \def\url#1{\texttt{#1}}\fi
\expandafter\ifx\csname urlprefix\endcsname\relax\def\urlprefix{URL }\fi

\bibitem{rbinsmp2017}
J.~W. Banks, W.~D. Henshaw, D.~W. Schwendeman, Q.~Tang, A stable partitioned
  {FSI} algorithm for rigid bodies and incompressible flow. {Part I}: Model
  problem analysis, J. Comput. Phys. 343 (2017) 432--468.

\bibitem{rbins2017}
J.~W. Banks, W.~D. Henshaw, D.~W. Schwendeman, Q.~Tang, A stable partitioned
  {FSI} algorithm for rigid bodies and incompressible flow. {Part II}: General
  formulation, J. Comput. Phys. 343 (2017) 469--500.

\bibitem{mog2006}
W.~D. Henshaw, D.~W. Schwendeman, Moving overlapping grids with adaptive mesh
  refinement for high-speed reactive and non-reactive flow, J. Comput. Phys.
  216~(2) (2006) 744--779\citeCount{54}.

\bibitem{ICNS}
W.~D. Henshaw, A fourth-order accurate method for the incompressible
  {N}avier-{S}tokes equations on overlapping grids, J. Comput. Phys. 113~(1)
  (1994) 13--25\citeCount{154}.

\bibitem{splitStep2003}
W.~D. Henshaw, N.~A. Petersson, A split-step scheme for the incompressible
  {Navier-Stokes} equations, in: M.~M. Hafez (Ed.), Numerical Simulation of
  Incompressible Flows, World Scientific, 2003, pp. 108--125\citeCount{30}.

\bibitem{max2006b}
W.~D. Henshaw, A high-order accurate parallel solver for {Maxwell}'s equations
  on overlapping grids, SIAM J. Sci. Comput. 28~(5) (2006)
  1730--1765\citeCount{13}.

\bibitem{pog2008a}
W.~D. Henshaw, D.~W. Schwendeman, Parallel computation of three-dimensional
  flows using overlapping grids with adaptive mesh refinement, J. Comput. Phys.
  227~(16) (2008) 7469--7502\citeCount{45}.

\bibitem{cheng2004three}
R.~Cheng, Y.~G. Lai, K.~B. Chandran, Three-dimensional fluid-structure
  interaction simulation of bileaflet mechanical heart valve flow dynamics,
  Ann. Biomed. Eng. 32~(11) (2004) 1471--1483.

\bibitem{dumont2007comparison}
K.~Dumont, J.~Vierendeels, R.~Kaminsky, G.~Van~Nooten, P.~Verdonck,
  D.~Bluestein, Comparison of the hemodynamic and thrombogenic performance of
  two bileaflet mechanical heart valves using a {CFD/FSI} model, J. Biomech.
  Eng. 129~(4) (2007) 558--565.

\bibitem{tai2007numerical}
C.~H. Tai, K.~M. Liew, Y.~Zhao, Numerical simulation of {3D} fluid--structure
  interaction flow using an immersed object method with overlapping grids,
  Comput. Struct. 85~(11) (2007) 749--762.

\bibitem{borazjani2008curvilinear}
I.~Borazjani, L.~Ge, F.~Sotiropoulos, Curvilinear immersed boundary method for
  simulating fluid structure interaction with complex {3D} rigid bodies, J.
  Comput. Phys. 227~(16) (2008) 7587--7620.

\bibitem{nobili2008numerical}
M.~Nobili, U.~Morbiducci, R.~Ponzini, C.~Del~Gaudio, A.~Balducci, M.~Grigioni,
  F.~M. Montevecchi, A.~Redaelli, Numerical simulation of the dynamics of a
  bileaflet prosthetic heart valve using a fluid--structure interaction
  approach, J. Biomech. 41~(11) (2008) 2539--2550.

\bibitem{de2009direct}
M.~D. De~Tullio, A.~Cristallo, E.~Balaras, R.~Verzicco, Direct numerical
  simulation of the pulsatile flow through an aortic bileaflet mechanical heart
  valve, J. Fluid Mech. 622 (2009) 259--290.

\bibitem{borazjani2010high}
I.~Borazjani, L.~Ge, F.~Sotiropoulos, High-resolution fluid--structure
  interaction simulations of flow through a bi-leaflet mechanical heart valve
  in an anatomic aorta, Ann. Biomed. Eng. 38~(2) (2010) 326--344.

\bibitem{sotiropoulos2009review}
F.~Sotiropoulos, I.~Borazjani, A review of state-of-the-art numerical methods
  for simulating flow through mechanical heart valves, Med. Biol. Eng. Comput.
  47~(3) (2009) 245--256.

\bibitem{takashi1992}
N.~Takashi, T.~J. Hughes, An arbitrary {Lagrangian-Eulerian} finite element
  method for interaction of fluid and a rigid body, Comput. Method. Appl. Mech.
  Engrg. 95~(1) (1992) 115--138.

\bibitem{HuPatankarZhu2001}
H.~H. Hu, N.~A. Patankar, M.~Y. Zhu, Direct numerical simulations of
  fluid-solid systems using the arbitrary langrangian-eulerian technique, J.
  Comput. Phys. 169~(2) (2001) 427--462.

\bibitem{vierendeels2005analysis}
J.~Vierendeels, K.~Dumont, E.~Dick, P.~Verdonck, Analysis and stabilization of
  fluid-structure interaction algorithm for rigid-body motion, AIAA J. 43~(12)
  (2005) 2549--2557.

\bibitem{coquerelle2008vortex}
M.~Coquerelle, G.-H. Cottet, A vortex level set method for the two-way coupling
  of an incompressible fluid with colliding rigid bodies, J. Comput. Phys.
  227~(21) (2008) 9121--9137.

\bibitem{gibou2012efficient}
F.~Gibou, C.~Min, Efficient symmetric positive definite second-order accurate
  monolithic solver for fluid/solid interactions, J. Comput. Phys. 231~(8)
  (2012) 3246--3263.

\bibitem{glowinski1999distributed}
R.~Glowinski, T.-W. Pan, T.~I. Hesla, D.~D. Joseph, A distributed {Lagrange}
  multiplier/fictitious domain method for particulate flows, Int. J. Multiphase
  Flow 25~(5) (1999) 755--794.

\bibitem{costarelli2016embedded}
S.~D. Costarelli, L.~Garelli, M.~A. Cruchaga, M.~A. Storti, R.~Ausensi, S.~R.
  Idelsohn, An embedded strategy for the analysis of fluid structure
  interaction problems, Comput. Method. Appl. Mech. Engrg. 300 (2016) 106--128.

\bibitem{kajishima2002interaction}
T.~Kajishima, S.~Takiguchi, Interaction between particle clusters and
  particle-induced turbulence, Int. J. Heat Fluid Flow 23~(5) (2002) 639--646.

\bibitem{uhlmann2005immersed}
M.~Uhlmann, An immersed boundary method with direct forcing for the simulation
  of particulate flows, J. Comput. Phys. 209~(2) (2005) 448--476.

\bibitem{kim2006immersed}
D.~Kim, H.~Choi, Immersed boundary method for flow around an arbitrarily moving
  body, J. Comput. Phys. 212~(2) (2006) 662--680.

\bibitem{lee2008immersed}
T.-R. Lee, Y.-S. Chang, J.-B. Choi, D.~W. Kim, W.~K. Liu, Y.-J. Kim, Immersed
  finite element method for rigid body motions in the incompressible
  {Navier}--{Stokes} flow, Comput. Method. Appl. Mech. Engrg. 197~(25) (2008)
  2305--2316.

\bibitem{breugem2012second}
W.-P. Breugem, A second-order accurate immersed boundary method for fully
  resolved simulations of particle-laden flows, J. Comput. Phys. 231~(13)
  (2012) 4469--4498.

\bibitem{kempe2012improved}
T.~Kempe, J.~Fr{\"o}hlich, An improved immersed boundary method with direct
  forcing for the simulation of particle laden flows, J. Comput. Phys. 231~(9)
  (2012) 3663--3684.

\bibitem{yang2012simple}
J.~Yang, F.~Stern, A simple and efficient direct forcing immersed boundary
  framework for fluid--structure interactions, J. Comput. Phys. 231~(15) (2012)
  5029--5061.

\bibitem{bhalla2013unified}
A.~P.~S. Bhalla, R.~Bale, B.~E. Griffith, N.~A. Patankar, A unified
  mathematical framework and an adaptive numerical method for fluid--structure
  interaction with rigid, deforming, and elastic bodies, J. Comput. Phys. 250
  (2013) 446--476.

\bibitem{yang2015non}
J.~Yang, F.~Stern, A non-iterative direct forcing immersed boundary method for
  strongly-coupled fluid--solid interactions, J. Comput. Phys. 295 (2015)
  779--804.

\bibitem{wang2015strongly}
C.~Wang, J.~D. Eldredge, Strongly coupled dynamics of fluids and rigid-body
  systems with the immersed boundary projection method, J. Comput. Phys. 295
  (2015) 87--113.

\bibitem{kim2016penalty}
Y.~Kim, C.~S. Peskin, A penalty immersed boundary method for a rigid body in
  fluid, Phys. Fluids 28~(3) (2016) 033603.

\bibitem{lacis2016stable}
U.~L{\=a}cis, K.~Taira, S.~Bagheri, A stable fluid--structure-interaction
  solver for low-density rigid bodies using the immersed boundary projection
  method, J. Comput. Phys. 305 (2016) 300--318.

\bibitem{corona2017integral}
E.~Corona, L.~Greengard, M.~Rachh, S.~Veerapaneni, An integral equation
  formulation for rigid bodies in stokes flow in three dimensions, J. Comput.
  Phys. 332 (2017) 504--519.

\bibitem{Saye2017a}
R.~Saye, Implicit mesh discontinuous {Galerkin} methods and interfacial gauge
  methods for high-order accurate interface dynamics, with applications to
  surface tension dynamics, rigid body fluid–structure interaction, and free
  surface flow: {Part I}, J. Comput. Phys. 344 (2017) 647--682.

\bibitem{Saye2017b}
R.~Saye, Implicit mesh discontinuous {Galerkin} methods and interfacial gauge
  methods for high-order accurate interface dynamics, with applications to
  surface tension dynamics, rigid body fluid–structure interaction, and free
  surface flow: {Part II}, J. Comput. Phys. 344 (2017) 683--723.

\bibitem{Koblitz2016}
A.~R. Koblitz, S.~Lovett, N.~Nikiforakis, W.~D. Henshaw, Direct numerical
  simulation of particulate flows with an overset grid method, J. Comput. Phys.
  343 (2017) 414--431.

\bibitem{kadapa2017}
C.~Kadapa, W.~G. Dettmer, D.~Peri{\'c}, A stabilised immersed boundary method
  on hierarchical b-spline grids for fluid-rigid body interaction with
  solid--solid contact, Comput. Method. Appl. Mech. Engrg. 318 (2017) 242--269.

\bibitem{schwarz2015temporal}
S.~Schwarz, T.~Kempe, J.~Fr{\"o}hlich, A temporal discretization scheme to
  compute the motion of light particles in viscous flows by an immersed
  boundary method, J. Comput. Phys. 281 (2015) 591--613.

\bibitem{tschisgale2017non}
S.~Tschisgale, T.~Kempe, J.~Fr{\"o}hlich, A non-iterative immersed boundary
  method for spherical particles of arbitrary density ratio, J. Comput. Phys.
  339 (2017) 432--452.

\bibitem{lrb2013}
J.~W. Banks, W.~D. Henshaw, B.~Sj{\"o}green, A stable {FSI} algorithm for light
  rigid bodies in compressible flow, J. Comput. Phys. 245 (2013)
  399--430\citeCount{1}.

\bibitem{fsi2012}
J.~W. Banks, W.~D. Henshaw, D.~W. Schwendeman, Deforming composite grids for
  solving fluid structure problems, J. Comput. Phys. 231~(9) (2012)
  3518--3547\citeCount{5}.

\bibitem{flunsi2016}
J.~W. Banks, W.~D. Henshaw, A.~Kapila, D.~W. Schwendeman, An added-mass
  partitioned algorithm for fluid-structure interactions of compressible fluids
  and nonlinear solids, J. Comput. Phys. 305 (2016) 1037--1064\citeCount{0}.

\bibitem{fib2014}
J.~W. Banks, W.~D. Henshaw, D.~W. Schwendeman, An analysis of a new stable
  partitioned algorithm for {FSI} problems. {Part I}: Incompressible flow and
  elastic solids, J. Comput. Phys. 269 (2014) 108--137\citeCount{2}.

\bibitem{fis2014}
J.~W. Banks, W.~D. Henshaw, D.~W. Schwendeman, An analysis of a new stable
  partitioned algorithm for {FSI} problems. {Part II}: Incompressible flow and
  structural shells, J. Comput. Phys. 268 (2014) 399--416\citeCount{2}.

\bibitem{beamins2016}
L.~Li, W.~D. Henshaw, J.~W. Banks, D.~W. Schwendeman, G.~A. Main, A stable
  partitioned {FSI} algorithm for incompressible flow and deforming beams, J.
  Comput. Phys. 312 (2016) 272--306.

\bibitem{yoon2018}
G.~Yoon, C.~Min, S.~Kim, A stable and convergent {Hodge} decomposition method
  for fluid--solid interaction, J. Sci. Comput. 76~(2) (2018) 727--758.

\bibitem{Petersson00}
N.~A. Petersson, Stability of pressure boundary conditions for {S}tokes and
  {N}avier-{S}tokes equations, J. Comput. Phys. 172~(1) (2001) 40--70.

\bibitem{YihFluidMechanics}
C.-S. Yih, Fluid Mechanics, West River Press, 1977.

\bibitem{ogen}
W.~D. Henshaw, Ogen: An overlapping grid generator for {O}verture, Research
  Report UCRL-MA-132237, Lawrence Livermore National Laboratory (1998).

\bibitem{ChanHybridMesh2009}
W.~M. Chan, Enhancements to the hybrid mesh approach to surface loads
  integration on overset structured grids, AIAA paper 2009-3990.

\bibitem{CGCN}
G.~S. Chesshire, W.~D. Henshaw, Conservation on composite overlapping grids,
  {IBM} Research Report RC 16531, IBM Research Division, Yorktown Heights, NY
  (1991).

\bibitem{CGCN94}
G.~S. Chesshire, W.~D. Henshaw, A scheme for conservative interpolation on
  overlapping grids, SIAM J. Sci. Comput. 15~(4) (1994) 819--845\citeCount{60}.

\bibitem{CginsUserGuide}
W.~D. Henshaw, {Cgins} user guide: An {O}verture solver for the incompressible
  {N}avier-{S}tokes equations on composite overlapping grids, Software Manual
  LLNL-SM-455851, Lawrence Livermore National Laboratory (2010).

\bibitem{petsc-user-ref}
S.~Balay, S.~Abhyankar, M.~F. Adams, J.~Brown, P.~Brune, K.~Buschelman,
  L.~Dalcin, V.~Eijkhout, W.~D. Gropp, D.~Kaushik, M.~G. Knepley, L.~C.
  McInnes, K.~Rupp, B.~F. Smith, S.~Zampini, H.~Zhang, H.~Zhang, {PETS}c users
  manual, Tech. Rep. ANL-95/11 - Revision 3.7, Argonne National Laboratory
  (2016).

\bibitem{Cate2002}
A.~ten Cate, C.~H. Nieuwstad, J.~J. Derksen, H.~E. A.~V. den Akker, Particle
  imaging velocimetry experiments and lattice-boltzmann simulations on a single
  sphere settling under gravity, Physics of Fluids 14~(11) (2002) 4012--4025.

\bibitem{mordant2000velocity}
N.~Mordant, J.-F. Pinton, Velocity measurement of a settling sphere, Eur. Phys.
  J. B, Condens. Matter Complex Phys. 18~(2) (2000) 343--352.

\bibitem{yang2014sharp}
J.~Yang, F.~Stern, A sharp interface direct forcing immersed boundary approach
  for fully resolved simulations of particulate flows, J. Fluids Eng. 136~(4)
  (2014) 040904.

\bibitem{cadProjection02}
W.~D. Henshaw, An algorithm for projecting points onto a patched {CAD} model,
  Engineering with Computers 18 (2002) 265--273\citeCount{17}.

\bibitem{CGMG}
W.~D. Henshaw, G.~S. Chesshire, Multigrid on composite meshes, SIAM J. Sci.
  Stat. Comput. 8~(6) (1987) 914--923\citeCount{34}.

\bibitem{automg}
W.~D. Henshaw, On multigrid for overlapping grids, SIAM J. Sci. Comput. 26~(5)
  (2005) 1547--1572\citeCount{16}.

\bibitem{glowinski2000distributed}
R.~Glowinski, T.-W. Pan, T.~I. Hesla, D.~D. Joseph, J.~P\'eriaux, A distributed
  {Lagrange} multiplier/fictitious domain method for the simulation of flow
  around moving rigid bodies: application to particulate flow, Comput. Method.
  Appl. Mech. Engrg. 184~(2) (2000) 241--267.

\end{thebibliography}

\end{document}